%% file: plapSSL_arxiv_v9.tex
\documentclass[letterpaper,11pt]{article}

\usepackage{algpseudocode}
\usepackage{algorithm}
\usepackage{amsfonts}
\usepackage[leqno]{amsmath}
\usepackage{amssymb}
\usepackage{amstext}
\usepackage{amsthm}
\usepackage[toc,page]{appendix}
\usepackage{authblk}
\usepackage{booktabs}
\usepackage[font=small]{caption}
\usepackage{cite}
\usepackage{datetime}
\usepackage{enumitem}
\usepackage{etoolbox}
\usepackage[T1]{fontenc}
\usepackage{framed}
\usepackage{fullpage}
\usepackage[margin=1.2in]{geometry}
\usepackage{graphics}
\usepackage{graphicx}
\usepackage{hyperref}
\usepackage[utf8]{inputenc}
\usepackage{mathrsfs}
\usepackage[framemethod=TikZ]{mdframed}
\usepackage{microtype}
\usepackage{multirow}
\usepackage{nicefrac}
\usepackage{pdflscape}
\usepackage{pgfplots}
\usepackage{rotating}
\usepackage{setspace}
\usepackage[position=b]{subcaption}
\usepackage{tikz}
\usepackage{times}
\usepackage{units}
\usepackage{url}
\usepackage{xcolor}
\usepackage{dsfont}

\providecommand{\keywords}[1]{\textbf{\textit{\noindent  Keywords and phrases. }} #1}
\providecommand{\subjclass}[1]{\textbf{\textit{\noindent Mathematics Subject Classification. }} #1}

\usetikzlibrary{shapes,arrows,plotmarks,matrix,positioning,fit,calc,3d,patterns,decorations.pathreplacing}
\usepgfplotslibrary{groupplots,fillbetween,decorations.softclip}
\pgfplotsset{compat=newest}

\definecolor{darkred}{rgb}{0.75,0,0}
\definecolor{darkgreen}{rgb}{0.1,0.6,0.1}
\definecolor{darkblue}{rgb}{0.1,0.1,0.6}

\newmdenv[
roundcorner=10pt,
backgroundcolor=gray!10,
linecolor=gray!10,
tikzsetting={draw=black,line width=3pt,dashed,dash pattern= on 10pt off 3pt},
outerlinewidth=0pt,
innerlinewidth=0pt,
]{figbox}

\newlength\figureheight
\newlength\figurewidth

\DeclareMathOperator*{\argmax}{\mathrm{argmax}}
\DeclareMathOperator*{\argmin}{\mathrm{argmin}}
\def\dd{\mathrm{d}}
\newcommand{\dist}{\mathrm{dist}}

\def\eps{\varepsilon}

\DeclareMathOperator*{\esssup}{\mathrm{ess}\,\mathrm{sup}}
\newcommand{\err}{\mathrm{err}}

\newcommand{\Id}{\mathrm{Id}}
\newcommand{\Lip}{\mathrm{Lip}}
\newcommand{\osc}{\mathrm{osc}}
\newcommand{\spt}{\mathrm{spt}}
\newcommand{\supp}{\mathrm{supp}}
\newcommand{\Vol}{\mathrm{Vol}}
\DeclareMathOperator*{\Glim}{\Gamma\text{-}\lim}
\DeclareMathOperator*{\weakto}{\rightharpoonup}

\def\l{\left(}
\def\r{\right)}
\def\la{\left|}
\def\ra{\right|}
\def\lb{\left\{}
\def\rb{\right\}}
\def\rd{\right.}
\def\ls{\left[}
\def\rs{\right]}

\newcommand{\bbI}{\mathbb{I}}
\newcommand{\bbN}{\mathbb{N}}
\newcommand{\bbP}{\mathbb{P}}
\newcommand{\bbR}{\mathbb{R}}
\newcommand{\cA}{\mathcal{A}}
\newcommand{\cC}{\mathcal{C}}

\newcommand{\cE}{\mathcal{E}}
\newcommand{\cF}{\mathcal{F}}
\newcommand{\cH}{\mathcal{H}}

\newcommand{\cP}{\mathcal{P}}
\newcommand{\cS}{\mathcal{S}}
\newcommand{\cX}{\mathcal{X}}

\def\Xint#1{\mathchoice
{\XXint\displaystyle\textstyle{#1}}%
{\XXint\textstyle\scriptstyle{#1}}%
{\XXint\scriptstyle\scriptscriptstyle{#1}}%
{\XXint\scriptscriptstyle\scriptscriptstyle{#1}}%
\!\int}
\def\XXint#1#2#3{{\setbox0=\hbox{$#1{#2#3}{\int}$ }
\vcenter{\hbox{$#2#3$ }}\kern-.6\wd0}}

\def\dashint{\Xint-}

\newcommand{\R}{\mathbb{R}}

\newcommand{\M}{\mathcal{M}}

\newcommand{\te}{\textrm}
\newcommand{\tacka}{\,\cdot\,}
\newcommand{\Ecp}{\cE_{\infty}^{(p)}}
\newcommand{\Ecpc}{\cE_{\infty,con}^{(p)}}

\definecolor{mygreen}{rgb}{0.1,0.75,0.2}

\newenvironment{listi}
  {\begin{list}
 {\textbf{(A\arabic{broj})}}
{ \usecounter{broj}}
  \setlength{\itemindent}{-1pt}
  \addtolength{\itemsep}{-3pt}
   \setlength{\labelwidth}{30pt}}
{   \end{list} }

\newtheorem{theorem}{Theorem}[section]
\newtheorem{lemma}[theorem]{Lemma}
\newtheorem{proposition}[theorem]{Proposition}

\newtheorem{mydef}[theorem]{Definition}
\theoremstyle{remark}
\newtheorem{remark}[theorem]{Remark}

\makeatletter
\newcommand{\leqnomode}{\tagsleft@true}
\newcommand{\reqnomode}{\tagsleft@false}

\newcounter{num}
\def\dsl#1{\textbf{\textcolor{darkred}{#1}}}
\def\mt#1{\textbf{\textcolor{darkblue}{#1}}}

\title{Analysis of $p$-Laplacian Regularization in Semi-Supervised Learning}
\author[1]{Dejan Slep\v{c}ev}
\author[2]{Matthew Thorpe}
\affil[1]{Department of Mathematical Sciences,\protect\\ Carnegie Mellon University,\protect\\ Pittsburgh, PA 15213, USA \vspace{\baselineskip}}
\affil[2]{Department of Applied Mathematics and Theoretical Physics,\protect\\ University of Cambridge,\protect\\ Cambridge, CB3 0WA, UK}
\date{October 2017}

\begin{document}

\maketitle
\newcounter{broj}

\begin{abstract}
We investigate a family of  regression problems in a semi-supervised setting. 
The task is to assign real-valued  labels to a set of $n$ sample points, provided a small training subset of $N$ labeled points. 
A goal of semi-supervised learning is to take advantage of the (geometric) structure provided by the large number of unlabeled data when assigning labels. 
We consider random  geometric graphs, with connection radius  $\eps(n)$, to represent the geometry of the data set.  Functionals  which  model the task
reward the regularity of the estimator function and impose or reward the agreement with the training data. Here we consider the discrete $p$-Laplacian regularization.

 We investigate asymptotic behavior when the number of unlabeled points increases, while the number of training  points remains fixed. We uncover a delicate interplay between the regularizing nature of the functionals considered and the nonlocality inherent to the graph constructions.
We rigorously obtain  almost optimal ranges on the scaling of $\eps(n)$ for the asymptotic consistency to hold. We prove that the minimizers of the discrete functionals in random setting converge uniformly to the desired continuum limit. 
Furthermore we discover that for the standard model used there is a restrictive upper bound on how quickly $\eps(n)$ must converge to zero as  $n \to \infty$.  We introduce a new model which is as simple as the original model, but overcomes this restriction. 
\end{abstract}

\noindent
\section{Introduction \label{sec:Intro}}

Due to its applicability across a large spectrum of problems semi-supervised learning (SSL) is an important tool in data analysis. It deals with situations when one has access to relatively few labeled points but potentially a large number of unlabeled data. We assume that we are given $N$ labeled points 
$\{ (x_i, y_i) \: : \: i=1, \dots, N,\; x_i \in \R^d, \; y_i \in \R\}$ and $n-N$  points $x_i$, $\,i=N+1, \dots, n$ drawn from a fixed, but unknown measure, $\mu$ supported in a compact subset of $\R^d$. The goal is to assign labels to the unlabeled points, while taking advantage of the information provided by the unlabeled points when designing the estimator. In particular the unlabeled points carry information on the structure of $\mu$, such as the geometry of its support, which can lead to better estimators. To access the information on $\mu$ in a way that carries over to high dimensions, we build a  graph whose vertices are data points and connect them if they are close enough, that is if they are within some distance $\eps>0$. 
More generally the edge weights are prescribed by using a decreasing function $\eta :[0, \infty) \to [0, \infty)$ with $\lim_{r \to \infty} \eta(r) =0$. For fixed scale $\eps>0$ we set the weights to be 
\[ W_{ij} = \eta_\eps (|x_i - x_j|) \]
where  $\eta_\eps = \frac{1}{\eps^d} \eta(\cdot/\eps)$.

\begin{figure}[th!]
\setlength\figureheight{0.35\textwidth}
\setlength\figurewidth{0.43\textwidth}
\begin{subfigure}[t]{0.47\textwidth}
\centering
\scriptsize
\input{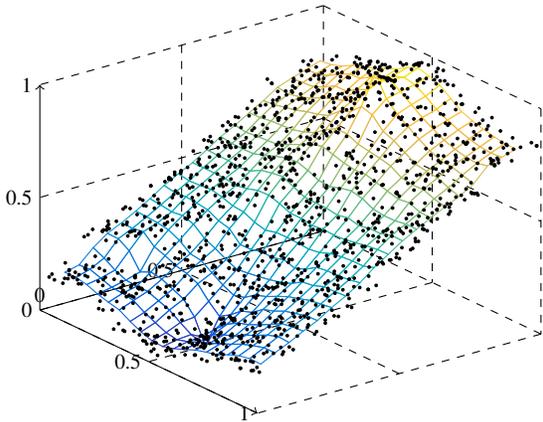}
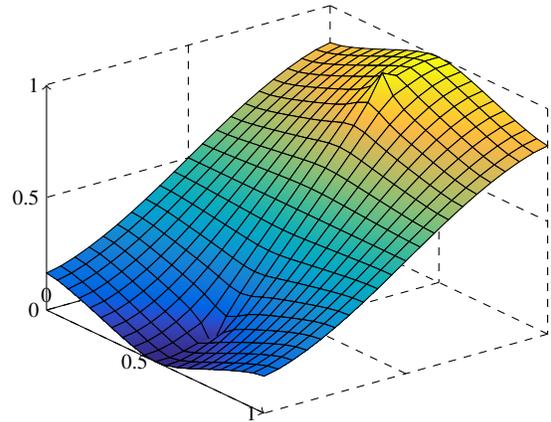
\caption{Minimizer of \eqref{DF} under constraint \eqref{CON} for $\eps= 0.058$ and $\eta = \mathds{1}_{[0,1]}$ and $n=1280$.
The grid is to aid visualization.
}
\end{subfigure}
\hspace*{0.04\textwidth}
\setlength\figureheight{0.35\textwidth}
\setlength\figurewidth{0.43\textwidth}
\begin{subfigure}[t]{0.47\textwidth}
\centering
\scriptsize
\input{ContSol_2D_p4.tikz}
\caption{
Minimizer of the continuum functional \eqref{CF} under  constraint \eqref{CON}. 
}
\end{subfigure}
\caption{
2D numerical experiment for measure $\mu$ with density one on $[0,1]^2$, training data $x_1 = (0.2,0.5)$ and $x_2 = (0.8, 0.5)$ and labels $y_1=0$ and $y_2 = 1$, and  $p=4$.}
\label{fig:intro}
\end{figure}

The regression problem is to find an estimator $u:\Omega_n:=\{x_i \::\: i=1,
\dots, n \} \to \R$ which agrees with preassigned labels.
To solve the regression problem  one considers objective functions which penalize the lack of smoothness of $u$ and take the structure of the graph into account. In particular here we consider the functionals which generalize the graph Laplacian, namely the graph $p$-Laplacian. 
A particular objective function we consider is
\begin{equation} \label{DF}
\cE_n^{(p)}(f) = \frac{1}{\eps_n^p n^2} \sum_{i,j=1}^n W_{ij} |f(x_i) - f(x_j)|^p. 
\end{equation}
We consider minimizing $\cE_n^{(p)}(f)$ under 
 the constraint that 
 \begin{equation} \label{CON}
f(x_i) = y_i\; \te{ for all } \; i=1, \dots, N. 
\end{equation}
A numerically computed example of the minimizer of the functional is shown  on Figure \ref{fig:intro}(a).

We investigate the asymptotic behavior in the limit when the number of unlabeled data goes to infinity,  which is consistent with semi-supervised learning paradigm of having few labeled points and an abundance of unlabeled data.
 As $n \to \infty$, $\eps(n) \to 0$ to increase the resolution and limit the computational cost. Namely as $\eps(n)$ is the length scale over  which the  information on $\mu$ is averaged, taking $\eps(n)$ to zero insures that the finer scales of $\mu$ are resolved as more data points become available.

We assume that  $\mu$ has density $\rho$ which has a positive lower bound on an open set $\Omega$ and is zero otherwise. While in this paper we consider data which are distributed in the set of full dimension,  we remark that there are no essential difficulties to extend the results to manifold setting, namely one where $\mu$ is a measure supported on compact manifold $\M$ of dimension $d$ embedded in $\R^D$. such extension has already been done for related problems concerning the graph laplacian \cite{GGHS}, where the modification of  background results (such as optimal transportation estimates) has been carried out.

The continuum limiting problem corresponds to minimizing 
\begin{equation} \label{CF}
\Ecp(f) = \sigma \int_\Omega \la \nabla f(x) \ra^p \rho^2(x) \, \dd x
\end{equation}
where  $\sigma$ is a constant that depends on $\eta$, subject to constraint that $f(x_i) = y_i$ for $i=1, \dots, N$. 
A numerically computed minimizer of the functional is shown in Figure \ref{fig:intro}(b).
Finiteness of $\Ecp(f)$  implies that $f$ lies in the Sobolev space $W^{1,p}(\Omega)$.
For the constraints to make sense it is needed that pointwise evaluation of functions is well defined, which is the case only if $p>d$, when the Sobolev embedding ensures that functions in $W^{1,p}$ are continuous. 
 When $p\leq d$  and $d >1$,  one cannot expect to be able to impose pointwise data.
Indeed spikes were in observed  in discrete models with graph-Laplacian-based regularizations (that is for $p=2$) by Nadler, Srebro, and Zhou in  \cite{NaSrZh09} who also argued that they arise since there exist functions with arbitrarily small energy $\Ecp(f)$,  for $p=2$, which agree with labels on the training set.
In \cite{elalaoui16} El Alaoui, Cheng, Ramdas,  Wainwright,  and Jordan go a step further and 
 suggest $p=d$ as the transition point between the regime where spikes appear and where solutions are ``smooth''. 
They argue, based on Sobolev embedding  theorem, that for $p\leq d$ the minimizers of  $\cE_n^{(p)}(f) $   can develop spikes as $n \to \infty$, while for $p>d$ they should not develop spikes (the authors consider $p \geq d+1$, but the same argument applies for $p>d$).
The authors also argue that for data purposes taking $p>d$ and close to $d$ is optimal since as $p \to \infty$ the solution forgets the information provided by the unlabeled points and only depends on the labeled ones. 

Our initial goal was to verify the conclusions of \cite{elalaoui16}. More precisely to show that minimizers of  $\cE_n^{(p)}(f) $, constrained to agree with the provided labels, converge, in the appropriate topology, to minimizers of $\Ecp(f)$, which also respect the labels, as $n \to \infty$ when $p>d$, and that they develop spikes when $p\leq d$. 
However we discovered an additional phenomenon, namely that the undesirable  spikes in the minimizers to graph $p$-Laplacian can occur even when $p>d$. 

Namely  \cite{elalaoui16} shows pointwise convergence of the form  
\[ \displaystyle{\lim_{\eps \to 0} \, \lim_{n\to \infty} \cE_n^{(p)}(f) = \Ecp(f)}, \]
 when $f$ is smooth enough. However considering a fixed function $f$ is not sufficient to conclude that the constrained  minimizers of $\cE_n^{(p)}$ converge to constrained minimizers of $ \Ecp$. In fact answering that question requires a set of tools from applied analysis which we discuss below. 
We show, roughly speaking, that for $d \geq 3$
 the convergence of minimizers holds if and only if 
\begin{equation} \label {eps_both_bounds}
\left( \frac{1}{n}\right)^{\frac{1}{p}} \gg \eps_n \gg \left(\frac{\log n}{n}\right)^\frac{1}{d} \quad \te{ as } n \to \infty. 
\end{equation}
The lower bound above is related to the connectivity of the graph constructed and was well understood, \cite{garciatrillos15aAAA, garciatrillos16}.
Our lower bounds for $d=1,2$ contain additional correction terms and are not optimal.
Our upper bound implies that the models  are in fact not consistent for a large family of scalings of $\eps$ on $n$ that were thus far thought to ensure consistency (namely for $1 \gg \eps_n \gg n^{-1/p}$).  Our work indicates that careful analytical approaches are needed and are in fact  capable of providing precise information on asymptotic consistency of algorithms. 



In the ``ill-posed'' regime $\eps_n^p n\to \infty$, under the usual connectivity requirement (which when $d \geq 3$ reads $\eps_n^{d} \frac{n}{\log n} \to \infty$), we are still able to establish the asymptotic behavior of algorithms.
Namely we show that minimizers of $\cE_n^{(p)}(f)$ with constraints converge, along subsequences, as $n \to \infty$ and $\eps_n \to 0$ to a minimizer of $\Ecp(f)$ without constraints. That is, the labels are forgotten in the limit as $n \to \infty$. 
This explains why, for large $n$, minimizers of $\cE_n^{(p)}$ are `spikey'.
The need to consider subsequences in the limit is due to the fact that minimizers of $\Ecp(f)$ without constraints are nonunique.

While the degeneracy of the problem when $p\leq d$ was known, \cite{elalaoui16}, we believe that degeneracy when $p>d$ and $\eps_n^p n\to \infty$ is a new and at first surprising result.
The heuristic explanation for the appearance of spikes is that the discrete $p$-Laplacian does not share the regularizing properties of the continuum $p$-Laplacian. Namely the discrete $p$-Laplacian still involves averaging over the length scale $\eps$ and thus more closely resembles an integral operator (one in \eqref{eq:Proofs:Compact:EnNL} to be precise). This allows 
high-frequency irregularities to form, without paying a high price in the energy. 
In particular, if we consider one labeled point taking the value 1, say $f_n(x_1)=1$, while $f_n(x_i) = 0$ for all $i\geq 2$ then
\[ \cE_n^{(p)}(f_n) = \frac{2}{\eps_n^p n^2} \sum_{j=2}^n \frac{1}{\eps_n^d} \eta\l\frac{|x_1-x_j|}{\eps_n}\r = \frac{2}{\eps_n^p n}\,  \eta_{\eps_n} \ast \mu_n (x_1) \to 0 \]
as $n \to \infty$, when $\eps_n^p n\to \infty$.
Note that $f_n$ exhibits degeneracy while $\cE_n^{(p)}(f_n)\to 0$. 
\medskip

In addition to the constrained problem above we also consider the problem where the agreement with the labels provided is imposed through a penalty term. Our results and analysis are analogous. 

Using the insights of our analysis, we define a new model which is quite similar to the original one, but for which the asymptotic consistency holds with only upper bound requirement being that $\eps_n \to 0$ as $n \to \infty$. 

To prove our results we use the tools of calculus of variations and optimal transportation. In particular we use the setup for convergence of objective functionals defined on graphs to their continuum limits 
developed in \cite{garciatrillos16}. This includes the definition of the proper topology ($TL^p$) to compare functionals defined on finite discrete objects (graphs) with their continuum limits. However the $TL^p$ topology, which is an extension of the $L^p$ topology, is not strong enough to ensure that the labels are preserved in the limit. For this reason we also need to consider a stronger topology, namely the one of uniform convergence. Proving the needed local regularity results for the discrete $p$-Laplacian (Lemma \ref{lem:Proofs:Compact:Osc}) and the compactness results needed to ensure the locally uniform convergence are the main technical contributions of the paper. We note that to the best of our knowledge, our results are the first where one proves (locally) uniform convergence of minimizers of nonlinear functionals in random discrete setting to the minimizers of the corresponding continuum functional. 

We note that our results on asymptotic behavior of minimizers do not provide any 
error estimates for finite $n$ and do not provide precise guidance on what $\eps$ would lead to best approximation.  In Section \ref{sec:Ex}, we numerically investigate prototypical examples in one and two dimensions to shed some light on these issues. We numerically observe the the predicted critical scalings for $\eps_n$ 
given in \eqref{eps_both_bounds}. We also numerically compare the results with our improved model~\eqref{eq:cF}. In investigating how precisely the observed error depends on $\eps$ we find that the error is smallest when $\eps$ is quite close to the connectivity radius on the graph. This is interesting and at first surprising. Rigorously explaining the  phenomenon is a in our opinion an valuable open problem. 

The paper is organized as follows.
We complete the introduction with a review on related works.
In Section~\ref{sec:MainRes} we give a precise description of the problem with the assumptions and  state the main results.
Section \ref{sec:Back}  contains a brief overview of background results we use.
This includes a description of the $TL^p$ topology, which we use for discrete-to-continuum convergence, and a short overview of $\Gamma$-convergence and optimal transportation.
Section~\ref{sec:Proofs} contains the proofs of the main results given in Section~\ref{sec:MainRes}.
In Section~\ref{sec:IM} we present an improved model that, while similar to the constrained problem
for $\cE_n^{(p)}(f)$, is asymptotically consistent with the desired limiting problem even when $\eps_n \to 0$ slowly as $n \to \infty$.
We conclude the paper with 1D and 2D numerical experiments in Section \ref{sec:Ex}.

\subsection{Discussion of Related Works \label{subsec:Intro:LitRev}}

The approach to semi-supervised learning using a weighted graph to represent the geometry of the unlabeled data and Laplacian based regularization was proposed by Zhu, Ghahramani, and Lafferty in \cite{ZhGhLa03}. It fits in the general theme of graph-Laplacian based approaches to machine learning tasks such as clustering, which are reviewed in  \cite{vonluxburg07}. See also \cite{BLSZ17} for a recent application to semi-supervised learning. 
Zhou and Sch\"{o}lkopf \cite{ZhoSch05} generalized the regularizers of  \cite{ZhGhLa03} to include a version of the graph $p$-Laplacian. 
The $p$-Laplacian regularization has also been used by B\"uhler and Hein  in clustering problems \cite{buhler09}, where values of $p$ close to $1$ are of particular interest due to connections with graph cuts. Graph based $p$-Laplacian regularization has found further applications in semi-supervised learning and image processing 
\cite{elmoataz12, elmoataz15, elmoataz17}. These papers also make the connection to the $\infty$-Laplacian, which is closely related to minimal Lipschitz extensions \cite{CrEvGa01}.

While the approach of \cite{ZhGhLa03}  has found many applications it was pointed out by 
Nadler, Srebro and Zhou \cite{NaSrZh09} that the estimator degenerates and becomes uninformative in $d \geq 2$, when the number of unlabeled data points $n \to \infty$.  
Almagir and von Luxburg \cite{alamgir11} explored the  $p$-resistances, the resulting distance on graphs, and connections to the $p$-Laplace regularization. Based on their analysis they suggested that $p=d$ should be a good choice to prevent degeneracy in the $n \to \infty$ limit.
El Alaoui, Cheng, Ramdas,  Wainwright,  and Jordan \cite{elalaoui16} show that for $p \leq d$ the 
problem degenerates as $n \to \infty$ and spikes can occur. 
They argue that regularizations with high $p \geq d+1$ are sufficient to prevent the appearance of spikes as $n \to \infty$, and lead to a well-posed problem in the limit. Here we make part of their claims rigorous, 
namely that if $p>d$ then the asymptotic consistency holds only if $\eps_n$ converges to zero sufficiently fast ($\eps_n n^p \to 0$ as $n \to \infty$). If $p>d$ and $\eps_n n^p \to \infty$ as $n \to \infty$ we prove that the problem still degenerates as $n \to \infty$ and that spikes occur.
We also introduce a modification to the discrete problem (by modifying  how the agreement with the assigned labels is imposed) which is well posed when $p>d$ without the need for $\eps_n$ to converge to $0$ quickly. 

There are other ways to regularize the SSL regression problems which ensure that no spikes occur. Namely Belkin and Niyogi \cite{BelNiy03SSL, BelNiy04SSL} consider estimators which are required to lie in the space spanned by a fixed number of eigenvectors of the graph Laplacian. Due to the smoothness of low eigenvectors of the Laplacian this prevents the formation of spikes.
 One can think of this approach in energy based setting where infinite penalty has been imposed on high frequencies. A softer, but still linear, way to do this is to consider (fractional) powers of the graph Laplacian, namely the regularity term 
$J_n(u) = \langle cL_n^\alpha f,f \rangle$ where $L_n$ is the graph Laplacian, and $\alpha>0$. 
This regularization was studied by Belkin and Zhou \cite{zhou11} who argue, again via regularity obtained by Sobolev embedding theorems, that taking $\alpha>\frac{d}{2}$ prevents spikes. However Dunlop, Stuart, and the authors have discovered that a similar phenomenon to one described in this paper. Namely 
even when $\alpha>\frac{d}{2}$ the limit may be degenerate, and spikes can occur, if $\eps_n$ converges to zero slowly, namely if $\eps_n^{2\alpha} n \to \infty$ as $n \to \infty$. 
\medskip

Our results fall in the class of asymptotic consistency results in machine learning. In general one is interested  the asymptotic behavior of an objective function posed on a random sample of $n$ points, and which also depends on a parameter $\eps$, $E_{n, \eps}(f_n)$ where $f_n$ is a real valued function defined at sample points. The limit is considered as $n \to \infty$ while $\eps_n \to 0$ at appropriate rate. The limiting problem is described by a continuum functional $E_\infty(f)$ which acts on real valued functions supported on domains or manifolds. Also relevant is the (nonlocal) continuum problem, $E_{\infty, \eps}(f)$ which describes the limit $n \to \infty$ while $\eps>0$ is kept fixed. 

The type of consistency that is needed for the conclusions, and the one we consider, is \emph{variational consistency}, namely that minimizers of  $E_{n, \eps_n}(f_n)$ converge to minimizers of $E_\infty(f)$ as 
$n \to \infty$ while $\eps_n \to 0$ at an appropriate rate. Proving such results includes choosing the right topology to compare the functions on discrete domain $f_n$ with those on the continuum domain $f$.

Many works in the literature are interested in a simpler notion of convergence, namely that for a fixed, sufficiently smooth, continuum function $f$ it holds that $E_{n, \eps_n}(f) \to E_\infty(f)$ as  $n \to \infty$ while $\eps_n \to 0$ at an appropriate rate, where by $E_{n, \eps_n}(f)$ we mean that the discrete functional is evaluated at the restriction of $f$ to the data points. We call this notion of convergence \emph{pointwise convergence}. A somewhat weaker notion of convergence is  what we here call \emph{iterated pointwise convergence}, namely considering $\lim_{\eps \to 0} \lim_{n \to\infty} E_{n,\eps}(f)$. Also relevant for the problems based on linear operators (namely the graph Laplacian) is \emph{spectral convergence} which asks for the eigenvalues and eigenvectors of the discrete operator to converge to eigenvalues and eigenfunction of the continuum one. This notion of the convergence is typically sufficient for the kind of conclusions we are investigating (however our problems are nonlinear).

 Pointwise (and similar notions of) convergence of graph Laplacians was studied by Belkin and Niyogi \cite{belkin07},   
Coifman and Lafon \cite{coifman1},
Gin\'e and Koltchinskii \cite{GK}, Hein, Audibert and von Luxburg \cite{hein_audi_vlux05}, 
Hein \cite{hein06},
Singer \cite{singer06}, and Ting, Huang, and  Jordan \cite{THJ}.
Spectral convergence was studied in the works of  Belkin and Niyogi \cite{belkin07} on the convergence of Laplacian eigenmaps, von Luxburg, Belkin,  and Bousquet \cite{vonluxburg08} and Pelletier and Pudlo \cite{PelPud11} on graph Laplacians, and of  Singer and Wu \cite{SinWu13} on the connection graph Laplacian. In these works on spectral convergence either $\eps$ remains fixed as $n \to \infty$ or $\eps(n) \to 0$ at an unspecified rate. The precise and almost optimal rates were obtained in \cite{garciatrillos15aAAA} using variational methods. 
Further problems involve obtaining error estimates between discrete and continuum objects. 
Laplacians on discretized manifolds was studied by Burago, Ivanov and Kurylev  \cite{BIK} who obtain precise error estimates for eigenvalues and eigenvectors. Related results on approximating elliptic equations on point clouds have been obtained by Li and Shi \cite{LiShi15}, and  Li, Shi, Sun, \cite{LiShiSun17}. Error bounds  for the spectral convergence of graph Laplacians have been considered by Wang \cite{wang14AAA}  and 
Garc\'ia Trillos, Gerlach, Hein and  one of the authors \cite{GGHS}.
Regarding graph $p$-Laplacians, the authors of \cite{elalaoui16} obtain iterated pointwise convergence of graph $p$-Laplacians to the continuum $p$-Laplacian. Finally we mention that for a different type of problems, namely for  nondominated sorting, of Calder, Esedo\=glu, and Hero \cite{CEH14} have obtained uniform convergence of discrete solutions to the solution of a continuum Hamilton-Jacobi equation.

%

To obtain the results on variational convergence of $\cE_n^{(p)}$ to $\Ecp$ needed to fully explain the asymptotics of discrete regression problems we combine tools of calculus of variations (in particular $\Gamma$-convergence)  and optimal transportation. This approach to asymptotics of problems posed on discrete random samples was developed by Garc\'ia-Trillos and one of the authors \cite{garciatrillos16,garciatrillos15aAAA}.
In \cite{garciatrillos16} they introduce the $TL^p$ topology for comparing the functions defined on the discrete sets to the ones defined in the continuum, and apply the approach to asymptotics of graph-cut based objective functions. We refer to this paper for a description of the rich background of the works that underpin the approach.
In \cite{garciatrillos15aAAA} the authors apply the approach to convergence of graph Laplacian based functionals.
Consistency of $k$-means clustering for paths with regularization was recently studied by Theil, Johansen and Cade, and one of the authors \cite{thorpe15}, using a similar viewpoint.
This technical setup has recently been used and extended to studies on modularity based clustering 
\cite{davis16AAA}, 
 Cheeger and ratio cuts \cite{garciatrillos15c}, 
neighborhood graph constructions for graph cut based clustering \cite{trillos16knn}, and 
classification problems \cite{thorpe17AAA, garciatrillos16aAAA}.

An alternative approach to related regression problems was developed by Fefferman and collaborators, Israel, Klartag and Luli, who look  for a function of sufficient regularity,  that extends a function $f^\dagger:E\to \bbR$ to the whole of $\bbR^d$ in such a way as to minimize the norm of the extension.
They show that appropriate  extensions exist and finding efficient constructions for $f$, for $C^m$ regularity~\cite{fefferman09,fefferman09a,fefferman09b}, and for Sobolev regularity~\cite{fefferman16,fefferman16a,fefferman16b}.
In the context of machine learning this is a supervised learning problem and
 only makes use of the labeled data. In our context the problem is independent of $\{x_i\}_{i=N+1}^n$ and does not use the geometry of the unlabeled data.

\section{Setting and Main Results \label{sec:MainRes}}

Let $\Omega$ be an open, bounded domain in $\R^d$.
Let $\{ (x_i, y_i) \::\: i=1,\dots,N\}$  with $x_i \in \Omega$ and $y_i \in \R$ be a collection of distinct labeled points. Throughout the paper we consider $N$ to be fixed. Considering a model where $N$ grows is an interesting problem, which we do not address here. 
We consider $\mu$ to be the measure representing the distribution of data.  We assume that $\supp \mu =\overline \Omega$ and that $\mu$ has density $\rho$ with respect to Lebesgue measure. We assume that $\rho$ is continuous and is bounded above and below by positive constants on $\Omega$. 

We assume that unlabeled data, $\{x_i\}_{i=N+1, \dots}$ are given  by a sequence of iid samples of measure $\mu$. The empirical measure induced by data points is given by $\mu_n=\frac{1}{n} \sum_{i=1}^n \delta_{x_i}$.
Let $G_n=(\Omega_n,E_n,W_n)$ be a graph with vertices $\Omega_n = \{x_i \::\: i=1, \dots,n\}$, edges $E_n = \{e_{ij}\}_{i,j=1}^n$ and edge weights $W_n=\{W_{ij}\}_{i,j=1}^n$.
For notational simplicity we will set $W_{ij}=0$ if there is no edge between $x_i$ and $x_j$.

We assume the following structure on edge weights
\begin{equation} \label{eq:MainRes:Wij}
W_{ij} = \eta_\eps(|x_i-x_j|)
\end{equation}
where $\eta_\eps(|x|) = \frac{1}{\eps^d} \eta\l \frac{|x|}{\eps} \r$, $\eta:[0,\infty)\to [0,\infty)$ is a nonincreasing kernel  and $\eps=\eps_n$ is a scaling parameter depending on $n$.
For example if $\eta(|x|) = \bbI_{|x|\leq 1}$ then $\eta_\eps(|x|)$ is $\frac{1}{\eps^d}$ if $|x|\leq \eps$ and 0 otherwise.
In this case vertices are only connected if they are closer than $\eps$.

We consider two models: one where the agreement of the response with the training variables is imposed as a constraint and the other where it is imposed via a penalty. We call these models \emph{constrained} and \emph{penalized}.

In the constrained model we construct our estimator as the minimizer of
\begin{equation} \label{eq:MainRes:En}
\cE_n^{(p)}(f) = \frac{1}{\eps_n^p} \frac{1}{n^2} \sum_{i,j=1}^n W_{ij} |f(x_i) - f(x_j) |^p
\end{equation}
among  $\{f : \Omega_n \to \R\}$ which satisfy the constraint $f(x_i) = y_i$ for all $i=1, \dots, N$.

For technical reasons it is convenient to define the functional over all $f$ and impose the constraint in the following way
\begin{equation} \label{eq:MainRes:Encon}
\cE_{n,con}^{(p)}(f) = 
\begin{cases}
\frac{1}{\eps_n^p} \frac{1}{n^2} \sum_{i,j=1}^n W_{ij} |f(x_i) - f(x_j) |^p & \te{if }  f(x_i) = y_i \, \te{ for } i=1,2,\dots, N \\
\infty & \te{else.}
\end{cases}
\end{equation}

We now turn to the penalized formulation. For $q > 0$ let
\[ R^{(q)}(f) = \sum_{i=1}^N |y_i-f(x_i)|^q. \]
We define the estimator as the minimizer of
\begin{equation} \label{eq:Spn}
 \cS_n^{(p)}(f) = \cE_n^{(p)}(f) + \lambda  R^{(q)}(f) 
\end{equation}
where $\lambda>0$ is a tunable parameter.
\medskip

We now introduce the continuum functionals that describe the limiting problems as $n \to \infty$.
Let 
\begin{equation} \label{eq:MainRes:Einfty}
\Ecp(f) = \begin{cases}
\sigma_\eta \int_\Omega \la \nabla f(x) \ra^p \rho^2(x) \, \dd x \;\;\, & \te{if } f \in W^{1,p}(\Omega), \\
\infty & \te{else}.
\end{cases}
\end{equation}
For $p>d$, Sobolev functions $f \in W^{1,p}$ are continuous and we can define
\begin{equation} \label{eq:MainRes:Einftycon}
\Ecpc(f) = \begin{cases}
\Ecp(f) & \te{if } f \in W^{1,p}(\Omega) \te{ and }
f(x_i) = y_i \te{ for } i = 1, \dots, N \\
\infty & \te{else}.
\end{cases}
\end{equation}
The constant  $\sigma_\eta$ above is defined, using   $e_1 = [1,0, \dots, 0]^T$, by 
\[ \sigma_\eta = \int_{\bbR^d} \eta(|x|) \, |x \cdot e_1|^p \; \dd x. \]

To describe the limit of the  penalized  model the large data limit we introduce
\begin{equation} \label{eq:Spninf}
 \cS_\infty^{(p)}(f) =  \Ecp(f) + \lambda  R^{(q)}(f).
\end{equation}

We note that functionals~\eqref{eq:MainRes:Einftycon} and~\eqref{eq:Spninf} are lower semi-continuous with respect to the $L^p$ norm.
In addition, coercivity of both functionals follows from Sobolev embeddings.
Coercivity and lower semi-continuity imply existence of minimizers, e.g.~\cite[Theorem 3.6]{fonseca07}.
Strict convexity implies that the minimizers are unique.

We are interested in asymptotic behavior of minimizers $f_n$ of the discrete models, say $\cE_{n,con}^{(p)}$. 
We say that $\cE_{n,con}^{(p)}$ is \emph{asymptotically consistent} with $\Ecpc$ if the minimizers $f_n$ of 
$\cE_{n,con}^{(p)}$ converge as $n \to \infty$ to a minimizer of $\Ecpc$. 
One should note topology of  the convergence $f_n\to f_\infty$ is not at this stage clear.

We observe that since $f_n : \Omega_n \to \R$, while $f:\Omega \to \R$ this issue is nontrivial. 
We use the $TL^p$ topology introduced in \cite{garciatrillos16}  precisely  to compare functions defined on different domains in a topology consistent with $L^p$ convergence.
We define the convergence rigorously in Section~\ref{sec:Back}. 

Another issue is the rate at which $\eps_n$ is allowed to converge to zero.
If $\eps_n \to 0$ too quickly then the graph becomes disconnected and hence it does not capture the geometry of $\Omega$ properly. The connectivity threshold \cite{penrose03} is $\eps_n \sim \left(\frac{\log n}{ n}\right)^{\frac{1}{d}} $. We require (when $d \geq 3$) $\eps_n \gg \left(\frac{\log n}{ n}\right)^{\frac{1}{d}} $ which means that our lower bound  is almost optimal. 
We discovered that if $\eps_n \to 0$ too slowly the discrete functional $\cE_{n,con}^{(p)}$ lacks sufficient regularity for the constraints to be preserved in the limit. The optimal upper bound on $\eps_n$ is discussed in Theorem \ref{thm:MainRes:Conv}.

\medskip

We now state our assumptions needed for  the main results.
\begin{listi}
\item $\Omega\subset\bbR^d$ is open, connected, bounded and with Lipschitz boundary;
\item The probability measure $\mu\in \cP(\Omega)$ has continuous density $\rho$ which is bounded above and below by strictly positive constants in $\Omega$;
\item There exists $N$ labeled points: $(x_i, y_i) \in \Omega \times \R$ for $i=1, \dots, N$; 
\item For $i >N$ the data points $x_i$, are iid samples of $ \mu$;
\item  Let $\eps_n$  be a sequence converging to $0$ satisfying the lower bound 
\[ \eps_n \gg 
\begin{cases}
\sqrt{\frac{\log\log n}{n}} & \text{if } d = 1 \medskip \\
\frac{(\log n)^{\frac{3}{4}}}{\sqrt{n}} & \text{if } d=2 \medskip \\ 
\left(\frac{\log n}{n}\right)^{\frac{1}{d}} & \text{if } d\geq 3;
\end{cases} \]
\item The kernel profile $\eta:[0, \infty) \to [0, \infty)$ is non-increasing;
\item $\eta$ is positive and continuous at $x=0$; 
\item The integral $\int_0^\infty \eta(t) |t|^{p+d} \, \dd t$ is finite (equivalently $\sigma_\eta=\int_{\bbR^d} \eta(|w|) |w\cdot e_1|^p \, \dd w<\infty$).
\end{listi}

The first main result of the paper is the following theorem. Its proof is presented in Section~\ref{sec:Proofs}.
\begin{theorem}[Consistency of the constrained model]
\label{thm:MainRes:Conv}
Let $p > 1$. Assume $\Omega$, $\mu$, $\eta$, and $x_i$ satisfy the assumptions \textbf{(A1) - (A8)}. Let graph weights $W_{ij}$ be given by~\eqref{eq:MainRes:Wij}. 
Let $f_n$ be a sequence of minimizers of $\cE_{n,con}^{(p)}$ defined in \eqref{eq:MainRes:Einftycon}. 
Then, almost surely, the sequence  $(\mu_n,f_n)$ is precompact in the $TL^p$ metric. The $TL^p$ limit of any convergent subsequence, $(\mu_{n_m}, f_{n_m})$, is of the form $(\mu,f)$ where $f \in W^{1,p}(\Omega)$. Furthermore,
\begin{itemize}
\item[(i)] if $\,n \eps_n^p \to 0$ as $n \to \infty$ then $f$ is continuous and 
\begin{itemize}
\item[(a)] $f_{n_m}$ converges locally uniformly to $f$, meaning that for any $\Omega' \subset\subset \Omega$
\[  \lim_{m \to \infty} \max_{\{k \leq n_m \::\: x_k \in \Omega'\}} |f(x_k) - f_{n_m}(x_k)| =0, \]
\item[(b)] $f$ is a minimizer of $\Ecpc$ defined in \eqref{eq:MainRes:Einftycon},
\item[(c)] the whole sequence $f_n$ converges to $f$ both in $TL^p$ and locally uniformly;
\end{itemize}
\item[(ii)] if $\,n \eps_n^p \to \infty$ as $n \to \infty$ then $f$ is a minimizer of $\Ecp$ defined in \eqref{eq:MainRes:Einfty}.
\end{itemize}
\end{theorem}
We note that  in case (i) assumption \textbf{(A5)} and $\,n \eps_n^p \to 0$ as $n \to \infty$ imply that 
$n^{-1/p} \gg \eps \gg n^{-1/d}$ which is only possible if $p>d$. Therefore in case (i) we always have that functions $f$ for which $\Ecp$ is finite are always continuous and thus it is possible to impose  pointwise values of $f$, as needed to define $\Ecpc$ in \eqref{eq:MainRes:Einftycon}.

The result (i) establishes the asymptotic consistency of the discrete constrained model with the 
constrained continuum weighted $p$-Laplacian model.

While the result (ii) looks similar its interpretation is different. It shows that the model ``forgets'' the constraints in the limit. Namely $\Ecp$ only has the gradient  term and no constraints! In particular its minimizers are constants over $\Omega$. What is happening is that $f_n$ develops narrow spikes near 
the labeled points $x_i$ and becomes nearly constant everywhere else. In the $TL^p$ limit the spikes disappear. 

This motivates referring to the scaling when $n^p \eps \to \infty$ as $n \to \infty$ as the 
 \emph{\textbf{degenerate}} regime.
On the other hand, we refer to the scaling of case (i)  as the \emph{\textbf{well-posed}} regime.

\medskip

The other main result is the convergence in the penalized model.
The proof is a straightforward extension of Theorem~\ref{thm:MainRes:Conv} in the special case $N=0$ (so that the constraint is not present). 
We include the proof in Section~\ref{subsec:MainRes:ConvSoft}.

\begin{proposition}
\label{prop:MainRes:ConvSoft}
Let $p > 1$. Assume $\Omega$, $\mu$, $\eta$, and $x_i$ satisfy the assumptions~\textbf{(A1)-(A8)}. Let graph weights $W_{ij}$ be given by~\eqref{eq:MainRes:Wij}. 
Let $f_n$ be a sequence of minimizers of $ \cS_n^{(p)}$ defined in \eqref{eq:Spn}. 
Then, almost surely, the sequence  $(\mu_n,f_n)$ is precompact in the $TL^p$ metric.
The $TL^p$ limit of any convergent subsequence, $(\mu_{n_m},f_{n_m})$, is of the form $(\mu,f)$ where $f \in W^{1,p}(\Omega)$. Furthermore,
\begin{itemize}
\item[(i)] if $\,n \eps_n^p \to 0$ as $n \to \infty$ then $f$ is continuous and 
\begin{itemize}
\item[(a)] $f_n$ converges locally uniformly to $f$, meaning that for any $\Omega' \subset\subset \Omega$
\[  \lim_{n \to \infty} \max_{\{k \leq n \::\: x_k \in \Omega'\}} |f(x_k) - f_n(x_k)| =0, \]
\item[(b)] $f$ is a minimizer of $ \cS_\infty^{(p)}$ defined in \eqref{eq:MainRes:Einftycon},
\item[(c)] the whole sequence $f_n$ converges to $f$ both in $TL^p$ and locally uniformly;
\end{itemize}
\item[(ii)] if $\,n \eps_n^p \to \infty$ as $n \to \infty$ then $f$ is a minimizer of $\Ecp$ defined in \eqref{eq:MainRes:Einfty}.
\end{itemize}
\end{proposition}

Again the result of (i) is a consistency result, while (ii) shows that the penalization of the labels is lost in the limit. 

\begin{remark}
\label{rem:MainRes:p=1}
The above results (Theorem~\ref{thm:MainRes:Conv} and Proposition~\ref{prop:MainRes:ConvSoft}) could also be extended to $p=1$, in which case the limiting functional $\cE^{(1)}_\infty$ would be a weighted $TV$ semi-norm $\cE_\infty^{(1)} = \sigma_\eta TV(\cdot;\rho)$ where
\[ TV(f;\rho) = \sup\left\{ \int_\Omega f\mathrm{div} \phi \, \mathrm{d} x \, : \, |\phi(x)| \leq \rho^2(x) \, \forall x\in \Omega, \phi\in C_c^\infty(\Omega;\mathbb{R}^d) \right\}. \]
A modification of the proofs contained here would prove the result, see also~\cite{garciatrillos16}.
\end{remark}

We recall that in Section~\ref{sec:IM} we propose an improved model that is well-posed when $p>d$ without requiring that $n \eps^p \to 0$.

\medskip

\section{Background Material \label{sec:Back}}


 
In an effort to make this paper more self-contained we briefly recall three key notions our work relies on. 
The first is $\Gamma$-convergence which is a notion of convergence of functionals developed for the analysis of sequences of variational problems. The second is the notion of optimal transportation, and  the third is the $TL^p$ space which we use to define the convergence of discrete functions to continuum functions.

\subsection{\texorpdfstring{$\Gamma$}{Gamma}--Convergence \label{subsec:Back:Gamma}}


$\Gamma$-convergence was introduced by De Giorgi in 1970's to study limits of variational problems.  We refer to~\cite{braides02,dalmaso93} for an in depth introduction to $\Gamma$-convergence.
Our application of $\Gamma$-convergence will be in a random setting.

\begin{mydef}[$\Gamma$-convergence]
\label{def:Back:Gamma:GamCon}
Let $(Z,d)$ be a metric space and $(\cX,\bbP)$ be a probability space.
For each $\omega\in \cX$ the functional $E_n^{(\omega)} :Z\to \bbR\cup\{\pm\infty\}$ is a random variable.
We say $E_n^{(\omega)}$ \textit{$\Gamma$-converge} almost surely on the domain $Z$ to $E_\infty :Z\to \bbR\cup\{\pm\infty\}$ with respect to $d$, and write $E_\infty = \Glim_{n \to \infty} E_n^{(\omega)}$, if there exists a set $\cX'\subset \cX$ with $\bbP(\cX') = 1$, such that for all $\omega\in \cX'$ and all $f\in Z$:
\begin{itemize}
\item[(i)] (liminf inequality) for every sequence $\{f_n\}_{n=1}^\infty$ converging to $f$
\[ E_\infty(f) \leq \liminf_{n\to \infty} E_n^{(\omega)}(f_n), \te{ and } \]
\item[(ii)] (recovery sequence) there exists a sequence $\{f_n\}_{n=1}^\infty$ converging to $f$ such that
\[ E_\infty(f) \geq \limsup_{n\to \infty} E_n^{(\omega)}(f_n). \]
\end{itemize}
\end{mydef}

For ease of notation we will suppress the dependence of $\omega$ on on our functionals, that is we apply the above definition to $E_n = \cE_n^{(p)}$. 
The almost sure statement in the above definition does not play a significant role in the proofs.
Basically it is enough to consider the set of realizations of $\{x_i\}_{i=1}^\infty$ such that the empirical measure converges $\text{weak}^*$. More precisely, we consider the set of realizations of $\{x_i\}_{i=1}^\infty$ such that the conclusions of Theorem~\ref{thm:Back:TransMapBound} hold. 
 
The fundamental result concerning $\Gamma$-convergence is the following convergence of minimizers result.
The proof can be found in~\cite[Theorem 1.21]{braides02} or~\cite[Theorem 7.23]{dalmaso93}.
 
\begin{theorem}[Convergence of Minimizers]
\label{thm:Back:Gamma:Conmin}
Let $(Z,d)$ be a metric space and $E_n: Z\to [0,\infty]$ be a sequence of functionals.
Let $f_n$ be a minimizing sequence for $E_n$.
If the set $\{f_n\}_{n=1}^\infty$ is precompact and $E_\infty = \Glim_n E_n$ where $E_\infty:Z\to[0,\infty]$ is not identically $+\infty$ then
\[ \min_Z E_\infty = \lim_{n\to \infty} \inf_Z E_n. \]
Furthermore any cluster point of $\{f_n\}_{n=1}^\infty$ is a minimizer of  $E_\infty$.
\end{theorem}
The theorem is also true if we replace minimizers with almost minimizers.
\medskip

We note that $\Gamma$-convergence is defined for functionals on a common metric space.
The next section overviews the metric space we use to analyze the asymptotics of our semi-supervised learning models, in particular it allows us to go from discrete to continuum.



\subsection{Optimal Transportation and Approximation of Measures} \label{sec:ot}

Here we recall the notion of optimal transportation between measures and the metric it introduces. Comprehensive treatment of the topic can be found in books of Villani \cite{villani09} and Santambrogio \cite{santambrogio}.

Given $\Omega$ is open and bounded, and probability measures $\mu$ and $\nu$ in $\cP(\overline \Omega)$ we define the set   $\Pi(\mu,\nu) $  of transportation maps, or couplings, between them to be  the set of probability measures on the product space $\pi \in \cP(\overline \Omega \times \overline \Omega)$ whose first marginal is $\mu$ and second marginal is $\nu$. 
We then define the $p$-optimal transportation distance (a.k.a. $p$-Wasserstein distance)
by
\[ d_p(\mu, \nu) = 
\begin{cases}
\displaystyle{\inf_{\pi\in\Pi(\mu,\nu) } \l \int_{\Omega\times \Omega} |x-y|^p \, \dd \pi(x,y)  \r^\frac{1}{p}}  \quad & \te{if } 1 \leq p < \infty \\
\displaystyle{ \inf_{\pi\in\Pi(\mu,\nu)}  \pi\text{-}\esssup_{(x,y)} |x-y|  } & \te{if } p=\infty.
\end{cases} \]
If $\mu$ has a density with respect to Lebesgue measure on $\Omega$, then the distance can be rewritten using transportation maps, $T:  \Omega \to \Omega$, instead of transportation plans,
\[ d_p(\mu, \nu) = 
\begin{cases}
\displaystyle{\inf_{\pi\in\Pi(\mu,\nu) } \l  \int_\Omega |x-T(x)|^p \, \dd \mu(x)    \r^\frac{1}{p}}  \quad & \te{if } 1 \leq p < \infty \\
\displaystyle{ \inf_{T_\sharp \mu = \nu}  \mu\text{-}\esssup_{x} |x-T(x)|  } & \te{if } p=\infty.
\end{cases} \]
where $T_\mu=\nu$ means that the push forward of measure $\mu$ by $T$ is measure $\nu$, namely that $T$ is Borel measurable  and such that for all $U \subset  \overline \Omega$, open, $\mu(T^{-1}(U)) = \nu(U)$. 

When $p< \infty$ the metric $d_p$ metrizes the weak convergence of measures.

Optimal transportation  plays an important role in comparing the discrete and continuum objects we study. In particular we use sharp estimates on the $\infty$-optimal transportation 
distance between a measure and the empirical measure of its sample. In the form below, for $d \geq 2$, they were established in \cite{garciatrillos15}, which extended the related results in  \cite{AKT, LeightonShor, ShorYukich, Talagrand}. For $d=1$ the estimates are simpler, and follow from the law of iterated logarithms. 

\begin{theorem}
\label{thm:Back:TransMapBound}
Let $\Omega\subset\bbR^d$ be open, connected and bounded with Lipschitz boundary.
Let $\mu$ be a probability measure on $\Omega$ with density (with respect to Lebesgue) $\rho$ which is bounded above and below by positive constants.
Let $x_1,x_2,\dots$ be a sequence of independent random variables with distribution $\mu$ and let $\mu_n$ be the empirical measure.
Then there exists a constants $C \geq c > 0$ such that almost surely there exists a sequence of transportation maps $\{T_n\}_{n=1}^\infty$ from $\mu$ to $\mu_n$ such that
\[ c \leq  \liminf_{n\to \infty} \frac{\|T_n-\Id\|_{L^\infty(\Omega)}}{\delta_n}  \leq \limsup_{n\to \infty} \frac{\|T_n-\Id\|_{L^\infty(\Omega)}}{\delta_n} \leq C \]
where
\[ \delta_n = \lb \begin{array}{ll} \sqrt{\frac{\log \log(n)}{n}} & \text{if } d=1 \\ \frac{(\log n)^{\frac{3}{4}}}{\sqrt{n}} & \text{if } d=2 \\ \frac{(\log n)^{\frac{1}{d}}}{n^{\frac{1}{d}}} & \text{if } d\geq 3. \end{array} \rd \]
\end{theorem}

\subsection{The \texorpdfstring{$TL^p$}{TLp} Space \label{subsec:Back:TLp}}



The discrete functionals we consider (e.g. $\cE_n^{(p)}$)  are defined for functions $f_n : \Omega_n \to \R$ where  $\Omega_n= \{x_i \::\: i=1, \dots, n\}$, while the limit functional 
$\Ecp$ acts on functions $f:\Omega \to \R$, where $\Omega$ is an open set. We can view 
$f_n$ as elements of $L^p(\mu_n)$ where $\mu_n$ is the empirical measure of the sample $\mu_n = \frac{1}{n} \sum_{i=1}^n \delta_{x_i}$. Likewise $f \in L^p(\mu)$ where $\mu$ is the measure with density $\rho$ out of which the points are sampled from. One would like how to compare $f$ and $f_n$ in a way that is consistent with $L^p$ topology. 
To do so we use the  $TL^p$ space was introduced in~\cite{garciatrillos16}, where it was used to study the continuum limit of the graph total variation (that is $\cE_n^{(1)}$).
Subsequent development of the $TL^p$ space has been carried out in~\cite{garciatrillos15aAAA, thorpe17bAAA, thorpe17cAAA}.

To compare the functions $f_n$ and $f$ above we need to take into account their domains, 
 or more precisely to account for  $\mu$ and $\mu_n$. 
For that purpose the space of configurations  is defined to be 
\[ TL^p(\Omega) = \lb(\mu,f) \, : \, \mu \in \cP(\overline \Omega), f \in L^p(\mu) \rb. \]
The metric on the space is 
\[ d_{TL^p}^p((\mu,f),(\nu,g)) = \inf \lb \int_{\Omega\times \Omega} |x-y|^p + |f(x) - g(y)|^p \, \dd \pi(x,y) \, : \, \pi\in\Pi(\mu,\nu) \rb \]
where $\Pi(\mu,\nu) $ the set of transportation plans defined in Section \ref{sec:ot}.
 We note that the minimizing $\pi$ exists and that  $TL^p$ space is a metric space,\cite{garciatrillos16}.

When $\mu$ has a density with respect to Lebesgue measure on $\Omega$, then the distance can be rewritten using transportation maps, $T$ instead of transportation plans,
\[ d_{TL^p}^p((\mu,f),(\nu,g)) = \inf \lb \int_\Omega |x-T(x)|^p + |f(x) - g(T(x))|^p \, \dd \mu(x) \, : \, T_{\#} \mu = \nu \rb. \]
This formula provides a clear interpretation of the distance in our setting. Namely 
to compare functions $f_n: \Omega_n\to \bbR$ we define a mapping $T_n:\Omega \to \Omega_n$ and compare the functions $\tilde{f}_n = f_n\circ T_n$ and $f$ in $L^p(\mu)$, while also accounting for the transport, namely the $|x - T_n(x)|^p$ term.

We remark that $TL^p(\overline \Omega)$ space is not complete and that its completion was discussed in \cite{garciatrillos16}. In the setting of this paper, since the corresponding measure is clear from context, we often say that $f_n$ converges in $TL^p$ to $f$ as a short way to say that $(\mu_n, f_n)$ converges in $TL^p$ to $(\mu,f)$.

\section{Regularity and Asymptotics of Discrete and Nonlocal Functionals \label{sec:Proofs}}

Here we present some of the key properties of the functionals involved that allow us to show the asymptotic consistency of Theorem \ref{thm:MainRes:Conv}. A fundamental new issue (compared to say \cite{garciatrillos15aAAA}) is that constraints in $\Ecp$ are imposed pointwise on a set of $\mu$ measure zero. [The reason that these constraints make sense is that for $p>d$ the finiteness of $\Ecp(f)$ implies that $f$ is continuous.] We note that the $TL^p$ convergence used in \cite{garciatrillos15aAAA} is not sufficient to imply that constraints are preserved. One needs a stronger convergence, like the uniform one. This raises the question on how to obtain the needed compactness of sequences $f_n$, that is how to show that uniform boundedness of $\cE_{n,  con}^{(p)}(f_n)$ implies the existence of a (locally) uniformly converging subsequence. Our approach combines discrete and continuum regularity results. Namely we obtain in  Lemma \ref{lem:Proofs:Compact:Osc} a local control of oscillation of $f_n$ over distances of order $\eps_n$. In Lemma \ref{lem:EntoNLBound} we show that discrete functionals $\cE_{n,  con}^{(p)}(f_n)$ control the values of the  associated nonlocal continuum  functionals $\cE_{\eps_n}^{(NL,p)}(\tilde f_n)$ (defined in~\eqref{eq:Proofs:Compact:EnNL} below) applied to an appropriate extrapolation $\tilde{f}_n$ of $f_n$. A simple but important point is that the discrete functionals at fixed $n$ are always closer to a nonlocal functional with nonlocality at scale $\eps_n$, than to the limiting functional. 
The issue is that these nonlocal functionals do not share the regularizing properties of the limiting functional. However we show in Lemma \ref{lem:Proofs:Compact:NLtoLBound} that control of the nonlocal energy is sufficient to provide regularity at scales larger than $\eps$. Combining these estimates is enough to imply the compactness with respect to (locally) uniform convergence, Lemma \ref{lem:Proofs:Compact:UnifConverg}.

\begin{lemma}[discrete regularity]
\label{lem:Proofs:Compact:Osc}
Let $p > 1$. Assume $\Omega$, $\mu$, $\eta$, and $x_i$ satisfy the assumptions \textbf{(A1) - (A8)}. Let graph weights $W_{ij}$ be given by~\eqref{eq:MainRes:Wij}. 
Let $\Omega_n=\{x_i\}_{i=1}^n$. 
For any $f_n : \Omega_n \to \R$, we define $\osc_\eps^{(n)}(f_n): \Omega_n \to \R$ by
\[ \osc_\eps^{(n)}(f_n)(x_i) =  \max_{z\in B(x_i,\eps)\cap \Omega_n} f_n(z) - \min_{z\in B(x_i,\eps)\cap \Omega_n} f_n(z). \]
For any $\alpha_0>0$, with probability one, there exist $n_0>0$ and $C>0$ (independent of $n$) such that for any $\alpha \geq \alpha_0$, all $n \geq n_0$, all $f_n : \Omega_n \to \R$, and all  $k\in \{1,2,\dots, n\}$
\[ \l \osc_{\alpha \eps_n}^{(n)}(f_n)(x_k) \r^p \leq C\alpha^p  n \eps^p_n \cE_n^{(p)}(f_n), \]
where  $\cE_n^{(p)}$ is defined by~\eqref{eq:MainRes:En}. 
\end{lemma}

\begin{proof}
Let $\tilde{\eta}(t) = a$ if $0\leq t<b$ and $\tilde{\eta}(t) = 0$  where $a$ and $b$ are chosen such that $\tilde{\eta}\leq \eta$.
We can furthermore choose $b$ so that $b \leq \alpha_0$.
For all $k\in \{1,\dots, n\}$  let
\begin{align*}
\bar{f}_n(x_k) & = \max_{z\in B(x_k,\frac{b\eps_n}{2})\cap\Omega_n} f_n(z), 
& \bar{x}_k & \in \argmax_{z\in B(x_k,\frac{b\eps_n}{2})\cap\Omega_n} f_n(z), \\
\underline{f}_n(x_k) & = \min_{z\in B(x_k,\frac{b\eps_n}{2})\cap\Omega_n} f_n(z), 
& \underline{x}_k & \in \argmin_{z\in B(x_k,\frac{b\eps_n}{2})\cap\Omega_n} f_n(z).
\end{align*}
Note that $\osc_{\frac{b\eps_n}{2}}^{(n)}(f_n)(x_k) = \bar{f}_n(x_k) - \underline{f}_n(x_k) $ and for all $x\in B\l x_k,\frac{b\eps_n}{2}\r \cap \Omega_n$ 
\begin{align*}
\text{(i)} & \quad \quad  \bar{f}_n(x_k)  - f_n(x)  \geq \frac12 \osc_{\frac{b\eps_n}{2}}(f_n)(x_k), \\
\text{or} \quad \text{(ii)} & \quad \quad f_n(x) - \underline{f}_n(x_k) \geq \frac12 \osc_{\frac{b\eps_n}{2}}(f_n)(x_k).
\end{align*}
Without a loss of generality we assume that (i) holds for at least half the points in $B\l x_k,\frac{b\eps_n}{2}\r\cap \Omega_n$.
Then,
\begin{align}
\begin{split}
\cE_n^{(p)}(f_n) & \geq \frac{1}{\eps^{p+d}_n n^2} \sum_{i,j=1}^n \tilde{\eta}\l \frac{|x_i-x_j|}{\eps_n}\r |f_n(x_i) - f_n(x_j)|^p \\
 & \geq \frac{1}{\eps^{p+d}_n n^2} \sum_{j: |x_j-\bar{x}_k|\leq b\eps_n} \tilde{\eta}\l \frac{|\bar{x}_k-x_j|}{\eps_n}\r |f_n(x_j) - f_n(\bar{x}_k)|^p \\
 & \geq \frac{a}{\eps^{p+d}_n n^2} \sum_{j: |x_j-x_k|\leq \frac{b\eps_n}{2}} |f_n(x_j) - f_n(\bar{x}_k)|^p, \;\; \quad \text{since } |x_k-\bar{x}_k|\leq \frac{b\eps_n}{2} \\
 & \geq \frac{a}{2^{p+1} \eps^{p+d}_n n^2} \l \osc_{\frac{b\eps_n}{2}}(f_n)(x_k) \r^p \#\lb j \, : \, |x_j-x_k|\leq \frac{b\eps_n}{2} \rb \\
 & = \frac{a}{2^{p+1} \eps^{p+d}_n n} \l \osc_{\frac{b\eps_n}{2}}(f_n)(x_k)  \r^p \mu_n \l B \l x_k, \frac{b \eps_n}{2} \r  \r.
\end{split} \label{eq:pf411}
\end{align}
where  $\mu_n = \frac{1}{n} \sum_{i=1}^n \delta_{x_i}$.
Now, for a transport map $T_n:\Omega\to \Omega_n$ from $\mu$ to $\mu_n$, satisfying the conclusions of Theorem~\ref{thm:Back:TransMapBound}, we have
\begin{align}
\begin{split}
\frac{1}{\eps_n^d }  \mu_n \l B\l x_k, \frac{b \eps_n}{2} \r \r & = \frac{1}{\eps_n^d} \int_\Omega \bbI_{\{|T_n(x)-x_k|\leq \frac{b\eps_n}{2}\}} \rho(x) \, \dd x \\
 & \geq \frac{\inf_{x\in \Omega} \rho}{\eps_n^d} \int_\Omega \bbI_{\{ |x-x_k|\leq \frac{b\eps_n}{2} - \|T_n-\Id\|_{L^\infty}\}} \, \dd x \\
 & = \big( \inf_{x\in \Omega} \rho(x) \big) \Vol\l B\l 0,\frac{b}{2}-\frac{\|T_n-\Id\|_{L^\infty}}{\eps_n} \r \r.
 \end{split} \label{eq:pf412}
\end{align}
We choose $n_0$ such that  for $n \geq n_0$ it holds that $\frac{\|T_n-\Id\|_{L^\infty}}{\eps_n}\leq \frac{b}{4}$. Combining \eqref{eq:pf411} and \eqref{eq:pf412} gives
\[ \l\osc_{\frac{b\eps_n}{2}}(f_n)(x_k) \r^p \leq \frac{2^{p+1} \eps_n^p n \cE_n^{(p)}(f_n)}{a \big( \inf_{x\in \Omega} \rho(x) \big) \Vol\l B\l 0, \frac{b}{4}\r \r} =: C_1 \eps_n^p n \cE_n^{(p)}(f_n). \]
For $\alpha>\alpha_0$, using $\alpha_0 \geq b$ and applying  the triangle inequality $\left \lfloor \frac{2 \alpha}{b}  \right \rfloor$ times,  we obtain
\[ \l \osc_{\alpha\eps_n}(f_n)(x_k) \r^p \leq C_1  \l \left \lfloor \frac{2 \alpha}{b} \right \rfloor +1 \r^p  \eps_n^p n\cE_n^{(p)}(f_n)
 \leq  C_1 \l \frac{ 3 \alpha}{b }\r^p   \eps_n^p n \cE_n^{(p)}(f_n) 
 \]
which completes the proof.
\end{proof}

\begin{lemma}[discrete to nonlocal control]
\label{lem:EntoNLBound}
Let $p\geq 1$. Assume $\Omega$, $\mu$, $\eta$, and $x_i$ satisfy \textbf{(A1) - (A8)}. Let graph weights $W_{ij}$ be given by~\eqref{eq:MainRes:Wij}. 
Let constants $a,b>0$ be such that for $\tilde{\eta}(|x|) = a$ for $|x|\leq b$ and $\tilde{\eta}(|x|) = 0$ otherwise it holds that $\tilde \eta \leq \eta$.
Let $T_n$ be a transport map satisfying the results of Theorem~\ref{thm:Back:TransMapBound} and let $\tilde{\eps}_n = \eps_n-\frac{2\|T_n-\Id\|_{L^\infty}}{b}$.
Then there exists constants $n_0>0$ and $C>0$ (independent of $n$ and $f_n$) such that for all $n\geq n_0$
\[ \cE_{\tilde{\eps}_n}^{(NL,p)}(f_n\circ T_n;  \tilde \eta) \leq C \cE_n^{(p)}(f_n ; \eta) \]
where $\cE_{\tilde{\eps}_n}^{(NL,p)}$ is defined by
\begin{equation} \label{eq:Proofs:Compact:EnNL} 
\cE_\eps^{(NL,p)}(f;\eta) = \frac{1}{\eps^p} \int_\Omega \int_\Omega \eta_\eps(|x-z|) |f(x)-f(z)|^p \, \dd x \, \dd z.
\end{equation}
\end{lemma}
\begin{proof}
Assume $\la\frac{x-z}{\tilde{\eps}_n}\ra < b$ then
\[ |T_n(x) - T_n(z)| \leq 2 \| T_n - \Id \|_{L^\infty} + | x-z| \leq 2 \| T_n - \Id \|_{L^\infty} + b \tilde{\eps}_n = b\eps_n. \]
So,
\[ \la \frac{x-z}{\tilde{\eps}_n} \ra < b \Rightarrow \la \frac{T_n(x) - T_n(z)}{\eps_n} \ra \leq b \]
and therefore
\[ \tilde{\eta}\l \frac{|x-z|}{\tilde{\eps}_n} \r \leq \tilde{\eta}\l \frac{|T_n(x) - T_n(z)|}{\eps_n} \r \leq \eta\l\frac{|T_n(x) - T_n(z)|}{\eps_n} \r. \]
Now,
\begin{align*}
\cE_{\tilde{\eps}_n}^{(NL,p)}(f_n\circ T_n) & \leq \frac{\eps_n^d}{\tilde{\eps}_n^{d+p}} \int_{\Omega^2} \eta_{\eps_n}(|T_n(x)-T_n(z)|) \la f_n(T_n(x)) - f_n(T_n(z)) \ra^p \, \dd x \, \dd z \\
 & = \frac{\eps_n^{d+p}}{\big(\inf_{x\in \Omega} \rho^2(x) \big) \tilde{\eps}_n^{d+p}} \cE_n^{(p)}(f_n).
\end{align*}
Since $\frac{\eps_n}{\tilde{\eps}_n}\to 1$ we are done.
\end{proof}

In the next lemma we show that  that boundedness of non-local energies implies regularity at scales greater $\eps$. This allows us to relate non-local bounds to local bounds after mollification. 
\begin{lemma}[nonlocal to averaged local]
\label{lem:Proofs:Compact:NLtoLBound}
Assume $\Omega\subset \bbR^d$ is open and bounded and $p\geq 1$.
Assume that $\eta:[0,\infty)\to [0,\infty)$ is non-increasing, $\eta(0)>0$ 
and $\eta$ is continuous near $0$.
Then there exists a constant $C\geq 1$ and a mollifier $J$ with $\supp(J)\subseteq \overline{B(0,1)}$ such that for all $\eps>0$, $f\in L^p(\Omega)$, $\Omega'\subset\subset \Omega$ with $\dist(\Omega',\partial \Omega)> \eps$ it holds that
\[ \cE_\infty^{(p)}( J_\eps * f;\Omega') \leq C\cE_\eps^{(NL,p)}(f). \]
where $\cE_\infty^{(p)}$ is defined by~\eqref{eq:MainRes:Einfty} and $\cE_\eps^{(NL,p)}$ is defined by~\eqref{eq:Proofs:Compact:EnNL}.
\end{lemma}
\begin{proof}
Let $J$ be a radially symmetric mollifier supported in $B(0,1)$ and such that for some $\beta>0$, $J \leq \beta \eta$ and 
$|\nabla J| \leq \beta \eta$.
Without loss of generality we can assume $\supp(\eta)\subset\overline{B(0,1)}$. 
Let $J_\eps(\tacka) = J(\tacka /\eps )/\eps^d$ and
let $g_\eps = J_\eps * f$.
For arbitrary  $x\in \Omega$ with $\dist(x,\partial \Omega)> {\eps}$ we have
\begin{align*}
|\nabla {g_\eps}(x) |  & = \la \int_\Omega \nabla J_{{\eps}} \l x-z\r f(z) \, \dd z \ra \\
& =  \la  \int_\Omega \nabla J_{{\eps}}\l x-z \r \l f(z) - f(x) \r \, \dd z - \int_{\bbR^d\setminus \Omega} \nabla J_{{\eps}} \l x-z \r f(x) \, \dd z \ra \\
& \leq \frac{\beta}{{\eps}^{d+1}} \int_\Omega {\eta}\l \frac{x-z}{{\eps}} \r \la f(z) - f(x) \ra \, \dd z + \frac{1}{{\eps^{d+1}}} \int_{\bbR^d\setminus \Omega} \la (\nabla J) \l \frac{x-z}{\eps}\r \ra \la f(x) \ra \, \dd z.
\end{align*}
where the second line follows from $ \int_{\bbR^d} \nabla J(w) \, \dd w = 0$.
For the second term we have
\begin{align*}
\frac{1}{{\eps}^{d+1}} \int_{\bbR^d\setminus \Omega} \la \nabla J \l\frac{x-z}{{\eps}}\r \ra |f(x)| \, \dd z & =0 
\end{align*}
since for all $z \in \R^d\setminus \Omega$ and $x \in \Omega$ with $\dist(x,\partial \Omega)> {\eps}$ it follows that $|x-z| > \eps$ and thus $\nabla J \l\frac{x-z}{{\eps}}\r =0$.
Therefore,
\begin{align*}
|\nabla {g_\eps}(x) |^p & \leq \beta^p \l \int_\Omega \frac{1}{{\eps}} {\eta}_{{\eps}}(x-z) \la f(z) - f(x) \ra \, \dd z \r^p \\
 & \leq \gamma_{{\eta}}^{p-1} \beta^p \int_\Omega {\eta}_{{\eps}}(x-z) \frac{\la f(z) - f(x) \ra^p}{{\eps}^p} \, \dd z 
\end{align*}
by Jensen's inequality and where $\gamma_{{\eta}} = \int_{B(0,1)} {\eta}(w) \, \dd w$.
Hence,
\begin{align*}
\int_{\Omega'} \la \nabla {g_\eps}(x) \ra^p \, \dd x & \leq \gamma_{{\eta}}^{p-1}\beta^p \int_\Omega \int_\Omega {\eta}_{{\eps}}(|x-z|) \la \frac{f(z) - f(x)}{{\eps}^p} \ra^p \, \dd z \, \dd x \\
 & \leq \gamma_{{\eta}}^{p-1}\beta^p \cE_{{\eps}}^{(NL,p)}(f)
\end{align*}
which completes the proof.
\end{proof}

We  prove the compactness property for bounded sequences.
The convergence of a subsequence is a consequence of being able to bound $\tilde{g}_n=J_{\eps_n}\ast (f_n\circ T_n)$ in $W^{1,p}$ (hence the sequence $\{\tilde{g}_n\}_{n}$ is precompact in $L^p(\mu)$) and show $\|f_n\circ T_n - \tilde{g}_n\|_{L^p}\to 0$.

\begin{proposition}[compactness]
\label{prop:Proofs:Compact:Compact}
Consider the assumptions and the graph construction of Lemma \ref{lem:Proofs:Compact:Osc}.
Then with probability one, any sequence $f_n:\Omega_n\to \bbR$ with $\sup_{n\in\bbN} \cE_{n}^{(p)}(f_n) < \infty$ and $\sup_{n\in\bbN} \|f_n\|_{L^\infty(\mu_n)} < \infty$ has a subsequence $f_{n_m}$ such that $(\mu_{n_m}, f_{n_m})$, converges in $TL^p$ to $(\mu,f)$ for some $f \in L^p(\mu)$.
\end{proposition}
\begin{proof}
Since $\cE_{n}^{(p)}(f_n) \geq C \cE_{n}^{(1)}(f_n)$ the compactness in $TL^1$ follows from 
Theorem 1.2 in \cite{garciatrillos16}. We note that from the proof of Theorem 1.2 it follows that there in fact exists a  subsequence $f_{n_m}$, and a sequence of transportation maps ${T_{n_m}}_\sharp \mu = \mu_{n_m}$ such that
\[ \lim_{m \to \infty} \| f - f_{n_m} \circ T_{n_m} \|_{L^1(\mu)} + \|T_{n_m} - \Id \|_{L^\infty(\mu)} =0. \]
Since $ \| f - f_{n_m} \circ T_{n_m} \|_{L^\infty(\mu)} \leq M <\infty$ for some $M \in \R$, the convergence of $f_{n_m}$ to $f$ in $TL^p$ follows by interpolation.
\end{proof}

%

\begin{lemma}[uniform convergence]
\label{lem:Proofs:Compact:UnifConverg}
Consider the assumptions and the graph construction of Lemma \ref{lem:Proofs:Compact:Osc}. Assume that 
$\eps_n^p n\to 0$ as $n \to \infty$, which, due to \textbf{(A5)}, implies that $p>d$. Furthermore assume that with probability one 
$(\mu_n, f_n) \to (\mu,f)$ in $TL^p$ metric as $n \to \infty$ and that $\sup_{n\in \bbN} \cE_n^{(p)}(f_n)<\infty$. 
Then $f\in C^{0,\gamma}(\Omega)$, with $\gamma=1-\frac{d}{p}>0$, and for all
$\Omega' \subset \subset \Omega$
\[  \max_{\{k \::\: x_k \in \Omega'\}}  |f(x_k) - f_{n}(x_k)| \to 0 \quad \quad \text{as } n\to \infty. \]
Moreover, if for all $k = 1, \dots, N$, $f_n(x_k) = y_k$  for all $n$, it follows that $f(x_k) = y_k$.
\end{lemma}

\begin{proof}
Find constants $a,b>0$ such that $\tilde{\eta}(t):= a$ if $|t|\leq b$ and $\tilde{\eta}(t):= 0$ if $|t|> b$ satisfies $\tilde{\eta}\leq \eta$.
Now we define $\tilde{f}_n=f_n\circ T_n$ where $T_n$ is the transportation map satisfying the conclusions on Theorem~\ref{thm:Back:TransMapBound}  and set $\tilde{\eps}_n = \eps_n-\frac{2\|T_n-\Id\|_{L^\infty}}{b}$. Then  for  $n$ sufficiently large $\tilde{\eps}_n> 0$, and $\frac{\eps_n}{\tilde{\eps_n}}\to 1$.
We note that if $|T_n(x)-T_n(z)|>b \eps_n$ then
\[ |x-z| \geq |T_n(x) - T_n(z)| - 2\|T_n-\Id\|_{L^\infty} > b \eps_n - 2\|T_n-\Id\|_{L^\infty} = \tilde{\eps}_n b. \]
Hence, $\tilde{\eta}\l\frac{|x-z|}{\tilde{\eps}_n}\r \leq \tilde{\eta}\l \frac{|T_n(x)-T_n(z)|}{\eps_n}\r$.
Let $\cE_{\tilde{\eps}}^{(NL,p)}$ be the non-local Dirichlet energy defined in \eqref{eq:Proofs:Compact:EnNL} with $\eps=\tilde{\eps}_n$ and $\eta=\tilde{\eta}$.
Then, by Lemma \ref{lem:EntoNLBound}
\begin{align*}
\cE_{\tilde{\eps}}^{(NL,p)}(\tilde{f}_n) & \leq 
C \cE_n^{(p)}(f_n).
\end{align*}
Hence, $\cE_{\tilde{\eps}}^{(NL,p)}(\tilde{f}_n)$ is bounded and therefore, by Lemma~\ref{lem:Proofs:Compact:NLtoLBound} we have that $\cE_{\infty}^{(p)}(J_{\tilde{\eps}_n}\ast f_n;\Omega')$ is bounded for every $\Omega'\subset\subset \Omega$.
One can easily show $\|J_{\tilde{\eps}_n}\ast f_n\|_{L^p(\Omega')} \leq \|\tilde{f}_n\|_{L^p}$ and therefore $J_{\tilde{\eps}_n} \ast \tilde{f}_n$ is locally bounded in $W^{1,p}$.
We also note that since $f_n \circ T_n$ converges to $f$ in $L^p(\mu)$ 
\begin{align*}
\| J_{\tilde{\eps}_n} \ast \tilde{f}_n - f\|_{L^p(\Omega')} \leq & \: \| J_{\tilde{\eps}_n} \ast \tilde{f}_n - J_{\tilde{\eps}_n} \ast f + J_{\tilde{\eps}_n} \ast f - f\|_{L^p(\Omega')} \\
\leq  & \: \| \tilde{f}_n -  f \|_{L^p(\Omega)} + \|J_{\tilde{\eps}_n} \ast f - f\|_{L^p(\Omega')} \to 0 \quad \te{ as } n \to \infty.
\end{align*}
Since $J_{\tilde{\eps}_n} \ast \tilde{f}_n \to f$ in $L^p(\Omega')$, by the compactness of the embedding of $W^{1,p}(\Omega')$ into $C^{0,\gamma}$ (Morrey's inequality), for $\gamma=1-\frac{d}{p}$, we have that 
\[ J_{\tilde{\eps}_n} \ast \tilde{f}_n  \to  f \quad \te{ uniformly on }\; \Omega' \;\te{ as } n \to \infty. \]
Therefore, for each $k\in \{1,\dots, N\}$, $J_{\tilde{\eps}_{n}} \ast \tilde{f}_{n}$ converges uniformly to $f$ on $B(0,\delta)$ for any $\delta$ such that $B(x_k,\delta)\subset \Omega$.
For any $x\in B(x_k,3\tilde{\eps}_n)\cap \Omega_n$ we have (for a constant $C$)
\[ |f_n(x_k) - f_n(x)| \leq \osc_{3\tilde{\eps}_n}(f_n)(x_k) \leq \osc_{4\eps_n}(f_n)(x_k) \leq \l 4^pC\cE_n^{(p)}(f_n) n \eps_n^p\r^{\frac{1}{p}} \to 0 \]
by Lemma~\ref{lem:Proofs:Compact:Osc}.
It follows that
\[  \max_{k=1, \dots, n} \max_{x\in B(x_k,3\tilde{\eps}_n)\cap\Omega_n} |f_n(x) - f_n(x_k)| \to 0. \]
To complete the proof we notice that for any $\Omega' \subset\subset  \Omega$
\begin{align*}
& \max_{\{k \::\: x_k \in \Omega'\}}  |f(x_k) - f_{n}(x_k)| \\
& \leq \max_{\{k \::\: x_k \in \Omega'\}}  |f(x_k) - J_{\tilde{\eps}_{n}} \ast \tilde{f}_{n}(x_k)| + |J_{\tilde{\eps}_{n}} \ast \tilde{f}_{n}(x_k) - f_n(x_k)| \\
& \leq \|f - J_{\tilde{\eps}_{n}} \ast \tilde{f}_{n} \|_{L^\infty(\Omega')} + \max_{\{k \::\: x_k \in \Omega'\}}  \int_{B(0,2\tilde{\eps}_{n})} J_{\tilde{\eps}_{n}}(x_k-x) \la f_{n}(T_n(x)) - f_{n}(x_k) \ra \, \dd x \\
& \leq \|f - J_{\tilde{\eps}_{n}} \ast \tilde{f}_{n} \|_{L^\infty(\Omega')} +  \max_{\{k \::\: x_k \in \Omega'\}}  \sup_{x\in B(x_k,3\tilde{\eps}_{n})\cap \Omega_{n}} \la f_{n}(x) - f_{n}(x_k) \ra 
\end{align*}
and the above converges to zero for all $x_k$.
\end{proof}

\subsection{Asymptotic Consistency via \texorpdfstring{$\Gamma$}{Gamma}--Convergence \label{subsec:Proofs:Gamma} }

We approach proving Theorem \ref{thm:MainRes:Conv} using $\Gamma$-convergence. Namely as pointed out in Section~\ref{subsec:Back:Gamma} convergence of minimizers follows from $\Gamma$-convergence and compactness. We use the general setup of \cite{garciatrillos16}. In particular we first establish in Lemma  \ref{lem:GammaNLL} that nonlocal functionals $\cE_{\eps_n}^{(NL,p)}$ $\,\Gamma$-converge to $\Ecp$. We then state and prove the $\Gamma$-convergence of  $\cE_{n,con}^{(p)}$ towards $\Ecp$ or $\Ecpc$ depending on how quickly $\eps_n \to 0$ as $n \to \infty$. Steps of proving this claim rely on Lemma 
\ref{lem:GammaNLL}.

\begin{lemma}[continuum nonlocal to local]
\label{lem:GammaNLL}
Let $p > 1$. Assume $\Omega$ satisfy the assumptions \textbf{(A1) - (A2)} and $\eta$ satisfies assumptions \textbf{(A6) - (A8)}. 
Then   $\cE_{\eps}^{(NL,p)} \,$, defined in \eqref{eq:Proofs:Compact:EnNL}, $\Gamma$-converges as $n \to \infty$ in $L^p(\Omega)$ to the functional $\Ecp$ defined in \eqref{eq:MainRes:Einfty}.
\end{lemma}
If $\rho$ is constant and $\Omega$ is convex this result is contained in the appendix to \cite{AB2}. For general $\Omega$ it follows from Theorem 8 in \cite{Ponce}. We remark that while the functional in \cite{Ponce} appears different the term $|x-y|^p$ which arises can be absorbed in the kernel. The results can be extended to general $\rho$ in a straightforward manner as has been done for $p=1$ in Section 4 of \cite{garciatrillos16} and has been remarked in Proposition 1.10 in \cite{garciatrillos15aAAA}.

\begin{theorem}[discrete to local $\Gamma$-convergence]
\label{thm:Proofs:Gamma:Gamma}
Let $p > 1$. Assume $\Omega$, $\mu$, $\eta$, $\eps_n$, and $x_i$ satisfy the assumptions \textbf{(A1) - (A8)}. Let graph weights $W_{ij}$ be given by~\eqref{eq:MainRes:Wij}. 
Let $M\geq \max_{i=1,\dots, N}|y_i|$.
Then with probability one  $\cE_{n,con}\,$, defined in \eqref{eq:MainRes:Encon}, $\Gamma$-converges as $n \to \infty$ in $TL^p$ metric  on the set $\{(\nu,g) \, : \, \nu \in \mathcal{P}(\Omega), \|g\|_{L^\infty(\nu)}\leq M\}$ to the functional 
\begin{equation*} 
\begin{cases}
\Ecpc \quad & \te{if } \lim_{n \to \infty} \,n \eps_n^p = 0 \\
\Ecp \quad & \te{if } \lim_{n \to \infty} \,n \eps_n^p = \infty \\
\end{cases}
\end{equation*}
where  $\Ecp$ is defined in \eqref{eq:MainRes:Einfty} and $\Ecpc$ is defined in \eqref{eq:MainRes:Einftycon}.
\end{theorem}
Restricting the space to the set of functions bounded by $M$ is really needed only for the case $\lim_{n \to \infty} \,n \eps_n^p = \infty$. It is required  since the functional $\Ecp$ is invariant under adding a constant and thus the  loss of constraints in the limit when $\lim_{n \to \infty} \,n \eps_n^p = \infty$ would lead to loss of compactness, without the restriction. We note that placing an upper bound on $f$ is not restrictive in practice since both discrete and continuum minimizers satisfy the bound. 

We prove the liminf inequalities  and the existence of a recovery sequence separately. Since $ \cE_{\infty}^{(p)} \leq \cE_{\infty,con}^{(p)}$ the liminf inequalities needed can be stated in the following way.
\begin{lemma}
\label{lem:Proofs:Gamma:Liminf}
Under the same conditions as Theorem~\ref{thm:Proofs:Gamma:Gamma}, with probability one, for any $f \in L^p$ with $\|f\|_{L^\infty(\mu)} \leq M$ and any sequence $f_n\to f$ in $TL^p$ with $\|f_n\|_{L^\infty(\mu_n)} \leq M$ we have
\begin{equation} \label{eq:temp_linf1}
\cE_{\infty}^{(p)}(f)  \leq \liminf_{n\to \infty} \cE_{n}^{(p)}(f_n) \leq \liminf_{n\to \infty} \cE_{n,con}^{(p)}(f_n). 
\end{equation}
Furthermore if $\lim_{n \to \infty} \,n \eps_n^p = 0 $ then
\begin{equation} \label{eq:temp_linf2}
\cE_{\infty,con}^{(p)}(f) \leq \liminf_{n\to \infty} \cE_{n,con}^{(p)}(f_n). 
\end{equation}
\end{lemma}
\begin{proof}
Let $f_n\to f$ in $TL^p$.
The first inequality of \eqref{eq:temp_linf1} follows from Lemma \ref{lem:GammaNLL} in the same way the analogous result is shown for $p=1$ in Section 5 of \cite{garciatrillos16}. The second inequality follows from definition of $ \cE_{n}^{(p)}$  and  $\cE_{n,con}^{(p)}$

When $\lim_{n \to \infty} \,n \eps_n^p = 0$ the  inequality \eqref{eq:temp_linf2}  is a consequence of  Lemma \ref{lem:Proofs:Compact:UnifConverg}.
\end{proof}

We now prove the existence of a recovery sequence. 
Since $ \cE_{\infty}^{(p)} \leq \cE_{\infty,con}^{(p)}$ we state it in the following way.
\begin{lemma}
\label{lem:Proofs:Gamma:Limsup}
Under the same conditions as Theorem~\ref{thm:Proofs:Gamma:Gamma}, with probability one, for any function $f \in L^p$, with $\|f\|_{L^\infty(\mu)}\leq M$ there exists a sequence $f_n$ satisfying $f_n\to f$ in $TL^p$ with $\|f_n\|_{L^\infty(\mu_n)}\leq M$ and
\begin{equation} \label{eq:Proofs:Gamma:Limsupcon}
\cE_{\infty,con}^{(p)}(f) \geq \limsup_{n\to \infty} \cE_{n,con}^{(p)}(f_n).
\end{equation}
Furthermore if $\lim_{n\to \infty} n \eps_n^p = \infty$ then
\begin{equation} \label{eq:Proofs:Gamma:Limsup}
\cE_{\infty}^{(p)}(f) \geq \limsup_{n\to \infty} \cE_{n,con}^{(p)}(f_n).
\end{equation}
\end{lemma}
\begin{proof}
The proof of the first inequality is a straightforward adaptation of the analogous result for $p=1$ in Section 5 of \cite{garciatrillos16}. The recovery sequence used is defined as a restriction of $f$ to $\Omega_n$: $f_n(x_i) =  f(x_i)$ for all $i=1,\dots, n$, and thus satisfies the constraints and $\|f_n\|_{L^\infty(\mu_n)}\leq M$. 

The same argument and recovery sequence construction can be used to show that 
with probability one, for any function $f \in L^p$, with $\|f\|_{L^\infty(\mu)}\leq M$ there exists a sequence $f_n$ satisfying $f_n\to f$ in $TL^p$ with $\|f_n\|_{L^\infty(\mu_n)}\leq M$ and
\begin{equation} \label{eq:Limsupuncon}
\cE_{\infty}^{(p)}(f) \geq \limsup_{n\to \infty} \cE_{n}^{(p)}(f_n).
\end{equation}

Let us now consider that case that $n \eps_n^p \to \infty$ as $n \to \infty$ and show the second inequality. 
Suppose $\cE_{\infty,con}^{(p)}(f)<\infty$ else the lemma is trivial. 
Let $f_n$ be the recovery sequence for \eqref{eq:Limsupuncon}.

We define $\hat f_n: \Omega_n \to \R$ by
\[  \hat f_n(x_i) = 
  \begin{cases} y_i \quad & \te{for } i=1, \dots, N, \\
   f_n(x_i) & \te{for } i=N+1,\dots, n.
   \end{cases}
\]
We note that $\hat f_n\to f$ in $TL^p$ with $\|\hat f_n\|_{L^\infty(\mu_n)}\leq M$.
To show \eqref{eq:Proofs:Gamma:Limsup} it suffices to show that

\begin{equation} \label{eq:temp_lsup1}
 \lim_{n \to \infty} \cE_{n}^{(p)}(f_n) - \cE_{n,con}^{(p)}(\hat f_n) = 0. 
\end{equation}
We may write,
\begin{align} \label{eq:temp_lsup2}
\begin{split}
\la \cE_{n}^{(p)}(f_n) - \cE_{n,con}^{(p)}(\hat f_n) \ra & \leq \frac{1}{\eps_n^p} \frac{2}{n^2} \sum_{i=1}^N \sum_{j=1}^n  \eta_{\eps_n} (|x_i - x_j|)  \la \, |f(x_i) - f(x_j) |^p - | y_i - f(x_j) |^p \ra \\
& \leq \frac{2^{p+1}M^p}{\eps_n^p n} \sum_{i=1}^N  \frac{1}{n} \sum_{j=1}^n  \eta_{\eps_n} (|x_i - x_j|) 
\end{split}
\end{align}
\noindent\emph{Step 1.} Let us consider first the case that $\eta(t) = a$ if $|t|<b$ and $\eta(t) = 0$ otherwise for some $a,b>0$. Then, using Theorem  \ref{thm:Back:TransMapBound}
\begin{align*}
 \frac{1}{n} \sum_{j=1}^n  \eta_{\eps_n} (|x_i - x_j|)  & \leq \frac{\eta(0)}{\eps^d} \mu_n (B(x_i, \eps b) \\
 & \leq \frac{\eta(0)}{\eps^d} \mu(B(x_i, \eps b + \| \Id - T_n \|_{L^\infty} )) \\
 & \leq \eta(0) \l \frac{\eps b + \| \Id - T_n \|_{L^\infty} }{\eps} \r^d \Vol(B(0,1)) \| \rho\|_{L^\infty} \leq C. 
\end{align*}
Combining this inequality with \eqref{eq:temp_lsup2} implies \eqref{eq:temp_lsup1}. 
%

\noindent\emph{Step 2.} Consider now general $\eta$ satisfying \textbf{(A6)-(A8)}. Let 
\[ \tilde \eta(t) = 
    \begin{cases}
        \eta(0) \quad & \te{if } |t| \leq 1 \\
        \eta(t)   & \te{otherwise.}
    \end{cases}
\]
Note that $\tilde \eta$ is radially nonincreasing, $\tilde \eta \geq \eta$, and that $\tilde \eta((|x|-1)_+) \leq \tilde \eta(|x|/2)$. Theorem  \ref{thm:Back:TransMapBound} implies that  for $n$ large $\| \Id - T_n \|_{L^\infty} \leq \eps_n$. Consequently
\begin{align*}
 \frac{1}{n} \sum_{j=1}^n  \eta_{\eps_n} (|x_i - x_j|)  & \leq  \frac{1}{n} \sum_{j=1}^n  \tilde\eta_{\eps_n} (|x_i - x_j|)  \\
& =  \frac{1}{\eps_n^d} \int_\Omega \tilde \eta  \l \frac{|x_i - T_n(y)|}{\eps_n} \r d \mu(y)  \\
& \leq  \frac{1}{\eps_n^d}  \int_\Omega \tilde \eta \l \frac{|x_i - y|}{2 \eps_n} \r d \mu(y) \leq C
\end{align*}
where the penultimate inequality follows from $\frac{|x_i-T_n(y)|}{\eps_n} \geq \l \frac{|x_i-y|-\|T_n-\Id\|_{L^\infty}}{\eps_n} \r_+ \geq \l \frac{|x_i-y|}{\eps_n} - 1\r_+$.
Again combining this estimate with \eqref{eq:temp_lsup2} implies \eqref{eq:temp_lsup1}. 
\end{proof}
\medskip

We now state the $\Gamma$-convergence result relevant for the penalized model $\cS_n^{(p)}$.
\begin{lemma}
\label{lem:Proofs:ConvSoft:Gamma&Compact}
Under the conditions of Proposition~\ref{prop:MainRes:ConvSoft} we have:
\begin{itemize}
\item (compactness) Any sequence $f_n:\Omega_n\to \R$ with $\sup_{n\in \bbN} \cS_n^{(p)}(f_n) +  \|f_n\|_{L^\infty(\mu_n)}<\infty$ has, with probability one, a subsequence $f_{n_m}$ such that there exists $f_\infty\in W^{1,p}$ with $f_{n_m}\to f_\infty$ in $TL^p$.
\item ($\Gamma$-convergence, well-posed regime) If $\eps_n^p n \to 0$ then, with probability one,  on the set $(\mu_n,f_n)$ with $\|f_n\|_{L^\infty(\mu_n)}\leq M$,
\[ \Glim_{n\to \infty} \l \cE_n^{(p)} + \lambda R^{(q)}\r = \Ecp + \lambda R^{(q)} \]
where the $\Gamma$-convergence is considered in $TL^p$ topology.
\item ($\Gamma$-convergence, degenerate regime) If $\eps_n^p n \to \infty$ then, with probability one, on the set $(\mu_n,f_n)$ with $\|f_n\|_{L^\infty(\mu_n)}\leq M$,
\[ \Glim_{n\to \infty} \l \cE_n^{(p)} + \lambda R^{(q)}\r = \Ecp, \]
where the $\Gamma$-convergence is considered in $TL^p$ topology.
\end{itemize}
\end{lemma}
\begin{proof}
The compactness follows directly from Proposition~\ref{prop:Proofs:Compact:Compact}.

When $\eps_n^p n \to 0$, for the liminf inequality assume $f_n\to f$ in $TL^p$ and $\liminf_{n\to \infty} \cE_n^{(p)}(f_n)<\infty$. Then by Lemma~\ref{lem:Proofs:Compact:UnifConverg} $f_n(x_k)\to f(x_k)$ for all $k\in \{1,\dots, N\}$ and hence $\lambda R^{(q)}(f_n) \to \lambda R^{(q)}(f)$. By \eqref{eq:temp_linf1} of Lemma~\ref{lem:Proofs:Gamma:Liminf} we have $\liminf_{n\to \infty} \l \cE_n^{(p)}(f_n) + \lambda R^{(q)}(f_n) \r \geq \Ecp(f) + \lambda R^{(q)}(f)$.
The limsup inequality follows in a similar manner from equation~\eqref{eq:Limsupuncon} and Lemma~\ref{lem:Proofs:Compact:UnifConverg}.


If $\eps_n^p n \to \infty$, then the liminf inequality follows from \eqref{eq:temp_linf1} of Lemma~\ref{lem:Proofs:Gamma:Liminf}, while, the limsup inequality follows directly from
\[ \limsup_{n\to \infty} \cE_n^{(p)}(f_n) + \lambda R^{(q)}(f_n) \leq \limsup_{n\to \infty} \cE^{(p)}_{n,con}(f_n) \leq \cE_\infty^{(p)}(f) \]
and Lemma \ref{lem:Proofs:Gamma:Limsup}.
\end{proof}

\subsection{Proofs of Theorem~\ref{thm:MainRes:Conv} \label{subsec:Proofs:ConvHard}  and Proposition 
 \ref{prop:MainRes:ConvSoft}} \label{subsec:MainRes:ConvSoft}

The $\Gamma$-convergence and compactness results above allow us to prove Theorem~\ref{thm:MainRes:Conv}.
It is a general result that  $\Gamma$-convergence and compactness  imply the convergence of minimizers
(as well as of almost minimizers) to a minimizer of the limiting problem, see~\cite[Theorem 1.21]{braides02} or Theorem~\ref{thm:Back:Gamma:Conmin}. 

\begin{proof}[Proof of Theorem~\ref{thm:MainRes:Conv}.]
Let $f_n$ be a minimizer of $\cE_{n,con}^{(p)}$.
Recall that $M \geq \|y\|_{L^\infty(\mu_n)}$. Note that if $\|f_n\|_{L^\infty(\mu_n)} > M$ then, since the graph is connected with high probability  $p_n$, such that $\sum_{n=1}^\infty (1-p_n) < \infty$, for $\hat f_n = (f_n \wedge M) \vee (-M)$ we have
$\cE_{n,con}^{(p)}(\hat f_n) < \cE_{n,con}^{(p)}(f_n)$ which contradicts the definition of $f_n$. 
Thus with high probability $\|f_n\|_{L^\infty}\leq M$ for each $n$, hence we can restrict the minimization to the set of $(f_n,\mu_n)$ such that $\|f_n\|_{L^\infty(\mu_n)}\leq M$.
This allows us to consider the setting of Theorem \ref{thm:Proofs:Gamma:Gamma}.

By compactness result of Proposition \ref{prop:Proofs:Compact:Compact} there exists a subsequence 
$f_{n_m}$ converging in $TL^p$ to $f \in L^p(\mu)$. 

To prove (i) assume that $n \eps_n^p \to 0$ as $n \to \infty$. 
The uniform convergence of statement (a) then follows from Lemma \ref{lem:Proofs:Compact:UnifConverg}. 
The $\Gamma$-convergence result of Theorem \ref{thm:Proofs:Gamma:Gamma}  implies
  that  $f $ minimizes  $\cE_{\infty,con}^{(p)}$. Since the minimizer of $\Ecpc$ is unique the convergence holds along the whole sequence, thus establishing statement (c). 

To prove (ii) assume that $n \eps_n^p \to 0$ as $n \to \infty$. 
Again, Theorem \ref{thm:Proofs:Gamma:Gamma}  implies  that  $f $ minimizes  $\cE_{\infty}^{(p)}$.
\end{proof}
\medskip


The results of the Proposition~\ref{prop:MainRes:ConvSoft} are proved by the same arguments; using 
Lemma \ref{lem:Proofs:ConvSoft:Gamma&Compact} instead of Theorem \ref{thm:Proofs:Gamma:Gamma}.
%

\section{Improved Model  \label{sec:IM}}

In Theorem \ref{thm:MainRes:Conv} we proved that  the model $\cE_{n,con}^{(p)}$, defined in \eqref{eq:MainRes:Encon}, is consistent as $n \to \infty$ and lower bounds \textbf{(A5)} hold,  only if 
\[ \frac{1}{n^p} \gg \eps_n. \]
This upper bound is  undesirable as it restricts the range of $\eps$ that can be used. Furthermore in a nonasymptotic regime, for large but fixed finite $n$, it provides no guidance to what $\eps$ are appropriate (small enough). Finally as our numerical experiments show, see Figures~\ref{fig:Ex:1D:ErrorVEpsp=1.5}(a) and~\ref{fig:Ex:1D:ErrorVEpsp=2}(a), the range of $\eps$ for which the limiting problem is approximated well can be quite narrow. This problem is particularly pronounced if $p>d$ is close to $d$, which is the regime identified in \cite{elalaoui16} as the most relevant for semi-supervised learning. 

It would be advantageous to have another model, asymptotical consistent with $\cE_{\infty,con}^{(p)}$ which would not require an upper bound on $\eps_n$ (other than $\eps_n$ to converge to zero) as $n \to \infty$. 
Here we introduce a new, related, model $\cF_{n,con}^{(p)}$ which has the desired properties, and whose minimizers can be computed with the same algorithms as those for $\cE_{n,con}^{(p)}$.


We define the set of functions which are constant near the labeled points:
\[ C_n^{(\delta)} = \{ f : \Omega_n \to \R \::\: f(x_k) = y_i \; \te{ whenever } |x_k - x_i|<\delta \te{ for } i=1, \dots, N\} \]
Let $L = \min\{ |x_i - x_j| \::\: i \neq j\}/2$ and $R_n = \min\{2 \eps_n, L\}$. 
The new functional is defined by
\begin{equation} \label{eq:cF}
\cF_{n,con}^{(p)}(f) = 
\begin{cases}
\frac{1}{\eps_n^p} \frac{1}{n^2} \sum_{i,j=1}^n W_{ij} |f(x_i) - f(x_j) |^p & \te{if }  f \in C_n^{(R_n)} \\
\infty & \te{else.}
\end{cases}
\end{equation}
We note that for $f \in C_n^{(R_n)}$,   $\, \cF_{n,con}^{(p)}(f) = \cE_{n,con}^{(p)}(f)$ and that 
$\, \cF_{n,con}^{(p)}(f) \geq  \cE_{n,con}^{(p)}(f)$ for all $f$.

For the asymptotic consistency we still need to require $p>d$, since only then is the limiting model
$\cE_{\infty,con}^{(p)}$ well defined. In Theorem \ref{thm:MainRes:Conv} this followed from the assumption $n \eps_n^p \to 0$ as $n \to \infty$. Since we no longer require the upper bound on $\eps_n$ we need to require $p>d$  explicitly.

\begin{theorem}[Consistency of the improved model]
\label{thm:MainRes:Conv_IM}
Let $p > d$. Assume $\Omega$, $\mu$, $\eta$, and $x_i$ satisfy the assumptions \textbf{(A1) - (A8)}. Let graph weights $W_{ij}$ be given by~\eqref{eq:MainRes:Wij}. 
Let $f_n$ be a sequence of minimizers of $\cF_{n,con}^{(p)}$ defined in \eqref{eq:cF}. 
Then, almost surely, the sequence  $(\mu_n,f_n)$ is precompact in the $TL^p$ metric. The $TL^p$ limit of any convergent subsequence, $(\mu_{n_m}, f_{n_m})$, is of the form $(\mu,f)$ where $f \in W^{1,p}(\Omega)$ 
is a minimizer of $\Ecpc$ defined in \eqref{eq:MainRes:Einftycon}.
\end{theorem}

Proof of the theorem is a straightforward modification of the proof of Theorem \ref{thm:MainRes:Conv}. 
It relies on the following $\Gamma$-convergence result.

\begin{theorem}[discrete to local $\Gamma$-convergence]
\label{thm:Proofs:Gamma:Gamma_IM}
Let $M\geq \max_{i=1,\dots, N}|y_i|$. Under the conditions of Theorem~\ref{thm:MainRes:Conv_IM}, 
with probability one  $\cF_{n,con}\,$ $\Gamma$-converges as $n \to \infty$ in $TL^p$ metric  on the set $\{(\nu,g) \, : \, \nu \in \mathcal{P}(\Omega), \|g\|_{L^\infty(\nu)}\leq M\}$ to the functional $\Ecpc$.
\end{theorem}

We note that 
statement \eqref{eq:temp_linf1} of Lemma \ref{lem:Proofs:Gamma:Liminf}, and Proposition \ref{prop:Proofs:Compact:Compact} hold for $\cF_{n,con}^{(p)}$ since $\cE_{n,con}^{(p)} \leq \cF_{n,con}^{(p)}$. 
We now turn to proving the liminf property and the existence of recovery sequence needed to show that 
$\cF_{n,con}^{(p)}$ $\Gamma$ converges in $TL^p$ topology to $\cE_{\infty,con}^{(p)}$.

\begin{lemma}
\label{lem:Proofs:Gamma:Liminf_IM}
Under the conditions of Theorem~\ref{thm:MainRes:Conv_IM}, with probability one, for any $f \in L^\infty(\mu)$ with $\|f\|_{L^\infty(\mu)} \leq M$ and any sequence $f_n\to f$ in $TL^p$ with $\|f_n\|_{L^\infty(\mu_n)} \leq M$ we have
\begin{equation} \label{eq:IM:liminf}
\cE_{\infty,con}^{(p)}(f) \leq \liminf_{n\to \infty} \cF_{n,con}^{(p)}(f_n). 
\end{equation}
\end{lemma}
\begin{proof}
Consider a sequence $f_n$, uniformly bounded in $L^\infty(\mu_n)$ and convergent in $TL^p$  and such that $\liminf_{n\to \infty} \cF_{n,con}^{(p)}(f_n) < \infty$. Without a loss of generality we assume $\lim_{n\to \infty} \cF_{n,con}^{(p)}(f_n) < \infty$. Note that in contrast to Lemma \ref{lem:Proofs:Gamma:Liminf} we no longer require 
$n \eps_n^p \to 0$ as $n \to \infty$. Therefore we can no longer use the uniform convergence of Lemma 
\ref{lem:Proofs:Compact:UnifConverg}.  

Nevertheless since for $n$ large $f_n = y_i$ on $B(x_i, 2 \eps)$
and $\| \Id - T_n\|_{L^\infty} < \eps$ we conclude that $\tilde f_n := f_n \circ T_n = y_i$ on $B(x_i, \eps)$
and consequently that for  $g_n := J_{\eps_n} \ast \tilde f_n$ it holds that $g_n(x_i) = y_i$. Furthermore note that $\| g_n\|_{L^\infty} \leq M$.  By bounds of Lemma \ref{lem:EntoNLBound} and Lemma \ref{lem:Proofs:Compact:NLtoLBound}, $g_n$ is uniformly bounded in $W^{1,p}(\Omega')$ for any $\Omega' \subset\subset \Omega$. 
Arguing as in the proof of Lemma \ref{lem:Proofs:Compact:UnifConverg} we conclude that $g_n \to f$ in $L^p(\Omega)$. 
Since $p>d$, $W^{1,p}$ is compactly embedded in the space of continuous functions. This implies that
 $g_n$ uniformly converges to $f$ on sets compactly contained in $\Omega$. Therefore $f(x_i) = y_i$ for all $i =1, \dots, N$. Combining this with statement \eqref{eq:temp_linf1} of Lemma \ref{lem:Proofs:Gamma:Liminf} yields \eqref{eq:IM:liminf}. 
\end{proof}

\begin{lemma}
\label{lem:Proofs:Gamma:Limsup_IM}
Under the conditions of Theorem~\ref{thm:MainRes:Conv_IM}, with probability one, for any $f \in L^\infty(\mu)$ with $\|f\|_{L^\infty(\mu)} \leq M$ there exists a sequence $f_n\to f$ in $TL^p$ with $\|f_n\|_{L^\infty(\mu_n)} \leq M$ such that
\begin{equation} \label{eq:temp_linf2_IM}
\cE_{\infty,con}^{(p)}(f) \geq \limsup_{n\to \infty} \cF_{n,con}^{(p)}(f_n). 
\end{equation}
\end{lemma}
\begin{proof}
Assume $\|f\|_{L^\infty(\mu)} \leq M$ and $\cE_{\infty,con}^{(p)}(f) < \infty$. Then $f \in W^{1, p}(\Omega)$ and  since $p>d$, $f$ is continuous. Furthermore $f(x_i) = y_i$ for all $i=1, \dots, N$. 

If there exists $\delta>0$ such that $f \in W^{1,p}(\Omega)$ satisfies $f(x) = y_i$ for all $x \in B(x_i,\delta)$ and $i=1,\dots, N$ then the proof of \eqref{eq:temp_linf2_IM} is the same as the proof of \eqref{eq:Proofs:Gamma:Limsupcon}. In particular one can use the restriction of $f$ to data points to construct a recovery sequence. 

To treat general $f$ in $W^{1,p}(\Omega)$ it suffices to find a sequence $g_n \in W^{1,p}(\Omega)$ satisfying the conditions above, namely such that 
$\| g_n \|_{L^{\infty}} \leq M$, $\,g_n(x) = y_i$ for all $x \in B(x_i,\delta_n)$ for a sequence $\delta_n \geq R$ converging to zero, which satisfies
\begin{equation} \label{eq:temp_IM1}
 \lim_{n \to \infty}  \cE_{\infty,con}^{(p)}(g_n) = \cE_{\infty,con}^{(p)}(f). 
\end{equation}
We construct the sequence in the following way. Let $\theta$ be a cut-off function supported in $B(0,2)$. That is assume $\theta : \R^d \to [0,1]$ is smooth, radially symmetric and nonincreasing such that $\theta = 1$ on $B(0,1)$, $\theta = 0$ outside of $B(0,2)$, and $|\nabla \theta|<2$. Define $\theta_\delta(z) = \theta(z/\delta)$. 

We first consider the case $N=1$. Let
\[ g_n(x) = (1-\theta_{\delta_n}(x -x_1)) f(x) + \theta_{\delta_n}(x -x_1) y_1. \]
Then 
\begin{align*}
\la \cE_{\infty,con}^{(p)}(g_n) - \cE_{\infty,con}^{(p)}(f) \ra \geq \sigma_n
\int_{\Omega} \la |\nabla g_n|^p - |\nabla f|^p \ra  \rho^2 \, \dd x & \leq 
\sigma_n   \int_{B(0, 2 \delta_n)}  \l |\nabla g_n|^p + |\nabla f|^p \r \rho^2 \, \dd x
\end{align*}
We estimate
\begin{align*}
   \int_{B(0, 2 \delta_n)} |\nabla g_n|^p \rho^2 \, \dd x & \leq 
   2^p  \int_{B(0, 2 \delta_n)} |(f(x_1) - f(x)) \nabla  \theta_{\delta_n}(x - x_1)|^p + |\nabla f(x)|^p \rho^2 \, \dd x
\end{align*}
Using that $f \in C^{0, 1-d/p}$ and furthermore, 
by the remark following Theorem~4 in Section~5.6.2 of~\cite{evans10} 
we obtain
\begin{align*}
\int_{B(0, 2 \delta_n)} \! |(f(x) - f(x_1)) \nabla  \theta_{\delta_n}(x - x_1)|^p \rho^2(x) \, \dd x & \leq C_1\delta_n^{p-d} \|\nabla f\|_{L^p(B(0,2\delta_n))}^p \| \nabla  \theta_{\delta_n}\|_{L^p(B(x_1, 2 \delta_n))}^p \\
& \leq C_1 \|\nabla f\|_{L^p(B(0,2\delta_n))}^p \|\nabla  \theta\|_{L^p(\bbR^d)}^p.
\end{align*}
Since $\lim_{ n \to \infty}  \int_{B(0, 4 \delta_n)}  |\nabla  f(x)|^p dx =0$, by combining the inequalities above we conclude that \eqref{eq:temp_IM1} holds.

Generalizing to $N>1$ is straightforward. 
\end{proof}

\section{Numerical Experiments \label{sec:Ex}}

The results of Theorem \ref{thm:MainRes:Conv} show that when $\eps_n^p n\to 0$ then the solutions to  the SSL problem~\eqref{eq:MainRes:Encon} converge to a solution of the continuum  constrained problem \eqref{eq:MainRes:Einfty}, while when $\eps_n^p n\to \infty$ they degenerate as $n \to \infty$. However, in practice, for finite $n$, this does not provide a precise guidance on what $\eps$ are appropriate.
We investigate, via numerical experiments in 1D and 2D, the affect of $\eps$ on solutions to~\eqref{eq:MainRes:Encon} in elementary examples.
We also numerically compare the results with our improved model~\eqref{eq:cF}.

\subsection{1D Numerical Experiments \label{subsec:Ex:1D}}

\begin{figure}[ht!]
\setlength\figureheight{0.35\textwidth}
\setlength\figurewidth{0.47\textwidth}
\centering
\scriptsize
\begin{subfigure}[t]{0.47\textwidth}
\centering
\scriptsize
\input{Error_n1280_Model1_1D_p1p5.tikz} 
\caption{
Error \eqref{err_def} for $n=1280$.
Black line is the mean error, dashed lines are the 10\% and 90\% quantiles.
Connectivity bound $\eps_{\mathrm{conn}}$ in blue, optimal $\eps_*^{(1.5)}$ in red and upper bound $\eps_{\mathrm{upper}}^{(1.5)}$ in orange.
Blue line is the percentage of connected graphs for a given $\eps$.
}
\end{subfigure}
\hspace*{0.04\textwidth}
\begin{subfigure}[t]{0.47\textwidth}
\centering
\scriptsize
\input{Realisations_epsilon0p022143_n1280_Model1_1D_p1p5.tikz}
\caption{
We plot the functions output by the algorithm corresponding to nine realizations of the data for $n=1280$ and $\eps=0.022$ (marked in yellow in Figure (a)). 
}
\end{subfigure}
\medskip

\begin{subfigure}[t]{0.47\textwidth}
\centering
\scriptsize
\input{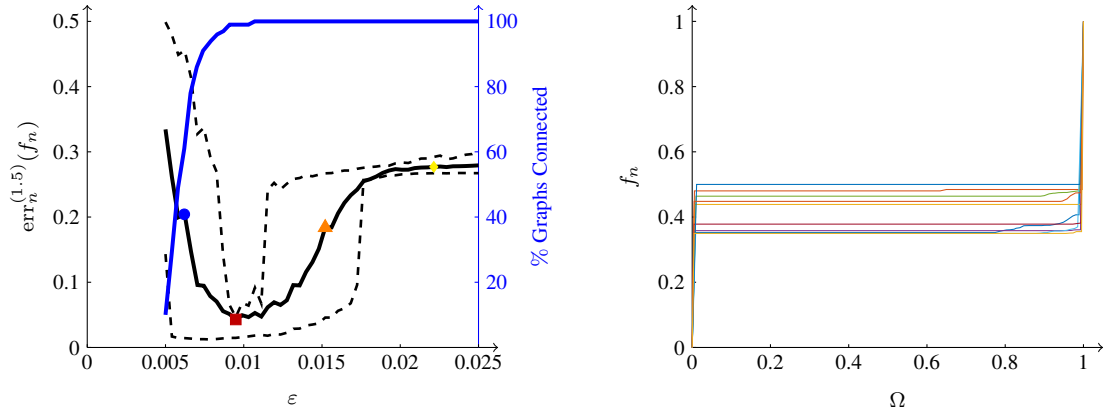}
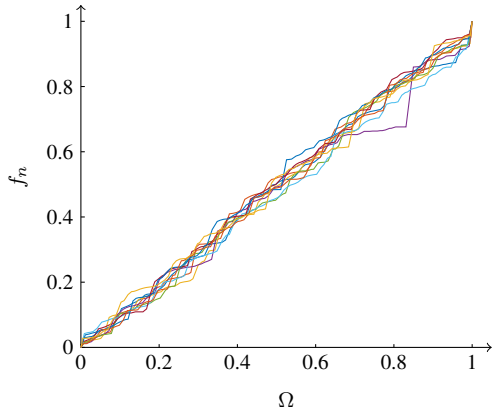
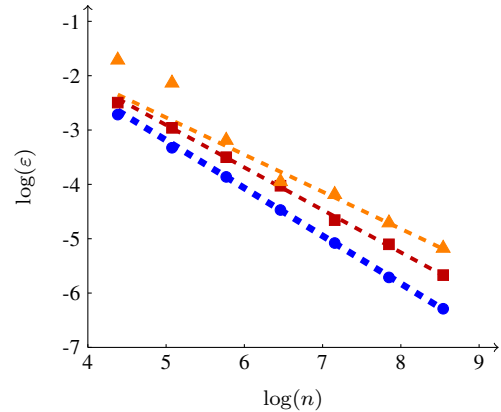
\caption{
The output of the algorithm, $f_n$, for nine realizations of the data for $n=1280$ and $\eps=\eps_{*}^{(1.5)}$.
}
\end{subfigure}
\hspace*{0.04\textwidth}
\begin{subfigure}[t]{0.47\textwidth}
\centering
\scriptsize
\input{ScalingInEpsilon_Model1_1D_p1p5.tikz}
\caption{
Orange triangles dots are $\eps_{\mathrm{upper}}^{(1.5)}$, red squares are  $\eps_*^{(1.5)}$, and blue dots are $\eps_{\mathrm{conn}}$. Lines show the linear fit over last 5 points.}
\end{subfigure}
\caption{
1D numerical experiments for \eqref{eq:MainRes:Encon} with $p=1.5$, averaged over 100 realizations.
\label{fig:Ex:1D:ErrorVEpsp=1.5}
}
\end{figure}

Let $\mu$ be the uniform measure on $[0,1]$ and consider $\eta$ defined by $\eta(t) = 1$ if $t\leq 1$ and $\eta(t) = 0$ otherwise.
We consider two different values of $p$: $p=1.5$ and $p=2$.
The training set is $\{(0,0),(1,1)\}$, that is we condition on functions $f_n$ taking the value 0 at $x_1=0$ and taking the value 1 at $x_2=1$ (so $N=2$).
We avoid using $p=1$ since any increasing function $f$ with $f(0)=0$ and $f(1)=1$ is a minimizer to the limiting problem.
For $p>1$ the solution to the constrained limiting problem is $f^\dagger(x) = x$ (note that this is independent of $p$).
Since $f^\dagger$ is continuous we can consider the following simple-to-compute notion of error:
\begin{equation} \label{err_def}
 \err^{(p)}_n(f_n) = \|f_n-f^\dagger \|_{L^p(\mu_n)}.
\end{equation}

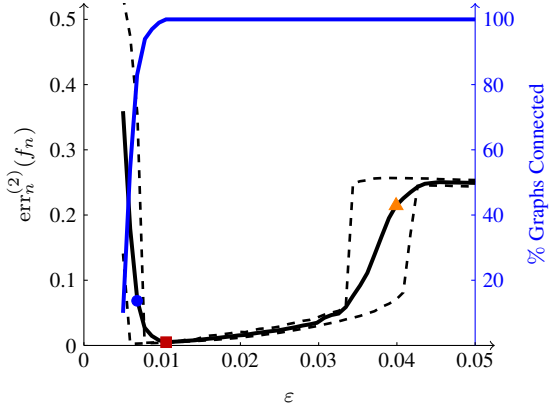
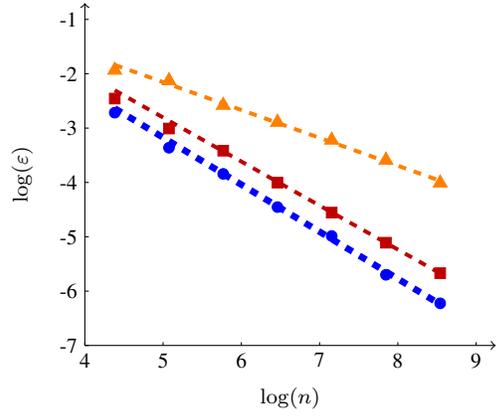
\begin{figure}[t!]
\setlength\figureheight{0.35\textwidth}
\setlength\figurewidth{0.47\textwidth}
\centering
\scriptsize

\begin{subfigure}[t]{0.47\textwidth}
\centering
\scriptsize
\input{Error_n1280_Model1_1D_p2.tikz}
\caption{
Error \eqref{err_def} for $n=1280$.
Black line is the mean error, dashed lines are the 10\% and 90\% quantiles.
Connectivity bound $\eps_{\mathrm{conn}}$ in blue, optimal $\eps_*^{(2)}$ in red and upper bound $\eps_{\mathrm{upper}}^{(2)}$ in orange.
Blue line is the percentage of connected graphs for a given $\eps$.
}
\end{subfigure}
\hspace*{0.04\textwidth}
%
\begin{subfigure}[t]{0.47\textwidth}
\centering
\scriptsize
\input{ScalingInEpsilon_Model1_1D_p2.tikz}
\caption{
Orange triangles dots are $\eps_{\mathrm{upper}}^{(2)}$, red squares are  $\eps_*^{(2)}$, and blue dots are $\eps_{\mathrm{conn}}$. Lines show the linear fit over last 5 points.
}
\end{subfigure}

\caption{
1D numerical experiments averaged over 100 realizations  for \eqref{eq:MainRes:Encon} with  $p=2$.
\label{fig:Ex:1D:ErrorVEpsp=2}
}
\end{figure}

To find minimizers of~\eqref{eq:MainRes:Encon} we use coordinate gradient descent.
The number of data points varies from $n = 80$ to $n=5120$. For each $n$, $\eps$ and $p$ we consider 100 different realizations of the random sample and plot the average results. 
When $\eps$ is too small the graph is disconnected and we should not expect informative solutions, when $\eps$ is large we expect discontinuities  to arise and  cause degeneracy.
In Figure~\ref{fig:Ex:1D:ErrorVEpsp=1.5}(a) and Figure~\ref{fig:Ex:1D:ErrorVEpsp=2}(a) we plot the error as a function of $\eps$ for fixed $n=1280$. 
We see clear regions where $\eps$ is too small and where $\eps$ is too large, with the intermediate range producing good estimators. Plots of minimizers for a particular $\eps$ in the ``large-$\eps$" region, show that they exhibit discontinuities, as expected. 

\begin{figure}[t!]
\centering
\scriptsize
\setlength\figureheight{0.35\textwidth}
\setlength\figurewidth{0.47\textwidth}
\begin{subfigure}[t]{0.47\textwidth}
\centering
\scriptsize
\input{SmoothEpsConnVError_n1280_Model1_1D_p1p5.tikz}
\caption{Error \eqref{err_def} for the output $f_n$ for $n=1280$ and $p=1.5$.
The solid line is the mean error, the dashed lines are the 10\% and 90\% quantiles.
}
\end{subfigure}
\hspace*{0.04\textwidth}
\begin{subfigure}[t]{0.47\textwidth}
\centering
\scriptsize
\input{SmoothEpsConnVError_n1280_Model1_1D_p2.tikz}
\caption{
Error \eqref{err_def} for the output $f_n$ for $n=1280$ and $p=2$.
The solid line is the mean error, the dashed lines are the 10\% and 90\% quantiles.
}
\end{subfigure}

\caption{
Error shifted by connectivity radius using the same results as in Figures~\ref{fig:Ex:1D:ErrorVEpsp=1.5} and~\ref{fig:Ex:1D:ErrorVEpsp=2}.
\label{fig:Ex:1D:ErrorVConnEps}
}
\end{figure}
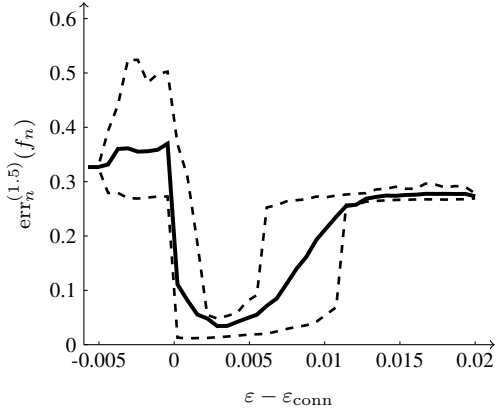
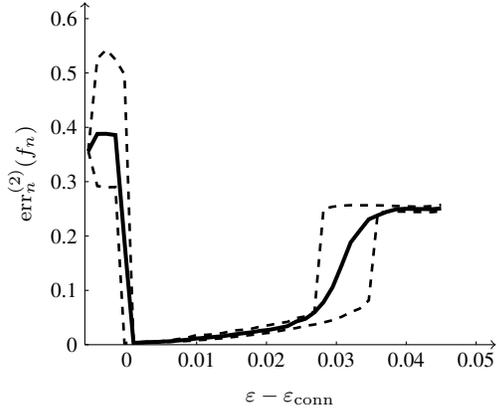

\begin{figure}[ht!]
\centering
\scriptsize
\setlength\figureheight{0.35\textwidth}
\setlength\figurewidth{0.47\textwidth}
\begin{subfigure}[t]{0.47\textwidth}
\centering
\scriptsize
\input{Error_n1280_Model2_1D_p2.tikz}
\caption{
Error \eqref{err_def} for $n=1280$.
Black line is the mean error, dashed lines are the 10\% and 90\% quantiles.
 Blue dot is the connectivity bound $\eps_{\mathrm{conn}}$. 
Blue line is the percentage of connected graphs  for given $\eps$.
}
\end{subfigure}
\hspace*{0.04\textwidth}
\begin{subfigure}[t]{0.47\textwidth}
\centering
\scriptsize
\input{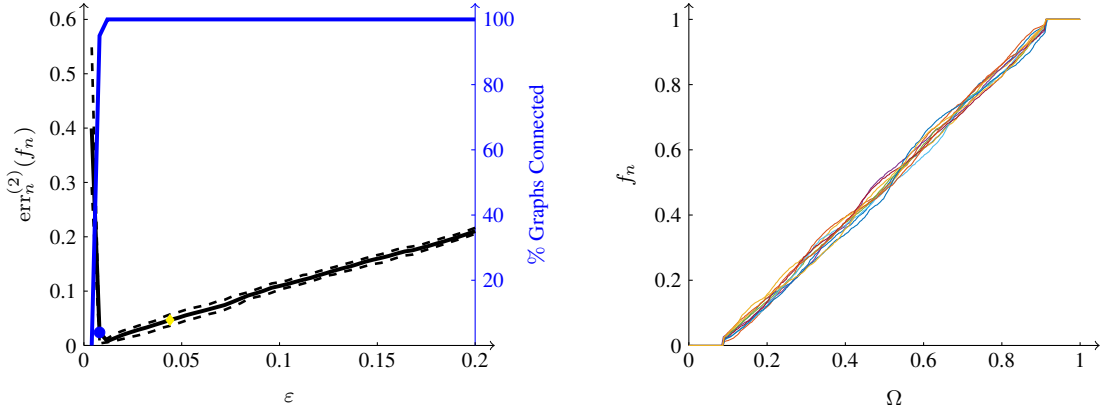}
\caption{
We plot the functions output from the algorithm corresponding to multiple realizations of the data for $n=1280$ and $\eps=0.045$ (marked in yellow in Figure (a)). 
}
\end{subfigure}

\caption{
1D numerical experiments averaged over 100 realizations for model \eqref{eq:cF} with $p=2$.
\label{fig:Ex:1D:ErrorVEpsModel2p=2}
}

\end{figure}

To measure how the transition point in $\eps$ where minimizers change behavior scale with $n$ we define the following:
\begin{itemize}
\item[(i)] Given a realization  $\{x_i^\omega\}_{i=1}^n$ let $\eps_{\mathrm{conn}}(n;\omega)$ be the connectivity radius for the particular realization, $\omega$ that is the smallest $\eps$ such that the graph with weights $W_{ij} =\eta_\eps(|x_i-x_j|)$ is connected. The value $\eps_{\mathrm{conn}}(n)= \frac{1}{M}\sum_{i=1}^M \eps_{\mathrm{conn}}(n;\omega_i)$ is the connectivity radius averaged over the realizations considered. We considered $M=100$ realizations.
\item[(ii)] $\eps_*^{(p)}(n)$ is the empirically best choice for $\eps$, namely the $\eps$ that minimizes $\err_n^{(p)}(f_n)$ where $f_n$ is the minimizer of~\eqref{eq:MainRes:Encon} with $\eps_n=\eps$; again averaged over $M=100$ realizations. 
\item[(iii)] $\eps_{\mathrm{upper}}^{(p)}(n)$ is the upper bound on $\eps$ for which the algorithm behaves well, which we  identify  as the maximizer of the second derivative of $-\err_n^{(p)}(f_n)$ with respect to $\eps$, among $\eps \geq \eps_*^{(p)}(n)$. While computing  $\eps_{\mathrm{upper}}^{(p)}(n)$ we smooth the error slightly so that the method is robust to small perturbations. As above the value is averaged over 100 realizations. 
\end{itemize}

All of these points are highlighted in Figure~\ref{fig:Ex:1D:ErrorVEpsp=1.5}(a)
and Figure~\ref{fig:Ex:1D:ErrorVEpsp=2}(a). 
In Figure~\ref{fig:Ex:1D:ErrorVEpsp=1.5}(d) and Figure~\ref{fig:Ex:1D:ErrorVEpsp=2}(d) we plot how these values of $\eps$ scale with $n$.
The best linear fit (based on five largest values of $n$) in the log-log domain gives the following scalings
\begin{align*}
\eps_*^{(1.5)} & \approx \frac{2.719}{n^{0.781}} & \eps_{\mathrm{upper}}^{(1.5)} & \approx \frac{1.905}{n^{0.683}} \\
\eps_*^{(2)} & \approx \frac{3.472}{n^{0.810}} & \eps_{\mathrm{upper}}^{(2)} & \approx \frac{1.507}{n^{0.513}} \\
\eps_{\mathrm{conn}} & \approx \frac{3.342}{n^{0.879}}.
\end{align*} 
We observe that asymptotic scaling established in Theorem \ref{thm:MainRes:Conv}  for $\eps_{\mathrm{upper}}^{(p)}$ is $\frac{1}{n^{0.5}}$ for $p=2$ and $\frac{1}{n^{0.667}}$ for $p=1.5$, which is very close to our numerical results.
The true scaling in the connectivity of the graph is $\frac{\log(n)}{n}$, our numerical results behave approximately as $\frac{1}{n^{0.879}}$. 
We note that if we use a linear fit in the log-log domain over the same range as considered above then $\frac{n}{\log(n)}$ is approximated by $\frac{1}{n^{0.859}}$.

\begin{figure}[ht!]
\setlength\figureheight{0.35\textwidth}
\setlength\figurewidth{0.47\textwidth}
\centering
\scriptsize
\begin{subfigure}[t]{0.47\textwidth}
\centering
\scriptsize
\input{Error_n1280_Model1_2D_p2.tikz} 
\caption{
Error \eqref{err_def2Ddeg}  for $n=1280$.
Black line is the mean error, dashed lines are the 10\% and 90\% quantiles.
 Blue dot is the connectivity bound $\eps_{\mathrm{conn}}$. 
Blue line is the percentage of connected graphs for given $\eps$.
}
\end{subfigure}
\hspace*{0.04\textwidth}
\setlength\figureheight{0.35\textwidth}
\setlength\figurewidth{0.43\textwidth}
\begin{subfigure}[t]{0.47\textwidth}
\centering
\scriptsize
\input{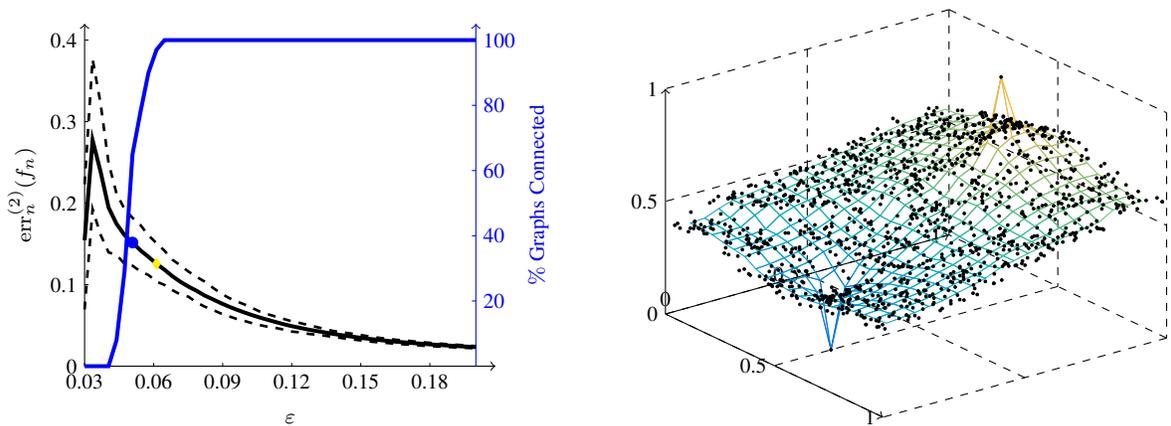}
\caption{
We plot an example of a function output from the algorithm corresponding to $n=1280$ and $\eps=0.06$ (marked in yellow in Figure (a)).
The grid is to aid visualisation.
}
\end{subfigure}
\caption{
2D numerical experiments averaged over 100 realizations for \eqref{eq:MainRes:Encon} with $p=2$.
\label{fig:Ex:2D:ErrorVEpsp=2}
}
\end{figure}

\begin{figure}[t!]
\setlength\figureheight{0.35\textwidth}
\setlength\figurewidth{0.47\textwidth}
\centering
\scriptsize
\begin{subfigure}[t]{0.47\textwidth}
\centering
\scriptsize
\input{Error_n1280_Model1_2D_p4.tikz} 
\caption{
Error \eqref{err_def2D} for $n=1280$.
Black line is the mean error, dashed lines are the 10\% and 90\% quantiles.
Blue dot is the connectivity bound $\eps_{\mathrm{conn}}$, red is $\eps_*^{(4)}$, and orange is the upper bound $\eps_{\mathrm{upper}}^{(4)}$. 
Blue line is the percentage of connected graphs for given $\eps$.
}
\end{subfigure}
\hspace*{0.04\textwidth}
%
\setlength\figureheight{0.35\textwidth}
\setlength\figurewidth{0.47\textwidth}
\begin{subfigure}[t]{0.47\textwidth}
\centering
\scriptsize
\input{ScalingInEpsilon_Model1_2D_p4.tikz}
\caption{
Orange triangles are $\eps_{\mathrm{upper}}^{(4)}$, red squares are  $\eps_*^{(4)}$, and blue dots are $\eps_{\mathrm{conn}}$. Lines show the linear fit over last 5 points.
}
\end{subfigure}
\caption{
2D numerical experiments averaged over 100 realizations for \eqref{eq:MainRes:Encon} with $p=4$.
\label{fig:Ex:2D:ErrorVEpsp=4}
}
\end{figure}
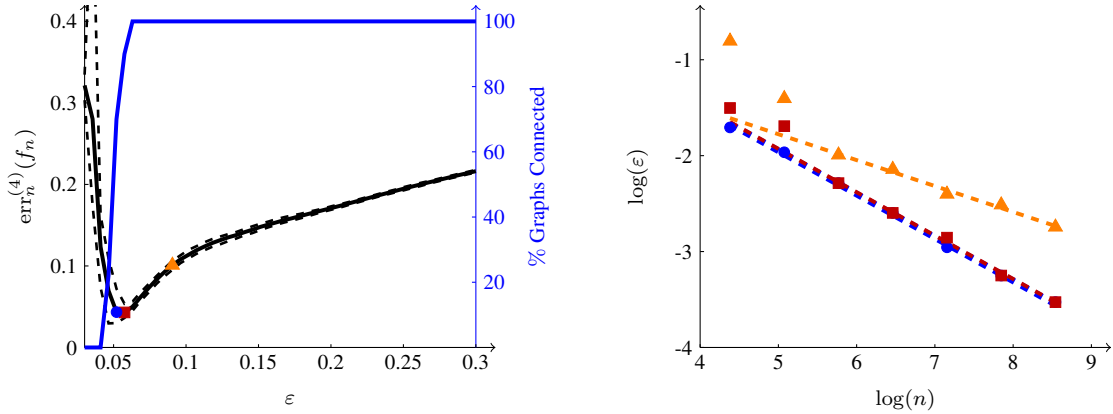

We observe that optimal choice $\eps_*^{(p)}$ is quite close to the connectivity radius $\eps - \eps_{\mathrm{conn}}(n)$. Choosing $\eps$ smaller than this results in a big error from a small number of realizations. 
To further investigate the proximity of the connectivity radius and the optimal choice of $\eps$  we plot in Figure~\ref{fig:Ex:1D:ErrorVConnEps} the error as the function of the size of $\eps$ relative to the connectivity radius.  More precisely 
we consider $\err_n^{(p)}(f_n; \eps, \omega)$, where $f_n$ is the minimizer of~\eqref{eq:MainRes:Encon} for given $\eps$ and realization $\omega$,  as a function of  $\eps - \eps_{\mathrm{conn}}(n;\omega)$ and then average over $M=100$ realizations. 
We observe that, for both $p=1.5$ and $p=2$, the error is the smallest when $\eps$ is quite close to the connectivity radius. The slight difference is that for $p=1.5$ there is a short interval beyond the connectivity radius where the error is still decreasing. 
\begin{remark} \label{rem:61} The close proximity of the optimal epsilon to the connectivity radius, both for the original model and the improved model (Figure \ref{fig:Ex:1D:ErrorVEpsModel2p=2}[a]) and both in 1D and 2D (Figure~\ref{fig:Ex:2D:ErrorVEpsp=4}[a]),  is not obvious since for $\eps$ small (i.e. relatively close to the connectivity radius) $\cE_n^{(p)}(f)$ is a poor approximation of $\Ecp(f)$, even for $f$ a fixed smooth function. Explaining the observed behavior of the error is an interesting open problem, that we believe should be approached from the viewpoint to stochastic homogenization.  
\end{remark} 

The improved model \eqref{eq:cF}, for which we show results in Figure~\ref{fig:Ex:1D:ErrorVEpsModel2p=2}, is far more robust to the choice of $\eps$.
We plot the error as a function of $\eps$ for $n=1280$ and we see a much larger range in the admissible choices of $\eps$.
To highlight the difference we plot in Figure~\ref{fig:Ex:1D:ErrorVEpsModel2p=2}(b) outputs from multiple realizations of the data under the same conditions as for Figure~\ref{fig:Ex:1D:ErrorVEpsp=2}(b), in particular we use the same choice of $\eps$. Note that the horizontal axis covers a much larger range on Figure~\ref{fig:Ex:1D:ErrorVEpsModel2p=2}(b).
The comparison shows that model \eqref{eq:MainRes:Encon} does not produce a reasonable output when $\eps \gtrsim 0.04$, while all outputs of \eqref{eq:cF} (when $\eps$ is larger than the connectivity radius) are close to the truth.

\subsection{2D Numerical Experiments \label{subsec:Ex:2D}}

Let $\mu$ be the uniform measure on $\Omega=[0,1]\times [0,1]$, and $\eta(t)=1$ if $|t|\leq 1$, $\eta(t)=0$ otherwise.
In 2D the critical value of $p$ is $p=2$, and we therefore choose to investigate $p=2$ and $p=4$.
The training set is $x_1 = (0.2,0.5)$, $x_2 = (0.8,0.5)$, with labels $y_1 = 0$, $y_2 = 1$.
In contrast to the 1D example the solution to the continuum problem~\eqref{eq:MainRes:Einftycon} (in the well-posed regime) depends on $p$ and furthermore cannot be solved analytically.
To estimate the solution we discretised~\eqref{eq:MainRes:Einftycon} on a uniform grid and ran a gradient descent algorithm to approximate the minimiser.
In the case when $p=4$ we plot our numerical approximation of the continuum minimizer to~\eqref{eq:MainRes:Einftycon} in Figure~\ref{fig:intro}(b).
For $p>2$ we define the error by 
\begin{equation} \label{err_def2D}
\err^{(p)}_n(f_n) = \|f_n-f^{p,\dagger} \|_{L^p(\mu_n)} 
\end{equation}
where $f^{p,\dagger}$  minimizes~\eqref{eq:MainRes:Einftycon}.
In the ill-posed case ($p\leq 2$) any constant function is a minimizer to the continuum problem, in which case we define the error as 
\begin{equation} \label{err_def2Ddeg}
\err^{(p)}_n(f_n) = \inf_{c\in \bbR} \|f_n-c \|_{L^p(\mu_n)}.
\end{equation}

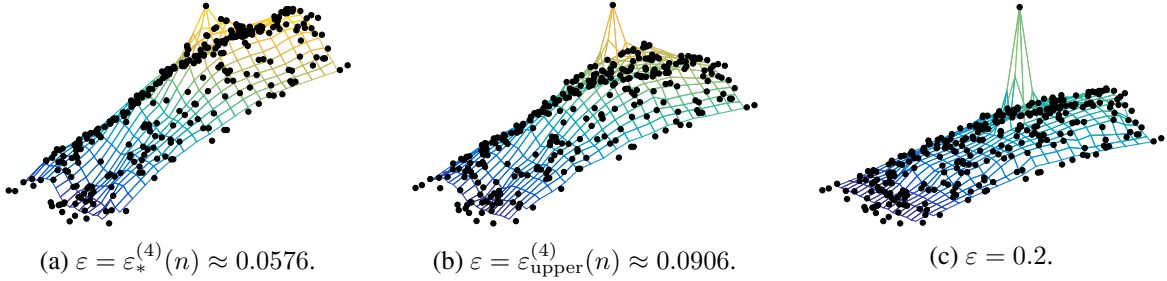
\begin{figure}[t!]
\setlength\figureheight{0.2442\textwidth}
\setlength\figurewidth{0.3\textwidth}
\centering
\scriptsize
\begin{subfigure}[t]{0.3\textwidth}
\centering
\scriptsize
\input{Realisation_Zoom_epsilon_0p0576_n1280_Model1_2D.tikz}
\caption{
$\eps = \eps_*^{(4)}(n) \approx 0.0576$.
}
\end{subfigure}
\hspace*{0.04\textwidth}
\begin{subfigure}[t]{0.3\textwidth}
\centering
\scriptsize
\input{Realisation_Zoom_epsilon_0p0906_n1280_Model1_2D.tikz}
\caption{
$\eps = \eps^{(4)}_{\mathrm{upper}}(n) \approx 0.0906$.
}
\end{subfigure}
\hspace*{0.04\textwidth}
\begin{subfigure}[t]{0.3\textwidth}
\centering
\scriptsize
\input{Realisation_Zoom_epsilon_0p2_n1280_Model1_2D.tikz}
\caption{
$\eps = 0.2$.
}
\end{subfigure}

%

\caption{
Realizations of~\eqref{eq:MainRes:Encon} with $p=4$ and $n=1280$ for a select choices of $\eps$. Only the part of the domain near labeled point $x_2$ is shown.
The grids are to aid visualization.
\label{fig:Ex:2D:p=4Spikes}
}
\end{figure}

To find minimizers of~\eqref{eq:MainRes:Encon} for $p=4$ we use coordinate gradient descent. 
For $p=2$ we use the method of~\cite{ZhGhLa03} that exactly solves the Euler Lagrange equation ($(L_nf)_i = 0$, where $L_n$ is the graph Laplacian, for $i>2$ and $f_1=0$, $f_2=1$).
The number of data points varies from 80 to 5120. We use 100 different realizations of the data $\{x_i\}_{i=1}^n$ for each $n$ and each choice of $\eps$.

Figure~\ref{fig:Ex:2D:ErrorVEpsp=2} shows the results in the ill-posed regime, for $p=2$ and $n=1280$.
We observe that the solutions form spikes in order to satisfy the constraints. 
Spikes are present for all $\eps$ beyond the connectivity threshold, and grow as $\eps$
increases (recall that the solution to the continuum problem is a constant and therefore the error decreasing indicates convergence to a constant solution with spikes around the two training data points).

In Figures~\ref{fig:intro} and \ref{fig:Ex:2D:ErrorVEpsp=4} show the results for $p=4$ which is in the well-posed range. Figure~\ref{fig:intro}(a) presents the numerically computed discrete minimizer for the optimal radius $\eps=\eps^{(4)}_*$. 
We observe in Figure~\ref{fig:Ex:2D:ErrorVEpsp=4}(a) that, similarly to 1D, for $\eps$ below the average connectivity radius, or $\eps$ is large, the error is high, and is the lowest for $\eps$ in between and close to connectivity threshold. In contrast to the 1D results we do notice that the transitions between the well-posed and ill-posed regime is gradual.

The numerical scaling for $\eps_*^{(p)}$, $\eps_{\mathrm{upper}}^{(p)}$, and $\eps_{\mathrm{conn}}$  (with the same definitions of quantities as in the 1D experiments in the previous subsection) we find is
\begin{align*}
\eps_*^{(4)} & \approx \frac{1.394}{n^{0.452}} & \eps_{\mathrm{upper}}^{(4)} & \approx \frac{0.654}{n^{0.270}} & \eps_{\mathrm{conn}} & \approx \frac{1.368}{n^{0.452}}.
\end{align*} 
The connectivity radius should scale according to $\sqrt{\frac{\log(n)}{n}}$, which is close to our observed rate of $n^{-0.452}$ (in fact when linear fitting $\log\eps$ to $\frac{\log n}{n}$ one obtains $\eps_{\mathrm{conn}} \approx 0.829 \l\frac{\log n}{n}\r^{0.526}$).
Our theoretical predictions give the scaling of the upper bound as $\eps_{\mathrm{upper}}^{(4)} \asymp n^{-0.25}$, close to our numerical rate of $n^{-0.270}$. 

\begin{figure}[t!]
\centering
\scriptsize
\setlength\figureheight{0.35\textwidth}
\setlength\figurewidth{0.47\textwidth}

\begin{subfigure}[b]{0.47\textwidth}
\centering
\scriptsize
\input{SmoothEpsConnVError_2DModel1_n1280_p4.tikz}
\caption{
Error \eqref{err_def2D} using the same results as in Figure~\ref{fig:Ex:2D:ErrorVEpsp=4}.
The solid line is the mean error, the dashed lines are the 10\% and 90\% quantiles.
}
\end{subfigure}
\hspace*{0.04\textwidth}
\begin{subfigure}[b]{0.47\textwidth}
\centering
\scriptsize
\includegraphics[width=\textwidth]{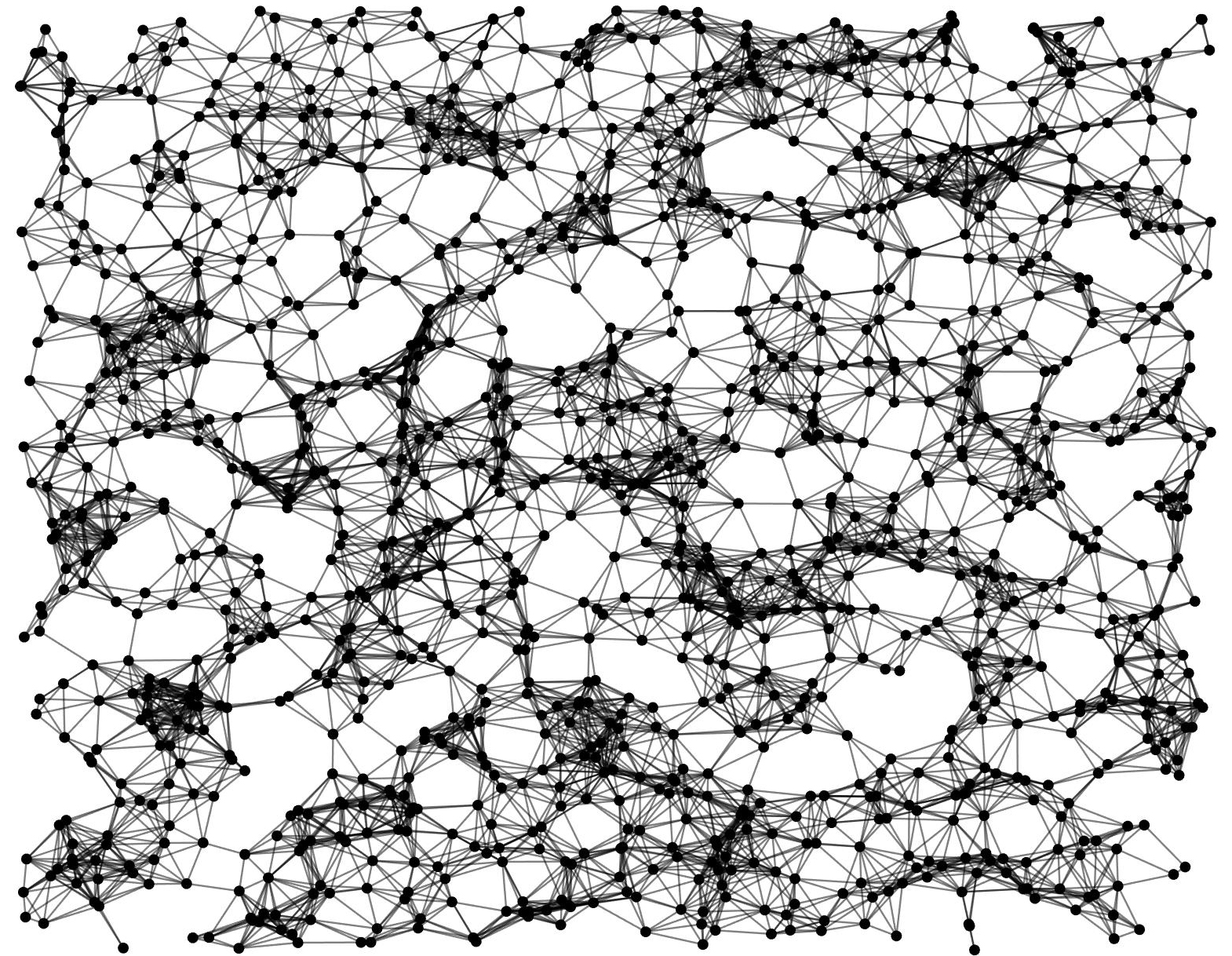}
\caption{
Example graph in 2D for $\eps=\eps_*^{(4)}(n)$ and $n=1280$.
\vspace{\baselineskip} 
}
\end{subfigure}

\caption{
Error dependency on the connectivity radius and the graph for optimal $\eps$, for $n=1280$ and $p=4$.
\label{fig:Ex:2D:ErrorVConnEps}
}
\end{figure}

\begin{figure}[ht!]
\centering
\scriptsize
\setlength\figureheight{0.35\textwidth}
\setlength\figurewidth{0.47\textwidth}
\begin{subfigure}[b]{0.47\textwidth}
\centering
\scriptsize
\input{Error_n1280_Model2_2D_p4.tikz}
\caption{
The black line is the mean error \eqref{err_def2D} for model \eqref{eq:MainRes:Encon}.
The error for the improved model \eqref{eq:cF} with 
constraint radius $R_n=2\eps$ is in yellow, $R_n=\eps$ is in orange, and $R_n=\eps/2$ is in red.
}
\end{subfigure}
\hspace*{0.04\textwidth}
\setlength\figureheight{0.35\textwidth}
\setlength\figurewidth{0.43\textwidth}
\begin{subfigure}[b]{0.47\textwidth}
\centering
\scriptsize
\input{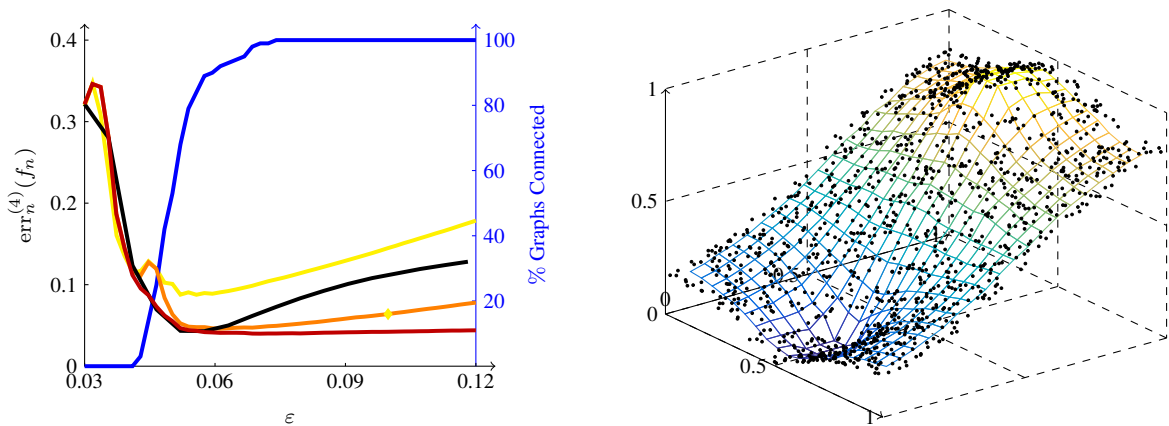}
\caption{
We plot an example of a function output from the algorithm corresponding to $n=1280$ and $\eps=0.1$ for the constraint set of size $\eps$ (marked in yellow in Figure (a)).
The grid is to aid visualisation. 
}
\end{subfigure}

\caption{
Experiments for improved  model~\eqref{eq:cF} with $n=1280$ and $p=4$ 
averaged over 100 realizations.
\label{fig:Ex:2D:ErrorVEpsModel2p=4}
}
\end{figure}

In Figure~\ref{fig:Ex:2D:p=4Spikes} we show instances of numerically computed minimizers of~\eqref{eq:MainRes:Encon} for increasing values of $\eps$. They show that the breakdown of the numerical approximation of the continuum solution (shown in Figure~\ref{fig:Ex:2D:ErrorVEpsModel2p=4}(c)) happens via development of spikes.

As in the 1D examples (shown on Figure \ref{fig:Ex:1D:ErrorVConnEps}) we investigate 
the proximity of the optimal radius $\eps_*^{(p)}(n;\omega)$ to the connectivity radius 
$\eps_{conn}(n;\omega)$, where $\omega$ is the sample considered.
In Figure~\ref{fig:Ex:2D:ErrorVConnEps} we plot the error, $\err_n^{(p)}(f_n;\omega)$, against $\eps-\eps_{\mathrm{conn}}(n;\omega)$ for $n=1280$ and $p=4$, averaged over $100$ samples.
The phenomena we observe is similar to the 1D case; the error is large and highly variable for $\eps$ below the connectivity radius.
There is a sharp transition to the well-posed regime, as soon as the graph is connected with the error increasing with $\eps$. As we explain in Remark \ref{rem:61} it is an intriguing and important open problem to explain why the error is the smallest for rather coarse graphs (Figure \ref{fig:Ex:2D:ErrorVConnEps}(b)).

Our theoretical result in Section~\ref{sec:IM} showed that minimizers of the improved model, \eqref{eq:cF}, converge as $n \to \infty$ to the correct solution if $1 \gg \eps_n  \gg \left(\ln n/n \right)^{1/d}$, regardless of how slowly $\eps_n \to 0$. Here we numerically  investigate two issues. One is how precisely does the error of the improved model depend on $\eps$ for fixed $n$. The other is to compare the observed error of the improved model to the original model. 
Recall that in the improved model we prove convergence when the labels are extended around the training set to balls of radius $2 \eps$. This is needed in our proof to ensure that spikes do not form. Here we numerically investigate if extending the labels to smaller balls is sufficient to prevent the spike formation. In particular in Figure~\ref{fig:Ex:2D:ErrorVEpsModel2p=4}(a) we display the error for fixed $n=1280$ and constraint ball radii $2 \eps$, $\eps$ and $\eps/2$. The numerics show that even radius $\eps/2$ is sufficient to prevent spike formation and that it allows for better approximation of the continuum solution. We also observe that fixing the labels on larger sets can
significantly impact the accuracy of approximation. This issue is less pronounced for larger values of $n$, where the connectivity radius is small compared to distances between the labeled points. 


\subsection*{Acknowledgements}

The authors thank Matt Dunlop and Andrew Stuart for enlightening exchanges.
The authors are grateful to Nicol\'as Garc\'ia Trillos for careful reading of the manuscript and insightful remarks.
This material is based on work supported by the National Science Foundation under the  grants CCT 1421502 and DMS 1516677.
The authors are also grateful to the Center for Nonlinear Analysis (CNA) for support.
MT is grateful for the support of the Cantab Capital Institute for the Mathematics of Information at the University of Cambridge.

\bibliographystyle{plain}
\bibliography{ref_plap}

\end{document}

%% file: ContSol_2D_p4.tikz
\definecolor{mycolor1}{rgb}{0,0.447,0.741}

\begin{tikzpicture}

\begin{axis}[%
width=\figurewidth,
height=\figureheight,
unbounded coords=jump,
clip=false,
view={-37.5}{30},
scale only axis,
every outer x axis line/.append style={darkgray!60!black},
every x tick label/.append style={font=\color{darkgray!60!black}},
xmin=-0,
xmax=1,
xmajorgrids,
every outer y axis line/.append style={darkgray!60!black},
every y tick label/.append style={font=\color{darkgray!60!black}},
ymin=0,
ymax=1,
ymajorgrids,
every outer z axis line/.append style={darkgray!60!black},
every z tick label/.append style={font=\color{darkgray!60!black}},
zmin=0,
zmax=1,
zmajorgrids,
hide axis,
grid style={solid},
]
\draw [->] (axis cs: 0,1,0) -- (axis cs: 1,1,0);
\draw [->] (axis cs: 0,1,0) -- (axis cs: 0,0,0);
\draw [->] (axis cs: 0,1,0) -- (axis cs: 0,1,1);
\draw [dashed,-] (axis cs: 1,0.5,0) -- (axis cs: 0,0.5,0) node[left] {0.5};
\draw [dashed,-] (axis cs: 1,0,0) -- (axis cs: 0,0,0) node[left] {1};
\draw [dashed,-] (axis cs: 1,0,0) -- (axis cs: 1,1,0) node[left] {1};
\draw [dashed,-] (axis cs: 0.5,0,0) -- (axis cs: 0.5,1,0) node[left] {0.5};
\draw [dashed,-] (axis cs: 1,1,1) -- (axis cs: 1,1,0); 
\draw [dashed,-] (axis cs: 0.5,1,1) -- (axis cs: 0.5,1,0); 
\draw [dashed,-] (axis cs: 1,1,0) -- (axis cs: 0,1,0) node[left] {0};
\draw [dashed,-] (axis cs: 0,0,0) -- (axis cs: 0,1,0) node[above] {0};
\draw [dashed,-] (axis cs: 1,1,0.5) -- (axis cs: 0,1,0.5) node[left] {0.5};
\draw [dashed,-] (axis cs: 1,1,1) -- (axis cs: 0,1,1) node[left] {1};
\draw [dashed,-] (axis cs: 1,1,1) -- (axis cs: 1,0,1); 
\draw [dashed,-] (axis cs: 1,1,0.5) -- (axis cs: 1,0,0.5); 
\draw [dashed,-] (axis cs: 1,0.5,1) -- (axis cs: 1,0.5,0); 
\draw [dashed,-] (axis cs: 1,0,1) -- (axis cs: 1,0,0); 

\addplot3[%
surf,
colormap={mymap}{[1pt] rgb(0pt)=(0.2081,0.1663,0.5292); rgb(1pt)=(0.211624,0.189781,0.577676); rgb(2pt)=(0.212252,0.213771,0.626971); rgb(3pt)=(0.2081,0.2386,0.677086); rgb(4pt)=(0.195905,0.264457,0.7279); rgb(5pt)=(0.170729,0.291938,0.779248); rgb(6pt)=(0.125271,0.324243,0.830271); rgb(7pt)=(0.0591333,0.359833,0.868333); rgb(8pt)=(0.0116952,0.38751,0.881957); rgb(9pt)=(0.00595714,0.408614,0.882843); rgb(10pt)=(0.0165143,0.4266,0.878633); rgb(11pt)=(0.0328524,0.443043,0.871957); rgb(12pt)=(0.0498143,0.458571,0.864057); rgb(13pt)=(0.0629333,0.47369,0.855438); rgb(14pt)=(0.0722667,0.488667,0.8467); rgb(15pt)=(0.0779429,0.503986,0.838371); rgb(16pt)=(0.0793476,0.520024,0.831181); rgb(17pt)=(0.0749429,0.537543,0.826271); rgb(18pt)=(0.0640571,0.556986,0.823957); rgb(19pt)=(0.0487714,0.577224,0.822829); rgb(20pt)=(0.0343429,0.596581,0.819852); rgb(21pt)=(0.0265,0.6137,0.8135); rgb(22pt)=(0.0238905,0.628662,0.803762); rgb(23pt)=(0.0230905,0.641786,0.791267); rgb(24pt)=(0.0227714,0.653486,0.776757); rgb(25pt)=(0.0266619,0.664195,0.760719); rgb(26pt)=(0.0383714,0.674271,0.743552); rgb(27pt)=(0.0589714,0.683757,0.725386); rgb(28pt)=(0.0843,0.692833,0.706167); rgb(29pt)=(0.113295,0.7015,0.685857); rgb(30pt)=(0.145271,0.709757,0.664629); rgb(31pt)=(0.180133,0.717657,0.642433); rgb(32pt)=(0.217829,0.725043,0.619262); rgb(33pt)=(0.258643,0.731714,0.595429); rgb(34pt)=(0.302171,0.737605,0.571186); rgb(35pt)=(0.348167,0.742433,0.547267); rgb(36pt)=(0.395257,0.7459,0.524443); rgb(37pt)=(0.44201,0.748081,0.503314); rgb(38pt)=(0.487124,0.749062,0.483976); rgb(39pt)=(0.530029,0.749114,0.466114); rgb(40pt)=(0.570857,0.748519,0.44939); rgb(41pt)=(0.609852,0.747314,0.433686); rgb(42pt)=(0.6473,0.7456,0.4188); rgb(43pt)=(0.683419,0.743476,0.404433); rgb(44pt)=(0.71841,0.741133,0.390476); rgb(45pt)=(0.752486,0.7384,0.376814); rgb(46pt)=(0.785843,0.735567,0.363271); rgb(47pt)=(0.818505,0.732733,0.34979); rgb(48pt)=(0.850657,0.7299,0.336029); rgb(49pt)=(0.882433,0.727433,0.3217); rgb(50pt)=(0.913933,0.725786,0.306276); rgb(51pt)=(0.944957,0.726114,0.288643); rgb(52pt)=(0.973895,0.731395,0.266648); rgb(53pt)=(0.993771,0.745457,0.240348); rgb(54pt)=(0.999043,0.765314,0.216414); rgb(55pt)=(0.995533,0.786057,0.196652); rgb(56pt)=(0.988,0.8066,0.179367); rgb(57pt)=(0.978857,0.827143,0.163314); rgb(58pt)=(0.9697,0.848138,0.147452); rgb(59pt)=(0.962586,0.870514,0.1309); rgb(60pt)=(0.958871,0.8949,0.113243); rgb(61pt)=(0.959824,0.921833,0.0948381); rgb(62pt)=(0.9661,0.951443,0.0755333); rgb(63pt)=(0.9763,0.9831,0.0538)},
shader=faceted interp,
faceted color=black,
mesh/rows=21]
table[row sep=crcr,header=false] {
0 0 0.164524389157533\\
0 0.05 0.155428287480189\\
0 0.1 0.139552323737823\\
0 0.15 0.120921473801775\\
0 0.2 0.10083388913566\\
0 0.25 0.07988581814176\\
0 0.3 0.0586544734647471\\
0 0.35 0.0380262985772315\\
0 0.4 0.0194655094840467\\
0 0.45 0.00546564800457887\\
0 0.5 4.82790605547743e-06\\
0 0.55 0.00544408625524704\\
0 0.6 0.0194523175235128\\
0 0.65 0.0380246018605653\\
0 0.7 0.0586589308018936\\
0 0.75 0.0798931531138446\\
0 0.8 0.100842734721718\\
0 0.85 0.120932258270878\\
0 0.9 0.139569361481521\\
0 0.95 0.155469637053011\\
0 1 0.164495863989141\\
0.05 0 0.178138152947402\\
0.05 0.05 0.169059378770813\\
0.05 0.1 0.152989064623178\\
0.05 0.15 0.134105548721012\\
0.05 0.2 0.113650413923588\\
0.05 0.25 0.09218701769275\\
0.05 0.3 0.0702211753187785\\
0.05 0.35 0.0484547211564335\\
0.05 0.4 0.0279153484163484\\
0.05 0.45 0.0102430029288555\\
0.05 0.5 6.03703491288344e-05\\
0.05 0.55 0.0102466898083719\\
0.05 0.6 0.0279212791729663\\
0.05 0.65 0.048458663431875\\
0.05 0.7 0.0702247100505089\\
0.05 0.75 0.0921903302446376\\
0.05 0.8 0.113652958172921\\
0.05 0.85 0.134107041577211\\
0.05 0.9 0.152989421242292\\
0.05 0.95 0.169057475662356\\
0.05 1 0.178102417679454\\
0.1 0 0.200007100387982\\
0.1 0.05 0.191659711067635\\
0.1 0.1 0.176023223732355\\
0.1 0.15 0.157323214504165\\
0.1 0.2 0.136855938641845\\
0.1 0.25 0.115182426530185\\
0.1 0.3 0.0927531606791931\\
0.1 0.35 0.0700991549973339\\
0.1 0.4 0.0477525551134435\\
0.1 0.45 0.0263662797794031\\
0.1 0.5 0.0110954770446266\\
0.1 0.55 0.0263021099613984\\
0.1 0.6 0.0477381183980303\\
0.1 0.65 0.0700929854496628\\
0.1 0.7 0.0927469743416953\\
0.1 0.75 0.115177256170145\\
0.1 0.8 0.136852145825679\\
0.1 0.85 0.157319968425095\\
0.1 0.9 0.176019549753227\\
0.1 0.95 0.191654701745569\\
0.1 1 0.199988979696443\\
0.15 0 0.226799787415329\\
0.15 0.05 0.219245699465908\\
0.15 0.1 0.204429510626382\\
0.15 0.15 0.186430529527188\\
0.15 0.2 0.166601131223144\\
0.15 0.25 0.145546819025606\\
0.15 0.3 0.123721155926359\\
0.15 0.35 0.101411580053452\\
0.15 0.4 0.0781137610731441\\
0.15 0.45 0.0521777994642257\\
0.15 0.5 0.0290440593529524\\
0.15 0.55 0.0518685273656458\\
0.15 0.6 0.0779867803382731\\
0.15 0.65 0.101367465875402\\
0.15 0.7 0.123704360654157\\
0.15 0.75 0.145536196035612\\
0.15 0.8 0.16659269092319\\
0.15 0.85 0.186423440396853\\
0.15 0.9 0.204422697505775\\
0.15 0.95 0.219238413154391\\
0.15 1 0.226786056476295\\
0.2 0 0.25757658133229\\
0.2 0.05 0.250904468245379\\
0.2 0.1 0.237322910585246\\
0.2 0.15 0.220631131632673\\
0.2 0.2 0.202217374690251\\
0.2 0.25 0.182730623561849\\
0.2 0.3 0.162563892956579\\
0.2 0.35 0.141520038489461\\
0.2 0.4 0.117244737030394\\
0.2 0.45 0.0801700981285665\\
0.2 0.5 5.6769941765632e-32\\
0.2 0.55 0.07870414626306\\
0.2 0.6 0.116775238064686\\
0.2 0.65 0.141390153109254\\
0.2 0.7 0.162526429090338\\
0.2 0.75 0.182716290588015\\
0.2 0.8 0.202207853368094\\
0.2 0.85 0.220622903416795\\
0.2 0.9 0.237314983079423\\
0.2 0.95 0.250896564143481\\
0.2 1 0.257565199488038\\
0.25 0 0.291910234154372\\
0.25 0.05 0.286215112797249\\
0.25 0.1 0.274272731490467\\
0.25 0.15 0.259481263548746\\
0.25 0.2 0.24320227831819\\
0.25 0.25 0.226072274970022\\
0.25 0.3 0.208358209380241\\
0.25 0.35 0.189466535026833\\
0.25 0.4 0.166488279763229\\
0.25 0.45 0.1354570108007\\
0.25 0.5 0.112934567173873\\
0.25 0.55 0.135351857200025\\
0.25 0.6 0.166084613021653\\
0.25 0.65 0.189285387585844\\
0.25 0.7 0.208301595880815\\
0.25 0.75 0.226054432410252\\
0.25 0.8 0.243193669379802\\
0.25 0.85 0.259474285812781\\
0.25 0.9 0.274265822201594\\
0.25 0.95 0.286208365803638\\
0.25 1 0.291901581165803\\
0.3 0 0.329408658631372\\
0.3 0.05 0.324770520349572\\
0.3 0.1 0.314817905240229\\
0.3 0.15 0.30243095851275\\
0.3 0.2 0.288852734903451\\
0.3 0.25 0.274648251093896\\
0.3 0.3 0.259986194079808\\
0.3 0.35 0.244322579778519\\
0.3 0.4 0.226315426702429\\
0.3 0.45 0.207466919339926\\
0.3 0.5 0.198097777547134\\
0.3 0.55 0.207422261715544\\
0.3 0.6 0.226165957816724\\
0.3 0.65 0.244192322704462\\
0.3 0.7 0.25993078984593\\
0.3 0.75 0.274629884718113\\
0.3 0.8 0.288845768288016\\
0.3 0.85 0.302426303015486\\
0.3 0.9 0.314813359762192\\
0.3 0.95 0.324766193128896\\
0.3 1 0.32940323777732\\
0.35 0 0.369586334826623\\
0.35 0.05 0.366064177048625\\
0.35 0.1 0.358376324348668\\
0.35 0.15 0.34878014324356\\
0.35 0.2 0.338304910741397\\
0.35 0.25 0.327405727221796\\
0.35 0.3 0.316212549598215\\
0.35 0.35 0.304474588562278\\
0.35 0.4 0.292056260217923\\
0.35 0.45 0.280999696675369\\
0.35 0.5 0.276258188345146\\
0.35 0.55 0.280958172374221\\
0.35 0.6 0.291986907255494\\
0.35 0.65 0.30440207359106\\
0.35 0.7 0.316172212115969\\
0.35 0.75 0.32739052058175\\
0.35 0.8 0.338299911525985\\
0.35 0.85 0.348777836374386\\
0.35 0.9 0.358374484626157\\
0.35 0.95 0.366062717363599\\
0.35 1 0.369584182555616\\
0.4 0 0.41185800161118\\
0.4 0.05 0.409491029036563\\
0.4 0.1 0.404264026047661\\
0.4 0.15 0.397728571115087\\
0.4 0.2 0.390619312057696\\
0.4 0.25 0.383256385796521\\
0.4 0.3 0.375747541163943\\
0.4 0.35 0.368054792752322\\
0.4 0.4 0.360417766942454\\
0.4 0.45 0.354204079259775\\
0.4 0.5 0.351711159016343\\
0.4 0.55 0.354174862627817\\
0.4 0.6 0.360375513372257\\
0.4 0.65 0.368013010691855\\
0.4 0.7 0.37572122090428\\
0.4 0.75 0.383245427014034\\
0.4 0.8 0.390616171987689\\
0.4 0.85 0.397728110361196\\
0.4 0.9 0.404264487825166\\
0.4 0.95 0.409492127094168\\
0.4 1 0.411858611954906\\
0.45 0 0.455564026760633\\
0.45 0.05 0.454375249278061\\
0.45 0.1 0.451732497554749\\
0.45 0.15 0.448425367073331\\
0.45 0.2 0.444836253793135\\
0.45 0.25 0.441131324463866\\
0.45 0.3 0.437375205458634\\
0.45 0.35 0.433591575132667\\
0.45 0.4 0.429970280760088\\
0.45 0.45 0.427152541063914\\
0.45 0.5 0.426053297400218\\
0.45 0.55 0.427131533359829\\
0.45 0.6 0.429940056764879\\
0.45 0.65 0.433563326746258\\
0.45 0.7 0.437356925951472\\
0.45 0.75 0.441123414092478\\
0.45 0.8 0.444834394050412\\
0.45 0.85 0.448426045371822\\
0.45 0.9 0.451734450155919\\
0.45 0.95 0.454378069238012\\
0.45 1 0.455566465390531\\
0.5 0 0.499998462027992\\
0.5 0.05 0.499998288634717\\
0.5 0.1 0.499998768212721\\
0.5 0.15 0.499999470397903\\
0.5 0.2 0.500000706203728\\
0.5 0.25 0.5000034292239\\
0.5 0.3 0.500007892485068\\
0.5 0.35 0.500012206794989\\
0.5 0.4 0.500013287321228\\
0.5 0.45 0.500009162165246\\
0.5 0.5 0.5\\
0.5 0.55 0.499990837834751\\
0.5 0.6 0.499986712678773\\
0.5 0.65 0.499987793205013\\
0.5 0.7 0.499992107514934\\
0.5 0.75 0.4999965707761\\
0.5 0.8 0.499999293796276\\
0.5 0.85 0.500000529602099\\
0.5 0.9 0.50000123178728\\
0.5 0.95 0.500001711365282\\
0.5 1 0.500001537972009\\
0.55 0 0.544433534609471\\
0.55 0.05 0.54562193076199\\
0.55 0.1 0.548265549844082\\
0.55 0.15 0.551573954628174\\
0.55 0.2 0.555165605949592\\
0.55 0.25 0.558876585907523\\
0.55 0.3 0.562643074048527\\
0.55 0.35 0.566436673253741\\
0.55 0.4 0.570059943235122\\
0.55 0.45 0.572868466640173\\
0.55 0.5 0.57394670259978\\
0.55 0.55 0.572847458936084\\
0.55 0.6 0.570029719239912\\
0.55 0.65 0.566408424867336\\
0.55 0.7 0.56262479454137\\
0.55 0.75 0.558868675536134\\
0.55 0.8 0.555163746206866\\
0.55 0.85 0.551574632926669\\
0.55 0.9 0.548267502445261\\
0.55 0.95 0.545624750721939\\
0.55 1 0.544435973239366\\
0.6 0 0.588141388045088\\
0.6 0.05 0.590507872905831\\
0.6 0.1 0.595735512174837\\
0.6 0.15 0.602271889638808\\
0.6 0.2 0.609383828012309\\
0.6 0.25 0.616754572985966\\
0.6 0.3 0.624278779095719\\
0.6 0.35 0.631986989308149\\
0.6 0.4 0.639624486627744\\
0.6 0.45 0.645825137372183\\
0.6 0.5 0.648288840983655\\
0.6 0.55 0.645795920740225\\
0.6 0.6 0.639582233057548\\
0.6 0.65 0.631945207247679\\
0.6 0.7 0.624252458836059\\
0.6 0.75 0.616743614203479\\
0.6 0.8 0.609380687942305\\
0.6 0.85 0.602271428884916\\
0.6 0.9 0.595735973952339\\
0.6 0.95 0.590508970963435\\
0.6 1 0.588141998388823\\
0.65 0 0.630415817444383\\
0.65 0.05 0.633937282636406\\
0.65 0.1 0.641625515373843\\
0.65 0.15 0.651222163625614\\
0.65 0.2 0.661700088474015\\
0.65 0.25 0.672609479418249\\
0.65 0.3 0.683827787884031\\
0.65 0.35 0.695597926408939\\
0.65 0.4 0.708013092744506\\
0.65 0.45 0.719041827625778\\
0.65 0.5 0.723741811654854\\
0.65 0.55 0.719000303324631\\
0.65 0.6 0.707943739782079\\
0.65 0.65 0.695525411437724\\
0.65 0.7 0.683787450401785\\
0.65 0.75 0.672594272778204\\
0.65 0.8 0.661695089258605\\
0.65 0.85 0.651219856756443\\
0.65 0.9 0.641623675651334\\
0.65 0.95 0.633935822951376\\
0.65 1 0.63041366517338\\
0.7 0 0.670596762222683\\
0.7 0.05 0.675233806871105\\
0.7 0.1 0.685186640237806\\
0.7 0.15 0.697573696984513\\
0.7 0.2 0.711154231711986\\
0.7 0.25 0.725370115281887\\
0.7 0.3 0.740069210154071\\
0.7 0.35 0.75580767729554\\
0.7 0.4 0.773834042183278\\
0.7 0.45 0.792577738284456\\
0.7 0.5 0.801902222452866\\
0.7 0.55 0.792533080660074\\
0.7 0.6 0.773684573297573\\
0.7 0.65 0.755677420221482\\
0.7 0.7 0.740013805920191\\
0.7 0.75 0.725351748906106\\
0.7 0.8 0.711147265096555\\
0.7 0.85 0.697569041487256\\
0.7 0.9 0.685182094759772\\
0.7 0.95 0.675229479650425\\
0.7 1 0.670591341368633\\
0.75 0 0.708098418834193\\
0.75 0.05 0.713791634196364\\
0.75 0.1 0.725734177798405\\
0.75 0.15 0.740525714187215\\
0.75 0.2 0.756806330620197\\
0.75 0.25 0.773945567589748\\
0.75 0.3 0.791698404119189\\
0.75 0.35 0.810714612414161\\
0.75 0.4 0.833915386978352\\
0.75 0.45 0.864648142799975\\
0.75 0.5 0.887065432826127\\
0.75 0.55 0.8645429891993\\
0.75 0.6 0.833511720236771\\
0.75 0.65 0.810533464973168\\
0.75 0.7 0.791641790619757\\
0.75 0.75 0.773927725029979\\
0.75 0.8 0.756797721681814\\
0.75 0.85 0.740518736451253\\
0.75 0.9 0.72572726850953\\
0.75 0.95 0.713784887202747\\
0.75 1 0.70808976584563\\
0.8 0 0.742434800511958\\
0.8 0.05 0.749103435856515\\
0.8 0.1 0.762685016920577\\
0.8 0.15 0.779377096583205\\
0.8 0.2 0.797792146631903\\
0.8 0.25 0.817283709411986\\
0.8 0.3 0.837473570909663\\
0.8 0.35 0.858609846890747\\
0.8 0.4 0.883224761935316\\
0.8 0.45 0.92129585373694\\
0.8 0.5 1\\
0.8 0.55 0.919829901871434\\
0.8 0.6 0.882755262969607\\
0.8 0.65 0.858479961510538\\
0.8 0.7 0.837436107043423\\
0.8 0.75 0.817269376438153\\
0.8 0.8 0.797782625309752\\
0.8 0.85 0.779368868367324\\
0.8 0.9 0.76267708941475\\
0.8 0.95 0.749095531754617\\
0.8 1 0.742423418667715\\
0.85 0 0.773213943523702\\
0.85 0.05 0.780761586845605\\
0.85 0.1 0.79557730249422\\
0.85 0.15 0.813576559603142\\
0.85 0.2 0.833407309076809\\
0.85 0.25 0.854463803964386\\
0.85 0.3 0.876295639345845\\
0.85 0.35 0.898632534124594\\
0.85 0.4 0.922013219661726\\
0.85 0.45 0.948131472634358\\
0.85 0.5 0.970955940647046\\
0.85 0.55 0.94782220053577\\
0.85 0.6 0.921886238926856\\
0.85 0.65 0.898588419946541\\
0.85 0.7 0.876278844073644\\
0.85 0.75 0.854453180974394\\
0.85 0.8 0.833398868776851\\
0.85 0.85 0.813569470472812\\
0.85 0.9 0.79557048937361\\
0.85 0.95 0.780754300534086\\
0.85 1 0.773200212584671\\
0.9 0 0.800011020303551\\
0.9 0.05 0.808345298254429\\
0.9 0.1 0.823980450246771\\
0.9 0.15 0.842680031574903\\
0.9 0.2 0.863147854174324\\
0.9 0.25 0.884822743829863\\
0.9 0.3 0.907253025658305\\
0.9 0.35 0.929907014550337\\
0.9 0.4 0.952261881601974\\
0.9 0.45 0.973697890038605\\
0.9 0.5 0.988904522955375\\
0.9 0.55 0.973633720220597\\
0.9 0.6 0.952247444886557\\
0.9 0.65 0.929900845002666\\
0.9 0.7 0.907246839320809\\
0.9 0.75 0.884817573469821\\
0.9 0.8 0.863144061358159\\
0.9 0.85 0.842676785495835\\
0.9 0.9 0.823976776267641\\
0.9 0.95 0.808340288932361\\
0.9 1 0.799992899612021\\
0.95 0 0.821897582320548\\
0.95 0.05 0.830942524337645\\
0.95 0.1 0.847010578757707\\
0.95 0.15 0.865892958422786\\
0.95 0.2 0.886347041827077\\
0.95 0.25 0.90780966975536\\
0.95 0.3 0.929775289949498\\
0.95 0.35 0.951541336568128\\
0.95 0.4 0.972078720827036\\
0.95 0.45 0.989753310191638\\
0.95 0.5 0.999939629650868\\
0.95 0.55 0.989756997071154\\
0.95 0.6 0.972084651583655\\
0.95 0.65 0.951545278843569\\
0.95 0.7 0.929778824681218\\
0.95 0.75 0.907812982307255\\
0.95 0.8 0.886349586076417\\
0.95 0.85 0.86589445127899\\
0.95 0.9 0.847010935376824\\
0.95 0.95 0.830940621229187\\
0.95 1 0.821861847052598\\
1 0 0.83550413601086\\
1 0.05 0.844530362946985\\
1 0.1 0.860430638518478\\
1 0.15 0.879067741729116\\
1 0.2 0.899157265278277\\
1 0.25 0.920106846886158\\
1 0.3 0.941341069198107\\
1 0.35 0.96197539813944\\
1 0.4 0.980547682476501\\
1 0.45 0.994555913744781\\
1 0.5 0.999995172093942\\
1 0.55 0.994534351995423\\
1 0.6 0.980534490515937\\
1 0.65 0.961973701422768\\
1 0.7 0.94134552653526\\
1 0.75 0.920114181858249\\
1 0.8 0.899166110864343\\
1 0.85 0.879078526198222\\
1 0.9 0.860447676262174\\
1 0.95 0.84457171251979\\
1 1 0.835475610842474\\
};

\end{axis}
\end{tikzpicture}%

%% file: Error_n1280_Model1_1D_p1p5.tikz
\begin{tikzpicture}
\def\xl{0}
\def\xu{0.025}
\def\yl{0}
\def\yu{0.5}

\begin{axis}[%
width=\figurewidth,
height=\figureheight,
scale only axis,
xmin={\xl-(\xu-\xl)/5},
xmax={\xu+(\xu-\xl)/5},
ymin={\yl-(\yu-\yl)/5},
ymax={\yu+(\yu-\yl)*0.05},
hide axis,
axis background/.style={fill=white!100}
]
\draw [->] (axis cs: \xl,\yl) -- (axis cs: {\xu+(\xu-\xl)*0.05},\yl);
\draw [->] (axis cs: \xl,\yl) -- (axis cs: \xl,{\yu+(\yu-\yl)*0.05});
\draw [blue,->] (axis cs: \xu,\yl) -- (axis cs: \xu,{\yu+(\yu-\yl)*0.05});
\node at (axis cs: {(\xu+(\xu-\xl)*0.05+\xl)/2},{\yl-0.16*(\yu-\yl)}) {$\eps$};
\node[rotate=90] at (axis cs: {\xl-0.16*(\xu-\xl)},{(\yu+(\yu-\yl)*0.05+\yl)/2}) {$\err_n^{(1.5)}(f_n)$};
\node[blue,rotate=90] at (axis cs: {\xu+0.16*(\xu-\xl)},{(\yu+(\yu-\yl)*0.05+\yl)/2}) {$\%$ Graphs Connected};
\foreach \yValue in {0,0.1,0.2,0.3,0.4,0.5} {
    \edef\temp{\noexpand\draw [-] ({\xl+(\xu-\xl)/100},\yValue) -- (\xl,\yValue) node[left] {\yValue};} 
    \temp
}
\foreach \yValue in {20,40,60,80,100} {
    \edef\temp{\noexpand\draw [blue,-] ({\xu-(\xu-\xl)/100},\yValue/200) -- (\xu,\yValue/200) node[right] {\yValue};} 
    \temp
}
\foreach \xValue in {0,0.005,0.01,0.015,0.02,0.025} {
	\edef\temp{\noexpand\draw [-] (\xValue,{\yl+(\yu-\yl)/100}) -- (\xValue,\yl) node[below] {\xValue};}    
    \temp
}
\addplot [
color=black,
solid,
mark options={solid},
line width=1.5pt,
forget plot
]
table[row sep=crcr]{
0.005 0.334422390385423\\
0.00540816326530612 0.261550604272811\\
0.00581632653061225 0.199581185503187\\
0.00622448979591837 0.204192532165609\\
0.00663265306122449 0.14676702697721\\
0.00704081632653061 0.0957602894517509\\
0.00744897959183673 0.0946119489053411\\
0.00785714285714286 0.0787884381699456\\
0.00826530612244898 0.0704959250754131\\
0.0086734693877551 0.0556984261184068\\
0.00908163265306123 0.0513394423497593\\
0.00948979591836735 0.0429589210780873\\
0.00989795918367347 0.0494264863852974\\
0.0103061224489796 0.0465515453074087\\
0.0107142857142857 0.0528960926095457\\
0.0111224489795918 0.0474396664261777\\
0.011530612244898 0.0616337542890265\\
0.0119387755102041 0.0692814251088099\\
0.0123469387755102 0.0651692132422151\\
0.0127551020408163 0.0723782851730873\\
0.0131632653061225 0.095634279564046\\
0.0135714285714286 0.0951133486200253\\
0.0139795918367347 0.115803678005684\\
0.0143877551020408 0.130789638187073\\
0.0147959183673469 0.151230090969197\\
0.0152040816326531 0.183873439742014\\
0.0156122448979592 0.184114529393456\\
0.0160204081632653 0.205556766745822\\
0.0164285714285714 0.221716183466977\\
0.0168367346938776 0.233421412766647\\
0.0172448979591837 0.242074287047287\\
0.0176530612244898 0.254896126621584\\
0.0180612244897959 0.257980592860026\\
0.018469387755102 0.262135758293352\\
0.0188775510204082 0.267679192283011\\
0.0192857142857143 0.271345582075571\\
0.0196938775510204 0.273483839241947\\
0.0201020408163265 0.272521491252163\\
0.0205102040816327 0.272599980436945\\
0.0209183673469388 0.275045084453463\\
0.0213265306122449 0.27586697030345\\
0.021734693877551 0.276142856733696\\
0.0221428571428571 0.276894989765294\\
0.0225510204081633 0.277469084353334\\
0.0229591836734694 0.276616128398288\\
0.0233673469387755 0.27782249944395\\
0.0237755102040816 0.277942581501519\\
0.0241836734693878 0.278225284846779\\
0.0245918367346939 0.278659023095876\\
0.025 0.279225248611443\\
};
\addplot [
color=black,
dashed,
line width=1pt,
forget plot
]
table[row sep=crcr]{
0.005 0.143181044231161\\
0.00540816326530612 0.0167700899908012\\
0.00581632653061225 0.0145545406120258\\
0.00622448979591837 0.0143164592710312\\
0.00663265306122449 0.0135918796237195\\
0.00704081632653061 0.0127024000517814\\
0.00744897959183673 0.0126995614049999\\
0.00785714285714286 0.0122895526214388\\
0.00826530612244898 0.0129676966090044\\
0.0086734693877551 0.0139257128014162\\
0.00908163265306123 0.0146577032165283\\
0.00948979591836735 0.0146947844296199\\
0.00989795918367347 0.0159409592120544\\
0.0103061224489796 0.0169398384477685\\
0.0107142857142857 0.0183055677981211\\
0.0111224489795918 0.0193044357367384\\
0.011530612244898 0.0180029736037405\\
0.0119387755102041 0.0206527459094304\\
0.0123469387755102 0.019858398686843\\
0.0127551020408163 0.0226067374112327\\
0.0131632653061225 0.0286143693943744\\
0.0135714285714286 0.0289583710392055\\
0.0139795918367347 0.0336394837628414\\
0.0143877551020408 0.0351211433553101\\
0.0147959183673469 0.0369447114934632\\
0.0152040816326531 0.045854946707788\\
0.0156122448979592 0.045826443390944\\
0.0160204081632653 0.0561933336279921\\
0.0164285714285714 0.0601177975453659\\
0.0168367346938776 0.0677838105484303\\
0.0172448979591837 0.0959282156240198\\
0.0176530612244898 0.256798859972989\\
0.0180612244897959 0.257243321543965\\
0.018469387755102 0.259882156815927\\
0.0188775510204082 0.26189958477093\\
0.0192857142857143 0.264537315192084\\
0.0196938775510204 0.264027408711686\\
0.0201020408163265 0.265471560243867\\
0.0205102040816327 0.265794294748347\\
0.0209183673469388 0.26665180699907\\
0.0213265306122449 0.266947850315657\\
0.021734693877551 0.266993261687476\\
0.0221428571428571 0.267189643265946\\
0.0225510204081633 0.267283688369189\\
0.0229591836734694 0.267036528303421\\
0.0233673469387755 0.267283688369189\\
0.0237755102040816 0.267524818096987\\
0.0241836734693878 0.267459179203948\\
0.0245918367346939 0.267283688369189\\
0.025 0.267283688369189\\
};
\addplot [
color=black,
dashed,
line width=1pt,
forget plot
]
table[row sep=crcr]{
0.005 0.498881006299626\\
0.00540816326530612 0.479129193922726\\
0.00581632653061225 0.448789222451921\\
0.00622448979591837 0.455580452165878\\
0.00663265306122449 0.413476375574355\\
0.00704081632653061 0.327441106071372\\
0.00744897959183673 0.336706696580374\\
0.00785714285714286 0.292949394969424\\
0.00826530612244898 0.269357680775747\\
0.0086734693877551 0.143454776853322\\
0.00908163265306123 0.0656078931874417\\
0.00948979591836735 0.0466243881379128\\
0.00989795918367347 0.0683631834733763\\
0.0103061224489796 0.06420962113295\\
0.0107142857142857 0.0916335068421281\\
0.0111224489795918 0.0674160996102911\\
0.011530612244898 0.241901819659948\\
0.0119387755102041 0.252424161938355\\
0.0123469387755102 0.242295602591958\\
0.0127551020408163 0.251125805022337\\
0.0131632653061225 0.258231581729732\\
0.0135714285714286 0.259846241260231\\
0.0139795918367347 0.264082474666096\\
0.0143877551020408 0.262464389045547\\
0.0147959183673469 0.264932903694642\\
0.0152040816326531 0.267033381243838\\
0.0156122448979592 0.267431862025498\\
0.0160204081632653 0.269642435149374\\
0.0164285714285714 0.271330451225099\\
0.0168367346938776 0.271327903581614\\
0.0172448979591837 0.274027523180341\\
0.0176530612244898 0.275577169894713\\
0.0180612244897959 0.276598721062431\\
0.018469387755102 0.277424155515926\\
0.0188775510204082 0.27888497618597\\
0.0192857142857143 0.278784285033725\\
0.0196938775510204 0.282278611194356\\
0.0201020408163265 0.283669793709653\\
0.0205102040816327 0.282853057136095\\
0.0209183673469388 0.286206813016263\\
0.0213265306122449 0.289680186036259\\
0.021734693877551 0.289851643563639\\
0.0221428571428571 0.292197814987491\\
0.0225510204081633 0.294230214255326\\
0.0229591836734694 0.288658197380923\\
0.0233673469387755 0.28993724426151\\
0.0237755102040816 0.292924851195895\\
0.0241836734693878 0.295556420459406\\
0.0245918367346939 0.296622666644807\\
0.025 0.301148134650579\\
};
\addplot [
color=orange,
mark size=3.0pt,
only marks,
mark=triangle*,
mark options={solid},
forget plot
]
table[row sep=crcr]{
0.0152040816326531 0.183873439742014\\
};
\addplot [
color=darkred,
mark size=2.0pt,
only marks,
mark=square*,
mark options={solid},
forget plot
]
table[row sep=crcr]{
0.00948979591836735 0.0429589210780873\\
};
\addplot [
color=blue,
mark size=2.0pt,
only marks,
mark=*,
mark options={solid},
forget plot
]
table[row sep=crcr]{
0.0062 0.2042\\
};
\addplot [
color=yellow,
mark size=2.0pt,
only marks,
mark=diamond*,
mark options={solid},
forget plot
]
table[row sep=crcr]{
0.022142857142857 0.276894989765294\\
};
\addplot [
color=blue,
solid,
mark options={solid},
line width=1.5pt,
forget plot
]
table[row sep=crcr]{
0.0050 0.0500\\
0.0054 0.1450\\
0.0058 0.2450\\
0.0062 0.3050\\
0.0066 0.3900\\
0.0070 0.4300\\
0.0074 0.4550\\
0.0079 0.4700\\
0.0083 0.4800\\
0.0087 0.4850\\
0.0091 0.4950\\
0.0095 0.4950\\
0.0099 0.4950\\
0.0103 0.4950\\
0.0107 0.5000\\
0.0111 0.5000\\
0.0115 0.5000\\
0.0119 0.5000\\
0.0123 0.5000\\
0.0128 0.5000\\
0.0132 0.5000\\
0.0136 0.5000\\
0.0140 0.5000\\
0.0144 0.5000\\
0.0148 0.5000\\
0.0152 0.5000\\
0.0156 0.5000\\
0.0160 0.5000\\
0.0164 0.5000\\
0.0168 0.5000\\
0.0172 0.5000\\
0.0177 0.5000\\
0.0181 0.5000\\
0.0185 0.5000\\
0.0189 0.5000\\
0.0193 0.5000\\
0.0197 0.5000\\
0.0201 0.5000\\
0.0205 0.5000\\
0.0209 0.5000\\
0.0213 0.5000\\
0.0217 0.5000\\
0.0221 0.5000\\
0.0226 0.5000\\
0.0230 0.5000\\
0.0234 0.5000\\
0.0238 0.5000\\
0.0242 0.5000\\
0.0246 0.5000\\
0.0250 0.5000\\
};
\end{axis}
\end{tikzpicture}%

%% file: Realisations_epsilon0p022143_n1280_Model1_1D_p1p5.tikz
\definecolor{mycolor1}{rgb}{0,0.447,0.741}
\definecolor{mycolor2}{rgb}{0.85,0.325,0.098}
\definecolor{mycolor3}{rgb}{0.929,0.694,0.125}
\definecolor{mycolor4}{rgb}{0.494,0.184,0.556}
\definecolor{mycolor5}{rgb}{0.466,0.674,0.188}
\definecolor{mycolor6}{rgb}{0.301,0.745,0.933}
\definecolor{mycolor7}{rgb}{0.635,0.078,0.184}

\begin{tikzpicture}
\def\xl{0}
\def\xu{1}
\def\yl{0}
\def\yu{1}

\begin{axis}[%
width=\figurewidth,
height=\figureheight,
scale only axis,
xmin={\xl-(\xu-\xl)/5},
xmax={\xu+(\xu-\xl)/5},
ymin={\yl-(\yu-\yl)/5},
ymax={\yu+(\yu-\yl)*0.05},
hide axis,
axis background/.style={fill=white!100}
]
\draw [->] (axis cs: \xl,\yl) -- (axis cs: {\xu+(\xu-\xl)*0.05},\yl);
\draw [->] (axis cs: \xl,\yl) -- (axis cs: \xl,{\yu+(\yu-\yl)*0.05});
\node at (axis cs: {(\xu+(\xu-\xl)*0.05+\xl)/2},{\yl-0.16*(\yu-\yl)}) {$\Omega$};
\node[rotate=90] at (axis cs: {\xl-0.16*(\xu-\xl)},{(\yu+(\yu-\yl)*0.05+\yl)/2}) {$f_n$};
\foreach \yValue in {0,0.2,0.4,0.6,0.8,1} {
    \edef\temp{\noexpand\draw [-] ({\xl+(\xu-\xl)/100},\yValue) -- (\xl,\yValue) node[left] {\yValue};} 
    \temp
}
\foreach \xValue in {0,0.2,0.4,0.6,0.8,1} {
	\edef\temp{\noexpand\draw [-] (\xValue,{\yl+(\yu-\yl)/100}) -- (\xValue,\yl) node[below] {\xValue};}    
    \temp
}
\addplot [
color=mycolor1,
solid,
forget plot
]
table[row sep=crcr]{
0 0\\
0.00739379296730858 0.352555462557445\\
0.015469936887653 0.352555462557445\\
0.0215054584944352 0.352555462557445\\
0.0280509988344291 0.352555462557445\\
0.0364928819372122 0.352555462557445\\
0.0455973142101893 0.352555462557445\\
0.0519138180681202 0.352555462557445\\
0.0584771525898904 0.352555462557445\\
0.0672899299668868 0.352555462557445\\
0.0751492214219124 0.352555462557445\\
0.0807756731697635 0.352555462557445\\
0.0910046439203794 0.352555462557445\\
0.10384756411149 0.352555462557445\\
0.109841901871349 0.352555462557445\\
0.117713748661732 0.352555462557445\\
0.124278170131476 0.352555462557445\\
0.131764421048501 0.352555462557445\\
0.139842873696654 0.352555462557445\\
0.146615176491381 0.352555462557445\\
0.153040093319723 0.352555462557445\\
0.162608878465381 0.352555462557445\\
0.168589015732343 0.352555462557445\\
0.173335129035222 0.352555462557445\\
0.180626844667006 0.352555462557445\\
0.184546470980641 0.352555462557445\\
0.19204309217505 0.352555462557445\\
0.200861070937291 0.352555462557445\\
0.205245781229856 0.352555462557445\\
0.21366684173522 0.352555462557445\\
0.224971149993058 0.352555462557445\\
0.232694500360531 0.352555462557445\\
0.239928697208167 0.352555462557445\\
0.248654978978848 0.352555462557445\\
0.258051261521408 0.352555462557445\\
0.266768048717129 0.352555462557445\\
0.276690089098999 0.352555462557445\\
0.282359177624182 0.352555462557445\\
0.287449930755672 0.352555462557445\\
0.293926362108711 0.352555462557445\\
0.302262230814328 0.352555462557445\\
0.309360864171716 0.352555462557445\\
0.316075912322564 0.352555462557445\\
0.328329286193228 0.352555462557445\\
0.33488600681161 0.352555462557445\\
0.341064886494271 0.352555462557445\\
0.348325869484179 0.352555462557445\\
0.352809739021552 0.352555462557445\\
0.359942987023561 0.352555462557445\\
0.369740432043976 0.352555462557445\\
0.381248969393156 0.352555462557445\\
0.391088816376553 0.352555462557445\\
0.402618198232193 0.352555462557445\\
0.41178780560957 0.352555462557445\\
0.417777273847238 0.352555462557445\\
0.423567572531339 0.352555462557445\\
0.435360258959498 0.352555462557445\\
0.441655139465557 0.352555462557445\\
0.448558500187295 0.352555462557445\\
0.453945321019544 0.352555462557445\\
0.460858722966701 0.352555462557445\\
0.468297687005314 0.352555462557445\\
0.475702347207985 0.352555462557445\\
0.486832147375197 0.352555462557445\\
0.490641026753033 0.352555462557445\\
0.498125221825437 0.352555462557445\\
0.504788377854927 0.352555462557445\\
0.509997218260962 0.352555462557445\\
0.51794313975113 0.352555462557445\\
0.522690682875737 0.352555462557445\\
0.529360804492995 0.352555462557445\\
0.534373244881685 0.352555462557445\\
0.538356046798203 0.352555462557445\\
0.547190304523309 0.352555462557445\\
0.555164669753827 0.352555462557445\\
0.564476471986747 0.352555462557445\\
0.573174800823797 0.352555462557445\\
0.578782856944427 0.352555462557445\\
0.583656810307454 0.352555462557445\\
0.590260502928193 0.352555462557445\\
0.596693603329232 0.352555462557445\\
0.606455524614159 0.352555462557445\\
0.612062796359086 0.352555462557445\\
0.622231949349152 0.352555462557445\\
0.632988091110383 0.352555462557445\\
0.640162190524043 0.352555462557445\\
0.651902736047255 0.352555462557445\\
0.657999176898195 0.352555462557445\\
0.665165495443073 0.352555462557445\\
0.674216768017507 0.352555462557445\\
0.684218591215579 0.352555462557445\\
0.691115822185419 0.352555462557445\\
0.697284572944579 0.352555462557445\\
0.70786724484914 0.352555462557445\\
0.716743398788843 0.352555462557445\\
0.724281291028591 0.352555462557445\\
0.731092418459908 0.352555462557445\\
0.740227824661933 0.352555462557445\\
0.74376444284129 0.352555462557445\\
0.752676189454188 0.352555462557445\\
0.758491336836409 0.352555462557445\\
0.765860263065131 0.352555462557445\\
0.774619082802172 0.352555462557445\\
0.787131225122028 0.356649126225008\\
0.796826469631431 0.358941542162193\\
0.803183685718723 0.360360896881072\\
0.812195237180006 0.36094060635619\\
0.822879852626606 0.362673728381346\\
0.828033750821027 0.365624234856089\\
0.83745130832171 0.367242642089106\\
0.847446327049683 0.373303008981538\\
0.860718832788239 0.373808093250062\\
0.866400348331296 0.373808093250062\\
0.870418175658732 0.373808093250062\\
0.875932585050108 0.373808093250062\\
0.884741938602735 0.373808093250062\\
0.893583828733885 0.373808093250062\\
0.900328239269656 0.374509311870687\\
0.90759653149143 0.374765898194186\\
0.914460467975331 0.374765898194186\\
0.925229682649938 0.37694814541043\\
0.932252667433979 0.378762804074589\\
0.939754551359623 0.380479839309553\\
0.948832323469715 0.38413836274849\\
0.957769289658373 0.390674506279421\\
0.972716686911837 0.403172776609847\\
0.980364397939405 0.406809326296243\\
0.988607655827051 0.406809326296243\\
1 1\\
};
\addplot [
color=mycolor2,
solid,
forget plot
]
table[row sep=crcr]{
0 0\\
0.00574170219256898 0.48016694662971\\
0.0119504369945436 0.48016694662971\\
0.0197893036144644 0.48016694662971\\
0.0275203386540629 0.48016694662971\\
0.0349114487151204 0.48016694662971\\
0.0445024259292603 0.48016694662971\\
0.0529965495788898 0.48016694662971\\
0.0593148754481776 0.48016694662971\\
0.0697661370856942 0.48016694662971\\
0.0755614928314223 0.48016694662971\\
0.0853424097556356 0.48016694662971\\
0.094875335377662 0.48016694662971\\
0.103476010159694 0.48016694662971\\
0.111662354786469 0.48016694662971\\
0.119021956285627 0.48016694662971\\
0.124634976232899 0.48016694662971\\
0.130692944862728 0.48016694662971\\
0.139517551623678 0.48016694662971\\
0.145931221073051 0.48016694662971\\
0.157939711056158 0.48016694662971\\
0.161815891408668 0.48016694662971\\
0.167777733763965 0.48016694662971\\
0.175998435385769 0.48016694662971\\
0.187064334466391 0.48016694662971\\
0.191156821919592 0.48016694662971\\
0.200697479993344 0.48016694662971\\
0.205995495786804 0.48016694662971\\
0.216404694845098 0.48016694662971\\
0.225027926871498 0.48016694662971\\
0.22829889473222 0.48016694662971\\
0.236500779829991 0.48016694662971\\
0.242754819542933 0.48016694662971\\
0.25061642594238 0.48016694662971\\
0.254340121030264 0.48016694662971\\
0.257831622953206 0.48016694662971\\
0.268730527205151 0.48016694662971\\
0.275143591604905 0.48016694662971\\
0.281033597442025 0.48016694662971\\
0.294417545267736 0.48016694662971\\
0.303055818149817 0.48016694662971\\
0.309779463154117 0.48016694662971\\
0.315839190663331 0.48016694662971\\
0.321889044748324 0.48016694662971\\
0.327735647371586 0.48016694662971\\
0.332657404600926 0.48016694662971\\
0.341350545277573 0.48016694662971\\
0.351502938309245 0.48016694662971\\
0.357167922887993 0.48016694662971\\
0.367813279783947 0.48016694662971\\
0.379504724926059 0.48016694662971\\
0.386210964541661 0.48016694662971\\
0.390258253115838 0.48016694662971\\
0.399076299306406 0.48016694662971\\
0.407257497033961 0.48016694662971\\
0.418965858008939 0.48016694662971\\
0.427893011783887 0.48016694662971\\
0.436831108024958 0.48016694662971\\
0.445517924787655 0.48016694662971\\
0.4523542779935 0.48016694662971\\
0.458250229930253 0.48016694662971\\
0.466692030731551 0.48016694662971\\
0.478551571947494 0.48016694662971\\
0.487497937020917 0.48016694662971\\
0.494305061564385 0.48016694662971\\
0.505318293512615 0.48016694662971\\
0.511687150940374 0.48016694662971\\
0.521683288883416 0.48016694662971\\
0.527617146484223 0.48016694662971\\
0.535360357405593 0.48016694662971\\
0.544613492471251 0.48016694662971\\
0.550567398142522 0.48016694662971\\
0.557973382429493 0.48016694662971\\
0.562695332404949 0.48016694662971\\
0.568253558572448 0.48016694662971\\
0.574818846782084 0.48016694662971\\
0.587740216790244 0.48016694662971\\
0.60225724630598 0.48016694662971\\
0.608897753444855 0.48016694662971\\
0.61770176165807 0.48016694662971\\
0.622859561478246 0.48016694662971\\
0.633604626727308 0.48016694662971\\
0.649603897630218 0.484407942892335\\
0.656112133155145 0.484407942892335\\
0.664818777773036 0.484407942892335\\
0.677668533845175 0.484407942892335\\
0.682775673676814 0.484407942892335\\
0.690642825471376 0.484407942892335\\
0.698798049009137 0.484407942892335\\
0.705944816935945 0.484407942892335\\
0.712908996541337 0.484407942892335\\
0.721080685309715 0.484407942892335\\
0.732878348357744 0.484407942892335\\
0.739358477240738 0.484407942892335\\
0.748948211088762 0.484407942892335\\
0.758264981926374 0.484407942892335\\
0.763412111788863 0.484407942892335\\
0.772877644798178 0.484407942892335\\
0.778525906797761 0.484407942892335\\
0.785561726726462 0.484407942892335\\
0.790609489110069 0.484407942892335\\
0.794991811540301 0.484407942892335\\
0.799940329894756 0.484407942892335\\
0.8100257299816 0.484407942892335\\
0.817850626456486 0.484407942892335\\
0.822175269656376 0.484407942892335\\
0.831155049198296 0.484407942892335\\
0.836806945264006 0.484407942892335\\
0.845055076099463 0.484407942892335\\
0.850345045099931 0.484407942892335\\
0.858919805395404 0.484407942892335\\
0.867152443975785 0.484407942892335\\
0.874696817034655 0.484407942892335\\
0.884766727199812 0.484407942892335\\
0.896265549267837 0.484407942892335\\
0.902121726831343 0.484407942892335\\
0.906460353840037 0.484407942892335\\
0.917279701266995 0.484407942892335\\
0.922427166264994 0.484407942892335\\
0.935025319747435 0.484407942892335\\
0.944463397835165 0.484407942892335\\
0.951646268514218 0.484407942892335\\
0.961341559304077 0.484407942892335\\
0.96972827838549 0.484407942892335\\
0.975491274760833 0.484407942892335\\
0.980030103485408 0.484407942892335\\
0.989025374018355 0.484407942892335\\
0.996804064219797 0.484407942892335\\
1 1\\
};
\addplot [
color=mycolor3,
solid,
forget plot
]
table[row sep=crcr]{
0 0\\
0.0058547626482045 0.438569906138374\\
0.0108289128661043 0.438569906138374\\
0.0195283723335676 0.438569906138374\\
0.02616583962539 0.438569906138374\\
0.0313865343371511 0.438569906138374\\
0.0387187144295068 0.438569906138374\\
0.0441529502512775 0.438569906138374\\
0.05090714161194 0.438569906138374\\
0.0603031013126979 0.438569906138374\\
0.0722482863876249 0.438569906138374\\
0.0804850877348134 0.438569906138374\\
0.0855292293319302 0.438569906138374\\
0.0916725758464556 0.438569906138374\\
0.0977822259464408 0.438569906138374\\
0.104092271505865 0.438569906138374\\
0.110094781082494 0.438569906138374\\
0.113470905665291 0.438569906138374\\
0.118595105034926 0.438569906138374\\
0.130216038167769 0.438569906138374\\
0.135045390337536 0.438569906138374\\
0.145446137931298 0.438569906138374\\
0.150694006207583 0.438569906138374\\
0.157442660747079 0.438569906138374\\
0.165847061154206 0.438569906138374\\
0.171326388413152 0.438569906138374\\
0.177818426331213 0.438569906138374\\
0.183233387548079 0.438569906138374\\
0.190775868658318 0.438569906138374\\
0.198941059784965 0.438569906138374\\
0.205214390179403 0.438569906138374\\
0.214478492503469 0.438569906138374\\
0.2204718937288 0.438569906138374\\
0.225628623040976 0.438569906138374\\
0.23365048915181 0.438569906138374\\
0.237536255564186 0.438569906138374\\
0.246093008005096 0.438569906138374\\
0.253428534948705 0.438569906138374\\
0.260978498065192 0.438569906138374\\
0.264350566760999 0.438569906138374\\
0.272360012841237 0.438569906138374\\
0.27892654363725 0.438569906138374\\
0.286374877912294 0.438569906138374\\
0.299818549206961 0.438569906138374\\
0.309569795470084 0.438569906138374\\
0.315782216088921 0.438569906138374\\
0.327707483155719 0.438569906138374\\
0.335341965148078 0.438569906138374\\
0.342493743320637 0.438569906138374\\
0.355836577586654 0.438569906138374\\
0.367840588461588 0.438569906138374\\
0.378816944104289 0.438569906138374\\
0.38597834171055 0.438569906138374\\
0.39305693338749 0.438569906138374\\
0.400756767814603 0.438569906138374\\
0.407452802331203 0.438569906138374\\
0.414113702966999 0.438569906138374\\
0.420166759224425 0.438569906138374\\
0.424589525032685 0.438569906138374\\
0.431578283033805 0.438569906138374\\
0.43727553949404 0.438569906138374\\
0.442094431090256 0.438569906138374\\
0.453053615291718 0.438569906138374\\
0.46203633139027 0.438569906138374\\
0.470866453517104 0.438569906138374\\
0.47688599854363 0.438569906138374\\
0.48503351845705 0.438569906138374\\
0.49351760700818 0.438569906138374\\
0.500815780581967 0.438569906138374\\
0.508927299735463 0.438569906138374\\
0.518767106879373 0.438569906138374\\
0.529698777479209 0.438569906138374\\
0.5358431438225 0.438569906138374\\
0.544014557729319 0.438569906138374\\
0.554236679683831 0.438569906138374\\
0.559513302617875 0.438569906138374\\
0.57232953716534 0.438569906138374\\
0.585208505701656 0.438569906138374\\
0.590943584404689 0.438569906138374\\
0.597851311386095 0.438569906138374\\
0.606094035307618 0.438569906138374\\
0.615566611217116 0.438569906138374\\
0.627243224270838 0.438569906138374\\
0.638653333470614 0.438569906138374\\
0.648859479046979 0.438569906138374\\
0.657936360549075 0.438569906138374\\
0.665924771227162 0.438569906138374\\
0.675784845352998 0.438569906138374\\
0.682365627225521 0.438569906138374\\
0.690501870095935 0.438569906138374\\
0.698745381123452 0.438569906138374\\
0.710181249815015 0.438569906138374\\
0.726145022248576 0.438569906138374\\
0.733036396329213 0.438569906138374\\
0.739571631222444 0.438569906138374\\
0.745171924731991 0.438569906138374\\
0.75294806996738 0.438569906138374\\
0.756536116361404 0.438569906138374\\
0.763641939248203 0.438569906138374\\
0.774600121685726 0.438569906138374\\
0.783080857373068 0.438569906138374\\
0.789481263141463 0.438569906138374\\
0.802441120999229 0.438569906138374\\
0.805947537323327 0.438569906138374\\
0.817051601024388 0.438569906138374\\
0.824081025732588 0.438569906138374\\
0.829199545456907 0.438569906138374\\
0.837799253790095 0.438569906138374\\
0.841452477034845 0.438569906138374\\
0.846123913874263 0.438569906138374\\
0.853410772956192 0.438569906138374\\
0.866608383008952 0.438569906138374\\
0.876343063972347 0.438569906138374\\
0.881891673073183 0.438569906138374\\
0.891349287617502 0.438569906138374\\
0.898225227165262 0.438569906138374\\
0.902838817030905 0.438569906138374\\
0.914091718738279 0.438569906138374\\
0.92282230963363 0.438569906138374\\
0.934408862662361 0.438569906138374\\
0.939821350994707 0.438569906138374\\
0.944961905274349 0.438569906138374\\
0.955960783954728 0.438569906138374\\
0.962274594641328 0.438569906138374\\
0.972621951985755 0.438569906138374\\
0.979067089686568 0.438569906138374\\
0.982907518102659 0.438569906138374\\
0.994594103525801 0.438569906138374\\
1 1\\
};
\addplot [
color=mycolor4,
solid,
forget plot
]
table[row sep=crcr]{
0 0\\
0.00695556425768351 0.358186920853873\\
0.0110007067489553 0.358186920853873\\
0.0187296012021148 0.358186920853873\\
0.027273838056825 0.358186920853873\\
0.0388728822312132 0.358186920853873\\
0.0464144065485507 0.358186920853873\\
0.0550800323053232 0.358186920853873\\
0.0610858978896993 0.358186920853873\\
0.0718335111784418 0.358186920853873\\
0.0829365227027477 0.358186920853873\\
0.0952249377767067 0.358186920853873\\
0.107300353999787 0.358186920853873\\
0.116090322900836 0.358186920853873\\
0.125127204843529 0.358186920853873\\
0.131898795610856 0.358186920853873\\
0.141255685747963 0.358186920853873\\
0.149348509348625 0.358186920853873\\
0.157769257667586 0.358186920853873\\
0.170124989756389 0.358186920853873\\
0.180292778930957 0.358186920853873\\
0.188777013920295 0.358186920853873\\
0.201816919431471 0.358186920853873\\
0.207447100533258 0.358186920853873\\
0.218973573588753 0.358186920853873\\
0.226866975930166 0.358186920853873\\
0.230503399918272 0.358186920853873\\
0.239459300067999 0.358186920853873\\
0.24529815148884 0.358186920853873\\
0.248963076370926 0.358186920853873\\
0.258343316177198 0.358186920853873\\
0.266057326637953 0.358186920853873\\
0.272539719857529 0.358186920853873\\
0.277675763377762 0.358186920853873\\
0.284512730542503 0.358186920853873\\
0.290873550005335 0.358186920853873\\
0.29966019704213 0.358186920853873\\
0.307672504806727 0.358186920853873\\
0.316623165149791 0.358186920853873\\
0.324947347731411 0.358186920853873\\
0.334112547793244 0.358186920853873\\
0.349215682474474 0.358186920853873\\
0.35878174140528 0.358186920853873\\
0.364169815809088 0.358186920853873\\
0.377354052796315 0.358186920853873\\
0.382159717413028 0.358186920853873\\
0.387675442381081 0.358186920853873\\
0.399651376251194 0.358186920853873\\
0.408256997407138 0.358186920853873\\
0.414559882142097 0.358186920853873\\
0.421425993062192 0.358186920853873\\
0.430012645276715 0.358186920853873\\
0.440522975225627 0.358186920853873\\
0.446474787180896 0.358186920853873\\
0.452745586810453 0.358186920853873\\
0.45979305005008 0.358186920853873\\
0.464508537799438 0.358186920853873\\
0.470091080953767 0.358186920853873\\
0.479273100516698 0.358186920853873\\
0.485327620278467 0.358186920853873\\
0.494682272975604 0.358186920853873\\
0.502720976769329 0.358186920853873\\
0.509341524258162 0.358186920853873\\
0.513596986544479 0.358186920853873\\
0.52685156365226 0.358186920853873\\
0.530236951418414 0.358186920853873\\
0.538343345725905 0.358186920853873\\
0.546173202796172 0.358186920853873\\
0.553354365858209 0.358186920853873\\
0.560200793075981 0.358186920853873\\
0.569404000426869 0.358186920853873\\
0.577421092231882 0.358186920853873\\
0.583270767107861 0.358186920853873\\
0.588424946235451 0.358186920853873\\
0.595647138403379 0.358186920853873\\
0.608134787561232 0.358186920853873\\
0.614841591730581 0.358186920853873\\
0.62143793881899 0.358186920853873\\
0.626920533954389 0.358186920853873\\
0.635916969904574 0.358186920853873\\
0.642929232203666 0.358186920853873\\
0.647157337412638 0.358186920853873\\
0.654689160864411 0.358186920853873\\
0.665422984683694 0.358186920853873\\
0.673525184114048 0.358186920853873\\
0.682718669768961 0.358186920853873\\
0.687981612171369 0.358186920853873\\
0.694873953000856 0.358186920853873\\
0.701176415076988 0.358186920853873\\
0.70596816239836 0.358186920853873\\
0.712953967534033 0.358186920853873\\
0.724129859575127 0.358186920853873\\
0.731693257315413 0.358186920853873\\
0.737586245908993 0.358186920853873\\
0.745526000025077 0.358186920853873\\
0.749819739291809 0.358186920853873\\
0.757594228420451 0.358186920853873\\
0.765218328373154 0.358186920853873\\
0.770078076725846 0.358186920853873\\
0.777991652347209 0.358186920853873\\
0.788954937795805 0.358186920853873\\
0.792915482946217 0.358186920853873\\
0.804570708396599 0.358186920853873\\
0.808845109949444 0.358186920853873\\
0.818226500516543 0.358186920853873\\
0.826132725362004 0.358186920853873\\
0.83085897922552 0.358186920853873\\
0.849848927708529 0.358186920853873\\
0.855248663446553 0.358186920853873\\
0.859434662387769 0.358186920853873\\
0.86449169472645 0.358186920853873\\
0.875395158374669 0.358186920853873\\
0.885974068117995 0.358186920853873\\
0.891633487444053 0.358186920853873\\
0.901515298343741 0.358186920853873\\
0.907208987391893 0.358186920853873\\
0.910829913613341 0.358186920853873\\
0.917789167915527 0.358186920853873\\
0.924910853827522 0.358186920853873\\
0.934354753455163 0.358186920853873\\
0.940199786594735 0.358186920853873\\
0.946040099626277 0.358186920853873\\
0.95119610429897 0.358186920853873\\
0.95652956585494 0.358186920853873\\
0.96934525323524 0.359111038577523\\
0.975343296152371 0.360353636604132\\
0.981534517468843 0.36076761540993\\
0.992530364873458 0.36076761540993\\
1 1\\
};
\addplot [
color=mycolor5,
solid,
forget plot
]
table[row sep=crcr]{
0 0\\
0.00986178771576918 0.46394547681857\\
0.0166906334208665 0.46394547681857\\
0.0262648977468932 0.46394547681857\\
0.0345415600685535 0.46394547681857\\
0.0459757500052422 0.46394547681857\\
0.0560832987508793 0.46394547681857\\
0.0623373007428145 0.46394547681857\\
0.0721495358660209 0.46394547681857\\
0.077794553191572 0.46394547681857\\
0.0907939622128378 0.46394547681857\\
0.100526576650534 0.46394547681857\\
0.109156077267806 0.46394547681857\\
0.114957040180962 0.46394547681857\\
0.119045048276651 0.46394547681857\\
0.129525634421168 0.46394547681857\\
0.135067910423582 0.46394547681857\\
0.145132167480987 0.46394547681857\\
0.153072761000881 0.46394547681857\\
0.158542406714321 0.46394547681857\\
0.167822298807575 0.46394547681857\\
0.173098489111572 0.46394547681857\\
0.178341254165829 0.46394547681857\\
0.18838546512271 0.46394547681857\\
0.197001916395814 0.46394547681857\\
0.202881563854272 0.46394547681857\\
0.215047242437538 0.46394547681857\\
0.225543631327761 0.46394547681857\\
0.230848393327049 0.46394547681857\\
0.236771529825916 0.46394547681857\\
0.240444089962755 0.46394547681857\\
0.2528290442544 0.46394547681857\\
0.258534700059398 0.46394547681857\\
0.269338971727638 0.46394547681857\\
0.277906989805809 0.46394547681857\\
0.285640033272589 0.46394547681857\\
0.294076142532642 0.46394547681857\\
0.297849448874277 0.46394547681857\\
0.303557971876334 0.46394547681857\\
0.309486452682348 0.46394547681857\\
0.318341928209682 0.46394547681857\\
0.327333913361938 0.46394547681857\\
0.331342498324984 0.46394547681857\\
0.337625715314739 0.46394547681857\\
0.345267496611504 0.46394547681857\\
0.354998855866891 0.46394547681857\\
0.361524343993699 0.46394547681857\\
0.371304619162649 0.46394547681857\\
0.380515121279633 0.46394547681857\\
0.385339749561173 0.46394547681857\\
0.392558304918604 0.46394547681857\\
0.399810667708561 0.46394547681857\\
0.40724416095544 0.46394547681857\\
0.4121970808357 0.46394547681857\\
0.421427424842169 0.46394547681857\\
0.427433239985109 0.46394547681857\\
0.4378149177291 0.46394547681857\\
0.443120793162966 0.46394547681857\\
0.447943944382515 0.46394547681857\\
0.456978542644923 0.46394547681857\\
0.463890445713749 0.46394547681857\\
0.471634107051513 0.46394547681857\\
0.478303470835529 0.46394547681857\\
0.485216068288887 0.46394547681857\\
0.491177108696888 0.46394547681857\\
0.497613427897007 0.46394547681857\\
0.507460569325481 0.46394547681857\\
0.515575262354125 0.46394547681857\\
0.523368327688094 0.46394547681857\\
0.534440972091594 0.46394547681857\\
0.54365921515684 0.46394547681857\\
0.546961784472202 0.46394547681857\\
0.55765902902255 0.46394547681857\\
0.566453803848534 0.46394547681857\\
0.570942184812716 0.46394547681857\\
0.575524330096429 0.46394547681857\\
0.58402586086839 0.46394547681857\\
0.588815116754685 0.46394547681857\\
0.594842607302659 0.46394547681857\\
0.606398017697934 0.46394547681857\\
0.613844348165026 0.46394547681857\\
0.618996137661399 0.46394547681857\\
0.623130368944373 0.46394547681857\\
0.631953062364708 0.46394547681857\\
0.644433437579771 0.46394547681857\\
0.650898352097782 0.46394547681857\\
0.658978721716531 0.46394547681857\\
0.671518870026312 0.46394547681857\\
0.684632259752288 0.46394547681857\\
0.696946156409948 0.46394547681857\\
0.705311868877952 0.46394547681857\\
0.712874651037404 0.46394547681857\\
0.722355091168786 0.46394547681857\\
0.729637310528381 0.46394547681857\\
0.7397203548151 0.46394547681857\\
0.747496456145909 0.46394547681857\\
0.756115242862567 0.46394547681857\\
0.766538342675892 0.46394547681857\\
0.772312682949824 0.46394547681857\\
0.780678279100129 0.46394547681857\\
0.787920816347971 0.46394547681857\\
0.792125352181546 0.46394547681857\\
0.803699684586806 0.46394547681857\\
0.812582409396442 0.46394547681857\\
0.815864420298515 0.46394547681857\\
0.82117157908853 0.46394547681857\\
0.830322889110113 0.46394547681857\\
0.835122264147304 0.46394547681857\\
0.840882519028383 0.46394547681857\\
0.850078347443398 0.46394547681857\\
0.857199181599525 0.46394547681857\\
0.863936460499755 0.46394547681857\\
0.869794217967066 0.46394547681857\\
0.876426819826354 0.46394547681857\\
0.884022982249579 0.46394547681857\\
0.890252188455537 0.46394547681857\\
0.900025348430105 0.465641600993734\\
0.912123230995892 0.47108698171794\\
0.923184357616962 0.473087331940739\\
0.926419592738078 0.473087331940739\\
0.936172939315727 0.47474987445252\\
0.948090355816309 0.475298591595027\\
0.953810167969595 0.475298591595027\\
0.957899623069952 0.475298591595027\\
0.968688959933072 0.477491546434857\\
0.977179434255531 0.47858461026304\\
0.985118929615211 0.479967652761622\\
0.990782734504455 0.479967652761622\\
1 1\\
};
\addplot [
color=mycolor6,
solid,
forget plot
]
table[row sep=crcr]{
0 0\\
0.00878259826751082 0.349756433112332\\
0.018438681163385 0.349756433112332\\
0.0235441367374577 0.349756433112332\\
0.0315148861746423 0.349756433112332\\
0.0414724576385124 0.349756433112332\\
0.0531912451696391 0.349756433112332\\
0.0594862105630167 0.349756433112332\\
0.0673429999530191 0.349756433112332\\
0.0754114680842328 0.349756433112332\\
0.0826843159360094 0.349756433112332\\
0.0870721899742319 0.349756433112332\\
0.0934874547640633 0.349756433112332\\
0.102281397671848 0.349756433112332\\
0.10895495344172 0.349756433112332\\
0.118338718995126 0.349756433112332\\
0.127272407484126 0.349756433112332\\
0.136095278164739 0.349756433112332\\
0.142359710075819 0.349756433112332\\
0.149931152500551 0.349756433112332\\
0.154795926622713 0.349756433112332\\
0.163562496361934 0.349756433112332\\
0.170126922914619 0.349756433112332\\
0.179881960593076 0.349756433112332\\
0.192036792382424 0.349756433112332\\
0.199051423623143 0.349756433112332\\
0.210646429198318 0.349756433112332\\
0.217612147534108 0.349756433112332\\
0.220788328793562 0.349756433112332\\
0.227040612066322 0.349756433112332\\
0.233748103393261 0.349756433112332\\
0.244064725313517 0.349756433112332\\
0.248585126631879 0.349756433112332\\
0.257966104363363 0.349756433112332\\
0.264527737177137 0.349756433112332\\
0.276922287923573 0.349756433112332\\
0.286132599253151 0.349756433112332\\
0.294089193308014 0.349756433112332\\
0.301830586266931 0.349756433112332\\
0.3099097715229 0.349756433112332\\
0.316700869653725 0.349756433112332\\
0.325811679935 0.349756433112332\\
0.333788845599071 0.349756433112332\\
0.341062941572823 0.349756433112332\\
0.347736820278584 0.349756433112332\\
0.361384295744923 0.349756433112332\\
0.371508401213146 0.349756433112332\\
0.379371231852286 0.349756433112332\\
0.390123300644636 0.349756433112332\\
0.399563750261836 0.349756433112332\\
0.412099491183804 0.349756433112332\\
0.416886851902908 0.349756433112332\\
0.424891916128684 0.349756433112332\\
0.432742764345648 0.349756433112332\\
0.439093412106258 0.349756433112332\\
0.449815319498308 0.349756433112332\\
0.456054866817163 0.349756433112332\\
0.462618013514618 0.349756433112332\\
0.470420208572676 0.349756433112332\\
0.478895480094127 0.349756433112332\\
0.485225686578618 0.349756433112332\\
0.492554612958395 0.349756433112332\\
0.500303595045534 0.349756433112332\\
0.507485423592285 0.349756433112332\\
0.518499620414872 0.349756433112332\\
0.525211134330072 0.349756433112332\\
0.534731462990774 0.349756433112332\\
0.541978671223474 0.349756433112332\\
0.547918622676396 0.349756433112332\\
0.555484963347549 0.349756433112332\\
0.564391747518116 0.349756433112332\\
0.572975030505984 0.349756433112332\\
0.579310067874056 0.349756433112332\\
0.586182537298315 0.349756433112332\\
0.593028882322638 0.349756433112332\\
0.597411692105503 0.349756433112332\\
0.602590578654587 0.349756433112332\\
0.613996064123526 0.349756433112332\\
0.620779188623008 0.349756433112332\\
0.630466726573518 0.349756433112332\\
0.638873050657911 0.349756433112332\\
0.652390343366094 0.349756433112332\\
0.658228972339929 0.349756433112332\\
0.66427257813615 0.349756433112332\\
0.668416052109664 0.349756433112332\\
0.674854612706409 0.349756433112332\\
0.681620368151351 0.349756433112332\\
0.686439293538524 0.349756433112332\\
0.690665382977072 0.349756433112332\\
0.696393207711487 0.349756433112332\\
0.703876980960916 0.349756433112332\\
0.708590693093133 0.349756433112332\\
0.71862410773334 0.349756433112332\\
0.727699537043441 0.349756433112332\\
0.734186993079291 0.349756433112332\\
0.739788086439014 0.349756433112332\\
0.745881897626861 0.349756433112332\\
0.753747489047502 0.349756433112332\\
0.757682494828312 0.349756433112332\\
0.76198810125462 0.349756433112332\\
0.771219554309787 0.349756433112332\\
0.780843620515211 0.349756433112332\\
0.786754773040756 0.349756433112332\\
0.794844809435464 0.349756433112332\\
0.799757848783946 0.349756433112332\\
0.80267839922833 0.349756433112332\\
0.810475364612377 0.349756433112332\\
0.820139856965194 0.349756433112332\\
0.828969112819342 0.349756433112332\\
0.835073544564339 0.349756433112332\\
0.844223298373762 0.349756433112332\\
0.852133713116623 0.349756433112332\\
0.856483714600891 0.349756433112332\\
0.865314633002204 0.349756433112332\\
0.873145191153507 0.349756433112332\\
0.878976505399485 0.349756433112332\\
0.889150015539043 0.350415631693659\\
0.899046993282365 0.352220635555401\\
0.906161433881844 0.352220635555401\\
0.916007440418559 0.355160566830784\\
0.925691872508028 0.357116196395358\\
0.934503619782305 0.357172576605609\\
0.940250751199984 0.357172576605609\\
0.9456016092407 0.357562890719071\\
0.957987753608357 0.359377475940632\\
0.963839957327442 0.36144250785078\\
0.977821929335811 0.365046185400056\\
0.98528508353859 0.369856110044191\\
1 1\\
};
\addplot [
color=mycolor7,
solid,
forget plot
]
table[row sep=crcr]{
0 0\\
0.00422704508733129 0.378247709625623\\
0.0103864601643953 0.378247709625623\\
0.018728881435616 0.378247709625623\\
0.0276673514438054 0.378247709625623\\
0.0345860956659497 0.378247709625623\\
0.045035267224203 0.378247709625623\\
0.0521739993679095 0.378247709625623\\
0.0605563201647787 0.378247709625623\\
0.0673841625408863 0.378247709625623\\
0.0770176081553625 0.378247709625623\\
0.0827458977093649 0.378247709625623\\
0.0937349846263056 0.378247709625623\\
0.107549559449934 0.378247709625623\\
0.113610903968387 0.378247709625623\\
0.123095424978023 0.378247709625623\\
0.128605976462078 0.378247709625623\\
0.139523216457043 0.378247709625623\\
0.143432249988541 0.378247709625623\\
0.151399595789309 0.378247709625623\\
0.15540702762754 0.378247709625623\\
0.158670634644433 0.378247709625623\\
0.167896817387199 0.378247709625623\\
0.17848560096295 0.378247709625623\\
0.189011735586803 0.378247709625623\\
0.19778636868347 0.378247709625623\\
0.206294295459452 0.378247709625623\\
0.212539573893064 0.378247709625623\\
0.223384449377847 0.378247709625623\\
0.227825888704068 0.378247709625623\\
0.234765092425956 0.378247709625623\\
0.2417172380665 0.378247709625623\\
0.247345047755949 0.378247709625623\\
0.253102094806007 0.378247709625623\\
0.258941670700964 0.378247709625623\\
0.265108484808473 0.378247709625623\\
0.272229781715116 0.378247709625623\\
0.283728915225057 0.378247709625623\\
0.290139568894122 0.378247709625623\\
0.298652733992684 0.378247709625623\\
0.307002249263618 0.378247709625623\\
0.313956394294539 0.378247709625623\\
0.321214346912999 0.378247709625623\\
0.328102734941354 0.378247709625623\\
0.340123803930224 0.378247709625623\\
0.353922099086338 0.378247709625623\\
0.361235202689332 0.378247709625623\\
0.375557278994231 0.378247709625623\\
0.384131631559767 0.378247709625623\\
0.391417571670751 0.378247709625623\\
0.398376500922326 0.378247709625623\\
0.403570225854418 0.378247709625623\\
0.409373812850909 0.378247709625623\\
0.415937742985064 0.378247709625623\\
0.420922474449694 0.378247709625623\\
0.433490864435749 0.378247709625623\\
0.441168758043186 0.378247709625623\\
0.446877405234202 0.378247709625623\\
0.454687493799705 0.378247709625623\\
0.461550770648974 0.378247709625623\\
0.47022210235652 0.378247709625623\\
0.476295978788845 0.378247709625623\\
0.48514264270274 0.378247709625623\\
0.490417390076174 0.378247709625623\\
0.497621935180134 0.378247709625623\\
0.504512395357882 0.378247709625623\\
0.513334979712267 0.378247709625623\\
0.520743823494456 0.378247709625623\\
0.526570036176312 0.378247709625623\\
0.534698380287823 0.378247709625623\\
0.544005820786544 0.378247709625623\\
0.550551426880991 0.378247709625623\\
0.558438061248971 0.378247709625623\\
0.562544854933735 0.378247709625623\\
0.569783239876163 0.378247709625623\\
0.575422801087775 0.378247709625623\\
0.583927771600345 0.378247709625623\\
0.589922183867048 0.378247709625623\\
0.597834706953294 0.378247709625623\\
0.608638574353353 0.378247709625623\\
0.617165165126688 0.378247709625623\\
0.624833904620096 0.378247709625623\\
0.631455370343877 0.378247709625623\\
0.636322388763163 0.378247709625623\\
0.643373787791784 0.378247709625623\\
0.649861784265167 0.378247709625623\\
0.661325373436521 0.378247709625623\\
0.671376219882024 0.378247709625623\\
0.679507564673468 0.378247709625623\\
0.686335601816529 0.378247709625623\\
0.691746000560573 0.378247709625623\\
0.704875871579148 0.378247709625623\\
0.714705015047212 0.378247709625623\\
0.720351404610696 0.378247709625623\\
0.725640218598017 0.378247709625623\\
0.729363150436854 0.378247709625623\\
0.734909691095842 0.378247709625623\\
0.738754599043202 0.378247709625623\\
0.750957406448896 0.378247709625623\\
0.759931797478678 0.378247709625623\\
0.766995924384589 0.378247709625623\\
0.775070810592543 0.378247709625623\\
0.781151350811032 0.378247709625623\\
0.794880204076521 0.378247709625623\\
0.801848108863006 0.378247709625623\\
0.808495618866521 0.378247709625623\\
0.816989041337801 0.378247709625623\\
0.826129544257165 0.378247709625623\\
0.832101556795537 0.378247709625623\\
0.837469152535256 0.378247709625623\\
0.844097498984227 0.378247709625623\\
0.850487032123738 0.378247709625623\\
0.859246236405992 0.378247709625623\\
0.866376558535802 0.378247709625623\\
0.873451612279857 0.378247709625623\\
0.885807424296944 0.378247709625623\\
0.89942445034986 0.378247709625623\\
0.912312851185602 0.378247709625623\\
0.922611127686741 0.378247709625623\\
0.933760978383936 0.378247709625623\\
0.940775683259667 0.378247709625623\\
0.947864343851411 0.378247709625623\\
0.952671703125545 0.378247709625623\\
0.962914956117956 0.378247709625623\\
0.968765137741681 0.378247709625623\\
0.97266578227288 0.378247709625623\\
0.983773786423042 0.380153074150003\\
0.993352884437895 0.381442708081974\\
1 1\\
};
\addplot [
color=mycolor1,
solid,
forget plot
]
table[row sep=crcr]{
0 0\\
0.011294987029619 0.5\\
0.0211531037888401 0.5\\
0.0288295366276068 0.5\\
0.0360708441833427 0.5\\
0.042648414439866 0.5\\
0.0523782689367551 0.5\\
0.0572750719489611 0.5\\
0.0628903921014132 0.5\\
0.0677686078286365 0.5\\
0.0718392810052769 0.5\\
0.0817448706889525 0.5\\
0.0911594020356025 0.5\\
0.0994271927287376 0.5\\
0.105808880237034 0.5\\
0.112871874049711 0.5\\
0.120548502640382 0.5\\
0.133155061835275 0.5\\
0.14109960598725 0.5\\
0.147375598248104 0.5\\
0.152636675206917 0.5\\
0.158038788095334 0.5\\
0.163399523177265 0.5\\
0.168892988281861 0.5\\
0.178227003502948 0.5\\
0.184397929007989 0.5\\
0.19794068687282 0.5\\
0.202311696882467 0.5\\
0.211392502961193 0.5\\
0.219962566589883 0.5\\
0.229221650899553 0.5\\
0.234004526711326 0.5\\
0.240197213168007 0.5\\
0.245078270803916 0.5\\
0.249712090600504 0.5\\
0.253989260815944 0.5\\
0.262679530224913 0.5\\
0.268374892409072 0.5\\
0.279748078347046 0.5\\
0.283295077989597 0.5\\
0.289409248651995 0.5\\
0.296765447288324 0.5\\
0.305035125560225 0.5\\
0.315648042318761 0.5\\
0.323009240562171 0.5\\
0.335590503131237 0.5\\
0.348298669868739 0.5\\
0.355278371235184 0.5\\
0.363923371554591 0.5\\
0.368780392849433 0.5\\
0.373210799070582 0.5\\
0.377574313745246 0.5\\
0.388387874505996 0.5\\
0.397852561236768 0.5\\
0.405773957749308 0.5\\
0.409498456535682 0.5\\
0.422513926686829 0.5\\
0.43042187146669 0.5\\
0.442033244578875 0.5\\
0.449300299301643 0.5\\
0.460353634583088 0.5\\
0.464009406826216 0.5\\
0.474110167151251 0.5\\
0.485246792582896 0.5\\
0.489447649185367 0.5\\
0.495481217221123 0.5\\
0.504488704127786 0.5\\
0.51086196325129 0.5\\
0.525565546482794 0.5\\
0.532451282671558 0.5\\
0.53851592254642 0.5\\
0.54978562045773 0.5\\
0.558595965122075 0.5\\
0.56510177350841 0.5\\
0.574187825369029 0.5\\
0.580223618054033 0.5\\
0.587920031424803 0.5\\
0.595264308402191 0.5\\
0.602708071951839 0.5\\
0.609448817976728 0.5\\
0.615038990082154 0.5\\
0.620169150049723 0.5\\
0.62635098523043 0.5\\
0.634241731382303 0.5\\
0.641427781070083 0.5\\
0.647523307828224 0.5\\
0.660852498666442 0.5\\
0.668023064298222 0.5\\
0.675402881777588 0.5\\
0.683000425315273 0.5\\
0.689274404909753 0.5\\
0.693978823931994 0.5\\
0.700556654370292 0.5\\
0.710531604674458 0.5\\
0.719285006884383 0.5\\
0.726103630375646 0.5\\
0.731946356963858 0.5\\
0.739033595575275 0.5\\
0.74926088489798 0.5\\
0.757203439472106 0.5\\
0.763828631027597 0.5\\
0.773530282245349 0.5\\
0.783421156777976 0.5\\
0.79012674728373 0.5\\
0.795198852864066 0.5\\
0.808332904192068 0.5\\
0.81365844270575 0.5\\
0.821713358693843 0.5\\
0.830653780616089 0.5\\
0.841276691637419 0.5\\
0.848365669587279 0.5\\
0.86257404865077 0.5\\
0.867571742228158 0.5\\
0.880797336039962 0.5\\
0.886210589927232 0.5\\
0.896553962169967 0.5\\
0.905001415641549 0.5\\
0.915361475886514 0.5\\
0.926374887050366 0.5\\
0.934004257897 0.5\\
0.939371330301829 0.5\\
0.945708784719455 0.5\\
0.955377485039781 0.5\\
0.964207311930758 0.5\\
0.975552048150591 0.5\\
0.983019687478002 0.5\\
0.986362386342506 0.5\\
0.991224369768716 0.5\\
1 1\\
};
\addplot [
color=mycolor2,
solid,
forget plot
]
table[row sep=crcr]{
0 0\\
0.00499352060421809 0.448002265422398\\
0.0120086375544922 0.448002265422398\\
0.0152974025717445 0.448002265422398\\
0.0272284148927051 0.448002265422398\\
0.0325709523900779 0.448002265422398\\
0.0405388609370051 0.448002265422398\\
0.0521733793399074 0.448002265422398\\
0.0578831136018548 0.448002265422398\\
0.0625459083787632 0.448002265422398\\
0.0684811010209129 0.448002265422398\\
0.0744213248946434 0.448002265422398\\
0.0838784599168959 0.448002265422398\\
0.0907127518299737 0.448002265422398\\
0.0985108173303282 0.448002265422398\\
0.108795027707178 0.448002265422398\\
0.119082474650533 0.448002265422398\\
0.127283294235563 0.448002265422398\\
0.133070057963871 0.448002265422398\\
0.140877339490477 0.448002265422398\\
0.145358161630818 0.448002265422398\\
0.159617543363774 0.448002265422398\\
0.169169210942981 0.448002265422398\\
0.171963152393228 0.448002265422398\\
0.177431530475393 0.448002265422398\\
0.18531228228864 0.448002265422398\\
0.197387525827264 0.448002265422398\\
0.204804187670949 0.448002265422398\\
0.213497499309273 0.448002265422398\\
0.221230797327294 0.448002265422398\\
0.234686398814455 0.448002265422398\\
0.240412998977191 0.448002265422398\\
0.256677570063666 0.448002265422398\\
0.262945774929435 0.448002265422398\\
0.27102128033027 0.448002265422398\\
0.27900169922768 0.448002265422398\\
0.284668674703848 0.448002265422398\\
0.296250568677223 0.448002265422398\\
0.301783669539542 0.448002265422398\\
0.307634692827078 0.448002265422398\\
0.315117585333607 0.448002265422398\\
0.323573771162409 0.448002265422398\\
0.328148965378829 0.448002265422398\\
0.33155139646845 0.448002265422398\\
0.343376454728428 0.448002265422398\\
0.351832490842227 0.448002265422398\\
0.359637605202008 0.448002265422398\\
0.368690703096094 0.448002265422398\\
0.376661699617454 0.448002265422398\\
0.38436442493577 0.448002265422398\\
0.394147519929894 0.448002265422398\\
0.401082397525651 0.448002265422398\\
0.407814377077994 0.448002265422398\\
0.416803863184487 0.448002265422398\\
0.429261320123467 0.448002265422398\\
0.44238526430419 0.448002265422398\\
0.450725157626697 0.448002265422398\\
0.455600868298831 0.448002265422398\\
0.462673531815445 0.448002265422398\\
0.470569426834885 0.448002265422398\\
0.478104596790864 0.448002265422398\\
0.488587907055134 0.448002265422398\\
0.494672509457533 0.448002265422398\\
0.50012118950913 0.448002265422398\\
0.507235404023284 0.448002265422398\\
0.513212931717505 0.448002265422398\\
0.522800623253514 0.448002265422398\\
0.534851781462504 0.448002265422398\\
0.541072117669535 0.448002265422398\\
0.547063393219765 0.448002265422398\\
0.550618080824591 0.448002265422398\\
0.555704805838894 0.448002265422398\\
0.562087752936838 0.448002265422398\\
0.569123848750593 0.448002265422398\\
0.575763475477931 0.448002265422398\\
0.586902812399022 0.448002265422398\\
0.597655828321079 0.448002265422398\\
0.600634607825969 0.448002265422398\\
0.608727408132907 0.448002265422398\\
0.612056733511065 0.448002265422398\\
0.624792504483967 0.448002265422398\\
0.632311616900499 0.448002265422398\\
0.638106898801005 0.448002265422398\\
0.647670708675211 0.448002265422398\\
0.655246393161734 0.448002265422398\\
0.666127695807293 0.448002265422398\\
0.674023144585851 0.448002265422398\\
0.68140770007578 0.448002265422398\\
0.690652156869272 0.448002265422398\\
0.695354569218326 0.448002265422398\\
0.706403628082067 0.448002265422398\\
0.718390226563091 0.448002265422398\\
0.727551174790331 0.448002265422398\\
0.735076721797279 0.448002265422398\\
0.742647386011579 0.448002265422398\\
0.75038333982695 0.448002265422398\\
0.761412393252192 0.448002265422398\\
0.764766419916328 0.448002265422398\\
0.772630300582891 0.448002265422398\\
0.788648312670732 0.448002265422398\\
0.794651713607655 0.448002265422398\\
0.801679696425824 0.448002265422398\\
0.805091196330018 0.448002265422398\\
0.815525257516667 0.448002265422398\\
0.822350745184223 0.448002265422398\\
0.827864884930394 0.448002265422398\\
0.832462520340597 0.448002265422398\\
0.841086482318294 0.448002265422398\\
0.850022057267701 0.448002265422398\\
0.855682137026202 0.448002265422398\\
0.860504783567549 0.448002265422398\\
0.866536698009438 0.448002265422398\\
0.875917884715509 0.448002265422398\\
0.882416986744464 0.448002265422398\\
0.886089522946847 0.448002265422398\\
0.893986969993536 0.448002265422398\\
0.905355433327463 0.448002265422398\\
0.911752562387437 0.448002265422398\\
0.924522819159473 0.448002265422398\\
0.928818492915338 0.448002265422398\\
0.933912861197256 0.448002265422398\\
0.940823611486507 0.448002265422398\\
0.945911556390755 0.44881931841162\\
0.958045068497593 0.455670681028281\\
0.966649019659788 0.465154599497079\\
0.979837079883804 0.473846743892405\\
0.986832803561379 0.474211578554334\\
0.993204950336548 0.474211578554334\\
1 1\\
};
\addplot [
color=mycolor3,
solid,
forget plot
]
table[row sep=crcr]{
0 0\\
0.00497497637116562 0.349756433112332\\
0.00984661077772231 0.349756433112332\\
0.0133773394681144 0.349756433112332\\
0.019041221571737 0.349756433112332\\
0.0231409445716745 0.349756433112332\\
0.0301502280027306 0.349756433112332\\
0.0391174022472687 0.349756433112332\\
0.0477466875258334 0.349756433112332\\
0.0570054967863128 0.349756433112332\\
0.0645811358944621 0.349756433112332\\
0.0685330133105946 0.349756433112332\\
0.0820990152652959 0.349756433112332\\
0.0894562286100145 0.349756433112332\\
0.102537744255566 0.349756433112332\\
0.11025640739743 0.349756433112332\\
0.120540101972681 0.349756433112332\\
0.129483820942727 0.349756433112332\\
0.135213301540134 0.349756433112332\\
0.142480502892372 0.349756433112332\\
0.151040905936066 0.349756433112332\\
0.156461223178824 0.349756433112332\\
0.163876857308721 0.349756433112332\\
0.16897475055232 0.349756433112332\\
0.174926555204572 0.349756433112332\\
0.181390095080658 0.349756433112332\\
0.187351089972685 0.349756433112332\\
0.196000104283152 0.349756433112332\\
0.204520146831781 0.349756433112332\\
0.213961742036706 0.349756433112332\\
0.223833213065366 0.349756433112332\\
0.235392609177027 0.349756433112332\\
0.243849743992681 0.349756433112332\\
0.252960414080877 0.349756433112332\\
0.258567483407977 0.349756433112332\\
0.265153367329844 0.349756433112332\\
0.268504657419517 0.349756433112332\\
0.276806087675452 0.349756433112332\\
0.2851987227546 0.349756433112332\\
0.298098765869008 0.349756433112332\\
0.303447536197692 0.349756433112332\\
0.315187496917858 0.349756433112332\\
0.321265642705109 0.349756433112332\\
0.330073342142676 0.349756433112332\\
0.335677268558179 0.349756433112332\\
0.341890901904959 0.349756433112332\\
0.345134810004594 0.349756433112332\\
0.357406964372535 0.349756433112332\\
0.364514910149972 0.349756433112332\\
0.371327062272523 0.349756433112332\\
0.374168272276154 0.349756433112332\\
0.384377289046401 0.349756433112332\\
0.396646810285164 0.349756433112332\\
0.401743638987335 0.349756433112332\\
0.407457160041521 0.349756433112332\\
0.419237137962321 0.349756433112332\\
0.42757770381846 0.349756433112332\\
0.4350189001216 0.349756433112332\\
0.442453369979222 0.349756433112332\\
0.452186231538485 0.349756433112332\\
0.458519191442692 0.349756433112332\\
0.46541771924758 0.349756433112332\\
0.471890444818962 0.349756433112332\\
0.487173402037365 0.349756433112332\\
0.498240164969678 0.349756433112332\\
0.507043485188985 0.349756433112332\\
0.513716704783173 0.349756433112332\\
0.520625297290398 0.349756433112332\\
0.532132742782327 0.349756433112332\\
0.540003030160306 0.349756433112332\\
0.544644329015756 0.349756433112332\\
0.554742610952711 0.349756433112332\\
0.560817744637076 0.349756433112332\\
0.569452422396106 0.349756433112332\\
0.577823348674442 0.349756433112332\\
0.585323686164029 0.349756433112332\\
0.594880176738806 0.349756433112332\\
0.60555593677969 0.349756433112332\\
0.611003987301917 0.349756433112332\\
0.618163231427581 0.349756433112332\\
0.624285987114688 0.349756433112332\\
0.629343032626208 0.349756433112332\\
0.642087781938242 0.349756433112332\\
0.647246077779821 0.349756433112332\\
0.652424174950889 0.349756433112332\\
0.659906518792492 0.349756433112332\\
0.664759121761405 0.349756433112332\\
0.670359316556265 0.349756433112332\\
0.677025279892209 0.349756433112332\\
0.684226533280522 0.349756433112332\\
0.699484115110521 0.349756433112332\\
0.711436397494711 0.349756433112332\\
0.717872236193634 0.349756433112332\\
0.728406536854225 0.349756433112332\\
0.737601802937662 0.349756433112332\\
0.748634860504002 0.349756433112332\\
0.75433158787469 0.349756433112332\\
0.761112401298204 0.349756433112332\\
0.766701847152416 0.349756433112332\\
0.774674682225999 0.349756433112332\\
0.781149350858924 0.349756433112332\\
0.79231615699356 0.349756433112332\\
0.796703485430832 0.349756433112332\\
0.807617064697133 0.349756433112332\\
0.817431297510114 0.349756433112332\\
0.824724168227048 0.349756433112332\\
0.829519674169076 0.349756433112332\\
0.841214524435255 0.349756433112332\\
0.846745097165574 0.349756433112332\\
0.854289498013312 0.349756433112332\\
0.85883759502499 0.349756433112332\\
0.869059713110833 0.349756433112332\\
0.878781324770371 0.349756433112332\\
0.887306524560119 0.349756433112332\\
0.90340008717167 0.349756433112332\\
0.911399976936214 0.349756433112332\\
0.918356743426828 0.349756433112332\\
0.928165880853685 0.349756433112332\\
0.934211584438851 0.349756433112332\\
0.939003026159948 0.349756433112332\\
0.945975231654029 0.349756433112332\\
0.954477232795886 0.349756433112332\\
0.959579569911414 0.350212719996605\\
0.965997977712348 0.350850102423352\\
0.97363189453544 0.351095328429521\\
0.981919295683996 0.355317274326161\\
0.99021163389081 0.355317274326161\\
0.995053435433604 0.355317274326161\\
1 1\\
};
\end{axis}
\end{tikzpicture}%

%% file: ScalingInEpsilon_Model1_1D_p1p5.tikz
%
%
%
%
\begin{tikzpicture}
\def\xl{4}
\def\xu{9}
\def\yl{-7}
\def\yu{-1}

\begin{axis}[%
width=\figurewidth,
height=\figureheight,
scale only axis,
xmin={\xl-(\xu-\xl)/5},
xmax={\xu+(\xu-\xl)/5},
ymin={\yl-(\yu-\yl)/5},
ymax={\yu+(\yu-\yl)*0.05},
hide axis,
axis background/.style={fill=white!100}
]
\draw [->] (axis cs: \xl,\yl) -- (axis cs: {\xu+(\xu-\xl)*0.05},\yl);
\draw [->] (axis cs: \xl,\yl) -- (axis cs: \xl,{\yu+(\yu-\yl)*0.05});
\node at (axis cs: {(\xu+(\xu-\xl)*0.05+\xl)/2},{\yl-0.16*(\yu-\yl)}) {$\log(n)$};
\node[rotate=90] at (axis cs: {\xl-0.16*(\xu-\xl)},{(\yu+(\yu-\yl)*0.05+\yl)/2}) {$\log(\eps)$};
\foreach \yValue in {-7,-6,-5,-4,-3,-2,-1} {
    \edef\temp{\noexpand\draw [-] ({\xl+(\xu-\xl)/100},\yValue) -- (\xl,\yValue) node[left] {\yValue};} 
    \temp
}
\foreach \xValue in {4,5,6,7,8,9} {
	\edef\temp{\noexpand\draw [-] (\xValue,{\yl+(\yu-\yl)/100}) -- (\xValue,\yl) node[below] {\xValue};}    
    \temp
}
\addplot [
color=blue,
solid,
mark size=2.0pt,
mark=*,
mark options={solid},
only marks,
forget plot
]
table[row sep=crcr]{
4.382026634673881 -2.716246968306389\\
5.075173815233827 -3.326506489060566\\
5.768320995793772 -3.862261495869082\\
6.461468176353717 -4.471638793363569\\
7.154615356913663 -5.079263800485352\\
7.847762537473608 -5.711979241860379\\
8.540909718033554 -6.288716070575914\\
};
\addplot [
color=darkred,
solid,
mark size=2.0pt,
mark=square*,
mark options={solid},
only marks,
forget plot
]
table[row sep=crcr]{
4.38202663467388 -2.49557560613757\\
5.07517381523383 -2.95965621708018\\
5.76832099579377 -3.4997782103346\\
6.46146817635372 -4.02535169073515\\
7.15461535691366 -4.65753817150541\\
7.84776253747361 -5.10248209058736\\
8.54090971803355 -5.66967686216969\\
};
\addplot [
color=orange,
solid,
mark size=3.0pt,
mark=triangle*,
mark options={solid},
only marks,
forget plot
]
table[row sep=crcr]{
4.38202663467388 -1.71253342146107\\
5.07517381523383 -2.13568800652572\\
5.76832099579377 -3.18872278669751\\
6.46146817635372 -3.94839064959902\\
7.15461535691366 -4.1861913587132\\
7.84776253747361 -4.70826569501507\\
8.54090971803355 -5.17555807090543\\
};
\addplot [
color=blue,
dashed,
line width=2.5pt,
forget plot
]
table[row sep=crcr]{
4.382026634673881 -2.645472041266671\\
5.075173815233827 -3.254797001057718\\
5.768320995793772 -3.864121960848765\\
6.461468176353717 -4.473446920639812\\
7.154615356913663 -5.082771880430860\\
7.847762537473608 -5.692096840221907\\
8.540909718033554 -6.301421800012954\\
};
\addplot [
color=darkred,
dashed,
line width=1.5pt,
forget plot
]
table[row sep=crcr]{
4.38202663467388 -2.42419432365749\\
5.07517381523383 -2.96588709400973\\
5.76832099579377 -3.50757986436197\\
6.46146817635372 -4.04927263471421\\
7.15461535691366 -4.59096540506644\\
7.84776253747361 -5.13265817541868\\
8.54090971803355 -5.67435094577092\\
};
\addplot [
color=orange,
dashed,
line width=1.5pt,
forget plot
]
table[row sep=crcr]{
4.38202663467388 -2.3480074666533\\
5.07517381523383 -2.82136202803649\\
5.76832099579377 -3.29471658941967\\
6.46146817635372 -3.76807115080286\\
7.15461535691366 -4.24142571218605\\
7.84776253747361 -4.71478027356923\\
8.54090971803355 -5.18813483495242\\
};
\end{axis}
\end{tikzpicture}%

%% file: Error_n1280_Model1_1D_p2.tikz
\begin{tikzpicture}
\def\xl{0}
\def\xu{0.05}
\def\yl{0}
\def\yu{0.5}

\begin{axis}[%
width=\figurewidth,
height=\figureheight,
scale only axis,
xmin={\xl-(\xu-\xl)/5},
xmax={\xu+(\xu-\xl)/5},
ymin={\yl-(\yu-\yl)/5},
ymax={\yu+(\yu-\yl)*0.05},
hide axis,
axis background/.style={fill=white!100}
]
\draw [->] (axis cs: \xl,\yl) -- (axis cs: {\xu+(\xu-\xl)*0.05},\yl);
\draw [->] (axis cs: \xl,\yl) -- (axis cs: \xl,{\yu+(\yu-\yl)*0.05});
\draw [blue,->] (axis cs: \xu,\yl) -- (axis cs: \xu,{\yu+(\yu-\yl)*0.05});
\node at (axis cs: {(\xu+(\xu-\xl)*0.05+\xl)/2},{\yl-0.16*(\yu-\yl)}) {$\eps$};
\node[rotate=90] at (axis cs: {\xl-0.16*(\xu-\xl)},{(\yu+(\yu-\yl)*0.05+\yl)/2}) {$\err_n^{(2)}(f_n)$};
\node[blue,rotate=90] at (axis cs: {\xu+0.16*(\xu-\xl)},{(\yu+(\yu-\yl)*0.05+\yl)/2}) {$\%$ Graphs Connected};
\foreach \yValue in {0,0.1,0.2,0.3,0.4,0.5} {
    \edef\temp{\noexpand\draw [-] ({\xl+(\xu-\xl)/100},\yValue) -- (\xl,\yValue) node[left] {\yValue};} 
    \temp
}
\foreach \yValue in {20,40,60,80,100} {
    \edef\temp{\noexpand\draw [blue,-] ({\xu-(\xu-\xl)/100},\yValue/200) -- (\xu,\yValue/200) node[right] {\yValue};} 
    \temp
}
\foreach \xValue in {0,0.01,0.02,0.03,0.04,0.05} {
	\edef\temp{\noexpand\draw [-] (\xValue,{\yl+(\yu-\yl)/100}) -- (\xValue,\yl) node[below] {\xValue};}    
    \temp
}
\addplot [
color=black,
solid,
mark options={solid},
line width=1.5pt,
forget plot
]
table[row sep=crcr]{
0.005 0.35923359453563\\
0.00591836734693878 0.178874778793678\\
0.00683673469387755 0.0685402642725013\\
0.00775510204081633 0.0274101170924269\\
0.0086734693877551 0.0147094194341576\\
0.00959183673469388 0.00798070820203186\\
0.0105102040816327 0.00499334657114678\\
0.0114285714285714 0.00568347952955347\\
0.0123469387755102 0.00626294767045227\\
0.013265306122449 0.00726331802510508\\
0.0141836734693878 0.00791386708080982\\
0.0151020408163265 0.00923822080820862\\
0.0160204081632653 0.0103816663724375\\
0.0169387755102041 0.0113839385244911\\
0.0178571428571429 0.0126669898810991\\
0.0187755102040816 0.0138912094059337\\
0.0196938775510204 0.0153690173682692\\
0.0206122448979592 0.0161297759014887\\
0.021530612244898 0.0180455605545357\\
0.0224489795918367 0.01963624985336\\
0.0233673469387755 0.0211198896964078\\
0.0242857142857143 0.0223518099933888\\
0.0252040816326531 0.0239612371159856\\
0.0261224489795918 0.0260240489845877\\
0.0270408163265306 0.0280872301276236\\
0.0279591836734694 0.0305431682300739\\
0.0288775510204082 0.0323942003624389\\
0.0297959183673469 0.034799732108701\\
0.0307142857142857 0.0423229175437825\\
0.0316326530612245 0.0466411464524435\\
0.0325510204081633 0.0489141291654087\\
0.033469387755102 0.0588374649722734\\
0.0343877551020408 0.0759084734443977\\
0.0353061224489796 0.0921620746734313\\
0.0362244897959184 0.111207505354159\\
0.0371428571428571 0.139096570280744\\
0.0380612244897959 0.167615157703657\\
0.0389795918367347 0.195225252247171\\
0.0398979591836735 0.214408857555729\\
0.0408163265306122 0.22553155640841\\
0.041734693877551 0.234911056613212\\
0.0426530612244898 0.243204694183776\\
0.0435714285714286 0.247922992689848\\
0.0444897959183674 0.24920342178815\\
0.0454081632653061 0.250411408692037\\
0.0463265306122449 0.25014871977874\\
0.0472448979591837 0.249851250721562\\
0.0481632653061225 0.249680299975588\\
0.0490816326530612 0.249471037981306\\
0.05 0.249257914652189\\
};
\addplot [
color=black,
dashed,
line width=1pt,
forget plot
]
table[row sep=crcr]{
0.005 0.14092784071431\\
0.00591836734693878 0.002481874372143\\
0.00683673469387755 0.00284747530012385\\
0.00775510204081633 0.00324643473715757\\
0.0086734693877551 0.00357827582830885\\
0.00959183673469388 0.00394907121563945\\
0.0105102040816327 0.00442919243935562\\
0.0114285714285714 0.00490769829798776\\
0.0123469387755102 0.00542959773669653\\
0.013265306122449 0.00601774621843581\\
0.0141836734693878 0.00664149161070162\\
0.0151020408163265 0.00735805238133553\\
0.0160204081632653 0.00827380785655626\\
0.0169387755102041 0.00885829418271865\\
0.0178571428571429 0.00978007659102725\\
0.0187755102040816 0.0106737580841098\\
0.0196938775510204 0.011595870014359\\
0.0206122448979592 0.0129835384040242\\
0.021530612244898 0.0142770404330864\\
0.0224489795918367 0.0154690536666754\\
0.0233673469387755 0.0170844906233164\\
0.0242857142857143 0.0181928102465016\\
0.0252040816326531 0.0198971191720192\\
0.0261224489795918 0.0218656371606336\\
0.0270408163265306 0.023655101635978\\
0.0279591836734694 0.0255626758485244\\
0.0288775510204082 0.0277044700993584\\
0.0297959183673469 0.0297047674367514\\
0.0307142857142857 0.0322354701751456\\
0.0316326530612245 0.0345657418452605\\
0.0325510204081633 0.0369205758662486\\
0.033469387755102 0.0399628549583085\\
0.0343877551020408 0.042470774184445\\
0.0353061224489796 0.0453086734652116\\
0.0362244897959184 0.0490161707929409\\
0.0371428571428571 0.0525641218224009\\
0.0380612244897959 0.0592031872806503\\
0.0389795918367347 0.0638993575445224\\
0.0398979591836735 0.0703129171580246\\
0.0408163265306122 0.0812604902632812\\
0.041734693877551 0.170818018482218\\
0.0426530612244898 0.244631478776035\\
0.0435714285714286 0.245605272145196\\
0.0444897959183674 0.245230193671353\\
0.0454081632653061 0.245104736800309\\
0.0463265306122449 0.245098662083886\\
0.0472448979591837 0.244454935721485\\
0.0481632653061225 0.244229659112379\\
0.0490816326530612 0.243955978020833\\
0.05 0.243928552511108\\
};
\addplot [
color=black,
dashed,
line width=1pt,
forget plot
]
table[row sep=crcr]{
0.005 0.534544642632926\\
0.00591836734693878 0.47443680346915\\
0.00683673469387755 0.354103390532575\\
0.00775510204081633 0.00401484245952688\\
0.0086734693877551 0.00446272452119644\\
0.00959183673469388 0.00503825335675189\\
0.0105102040816327 0.00556519223293461\\
0.0114285714285714 0.00626510984968964\\
0.0123469387755102 0.0069450297424329\\
0.013265306122449 0.00953673866683511\\
0.0141836734693878 0.00911957182632583\\
0.0151020408163265 0.0124937593300707\\
0.0160204081632653 0.0149923708114992\\
0.0169387755102041 0.0158590264170441\\
0.0178571428571429 0.0173998291128343\\
0.0187755102040816 0.0187609293460301\\
0.0196938775510204 0.021006532113933\\
0.0206122448979592 0.0203607955333577\\
0.021530612244898 0.0232430936765233\\
0.0224489795918367 0.0255067668338743\\
0.0233673469387755 0.0264867639354141\\
0.0242857142857143 0.0274869667169244\\
0.0252040816326531 0.0294736293597386\\
0.0261224489795918 0.032748849377639\\
0.0270408163265306 0.0350379051243412\\
0.0279591836734694 0.0376869057059823\\
0.0288775510204082 0.0392582800454824\\
0.0297959183673469 0.042720946914725\\
0.0307142857142857 0.0471970218162519\\
0.0316326530612245 0.0502808910187421\\
0.0325510204081633 0.0531071165168538\\
0.033469387755102 0.0594112152895213\\
0.0343877551020408 0.248456605097314\\
0.0353061224489796 0.251805413087787\\
0.0362244897959184 0.253420698816721\\
0.0371428571428571 0.256060479161836\\
0.0380612244897959 0.256635339365989\\
0.0389795918367347 0.256816008190549\\
0.0398979591836735 0.256648695440413\\
0.0408163265306122 0.256599218981885\\
0.041734693877551 0.256355507968854\\
0.0426530612244898 0.256153431441457\\
0.0435714285714286 0.256030214574561\\
0.0444897959183674 0.255878767872122\\
0.0454081632653061 0.25499473465156\\
0.0463265306122449 0.254746717608536\\
0.0472448979591837 0.254448551628818\\
0.0481632653061225 0.254060354796291\\
0.0490816326530612 0.253855234858472\\
0.05 0.25368443750695\\
};
\addplot [
color=orange,
mark size=3.0pt,
only marks,
mark=triangle*,
mark options={solid},
forget plot
]
table[row sep=crcr]{
0.0398979591836735 0.214408857555729\\
};
\addplot [
color=blue,
mark size=2.0pt,
only marks,
mark=*,
mark options={solid},
forget plot
]
table[row sep=crcr]{
0.0068 0.0685\\
};
\addplot [
color=darkred,
mark size=2.0pt,
only marks,
mark=square*,
mark options={solid},
forget plot
]
table[row sep=crcr]{
0.0105102040816327 0.00499334657114678\\
};
\addplot [
color=blue,
solid,
mark options={solid},
line width=1.5pt,
forget plot
]
table[row sep=crcr]{
0.0050 0.0500\\
0.0059 0.2700\\
0.0068 0.4150\\
0.0078 0.4700\\
0.0087 0.4850\\
0.0096 0.4950\\
0.0105 0.5000\\
0.0114 0.5000\\
0.0123 0.5000\\
0.0133 0.5000\\
0.0142 0.5000\\
0.0151 0.5000\\
0.0160 0.5000\\
0.0169 0.5000\\
0.0179 0.5000\\
0.0188 0.5000\\
0.0197 0.5000\\
0.0206 0.5000\\
0.0215 0.5000\\
0.0224 0.5000\\
0.0234 0.5000\\
0.0243 0.5000\\
0.0252 0.5000\\
0.0261 0.5000\\
0.0270 0.5000\\
0.0280 0.5000\\
0.0289 0.5000\\
0.0298 0.5000\\
0.0307 0.5000\\
0.0316 0.5000\\
0.0326 0.5000\\
0.0335 0.5000\\
0.0344 0.5000\\
0.0353 0.5000\\
0.0362 0.5000\\
0.0371 0.5000\\
0.0381 0.5000\\
0.0390 0.5000\\
0.0399 0.5000\\
0.0408 0.5000\\
0.0417 0.5000\\
0.0427 0.5000\\
0.0436 0.5000\\
0.0445 0.5000\\
0.0454 0.5000\\
0.0463 0.5000\\
0.0472 0.5000\\
0.0482 0.5000\\
0.0491 0.5000\\
0.0500 0.5000\\
};
\end{axis}
\end{tikzpicture}%

%% file: ScalingInEpsilon_Model1_1D_p2.tikz
%
%
%
%
\begin{tikzpicture}
\def\xl{4}
\def\xu{9}
\def\yl{-7}
\def\yu{-1}

\begin{axis}[%
width=\figurewidth,
height=\figureheight,
scale only axis,
xmin={\xl-(\xu-\xl)/5},
xmax={\xu+(\xu-\xl)/5},
ymin={\yl-(\yu-\yl)/5},
ymax={\yu+(\yu-\yl)*0.05},
hide axis,
axis background/.style={fill=white!100}
]
\draw [->] (axis cs: \xl,\yl) -- (axis cs: {\xu+(\xu-\xl)*0.05},\yl);
\draw [->] (axis cs: \xl,\yl) -- (axis cs: \xl,{\yu+(\yu-\yl)*0.05});
\node at (axis cs: {(\xu+(\xu-\xl)*0.05+\xl)/2},{\yl-0.16*(\yu-\yl)}) {$\log(n)$};
\node[rotate=90] at (axis cs: {\xl-0.16*(\xu-\xl)},{(\yu+(\yu-\yl)*0.05+\yl)/2}) {$\log(\eps)$};
\foreach \yValue in {-7,-6,-5,-4,-3,-2,-1} {
    \edef\temp{\noexpand\draw [-] ({\xl+(\xu-\xl)/100},\yValue) -- (\xl,\yValue) node[left] {\yValue};} 
    \temp
}
\foreach \xValue in {4,5,6,7,8,9} {
	\edef\temp{\noexpand\draw [-] (\xValue,{\yl+(\yu-\yl)/100}) -- (\xValue,\yl) node[below] {\xValue};}    
    \temp
}
\addplot [
color=orange,
solid,
mark size=3.0pt,
mark=triangle*,
mark options={solid},
only marks,
forget plot
]
table[row sep=crcr]{
4.38202663467388 -1.934546390405\\
5.07517381523383 -2.12708950127049\\
5.76832099579377 -2.58078842149128\\
6.46146817635372 -2.89302834472473\\
7.15461535691366 -3.2214301046733\\
7.84776253747361 -3.59023532048985\\
8.54090971803355 -4.00947834157886\\
};
\addplot [
color=blue,
solid,
mark size=2.0pt,
mark=*,
mark options={solid},
only marks,
forget plot
]
table[row sep=crcr]{
4.3820 -2.7162\\
5.0752 -3.3612\\
5.7683 -3.8430\\
6.4615 -4.4539\\
7.1546 -4.9854\\
7.8478 -5.6997\\
8.5409 -6.2249\\
};
\addplot [
color=darkred,
solid,
mark size=2.0pt,
mark=square*,
mark options={solid},
only marks,
forget plot
]
table[row sep=crcr]{
4.38202663467388 -2.4567357728213\\
5.07517381523383 -3.00805275794203\\
5.76832099579377 -3.41558611911426\\
6.46146817635372 -4.00275185881791\\
7.15461535691366 -4.55540867642903\\
7.84776253747361 -5.1092161227687\\
8.54090971803355 -5.66967686216969\\
};
\addplot [
color=blue,
dashed,
line width=2.5pt,
forget plot
]
table[row sep=crcr]{
4.3820 -2.6376\\
5.0752 -3.2386\\
5.7683 -3.8395\\
6.4615 -4.4405\\
7.1546 -5.0414\\
7.8478 -5.6423\\
8.5409 -6.2433\\
};
\addplot [
color=darkred,
dashed,
line width=1.5pt,
forget plot
]
table[row sep=crcr]{
4.38202663467388 -2.30466962783525\\
5.07517381523383 -2.86613420284142\\
5.76832099579377 -3.42759877784758\\
6.46146817635372 -3.98906335285375\\
7.15461535691366 -4.55052792785992\\
7.84776253747361 -5.11199250286608\\
8.54090971803355 -5.67345707787225\\
};
\addplot [
color=orange,
dashed,
line width=1.5pt,
forget plot
]
table[row sep=crcr]{
4.38202663467388 -1.83715738021549\\
5.07517381523383 -2.19261606180952\\
5.76832099579377 -2.54807474340355\\
6.46146817635372 -2.90353342499758\\
7.15461535691366 -3.2589921065916\\
7.84776253747361 -3.61445078818563\\
8.54090971803355 -3.96990946977966\\
};
\end{axis}
\end{tikzpicture}%

%% file: SmoothEpsConnVError_n1280_Model1_1D_p1p5.tikz
\begin{tikzpicture}
\def\xl{-0.006}
\def\xu{0.02}
\def\yl{0}
\def\yu{0.6}

\begin{axis}[%
width=\figurewidth,
height=\figureheight,
scale only axis,
xmin={\xl-(\xu-\xl)/5},
xmax={\xu+(\xu-\xl)/5},
ymin={\yl-(\yu-\yl)/5},
ymax={\yu+(\yu-\yl)*0.05},
hide axis,
axis background/.style={fill=white!100}
]
\draw [->] (axis cs: \xl,\yl) -- (axis cs: {\xu+(\xu-\xl)*0.05},\yl);
\draw [->] (axis cs: \xl,\yl) -- (axis cs: \xl,{\yu+(\yu-\yl)*0.05});
\node at (axis cs: {(\xu+(\xu-\xl)*0.05+\xl)/2},{\yl-0.16*(\yu-\yl)}) {$\eps-\eps_{\mathrm{conn}}$};
\node[rotate=90] at (axis cs: {\xl-0.16*(\xu-\xl)},{(\yu+(\yu-\yl)*0.05+\yl)/2}) {$\err_n^{(1.5)}(f_n)$};
\foreach \yValue in {0,0.1,0.2,0.3,0.4,0.5,0.6} {
    \edef\temp{\noexpand\draw [-] ({\xl+(\xu-\xl)/100},\yValue) -- (\xl,\yValue) node[left] {\yValue};} 
    \temp
}
\foreach \xValue in {-0.005,0,0.005,0.01,0.015,0.02} {
	\edef\temp{\noexpand\draw [-] (\xValue,{\yl+(\yu-\yl)/100}) -- (\xValue,\yl) node[below] {\xValue};}    
    \temp
}
\addplot [
color=black,
solid,
mark options={solid},
line width=1.5pt,
forget plot
]
table[row sep=crcr]{
-0.00571428571428571 0.326868611257633\\
-0.00505494505494506 0.326868611257633\\
-0.0043956043956044 0.331560483160871\\
-0.00373626373626374 0.360169926381856\\
-0.00307692307692308 0.36134148581714\\
-0.00241758241758242 0.355222133709712\\
-0.00175824175824176 0.355978313989425\\
-0.0010989010989011 0.358641846789478\\
-0.00043956043956044 0.370017885507593\\
0.000219780219780219 0.110894895327264\\
0.000879120879120878 0.0818106730240808\\
0.00153846153846154 0.0554103233971518\\
0.0021978021978022 0.0482084196371768\\
0.00285714285714286 0.0343149150954402\\
0.00351648351648352 0.0345205160756685\\
0.00417582417582418 0.0417651847578749\\
0.00483516483516484 0.0486600185539101\\
0.00549450549450549 0.0552184805220445\\
0.00615384615384615 0.0715527711209715\\
0.00681318681318681 0.0849940364169783\\
0.00747252747252747 0.112109802574128\\
0.00813186813186813 0.139864948261156\\
0.00879120879120879 0.162006345451644\\
0.00945054945054945 0.193136599273604\\
0.0101098901098901 0.214418952194632\\
0.0107692307692308 0.235593455469308\\
0.0114285714285714 0.255249608676502\\
0.0120879120879121 0.257538615797715\\
0.0127472527472527 0.268657479293963\\
0.0134065934065934 0.271868235636371\\
0.0140659340659341 0.274429054329616\\
0.0147252747252747 0.273963899808848\\
0.0153846153846154 0.275546487886914\\
0.016043956043956 0.276007621253123\\
0.0167032967032967 0.27762491321691\\
0.0173626373626374 0.277510132527514\\
0.018021978021978 0.277612932253669\\
0.0186813186813187 0.27735537216751\\
0.0193406593406593 0.27749290313727\\
0.02 0.273106982119822\\
};
\addplot [
color=black,
dashed,
line width=1pt,
forget plot
]
table[row sep=crcr]{
-0.00571428571428571 0.326868611257633\\
-0.00505494505494506 0.326868611257633\\
-0.0043956043956044 0.393213317643105\\
-0.00373626373626374 0.441306384141572\\
-0.00307692307692308 0.523838030641129\\
-0.00241758241758242 0.524123517366647\\
-0.00175824175824176 0.482274227364378\\
-0.0010989010989011 0.498144035427002\\
-0.00043956043956044 0.502540569466929\\
0.000219780219780219 0.36928649406381\\
0.000879120879120878 0.309410161368054\\
0.00153846153846154 0.187374217495101\\
0.0021978021978022 0.0546343933675904\\
0.00285714285714286 0.0482827858629013\\
0.00351648351648352 0.0529909439067465\\
0.00417582417582418 0.0584267814076059\\
0.00483516483516484 0.079941306978877\\
0.00549450549450549 0.0912963514507296\\
0.00615384615384615 0.252626133584753\\
0.00681318681318681 0.256725702624261\\
0.00747252747252747 0.263876825218553\\
0.00813186813186813 0.265840900976959\\
0.00879120879120879 0.267734949217783\\
0.00945054945054945 0.272632866541361\\
0.0101098901098901 0.271360193712037\\
0.0107692307692308 0.274542012559957\\
0.0114285714285714 0.276374343170981\\
0.0120879120879121 0.277839631617329\\
0.0127472527472527 0.278899308084108\\
0.0134065934065934 0.285597228838733\\
0.0140659340659341 0.286872745738879\\
0.0147252747252747 0.290638025778389\\
0.0153846153846154 0.286832204968464\\
0.016043956043956 0.287060273800372\\
0.0167032967032967 0.296467808012867\\
0.0173626373626374 0.294583376944975\\
0.018021978021978 0.289385722344079\\
0.0186813186813187 0.291879758911964\\
0.0193406593406593 0.290459198629167\\
0.02 0.278337174168239\\
};
\addplot [
color=black,
dashed,
line width=1pt,
forget plot
]
table[row sep=crcr]{
-0.00571428571428571 0.326868611257633\\
-0.00505494505494506 0.326868611257633\\
-0.0043956043956044 0.279291392485115\\
-0.00373626373626374 0.279291392485115\\
-0.00307692307692308 0.269669785938528\\
-0.00241758241758242 0.26894467218609\\
-0.00175824175824176 0.270721210781692\\
-0.0010989010989011 0.272621453167016\\
-0.00043956043956044 0.272686482452663\\
0.000219780219780219 0.0129209436646206\\
0.000879120879120878 0.0117463486031806\\
0.00153846153846154 0.0121508462232104\\
0.0021978021978022 0.0124240359477523\\
0.00285714285714286 0.013766228541326\\
0.00351648351648352 0.0149370523246561\\
0.00417582417582418 0.0166049373015936\\
0.00483516483516484 0.0182233971802154\\
0.00549450549450549 0.0192281259507584\\
0.00615384615384615 0.0201276117507191\\
0.00681318681318681 0.0244461950886395\\
0.00747252747252747 0.0290123482271903\\
0.00813186813186813 0.0327860076584147\\
0.00879120879120879 0.0357808123851624\\
0.00945054945054945 0.0422510436103017\\
0.0101098901098901 0.055536475318746\\
0.0107692307692308 0.0690947323235817\\
0.0114285714285714 0.257393693680829\\
0.0120879120879121 0.257644982418467\\
0.0127472527472527 0.263181636569473\\
0.0134065934065934 0.264534923053512\\
0.0140659340659341 0.265549869980384\\
0.0147252747252747 0.266346166206832\\
0.0153846153846154 0.266763195653944\\
0.016043956043956 0.267189955815156\\
0.0167032967032967 0.267413981840339\\
0.0173626373626374 0.267400854061731\\
0.018021978021978 0.266959580775504\\
0.0186813186813187 0.267533823995222\\
0.0193406593406593 0.267655911492533\\
0.02 0.269301584779951\\
};
\end{axis}
\end{tikzpicture}%

%% file: SmoothEpsConnVError_n1280_Model1_1D_p2.tikz
\begin{tikzpicture}
\def\xl{-0.006}
\def\xu{0.05}
\def\yl{0}
\def\yu{0.6}

\begin{axis}[%
width=\figurewidth,
height=\figureheight,
scale only axis,
xmin={\xl-(\xu-\xl)/5},
xmax={\xu+(\xu-\xl)/5},
ymin={\yl-(\yu-\yl)/5},
ymax={\yu+(\yu-\yl)*0.05},
hide axis,
axis background/.style={fill=white!100}
]
\draw [->] (axis cs: \xl,\yl) -- (axis cs: {\xu+(\xu-\xl)*0.05},\yl);
\draw [->] (axis cs: \xl,\yl) -- (axis cs: \xl,{\yu+(\yu-\yl)*0.05});
\node at (axis cs: {(\xu+(\xu-\xl)*0.05+\xl)/2},{\yl-0.16*(\yu-\yl)}) {$\eps-\eps_{\mathrm{conn}}$};
\node[rotate=90] at (axis cs: {\xl-0.16*(\xu-\xl)},{(\yu+(\yu-\yl)*0.05+\yl)/2}) {$\err_n^{(2)}(f_n)$};
\foreach \yValue in {0,0.1,0.2,0.3,0.4,0.5,0.6} {
    \edef\temp{\noexpand\draw [-] ({\xl+(\xu-\xl)/100},\yValue) -- (\xl,\yValue) node[left] {\yValue};} 
    \temp
}
\foreach \xValue in {0,0.01,0.02,0.03,0.04,0.05} {
	\edef\temp{\noexpand\draw [-] (\xValue,{\yl+(\yu-\yl)/100}) -- (\xValue,\yl) node[below] {\xValue};}    
    \temp
}
\addplot [
color=black,
solid,
mark options={solid},
line width=1.5pt,
forget plot
]
table[row sep=crcr]{
-0.00551020408163265 0.356319625051488\\
-0.00421507064364207 0.387910870445474\\
-0.00291993720565149 0.388039036433084\\
-0.00162480376766091 0.385744709566484\\
-0.00032967032967033 0.190539552818\\
0.000965463108320251 0.00343092335555912\\
0.00226059654631083 0.00410208807478225\\
0.00355572998430141 0.00494393023288898\\
0.00485086342229199 0.00540133270028949\\
0.00614599686028258 0.00637405346970963\\
0.00744113029827316 0.00769803585544885\\
0.00873626373626374 0.00927470118978609\\
0.0100313971742543 0.0114937151406981\\
0.0113265306122449 0.0127872451260853\\
0.0126216640502355 0.0144408594858558\\
0.0139167974882261 0.0159178631671229\\
0.0152119309262166 0.0186275332149004\\
0.0165070643642072 0.0207853736238825\\
0.0178021978021978 0.0227615907109629\\
0.0190973312401884 0.0255587562581845\\
0.020392464678179 0.0273516239732369\\
0.0216875981161696 0.0305996813223965\\
0.0229827315541601 0.034019402797075\\
0.0242778649921507 0.0427439595724703\\
0.0255729984301413 0.0476341792809506\\
0.0268681318681319 0.0595141280244367\\
0.0281632653061225 0.0782373975682226\\
0.029458398744113 0.105790151516745\\
0.0307535321821036 0.145941996459186\\
0.0320486656200942 0.188157841880998\\
0.0333437990580848 0.209789988506565\\
0.0346389324960754 0.230731473671727\\
0.0359340659340659 0.238080366168651\\
0.0372291993720565 0.24460108411022\\
0.0385243328100471 0.248969525132125\\
0.0398194662480377 0.250104482283101\\
0.0411145996860283 0.249741656123378\\
0.0424097331240188 0.249626894119501\\
0.0437048665620094 0.249362975957265\\
0.045 0.250781572022389\\
};
\addplot [
color=black,
dashed,
line width=1pt,
forget plot
]
table[row sep=crcr]{
-0.00551020408163265 0.356319625051488\\
-0.00421507064364207 0.526369430147787\\
-0.00291993720565149 0.54299052921176\\
-0.00162480376766091 0.52479380307027\\
-0.00032967032967033 0.497354488951177\\
0.000965463108320251 0.00428602554820995\\
0.00226059654631083 0.00510182859846033\\
0.00355572998430141 0.006022893865252\\
0.00485086342229199 0.00650568125465978\\
0.00614599686028258 0.00796556072090877\\
0.00744113029827316 0.0106003375118898\\
0.00873626373626374 0.0132520405978131\\
0.0100313971742543 0.0178399777117774\\
0.0113265306122449 0.0181310393187651\\
0.0126216640502355 0.0193254910887897\\
0.0139167974882261 0.0219437267057745\\
0.0152119309262166 0.0260273558057991\\
0.0165070643642072 0.0272799328059946\\
0.0178021978021978 0.0293168550707154\\
0.0190973312401884 0.0332242483916107\\
0.020392464678179 0.03604676644843\\
0.0216875981161696 0.0393771717141999\\
0.0229827315541601 0.0419207889379334\\
0.0242778649921507 0.0490473119520194\\
0.0255729984301413 0.0521281999583803\\
0.0268681318681319 0.0612544488180817\\
0.0281632653061225 0.249733411646164\\
0.029458398744113 0.253848837226081\\
0.0307535321821036 0.255800845916813\\
0.0320486656200942 0.256687487892756\\
0.0333437990580848 0.256508781969098\\
0.0346389324960754 0.256557440664151\\
0.0359340659340659 0.256386071793535\\
0.0372291993720565 0.256075161189611\\
0.0385243328100471 0.255112552470052\\
0.0398194662480377 0.25470086728225\\
0.0411145996860283 0.25412871008195\\
0.0424097331240188 0.254143021778791\\
0.0437048665620094 0.255021821904218\\
0.045 0.257188247599321\\
};
\addplot [
color=black,
dashed,
line width=1pt,
forget plot
]
table[row sep=crcr]{
-0.00551020408163265 0.356319625051488\\
-0.00421507064364207 0.291961201073158\\
-0.00291993720565149 0.289371556467731\\
-0.00162480376766091 0.290225590858739\\
-0.00032967032967033 0.00250534931281905\\
0.000965463108320251 0.00281671364106396\\
0.00226059654631083 0.00329809837815503\\
0.00355572998430141 0.0039682183022791\\
0.00485086342229199 0.00448071685428614\\
0.00614599686028258 0.0050658748560159\\
0.00744113029827316 0.0059801804874942\\
0.00873626373626374 0.00686825099945454\\
0.0100313971742543 0.00830552585451963\\
0.0113265306122449 0.00916731312814941\\
0.0126216640502355 0.0108280661244893\\
0.0139167974882261 0.0119579621204758\\
0.0152119309262166 0.0141156864247625\\
0.0165070643642072 0.0159923999051846\\
0.0178021978021978 0.0178974870369606\\
0.0190973312401884 0.0200402609898967\\
0.020392464678179 0.0218656371606336\\
0.0216875981161696 0.0248986420389301\\
0.0229827315541601 0.0280485736305431\\
0.0242778649921507 0.030741855041941\\
0.0255729984301413 0.034595260347265\\
0.0268681318681319 0.0371087849989068\\
0.0281632653061225 0.0409451918901744\\
0.029458398744113 0.0450783319582242\\
0.0307535321821036 0.0513494456514657\\
0.0320486656200942 0.0597553645641956\\
0.0333437990580848 0.0663738494300671\\
0.0346389324960754 0.0815114258092717\\
0.0359340659340659 0.243380271544802\\
0.0372291993720565 0.244758665895706\\
0.0385243328100471 0.244950610975333\\
0.0398194662480377 0.244877080113471\\
0.0411145996860283 0.244213427814701\\
0.0424097331240188 0.244202187742263\\
0.0437048665620094 0.24394865589929\\
0.045 0.245595456767387\\
};
\end{axis}
\end{tikzpicture}%

%% file: Error_n1280_Model2_1D_p2.tikz
\begin{tikzpicture}
\def\xl{0}
\def\xu{0.2}
\def\yl{0}
\def\yu{0.6}

\begin{axis}[%
width=\figurewidth,
height=\figureheight,
scale only axis,
xmin={\xl-(\xu-\xl)/5},
xmax={\xu+(\xu-\xl)/5},
ymin={\yl-(\yu-\yl)/5},
ymax={\yu+(\yu-\yl)*0.05},
hide axis,
axis background/.style={fill=white!100}
]
\draw [->] (axis cs: \xl,\yl) -- (axis cs: {\xu+(\xu-\xl)*0.05},\yl);
\draw [->] (axis cs: \xl,\yl) -- (axis cs: \xl,{\yu+(\yu-\yl)*0.05});
\draw [blue,->] (axis cs: \xu,\yl) -- (axis cs: \xu,{\yu+(\yu-\yl)*0.05});
\node at (axis cs: {(\xu+(\xu-\xl)*0.05+\xl)/2},{\yl-0.16*(\yu-\yl)}) {$\eps$};
\node[rotate=90] at (axis cs: {\xl-0.16*(\xu-\xl)},{(\yu+(\yu-\yl)*0.05+\yl)/2}) {$\err_n^{(2)}(f_n)$};
\node[blue,rotate=90] at (axis cs: {\xu+0.16*(\xu-\xl)},{(\yu+(\yu-\yl)*0.05+\yl)/2}) {$\%$ Graphs Connected};
\foreach \yValue in {0,0.1,0.2,0.3,0.4,0.5,0.6} {
    \edef\temp{\noexpand\draw [-] ({\xl+(\xu-\xl)/100},\yValue) -- (\xl,\yValue) node[left] {\yValue};} 
    \temp
}
\foreach \yValue in {20,40,60,80,100} {
    \edef\temp{\noexpand\draw [blue,-] ({\xu-(\xu-\xl)/100},\yValue*\yu/100) -- (\xu,\yValue*\yu/100) node[right] {\yValue};} 
    \temp
}
\foreach \xValue in {0,0.05,0.1,0.15,0.2} {
	\edef\temp{\noexpand\draw [-] (\xValue,{\yl+(\yu-\yl)/100}) -- (\xValue,\yl) node[below] {\xValue};}    
    \temp
}
\addplot [
color=black,
solid,
mark options={solid},
line width=1.5pt,
forget plot
]
table[row sep=crcr]{
0.004 0.399360175583313\\
0.008 0.0240297647037801\\
0.012 0.00800579426191366\\
0.016 0.0138214389390265\\
0.02 0.019767486508295\\
0.024 0.0245182281510002\\
0.028 0.0290724391226766\\
0.032 0.0330511253990308\\
0.036 0.0377251000317512\\
0.04 0.0421839866474933\\
0.044 0.0468434286850989\\
0.048 0.0507305899335914\\
0.052 0.0553775641794129\\
0.056 0.0592023736267586\\
0.06 0.0628177329235727\\
0.064 0.0666683954845118\\
0.068 0.0705871693145562\\
0.072 0.0744641728550719\\
0.076 0.0800934757864492\\
0.08 0.0860494383902287\\
0.084 0.0919681096931573\\
0.088 0.0961036366947307\\
0.092 0.10077181747395\\
0.096 0.106224673069998\\
0.1 0.109440769067757\\
0.104 0.112898286485043\\
0.108 0.116978117458009\\
0.112 0.12077936549673\\
0.116 0.124443700974651\\
0.12 0.129076480018226\\
0.124 0.133266451311748\\
0.128 0.137754610636614\\
0.132 0.141767636678411\\
0.136 0.145294445006412\\
0.14 0.149725302799491\\
0.144 0.154321840406989\\
0.148 0.157990413615735\\
0.152 0.161282045318845\\
0.156 0.164830490175289\\
0.16 0.169167762831431\\
0.164 0.172997847918422\\
0.168 0.175150878934648\\
0.172 0.179063562708379\\
0.176 0.18288728816897\\
0.18 0.18727222032827\\
0.184 0.191383988395447\\
0.188 0.195857224313684\\
0.192 0.200442399892372\\
0.196 0.205266658703943\\
0.2 0.210442023033041\\
};
\addplot [
color=black,
dashed,
line width=1pt,
forget plot
]
table[row sep=crcr]{
0.004 0.291547838909126\\
0.008 0.00376866332980507\\
0.012 0.00596934174314954\\
0.016 0.00868790292389969\\
0.02 0.0121266175353231\\
0.024 0.0156810620692378\\
0.028 0.01888548402947\\
0.032 0.0224699818473569\\
0.036 0.0271086990569988\\
0.04 0.0325688215341096\\
0.044 0.0361192290443366\\
0.048 0.0405679345650933\\
0.052 0.0455425076757271\\
0.056 0.0492204161766562\\
0.06 0.0533874802866064\\
0.064 0.0577743869810452\\
0.068 0.0609794112105881\\
0.072 0.0656819018639606\\
0.076 0.0716133892657272\\
0.08 0.0775355435061721\\
0.084 0.0855165408514899\\
0.088 0.0891116279310909\\
0.092 0.0921191664167059\\
0.096 0.0971599990648979\\
0.1 0.102712391738559\\
0.104 0.106516294825978\\
0.108 0.109713115256164\\
0.112 0.113831518318619\\
0.116 0.117897014167556\\
0.12 0.122486465874544\\
0.124 0.128230606529984\\
0.128 0.130733022118356\\
0.132 0.134962312629879\\
0.136 0.138378672981459\\
0.14 0.14185266277434\\
0.144 0.146195989767944\\
0.148 0.149887546269147\\
0.152 0.152081754506422\\
0.156 0.156499347236787\\
0.16 0.161387803705976\\
0.164 0.165472655070184\\
0.168 0.167262203173975\\
0.172 0.172508534400582\\
0.176 0.177108386613938\\
0.18 0.181586171559669\\
0.184 0.186237680261262\\
0.188 0.190946830386646\\
0.192 0.195558963045357\\
0.196 0.200283577014792\\
0.2 0.205436696851694\\
};
\addplot [
color=black,
dashed,
line width=1pt,
forget plot
]
table[row sep=crcr]{
0.004 0.548630276179698\\
0.008 0.00494584363282872\\
0.012 0.0129651819051006\\
0.016 0.0215801871397053\\
0.02 0.0273028660573232\\
0.024 0.03154959032979\\
0.028 0.0366151166703215\\
0.032 0.0406806532450794\\
0.036 0.0449735058612155\\
0.04 0.0516027836363204\\
0.044 0.0571134912902886\\
0.048 0.0609191421549057\\
0.052 0.0656419611932295\\
0.056 0.0692488873910527\\
0.06 0.0724645840564471\\
0.064 0.0757739363757367\\
0.068 0.0808852844745228\\
0.072 0.0841327468827202\\
0.076 0.0904903337983996\\
0.08 0.0941889412729805\\
0.084 0.0980238972178839\\
0.088 0.103122078052834\\
0.092 0.108363878088842\\
0.096 0.114106182041802\\
0.1 0.1164546315638\\
0.104 0.119593631411629\\
0.108 0.124569944841027\\
0.112 0.129636966591682\\
0.116 0.131487722914592\\
0.12 0.13652461127419\\
0.124 0.140692964243065\\
0.128 0.144234897026999\\
0.132 0.148577853029213\\
0.136 0.152309830605444\\
0.14 0.156181442305743\\
0.144 0.160766519567818\\
0.148 0.164698395403079\\
0.152 0.168015944667781\\
0.156 0.172336395766556\\
0.16 0.176844081527883\\
0.164 0.180082261048186\\
0.168 0.182267609829763\\
0.172 0.184845540055167\\
0.176 0.189514224599856\\
0.18 0.193877033359617\\
0.184 0.198465578267745\\
0.188 0.20156919060765\\
0.192 0.206229836428354\\
0.196 0.211155304587903\\
0.2 0.216462276743092\\
};
\addplot [
color=blue,
mark size=2.0pt,
only marks,
mark=*,
mark options={solid},
forget plot
]
table[row sep=crcr]{
0.008000000000000 0.024029764703780\\
};
\addplot [
color=yellow,
mark size=2.0pt,
only marks,
mark=diamond*,
mark options={solid},
forget plot
]
table[row sep=crcr]{
0.044 0.0468434286850989\\
};
\addplot [
color=blue,
solid,
mark options={solid},
line width=1.5pt,
forget plot
]
table[row sep=crcr]{
0.004000000000000 0\\
0.008000000000000 0.57000000000000\\
0.012000000000000 0.60000000000000\\
0.016000000000000 0.60000000000000\\
0.020000000000000 0.60000000000000\\
0.024000000000000 0.60000000000000\\
0.028000000000000 0.60000000000000\\
0.032000000000000 0.60000000000000\\
0.036000000000000 0.60000000000000\\
0.040000000000000 0.60000000000000\\
0.044000000000000 0.60000000000000\\
0.048000000000000 0.60000000000000\\
0.052000000000000 0.60000000000000\\
0.056000000000000 0.60000000000000\\
0.060000000000000 0.60000000000000\\
0.064000000000000 0.60000000000000\\
0.068000000000000 0.60000000000000\\
0.072000000000000 0.60000000000000\\
0.076000000000000 0.60000000000000\\
0.080000000000000 0.60000000000000\\
0.084000000000000 0.60000000000000\\
0.088000000000000 0.60000000000000\\
0.092000000000000 0.60000000000000\\
0.096000000000000 0.60000000000000\\
0.100000000000000 0.60000000000000\\
0.104000000000000 0.60000000000000\\
0.108000000000000 0.60000000000000\\
0.112000000000000 0.60000000000000\\
0.116000000000000 0.60000000000000\\
0.120000000000000 0.60000000000000\\
0.124000000000000 0.60000000000000\\
0.128000000000000 0.60000000000000\\
0.132000000000000 0.60000000000000\\
0.136000000000000 0.60000000000000\\
0.140000000000000 0.60000000000000\\
0.144000000000000 0.60000000000000\\
0.148000000000000 0.60000000000000\\
0.152000000000000 0.60000000000000\\
0.156000000000000 0.60000000000000\\
0.160000000000000 0.60000000000000\\
0.164000000000000 0.60000000000000\\
0.168000000000000 0.60000000000000\\
0.172000000000000 0.60000000000000\\
0.176000000000000 0.60000000000000\\
0.180000000000000 0.60000000000000\\
0.184000000000000 0.60000000000000\\
0.188000000000000 0.60000000000000\\
0.192000000000000 0.60000000000000\\
0.196000000000000 0.60000000000000\\
0.200000000000000 0.60000000000000\\
};
\end{axis}
\end{tikzpicture}%

%% file: Error_n1280_Model1_2D_p2.tikz
\begin{tikzpicture}
\def\xl{0.03}
\def\xu{0.2}
\def\yl{0}
\def\yu{0.4}

\begin{axis}[%
width=\figurewidth,
height=\figureheight,
scale only axis,
xmin={\xl-(\xu-\xl)/5},
xmax={\xu+(\xu-\xl)/5},
ymin={\yl-(\yu-\yl)/5},
ymax={\yu+(\yu-\yl)*0.05},
hide axis,
axis background/.style={fill=white!100}
]
\draw [->] (axis cs: \xl,\yl) -- (axis cs: {\xu+(\xu-\xl)*0.05},\yl);
\draw [->] (axis cs: \xl,\yl) -- (axis cs: \xl,{\yu+(\yu-\yl)*0.05});
\draw [blue,->] (axis cs: \xu,\yl) -- (axis cs: \xu,{\yu+(\yu-\yl)*0.05});
\node at (axis cs: {(\xu+(\xu-\xl)*0.05+\xl)/2},{\yl-0.16*(\yu-\yl)}) {$\eps$};
\node[rotate=90] at (axis cs: {\xl-0.16*(\xu-\xl)},{(\yu+(\yu-\yl)*0.05+\yl)/2}) {$\err_n^{(2)}(f_n)$};
\node[blue,rotate=90] at (axis cs: {\xu+0.16*(\xu-\xl)},{(\yu+(\yu-\yl)*0.05+\yl)/2}) {$\%$ Graphs Connected};
\foreach \yValue in {0,0.1,0.2,0.3,0.4} {
    \edef\temp{\noexpand\draw [-] ({\xl+(\xu-\xl)/100},\yValue) -- (\xl,\yValue) node[left] {\yValue};} 
    \temp
}
\foreach \yValue in {20,40,60,80,100} {
    \edef\temp{\noexpand\draw [blue,-] ({\xu-(\xu-\xl)/100},\yValue*\yu/100) -- (\xu,\yValue*\yu/100) node[right] {\yValue};} 
    \temp
}
\foreach \xValue in {0.03,0.06,0.09,0.12,0.15,0.18} {
	\edef\temp{\noexpand\draw [-] (\xValue,{\yl+(\yu-\yl)/100}) -- (\xValue,\yl) node[below] {\xValue};}    
    \temp
}
\addplot [
color=black,
solid,
mark options={solid},
line width=1.5pt,
forget plot
]
table[row sep=crcr]{
0.03 0.154473825733546\\
0.033469387755102 0.277075228599316\\
0.0369387755102041 0.238717081015403\\
0.0404081632653061 0.194698737998419\\
0.0438775510204082 0.174620732749205\\
0.0473469387755102 0.160920913408321\\
0.0508163265306122 0.151714594118157\\
0.0542857142857143 0.142082038224817\\
0.0577551020408163 0.134115623444389\\
0.0612244897959184 0.12560434241257\\
0.0646938775510204 0.118415956148628\\
0.0681632653061224 0.110327596480482\\
0.0716326530612245 0.10352290473955\\
0.0751020408163265 0.0976629385795629\\
0.0785714285714286 0.0926258919970774\\
0.0820408163265306 0.0868570801357794\\
0.0855102040816327 0.0822980838582127\\
0.0889795918367347 0.0776707917566752\\
0.0924489795918367 0.0733684943779934\\
0.0959183673469388 0.0694837288583903\\
0.0993877551020408 0.0659370615095196\\
0.102857142857143 0.0625808011102709\\
0.106326530612245 0.0595086070656042\\
0.109795918367347 0.0566684050652901\\
0.113265306122449 0.0539682946258416\\
0.116734693877551 0.0513340579955346\\
0.120204081632653 0.0489943484930732\\
0.123673469387755 0.0469726784830772\\
0.127142857142857 0.0449192007358224\\
0.130612244897959 0.0432306223984821\\
0.134081632653061 0.0415687901622146\\
0.137551020408163 0.0398723736859039\\
0.141020408163265 0.0383057314973846\\
0.144489795918367 0.0369013783489827\\
0.147959183673469 0.0355885480591165\\
0.151428571428571 0.0342664152191369\\
0.154897959183673 0.03311193707324\\
0.158367346938776 0.0320200962784642\\
0.161836734693878 0.0310219710739145\\
0.16530612244898 0.0300490504049847\\
0.168775510204082 0.0292586252846066\\
0.172244897959184 0.0284494245700793\\
0.175714285714286 0.0277081848845962\\
0.179183673469388 0.0269845675572336\\
0.18265306122449 0.0262731824070964\\
0.186122448979592 0.0256535591558541\\
0.189591836734694 0.0250521389504365\\
0.193061224489796 0.0244597859993872\\
0.196530612244898 0.0239715944112168\\
0.2 0.0234298966901685\\
};
\addplot [
color=black,
dashed,
line width=1pt,
forget plot
]
table[row sep=crcr]{
0.03 0.0695905710187722\\
0.033469387755102 0.192665743188714\\
0.0369387755102041 0.161352905066018\\
0.0404081632653061 0.139672519986906\\
0.0438775510204082 0.13635988671251\\
0.0473469387755102 0.125749515558758\\
0.0508163265306122 0.123283134459366\\
0.0542857142857143 0.116210604072617\\
0.0577551020408163 0.111860742677732\\
0.0612244897959184 0.103757253751892\\
0.0646938775510204 0.0996890115672994\\
0.0681632653061224 0.0932040876446983\\
0.0716326530612245 0.0879688127084623\\
0.0751020408163265 0.0852664567819008\\
0.0785714285714286 0.0784432042754419\\
0.0820408163265306 0.0732079743823847\\
0.0855102040816327 0.0693499347404843\\
0.0889795918367347 0.064946316541801\\
0.0924489795918367 0.0614451175403029\\
0.0959183673469388 0.0581436608558755\\
0.0993877551020408 0.0555528136163436\\
0.102857142857143 0.0527317709511996\\
0.106326530612245 0.050326088745001\\
0.109795918367347 0.0482740535149543\\
0.113265306122449 0.0467976579951572\\
0.116734693877551 0.0447488788107777\\
0.120204081632653 0.0424265835485494\\
0.123673469387755 0.0414961264779693\\
0.127142857142857 0.0396901045173159\\
0.130612244897959 0.0386087386979293\\
0.134081632653061 0.0374967626661161\\
0.137551020408163 0.0360732709886243\\
0.141020408163265 0.0347959321977711\\
0.144489795918367 0.0333999865618435\\
0.147959183673469 0.0322933110552715\\
0.151428571428571 0.0311569550487197\\
0.154897959183673 0.0299887164830934\\
0.158367346938776 0.0291941673029821\\
0.161836734693878 0.0283734010757324\\
0.16530612244898 0.0274146249677623\\
0.168775510204082 0.0267364289166666\\
0.172244897959184 0.0262805128619794\\
0.175714285714286 0.0253690460419458\\
0.179183673469388 0.0247641888536191\\
0.18265306122449 0.0242515667723063\\
0.186122448979592 0.0237669986100102\\
0.189591836734694 0.0232299955487561\\
0.193061224489796 0.0228580529829662\\
0.196530612244898 0.0223828518750508\\
0.2 0.0219723802110808\\
};
\addplot [
color=black,
dashed,
line width=1pt,
forget plot
]
table[row sep=crcr]{
0.03 0.223200464369571\\
0.033469387755102 0.376561391043664\\
0.0369387755102041 0.322384263172192\\
0.0404081632653061 0.249495969745129\\
0.0438775510204082 0.211698170177841\\
0.0473469387755102 0.192788658370948\\
0.0508163265306122 0.182042473988971\\
0.0542857142857143 0.17309574017666\\
0.0577551020408163 0.160497264656171\\
0.0612244897959184 0.150994531232244\\
0.0646938775510204 0.141710730772831\\
0.0681632653061224 0.130578464222958\\
0.0716326530612245 0.123721586947311\\
0.0751020408163265 0.116137078087014\\
0.0785714285714286 0.109727142123825\\
0.0820408163265306 0.103997256257226\\
0.0855102040816327 0.09668756177884\\
0.0889795918367347 0.0898437413208128\\
0.0924489795918367 0.0841125060151957\\
0.0959183673469388 0.0795787353976446\\
0.0993877551020408 0.0758619170836535\\
0.102857142857143 0.0713966938491349\\
0.106326530612245 0.0666761999226563\\
0.109795918367347 0.0635822323461074\\
0.113265306122449 0.0606137282382054\\
0.116734693877551 0.0577462881668123\\
0.120204081632653 0.0548740475512694\\
0.123673469387755 0.0528747398028701\\
0.127142857142857 0.0505349222828496\\
0.130612244897959 0.048218674099852\\
0.134081632653061 0.0457654342876038\\
0.137551020408163 0.0434736158733192\\
0.141020408163265 0.0417831884459227\\
0.144489795918367 0.040619052719844\\
0.147959183673469 0.0390754300386602\\
0.151428571428571 0.03754699059604\\
0.154897959183673 0.0360129046981047\\
0.158367346938776 0.0346272795426091\\
0.161836734693878 0.0335506111190413\\
0.16530612244898 0.0323901112451801\\
0.168775510204082 0.031613718076088\\
0.172244897959184 0.0306821631607694\\
0.175714285714286 0.0295982723150326\\
0.179183673469388 0.0288824074395526\\
0.18265306122449 0.0279890104797974\\
0.186122448979592 0.0273147985270478\\
0.189591836734694 0.0265260138083508\\
0.193061224489796 0.026041250114112\\
0.196530612244898 0.0252161711949748\\
0.2 0.0246985229586434\\
};
\addplot [
color=blue,
mark size=2.0pt,
only marks,
mark=*,
mark options={solid},
forget plot
]
table[row sep=crcr]{
0.0508163265306122 0.151714594118157\\
};
\addplot [
color=yellow,
mark size=2.0pt,
only marks,
mark=diamond*,
mark options={solid},
forget plot
]
table[row sep=crcr]{
0.061224489795918 0.125604342412570\\
};
\addplot [
color=blue,
solid,
mark options={solid},
line width=1.5pt,
forget plot
]
table[row sep=crcr]{
0.030000000000000 0\\
0.033469387755102 0\\
0.036938775510204 0\\
0.040408163265306 0\\
0.043877551020408 0.032000000000000\\
0.047346938775510 0.116000000000000\\
0.050816326530612 0.260000000000000\\
0.054285714285714 0.312000000000000\\
0.057755102040816 0.360000000000000\\
0.061224489795918 0.388000000000000\\
0.064693877551020 0.400000000000000\\
0.068163265306122 0.400000000000000\\
0.071632653061224 0.400000000000000\\
0.075102040816327 0.400000000000000\\
0.078571428571429 0.400000000000000\\
0.082040816326531 0.400000000000000\\
0.085510204081633 0.400000000000000\\
0.088979591836735 0.400000000000000\\
0.092448979591837 0.400000000000000\\
0.095918367346939 0.400000000000000\\
0.099387755102041 0.400000000000000\\
0.102857142857143 0.400000000000000\\
0.106326530612245 0.400000000000000\\
0.109795918367347 0.400000000000000\\
0.113265306122449 0.400000000000000\\
0.116734693877551 0.400000000000000\\
0.120204081632653 0.400000000000000\\
0.123673469387755 0.400000000000000\\
0.127142857142857 0.400000000000000\\
0.130612244897959 0.400000000000000\\
0.134081632653061 0.400000000000000\\
0.137551020408163 0.400000000000000\\
0.141020408163265 0.400000000000000\\
0.144489795918367 0.400000000000000\\
0.147959183673469 0.400000000000000\\
0.151428571428571 0.400000000000000\\
0.154897959183673 0.400000000000000\\
0.158367346938776 0.400000000000000\\
0.161836734693878 0.400000000000000\\
0.165306122448980 0.400000000000000\\
0.168775510204082 0.400000000000000\\
0.172244897959184 0.400000000000000\\
0.175714285714286 0.400000000000000\\
0.179183673469388 0.400000000000000\\
0.182653061224490 0.400000000000000\\
0.186122448979592 0.400000000000000\\
0.189591836734694 0.400000000000000\\
0.193061224489796 0.400000000000000\\
0.196530612244898 0.400000000000000\\
0.200000000000000 0.400000000000000\\
};
\end{axis}
\end{tikzpicture}%

%% file: Error_n1280_Model1_2D_p4.tikz
\begin{tikzpicture}
\def\xl{0.03}
\def\xu{0.3}
\def\yl{0}
\def\yu{0.4}

\begin{axis}[%
width=\figurewidth,
height=\figureheight,
scale only axis,
xmin={\xl-(\xu-\xl)/5},
xmax={\xu+(\xu-\xl)/5},
ymin={\yl-(\yu-\yl)/5},
ymax={\yu+(\yu-\yl)*0.05},
hide axis,
axis background/.style={fill=white!100}
]
\draw [->] (axis cs: \xl,\yl) -- (axis cs: {\xu+(\xu-\xl)*0.05},\yl);
\draw [->] (axis cs: \xl,\yl) -- (axis cs: \xl,{\yu+(\yu-\yl)*0.05});
\draw [blue,->] (axis cs: \xu,\yl) -- (axis cs: \xu,{\yu+(\yu-\yl)*0.05});
\node at (axis cs: {(\xu+(\xu-\xl)*0.05+\xl)/2},{\yl-0.16*(\yu-\yl)}) {$\eps$};
\node[rotate=90] at (axis cs: {\xl-0.16*(\xu-\xl)},{(\yu+(\yu-\yl)*0.05+\yl)/2}) {$\err_n^{(4)}(f_n)$};
\node[blue,rotate=90] at (axis cs: {\xu+0.16*(\xu-\xl)},{(\yu+(\yu-\yl)*0.05+\yl)/2}) {$\%$ Graphs Connected};
\foreach \yValue in {0,0.1,0.2,0.3,0.4,0.5,0.6} {
    \edef\temp{\noexpand\draw [-] ({\xl+(\xu-\xl)/100},\yValue) -- (\xl,\yValue) node[left] {\yValue};} 
    \temp
}
\foreach \yValue in {20,40,60,80,100} {
    \edef\temp{\noexpand\draw [blue,-] ({\xu-(\xu-\xl)/100},\yValue*\yu/100) -- (\xu,\yValue*\yu/100) node[right] {\yValue};} 
    \temp
}
\foreach \xValue in {0.05,0.1,0.15,0.2,0.25,0.3} {
	\edef\temp{\noexpand\draw [-] (\xValue,{\yl+(\yu-\yl)/100}) -- (\xValue,\yl) node[below] {\xValue};}    
    \temp
}
\addplot [
color=black,
solid,
mark options={solid},
line width=1.5pt,
forget plot
]
table[row sep=crcr]{
0.03 0.32147522070433\\
0.0355102040816326 0.280189256869023\\
0.0410204081632653 0.123528841683459\\
0.046530612244898 0.0695725337187247\\
0.0520408163265306 0.0433756116094781\\
0.0575510204081633 0.0429961638397036\\
0.0630612244897959 0.0499855877451332\\
0.0685714285714286 0.0616524757228678\\
0.0740816326530612 0.0727504834527598\\
0.0795918367346939 0.0832828853717174\\
0.0851020408163265 0.0925616995006072\\
0.0906122448979592 0.100936594101105\\
0.0961224489795918 0.108045071178075\\
0.101632653061224 0.113602824116318\\
0.107142857142857 0.119193575984676\\
0.11265306122449 0.123807340101691\\
0.118163265306122 0.128170088478517\\
0.123673469387755 0.131697757678649\\
0.129183673469388 0.135107613140361\\
0.13469387755102 0.137899127103548\\
0.140204081632653 0.141343010356724\\
0.145714285714286 0.144394622103667\\
0.151224489795918 0.147874220725178\\
0.156734693877551 0.150686530402943\\
0.162244897959184 0.153411831351363\\
0.167755102040816 0.156214502000634\\
0.173265306122449 0.158590166287783\\
0.178775510204082 0.161162940572627\\
0.184285714285714 0.163751708525589\\
0.189795918367347 0.166238603561372\\
0.19530612244898 0.168788067400726\\
0.200816326530612 0.171241582134136\\
0.206326530612245 0.173791039520445\\
0.211836734693878 0.176465766748792\\
0.21734693877551 0.179157944488205\\
0.222857142857143 0.181797597704467\\
0.228367346938776 0.184500287885861\\
0.233877551020408 0.187119848285201\\
0.239387755102041 0.189690015854103\\
0.244897959183673 0.192253700950701\\
0.250408163265306 0.194876066648544\\
0.255918367346939 0.197370224799306\\
0.261428571428571 0.199814288299062\\
0.266938775510204 0.202275585582573\\
0.272448979591837 0.204606790116668\\
0.277959183673469 0.207013927679539\\
0.283469387755102 0.209427742950064\\
0.288979591836735 0.211738230670111\\
0.294489795918367 0.214017796635208\\
0.3 0.216262721536863\\
};
\addplot [
color=black,
dashed,
line width=1pt,
forget plot
]
table[row sep=crcr]{
0.03 0.30336159149539\\
0.0355102040816326 0.1743856033045\\
0.0410204081632653 0.0734812072499518\\
0.046530612244898 0.0295418244545061\\
0.0520408163265306 0.0299609469719966\\
0.0575510204081633 0.0362701683722099\\
0.0630612244897959 0.0449945516973634\\
0.0685714285714286 0.0568788934359698\\
0.0740816326530612 0.0672427596138061\\
0.0795918367346939 0.0780006547827487\\
0.0851020408163265 0.0871930478082662\\
0.0906122448979592 0.0952826366998765\\
0.0961224489795918 0.103430356165293\\
0.101632653061224 0.109207109286168\\
0.107142857142857 0.114870044242266\\
0.11265306122449 0.119224267838555\\
0.118163265306122 0.123575721217886\\
0.123673469387755 0.126603978420375\\
0.129183673469388 0.13047971260313\\
0.13469387755102 0.133629722195529\\
0.140204081632653 0.137582670456729\\
0.145714285714286 0.140292050088396\\
0.151224489795918 0.143976688315098\\
0.156734693877551 0.147355494851444\\
0.162244897959184 0.150476112098023\\
0.167755102040816 0.153294438041075\\
0.173265306122449 0.155677160079418\\
0.178775510204082 0.158370159201094\\
0.184285714285714 0.161015692242127\\
0.189795918367347 0.164029444244661\\
0.19530612244898 0.166710302925609\\
0.200816326530612 0.169212036694196\\
0.206326530612245 0.171614689873313\\
0.211836734693878 0.174594017914497\\
0.21734693877551 0.177305365885956\\
0.222857142857143 0.18008558386336\\
0.228367346938776 0.182838953948359\\
0.233877551020408 0.185615151280202\\
0.239387755102041 0.187967051542028\\
0.244897959183673 0.190470300036404\\
0.250408163265306 0.192831982322447\\
0.255918367346939 0.195468146916186\\
0.261428571428571 0.197840941941256\\
0.266938775510204 0.200558897677694\\
0.272448979591837 0.203092185418785\\
0.277959183673469 0.205451562830728\\
0.283469387755102 0.207916580124759\\
0.288979591836735 0.210062377692017\\
0.294489795918367 0.212506763209595\\
0.3 0.214564057705266\\
};
\addplot [
color=black,
dashed,
line width=1pt,
forget plot
]
table[row sep=crcr]{
0.03 0.334654373054913\\
0.0355102040816326 0.565798046983202\\
0.0410204081632653 0.161955553790184\\
0.046530612244898 0.109115041924337\\
0.0520408163265306 0.0718015664461268\\
0.0575510204081633 0.0535501542624191\\
0.0630612244897959 0.0563293991067813\\
0.0685714285714286 0.0675257188305625\\
0.0740816326530612 0.0776300286748199\\
0.0795918367346939 0.0876190657263606\\
0.0851020408163265 0.0970872412652066\\
0.0906122448979592 0.106244753733623\\
0.0961224489795918 0.113023346889751\\
0.101632653061224 0.119658316694236\\
0.107142857142857 0.124179087696919\\
0.11265306122449 0.129008548466733\\
0.118163265306122 0.133036008755334\\
0.123673469387755 0.136118860152637\\
0.129183673469388 0.139090979903152\\
0.13469387755102 0.141674175384723\\
0.140204081632653 0.145202581219415\\
0.145714285714286 0.14859747365777\\
0.151224489795918 0.151738710263705\\
0.156734693877551 0.154398355132893\\
0.162244897959184 0.157484395483037\\
0.167755102040816 0.159839014074498\\
0.173265306122449 0.162044450526463\\
0.178775510204082 0.163865591563089\\
0.184285714285714 0.166348679525874\\
0.189795918367347 0.168490435302713\\
0.19530612244898 0.170766923808466\\
0.200816326530612 0.173432954774871\\
0.206326530612245 0.175693633602264\\
0.211836734693878 0.178230648854472\\
0.21734693877551 0.180860667045869\\
0.222857142857143 0.183519674236354\\
0.228367346938776 0.186104160071747\\
0.233877551020408 0.188849313032296\\
0.239387755102041 0.191445570897937\\
0.244897959183673 0.19390743542089\\
0.250408163265306 0.196791460175634\\
0.255918367346939 0.199206387636941\\
0.261428571428571 0.201711939568217\\
0.266938775510204 0.203894648947127\\
0.272448979591837 0.206214725084665\\
0.277959183673469 0.208632198066048\\
0.283469387755102 0.211007990888045\\
0.288979591836735 0.213330106018687\\
0.294489795918367 0.215507803513985\\
0.3 0.217800647377955\\
};
\addplot [
color=orange,
mark size=3.0pt,
only marks,
mark=triangle*,
mark options={solid},
forget plot
]
table[row sep=crcr]{
0.0906122448979592 0.100936594101105\\
};
\addplot [
color=darkred,
mark size=2.0pt,
only marks,
mark=square*,
mark options={solid},
forget plot
]
table[row sep=crcr]{
0.0575510204081633 0.0429961638397036\\
};
\addplot [
color=blue,
mark size=2.0pt,
only marks,
mark=*,
mark options={solid},
forget plot
]
table[row sep=crcr]{
0.0520408163265306 0.0433756116094781\\
};
\addplot [
color=blue,
solid,
mark options={solid},
line width=1.5pt,
forget plot
]
table[row sep=crcr]{
0.030000000000000 0\\
0.035510204081633 0\\
0.041020408163265 0\\
0.0465 0.0880\\
0.0520 0.2800\\
0.0575 0.3600\\
0.063061224489796 0.400000000000000\\
0.068571428571429 0.400000000000000\\
0.074081632653061 0.400000000000000\\
0.079591836734694 0.400000000000000\\
0.085102040816327 0.400000000000000\\
0.090612244897959 0.400000000000000\\
0.096122448979592 0.400000000000000\\
0.101632653061224 0.400000000000000\\
0.107142857142857 0.400000000000000\\
0.112653061224490 0.400000000000000\\
0.118163265306122 0.400000000000000\\
0.123673469387755 0.400000000000000\\
0.129183673469388 0.400000000000000\\
0.134693877551020 0.400000000000000\\
0.140204081632653 0.400000000000000\\
0.145714285714286 0.400000000000000\\
0.151224489795918 0.400000000000000\\
0.156734693877551 0.400000000000000\\
0.162244897959184 0.400000000000000\\
0.167755102040816 0.400000000000000\\
0.173265306122449 0.400000000000000\\
0.178775510204082 0.400000000000000\\
0.184285714285714 0.400000000000000\\
0.189795918367347 0.400000000000000\\
0.195306122448980 0.400000000000000\\
0.200816326530612 0.400000000000000\\
0.206326530612245 0.400000000000000\\
0.211836734693878 0.400000000000000\\
0.217346938775510 0.400000000000000\\
0.222857142857143 0.400000000000000\\
0.228367346938776 0.400000000000000\\
0.233877551020408 0.400000000000000\\
0.239387755102041 0.400000000000000\\
0.244897959183673 0.400000000000000\\
0.250408163265306 0.400000000000000\\
0.255918367346939 0.400000000000000\\
0.261428571428571 0.400000000000000\\
0.266938775510204 0.400000000000000\\
0.272448979591837 0.400000000000000\\
0.277959183673469 0.400000000000000\\
0.283469387755102 0.400000000000000\\
0.288979591836735 0.400000000000000\\
0.294489795918367 0.400000000000000\\
0.300000000000000 0.400000000000000\\
};
\end{axis}
\end{tikzpicture}%

%% file: ScalingInEpsilon_Model1_2D_p4.tikz
\begin{tikzpicture}
\def\xl{4}
\def\xu{9}
\def\yl{-4}
\def\yu{-0.6}

\begin{axis}[%
width=\figurewidth,
height=\figureheight,
scale only axis,
xmin={\xl-(\xu-\xl)/5},
xmax={\xu+(\xu-\xl)/5},
ymin={\yl-(\yu-\yl)/5},
ymax={\yu+(\yu-\yl)*0.05},
hide axis,
axis background/.style={fill=white!100}
]
\draw [->] (axis cs: \xl,\yl) -- (axis cs: {\xu+(\xu-\xl)*0.05},\yl);
\draw [->] (axis cs: \xl,\yl) -- (axis cs: \xl,{\yu+(\yu-\yl)*0.05});
\node at (axis cs: {(\xu+(\xu-\xl)*0.05+\xl)/2},{\yl-0.16*(\yu-\yl)}) {$\log(n)$};
\node[rotate=90] at (axis cs: {\xl-0.16*(\xu-\xl)},{(\yu+(\yu-\yl)*0.05+\yl)/2}) {$\log(\eps)$};
\foreach \yValue in {-1,-2,-3,-4} {
    \edef\temp{\noexpand\draw [-] ({\xl+(\xu-\xl)/100},\yValue) -- (\xl,\yValue) node[left] {\yValue};} 
    \temp
}
\foreach \xValue in {4,5,6,7,8,9} {
	\edef\temp{\noexpand\draw [-] (\xValue,{\yl+(\yu-\yl)/100}) -- (\xValue,\yl) node[below] {\xValue};}    
    \temp
}
\addplot [
color=orange,
solid,
mark size=3.0pt,
mark=triangle*,
mark options={solid},
only marks,
forget plot
]
table[row sep=crcr]{
4.382026634673881 -0.805333661288171\\
5.075173815233827 -1.404416768124039\\
5.76832099579377 -1.99270231056207\\
6.46146817635372 -2.14262044330137\\
7.15461535691366 -2.40116592166649\\
7.84776253747361 -2.51305420341153\\
8.54090971803355 -2.7444\\
};
\addplot [
color=blue,
solid,
mark size=2.0pt,
mark=*,
mark options={solid},
only marks,
forget plot
]
table[row sep=crcr]{
4.382026634673881 -1.705769021372532\\
5.075173815233827 -1.964656191876369\\
5.76832099579377 -2.28639040707407\\
6.46146817635372 -2.59709313051623\\
7.15461535691366 -2.95572693894029\\
7.84776253747361 -3.24996641193823\\
8.54090971803355 -3.52717718452272\\
};
\addplot [
color=darkred,
solid,
mark size=2.0pt,
mark=square*,
mark options={solid},
only marks,
forget plot
]
table[row sep=crcr]{
4.382026634673881 -1.503057508875529\\
5.075173815233827 -1.692375964036095\\
5.76832099579377 -2.28639040707407\\
6.46146817635372 -2.59709313051623\\
7.15461535691366 -2.8550834131606\\
7.84776253747361 -3.24996641193823\\
8.54090971803355 -3.5272\\
};
\addplot [
color=blue,
dashed,
line width=2.5pt,
forget plot
]
table[row sep=crcr]{
4.382026634673881 -1.669492080070583\\
5.075173815233827 -1.982936763702514\\
5.76832099579377 -2.29638144733444\\
6.46146817635372 -2.60982613096638\\
7.15461535691366 -2.92327081459831\\
7.84776253747361 -3.23671549823024\\
8.54090971803355 -3.55016018186217\\
};
\addplot [
color=darkred,
dashed,
line width=1.5pt,
forget plot
]
table[row sep=crcr]{
4.382026634673881 -1.649363374914646\\
5.075173815233827 -1.962808058546577\\
5.76832099579377 -2.27625274217851\\
6.46146817635372 -2.58969742581044\\
7.15461535691366 -2.90314210944237\\
7.84776253747361 -3.2165867930743\\
8.54090971803355 -3.53003147670623\\
};
\addplot [
color=orange,
dashed,
line width=1.5pt,
forget plot
]
table[row sep=crcr]{
4.3820 -1.6092\\
5.0752 -1.7966\\
5.7683 -1.9840\\
6.4615 -2.1714\\
7.1546 -2.3588\\
7.8478 -2.5462\\
8.5409 -2.7336\\
};
\end{axis}
\end{tikzpicture}%

%% file: Realisation_Zoom_epsilon_0p0576_n1280_Model1_2D.tikz
\definecolor{mycolor1}{rgb}{0,0.447,0.741}

\begin{tikzpicture}

\begin{axis}[%
width=\figurewidth,
height=\figureheight,
unbounded coords=jump,
clip=false,
view={-20}{10},
scale only axis,
every outer x axis line/.append style={darkgray!60!black},
every x tick label/.append style={font=\color{darkgray!60!black}},
xmin=0.5,
xmax=1,
xmajorgrids,
every outer y axis line/.append style={darkgray!60!black},
every y tick label/.append style={font=\color{darkgray!60!black}},
ymin=0.25,
ymax=0.75,
ymajorgrids,
every outer z axis line/.append style={darkgray!60!black},
every z tick label/.append style={font=\color{darkgray!60!black}},
zmin=0.4,
zmax=1,
zmajorgrids,
hide axis,
grid style={solid},
]

\addplot3[%
mesh,
colormap={mymap}{[1pt] rgb(0pt)=(0.2081,0.1663,0.5292); rgb(1pt)=(0.211624,0.189781,0.577676); rgb(2pt)=(0.212252,0.213771,0.626971); rgb(3pt)=(0.2081,0.2386,0.677086); rgb(4pt)=(0.195905,0.264457,0.7279); rgb(5pt)=(0.170729,0.291938,0.779248); rgb(6pt)=(0.125271,0.324243,0.830271); rgb(7pt)=(0.0591333,0.359833,0.868333); rgb(8pt)=(0.0116952,0.38751,0.881957); rgb(9pt)=(0.00595714,0.408614,0.882843); rgb(10pt)=(0.0165143,0.4266,0.878633); rgb(11pt)=(0.0328524,0.443043,0.871957); rgb(12pt)=(0.0498143,0.458571,0.864057); rgb(13pt)=(0.0629333,0.47369,0.855438); rgb(14pt)=(0.0722667,0.488667,0.8467); rgb(15pt)=(0.0779429,0.503986,0.838371); rgb(16pt)=(0.0793476,0.520024,0.831181); rgb(17pt)=(0.0749429,0.537543,0.826271); rgb(18pt)=(0.0640571,0.556986,0.823957); rgb(19pt)=(0.0487714,0.577224,0.822829); rgb(20pt)=(0.0343429,0.596581,0.819852); rgb(21pt)=(0.0265,0.6137,0.8135); rgb(22pt)=(0.0238905,0.628662,0.803762); rgb(23pt)=(0.0230905,0.641786,0.791267); rgb(24pt)=(0.0227714,0.653486,0.776757); rgb(25pt)=(0.0266619,0.664195,0.760719); rgb(26pt)=(0.0383714,0.674271,0.743552); rgb(27pt)=(0.0589714,0.683757,0.725386); rgb(28pt)=(0.0843,0.692833,0.706167); rgb(29pt)=(0.113295,0.7015,0.685857); rgb(30pt)=(0.145271,0.709757,0.664629); rgb(31pt)=(0.180133,0.717657,0.642433); rgb(32pt)=(0.217829,0.725043,0.619262); rgb(33pt)=(0.258643,0.731714,0.595429); rgb(34pt)=(0.302171,0.737605,0.571186); rgb(35pt)=(0.348167,0.742433,0.547267); rgb(36pt)=(0.395257,0.7459,0.524443); rgb(37pt)=(0.44201,0.748081,0.503314); rgb(38pt)=(0.487124,0.749062,0.483976); rgb(39pt)=(0.530029,0.749114,0.466114); rgb(40pt)=(0.570857,0.748519,0.44939); rgb(41pt)=(0.609852,0.747314,0.433686); rgb(42pt)=(0.6473,0.7456,0.4188); rgb(43pt)=(0.683419,0.743476,0.404433); rgb(44pt)=(0.71841,0.741133,0.390476); rgb(45pt)=(0.752486,0.7384,0.376814); rgb(46pt)=(0.785843,0.735567,0.363271); rgb(47pt)=(0.818505,0.732733,0.34979); rgb(48pt)=(0.850657,0.7299,0.336029); rgb(49pt)=(0.882433,0.727433,0.3217); rgb(50pt)=(0.913933,0.725786,0.306276); rgb(51pt)=(0.944957,0.726114,0.288643); rgb(52pt)=(0.973895,0.731395,0.266648); rgb(53pt)=(0.993771,0.745457,0.240348); rgb(54pt)=(0.999043,0.765314,0.216414); rgb(55pt)=(0.995533,0.786057,0.196652); rgb(56pt)=(0.988,0.8066,0.179367); rgb(57pt)=(0.978857,0.827143,0.163314); rgb(58pt)=(0.9697,0.848138,0.147452); rgb(59pt)=(0.962586,0.870514,0.1309); rgb(60pt)=(0.958871,0.8949,0.113243); rgb(61pt)=(0.959824,0.921833,0.0948381); rgb(62pt)=(0.9661,0.951443,0.0755333); rgb(63pt)=(0.9763,0.9831,0.0538)},
shader=flat,
mesh/rows=21]
table[row sep=crcr,header=false] {
0.5 0.25 NaN\\
0.5 0.275 NaN\\
0.5 0.3 NaN\\
0.5 0.325 NaN\\
0.5 0.35 NaN\\
0.5 0.375 NaN\\
0.5 0.4 NaN\\
0.5 0.425 NaN\\
0.5 0.45 NaN\\
0.5 0.475 NaN\\
0.5 0.5 NaN\\
0.5 0.525 NaN\\
0.5 0.55 NaN\\
0.5 0.575 NaN\\
0.5 0.6 NaN\\
0.5 0.625 NaN\\
0.5 0.65 NaN\\
0.5 0.675 NaN\\
0.5 0.7 NaN\\
0.5 0.725 NaN\\
0.5 0.75 NaN\\
0.525 0.25 NaN\\
0.525 0.275 NaN\\
0.525 0.3 0.50787932000803\\
0.525 0.325 0.512699133877458\\
0.525 0.35 0.508600907873804\\
0.525 0.375 0.510266531516891\\
0.525 0.4 0.505373976182989\\
0.525 0.425 0.504865170870924\\
0.525 0.45 0.509221690825499\\
0.525 0.475 0.504663956033215\\
0.525 0.5 0.508387491538436\\
0.525 0.525 0.520425432294147\\
0.525 0.55 0.538847995957221\\
0.525 0.575 0.554571648453545\\
0.525 0.6 0.552717580611094\\
0.525 0.625 0.549792803180708\\
0.525 0.65 0.539123700399622\\
0.525 0.675 0.539345093544173\\
0.525 0.7 0.528604585978598\\
0.525 0.725 0.51671532158265\\
0.525 0.75 NaN\\
0.55 0.25 NaN\\
0.55 0.275 0.529143810176602\\
0.55 0.3 0.536227033960117\\
0.55 0.325 0.551955107436949\\
0.55 0.35 0.52618696879765\\
0.55 0.375 0.533021918577683\\
0.55 0.4 0.533945126362574\\
0.55 0.425 0.538019448166253\\
0.55 0.45 0.534141863128681\\
0.55 0.475 0.533557766541141\\
0.55 0.5 0.537118828404025\\
0.55 0.525 0.547279073475805\\
0.55 0.55 0.55895556237264\\
0.55 0.575 0.574735996593107\\
0.55 0.6 0.582442231253984\\
0.55 0.625 0.579092911107269\\
0.55 0.65 0.574860067632671\\
0.55 0.675 0.574654562238363\\
0.55 0.7 0.564641357057541\\
0.55 0.725 0.54768041699839\\
0.55 0.75 NaN\\
0.575 0.25 NaN\\
0.575 0.275 0.562239738178179\\
0.575 0.3 0.570567050122982\\
0.575 0.325 0.58572192661572\\
0.575 0.35 0.58524381969426\\
0.575 0.375 0.582703690349799\\
0.575 0.4 0.570085317258173\\
0.575 0.425 0.565704372761515\\
0.575 0.45 0.559119708492739\\
0.575 0.475 0.55646398314167\\
0.575 0.5 0.551201624969455\\
0.575 0.525 0.567839277806028\\
0.575 0.55 0.585613144502053\\
0.575 0.575 0.59775838359394\\
0.575 0.6 0.606534544145391\\
0.575 0.625 0.603140260106618\\
0.575 0.65 0.604188737499302\\
0.575 0.675 0.60401126203046\\
0.575 0.7 0.595684362499914\\
0.575 0.725 0.577262654371058\\
0.575 0.75 NaN\\
0.6 0.25 NaN\\
0.6 0.275 0.596041872034898\\
0.6 0.3 0.604659099574175\\
0.6 0.325 0.612828848124433\\
0.6 0.35 0.615601832654681\\
0.6 0.375 0.613828920632211\\
0.6 0.4 0.605451695283381\\
0.6 0.425 0.597828203093585\\
0.6 0.45 0.589096592628369\\
0.6 0.475 0.586759233240659\\
0.6 0.5 0.593672401141747\\
0.6 0.525 0.602069006026432\\
0.6 0.55 0.614916955925502\\
0.6 0.575 0.627339852368099\\
0.6 0.6 0.630760262856639\\
0.6 0.625 0.630955395555313\\
0.6 0.65 0.624707104192331\\
0.6 0.675 0.619206185521495\\
0.6 0.7 0.609123850795229\\
0.6 0.725 0.597938503690977\\
0.6 0.75 NaN\\
0.625 0.25 NaN\\
0.625 0.275 0.62870859452134\\
0.625 0.3 0.637813342983841\\
0.625 0.325 0.642090581501042\\
0.625 0.35 0.643662408762164\\
0.625 0.375 0.643616800233759\\
0.625 0.4 0.647559136269626\\
0.625 0.425 0.647714954329294\\
0.625 0.45 0.638979464561966\\
0.625 0.475 0.634013649873741\\
0.625 0.5 0.641119674471849\\
0.625 0.525 0.648989522832957\\
0.625 0.55 0.655008941861897\\
0.625 0.575 0.661055693684829\\
0.625 0.6 0.654559118243629\\
0.625 0.625 0.651645812244302\\
0.625 0.65 0.650392778587202\\
0.625 0.675 0.641615054882288\\
0.625 0.7 0.627343887919243\\
0.625 0.725 0.611051831474\\
0.625 0.75 NaN\\
0.65 0.25 NaN\\
0.65 0.275 0.655691373481017\\
0.65 0.3 0.657640526126143\\
0.65 0.325 0.660836236185426\\
0.65 0.35 0.671161989077739\\
0.65 0.375 0.671054875097753\\
0.65 0.4 0.668533067668645\\
0.65 0.425 0.663382817523068\\
0.65 0.45 0.676701419143849\\
0.65 0.475 0.678462551707523\\
0.65 0.5 0.691437699815573\\
0.65 0.525 0.698293604577547\\
0.65 0.55 0.704862238579904\\
0.65 0.575 0.706801990789017\\
0.65 0.6 0.689909419379407\\
0.65 0.625 0.674003224755755\\
0.65 0.65 0.666429964010626\\
0.65 0.675 0.658550285965177\\
0.65 0.7 0.644461894649823\\
0.65 0.725 0.629312708228398\\
0.65 0.75 NaN\\
0.675 0.25 NaN\\
0.675 0.275 0.677304059035488\\
0.675 0.3 0.693174619231734\\
0.675 0.325 0.697002534637273\\
0.675 0.35 0.698279009266484\\
0.675 0.375 0.698006422991697\\
0.675 0.4 0.69929436275546\\
0.675 0.425 0.709555097965059\\
0.675 0.45 0.71107414037102\\
0.675 0.475 0.724266413171609\\
0.675 0.5 0.737312882394766\\
0.675 0.525 0.754614284204147\\
0.675 0.55 0.737697126328132\\
0.675 0.575 0.727218014752304\\
0.675 0.6 0.71322660326696\\
0.675 0.625 0.696416096105204\\
0.675 0.65 0.684185712258556\\
0.675 0.675 0.674561666926729\\
0.675 0.7 0.661637985513047\\
0.675 0.725 0.649631193329438\\
0.675 0.75 NaN\\
0.7 0.25 NaN\\
0.7 0.275 0.700529551954927\\
0.7 0.3 0.713860925408178\\
0.7 0.325 0.722694478985233\\
0.7 0.35 0.725868977226422\\
0.7 0.375 0.733624148111028\\
0.7 0.4 0.734007514127475\\
0.7 0.425 0.734011904194579\\
0.7 0.45 0.742519931442569\\
0.7 0.475 0.756221141204547\\
0.7 0.5 0.764214213586375\\
0.7 0.525 0.765267124564913\\
0.7 0.55 0.765543330692413\\
0.7 0.575 0.76214559581747\\
0.7 0.6 0.737501865263127\\
0.7 0.625 0.717661597068636\\
0.7 0.65 0.70172456856617\\
0.7 0.675 0.693110281003618\\
0.7 0.7 0.682620647876002\\
0.7 0.725 0.677728325861653\\
0.7 0.75 NaN\\
0.725 0.25 NaN\\
0.725 0.275 0.718260385246718\\
0.725 0.3 0.731248126341213\\
0.725 0.325 0.746778589214132\\
0.725 0.35 0.756274238428045\\
0.725 0.375 0.764055472475313\\
0.725 0.4 0.768838764243633\\
0.725 0.425 0.784513361793416\\
0.725 0.45 0.801710941350646\\
0.725 0.475 0.817682342826811\\
0.725 0.5 0.814216115825858\\
0.725 0.525 0.807952713435525\\
0.725 0.55 0.792581835025093\\
0.725 0.575 0.779249861899971\\
0.725 0.6 0.762237894170216\\
0.725 0.625 0.745470270683129\\
0.725 0.65 0.729086761322121\\
0.725 0.675 0.717699748012903\\
0.725 0.7 0.717822250601563\\
0.725 0.725 0.705361049363885\\
0.725 0.75 NaN\\
0.75 0.25 NaN\\
0.75 0.275 0.740460169494599\\
0.75 0.3 0.755993367935337\\
0.75 0.325 0.770969863020792\\
0.75 0.35 0.785647752799697\\
0.75 0.375 0.796476157144227\\
0.75 0.4 0.80497266123267\\
0.75 0.425 0.826485495709643\\
0.75 0.45 0.854972957822253\\
0.75 0.475 0.877161129602529\\
0.75 0.5 0.869192132372484\\
0.75 0.525 0.849601231248709\\
0.75 0.55 0.823069999367592\\
0.75 0.575 0.808983214853306\\
0.75 0.6 0.802184904500224\\
0.75 0.625 0.780253075355826\\
0.75 0.65 0.760419702945988\\
0.75 0.675 0.751106568353834\\
0.75 0.7 0.742828329015686\\
0.75 0.725 0.732993772866117\\
0.75 0.75 NaN\\
0.775 0.25 NaN\\
0.775 0.275 0.761641270783361\\
0.775 0.3 0.777226369650201\\
0.775 0.325 0.793190332478785\\
0.775 0.35 0.809677663317415\\
0.775 0.375 0.824268678773399\\
0.775 0.4 0.839945517155546\\
0.775 0.425 0.855883002514013\\
0.775 0.45 0.895862937070538\\
0.775 0.475 0.910272343593402\\
0.775 0.5 0.930654303622566\\
0.775 0.525 0.893088216088031\\
0.775 0.55 0.864072151290769\\
0.775 0.575 0.829844410360935\\
0.775 0.6 0.8237850151057\\
0.775 0.625 0.8060327892942\\
0.775 0.65 0.788687232918789\\
0.775 0.675 0.7739707652827\\
0.775 0.7 0.769117733043088\\
0.775 0.725 0.760626496368349\\
0.775 0.75 NaN\\
0.8 0.25 NaN\\
0.8 0.275 0.783062201086397\\
0.8 0.3 0.798044210881542\\
0.8 0.325 0.813410113304532\\
0.8 0.35 0.829797425817814\\
0.8 0.375 0.846161550505986\\
0.8 0.4 0.863543212268517\\
0.8 0.425 0.881229438522053\\
0.8 0.45 0.908924514823735\\
0.8 0.475 0.930743570999421\\
0.8 0.5 1\\
0.8 0.525 0.937843711886851\\
0.8 0.55 0.894147169562023\\
0.8 0.575 0.855395272397501\\
0.8 0.6 0.837214815922411\\
0.8 0.625 0.824082461832253\\
0.8 0.65 0.809816014684766\\
0.8 0.675 0.804610608103607\\
0.8 0.7 0.794162477041985\\
0.8 0.725 0.785126377810427\\
0.8 0.75 NaN\\
0.825 0.25 NaN\\
0.825 0.275 0.801491296101352\\
0.825 0.3 0.814673592233844\\
0.825 0.325 0.828257590546656\\
0.825 0.35 0.841181412046123\\
0.825 0.375 0.857106756032213\\
0.825 0.4 0.87714713614907\\
0.825 0.425 0.89423349319893\\
0.825 0.45 0.92204938047739\\
0.825 0.475 0.934665793809425\\
0.825 0.5 0.952424627469662\\
0.825 0.525 0.933278041416613\\
0.825 0.55 0.912684963083258\\
0.825 0.575 0.876526585340926\\
0.825 0.6 0.848136559415416\\
0.825 0.625 0.840464381713111\\
0.825 0.65 0.826401248326705\\
0.825 0.675 0.809918269391388\\
0.825 0.7 0.803122012306134\\
0.825 0.725 0.794985037876784\\
0.825 0.75 NaN\\
0.85 0.25 NaN\\
0.85 0.275 0.813830655309498\\
0.85 0.3 0.825324014387126\\
0.85 0.325 0.839601715372536\\
0.85 0.35 0.853388831512122\\
0.85 0.375 0.868013042718081\\
0.85 0.4 0.887697934612801\\
0.85 0.425 0.907515249219088\\
0.85 0.45 0.923136166317289\\
0.85 0.475 0.939003121681804\\
0.85 0.5 0.953389089758344\\
0.85 0.525 0.95130442334394\\
0.85 0.55 0.920382923405293\\
0.85 0.575 0.885455748554966\\
0.85 0.6 0.862605845488503\\
0.85 0.625 0.854352939672959\\
0.85 0.65 0.840749004969374\\
0.85 0.675 0.824108998519612\\
0.85 0.7 0.81089425971643\\
0.85 0.725 0.812201245363139\\
0.85 0.75 NaN\\
0.875 0.25 NaN\\
0.875 0.275 0.825247257473817\\
0.875 0.3 0.835808010665123\\
0.875 0.325 0.847909141808996\\
0.875 0.35 0.861566968689175\\
0.875 0.375 0.875949140198174\\
0.875 0.4 0.89577515929786\\
0.875 0.425 0.917595971412105\\
0.875 0.45 0.930204945223133\\
0.875 0.475 0.938685188812667\\
0.875 0.5 0.941082724118548\\
0.875 0.525 0.934781546987509\\
0.875 0.55 0.92635818773716\\
0.875 0.575 0.898007284943169\\
0.875 0.6 0.880841545326619\\
0.875 0.625 0.869955645912662\\
0.875 0.65 0.859844591904866\\
0.875 0.675 0.844423127708219\\
0.875 0.7 0.82883103315286\\
0.875 0.725 0.827751511459691\\
0.875 0.75 NaN\\
0.9 0.25 NaN\\
0.9 0.275 0.834899948321631\\
0.9 0.3 0.847759353716392\\
0.9 0.325 0.859331340641639\\
0.9 0.35 0.872313520615868\\
0.9 0.375 0.890371965750509\\
0.9 0.4 0.907851296859638\\
0.9 0.425 0.923889576809231\\
0.9 0.45 0.935878407899139\\
0.9 0.475 0.943166520753526\\
0.9 0.5 0.946375754488028\\
0.9 0.525 0.940315946130866\\
0.9 0.55 0.928240068198565\\
0.9 0.575 0.899604094006066\\
0.9 0.6 0.891953459255747\\
0.9 0.625 0.88487073451792\\
0.9 0.65 0.878095423072386\\
0.9 0.675 0.865816985395721\\
0.9 0.7 0.85008810907286\\
0.9 0.725 0.845358064770313\\
0.9 0.75 NaN\\
0.925 0.25 NaN\\
0.925 0.275 0.844454789201979\\
0.925 0.3 0.857223419502823\\
0.925 0.325 0.871721311702255\\
0.925 0.35 0.887936936256241\\
0.925 0.375 0.906870834225662\\
0.925 0.4 0.927524323487602\\
0.925 0.425 0.942662973898773\\
0.925 0.45 0.949536994644401\\
0.925 0.475 0.950551795848256\\
0.925 0.5 0.953119842078362\\
0.925 0.525 0.949674905465123\\
0.925 0.55 0.938410808265419\\
0.925 0.575 0.919643115480004\\
0.925 0.6 0.905681911177283\\
0.925 0.625 0.897931706547057\\
0.925 0.65 0.891472464523425\\
0.925 0.675 0.883086467348981\\
0.925 0.7 0.869438467492231\\
0.925 0.725 0.859554424776014\\
0.925 0.75 NaN\\
0.95 0.25 NaN\\
0.95 0.275 0.854009630082328\\
0.95 0.3 0.866947809165736\\
0.95 0.325 0.884168941066147\\
0.95 0.35 0.903109997572969\\
0.95 0.375 0.920361215672303\\
0.95 0.4 0.939398030278547\\
0.95 0.425 0.948566220492933\\
0.95 0.45 0.957341088436866\\
0.95 0.475 0.96486350612323\\
0.95 0.5 0.969562003005046\\
0.95 0.525 0.963161236799778\\
0.95 0.55 0.94841392254392\\
0.95 0.575 0.929800574019926\\
0.95 0.6 0.918684075137991\\
0.95 0.625 0.909947089312732\\
0.95 0.65 0.900757720072084\\
0.95 0.675 0.890227501187104\\
0.95 0.7 0.87647166049944\\
0.95 0.725 0.868986478016274\\
0.95 0.75 NaN\\
0.975 0.25 NaN\\
0.975 0.275 0.863564470962676\\
0.975 0.3 0.881750770502813\\
0.975 0.325 0.895452083719061\\
0.975 0.35 0.908443657489577\\
0.975 0.375 0.928983215608376\\
0.975 0.4 0.939460851315169\\
0.975 0.425 0.954007506859011\\
0.975 0.45 0.961293177283117\\
0.975 0.475 0.968771189272553\\
0.975 0.5 0.974983309861097\\
0.975 0.525 0.970686378950331\\
0.975 0.55 0.956719308034778\\
0.975 0.575 0.945542236769664\\
0.975 0.6 0.93301595586194\\
0.975 0.625 0.920301467085838\\
0.975 0.65 0.907396837516344\\
0.975 0.675 0.89471977244513\\
0.975 0.7 0.881622311875911\\
0.975 0.725 0.874342338067012\\
0.975 0.75 NaN\\
1 0.25 NaN\\
1 0.275 NaN\\
1 0.3 NaN\\
1 0.325 NaN\\
1 0.35 NaN\\
1 0.375 NaN\\
1 0.4 NaN\\
1 0.425 NaN\\
1 0.45 NaN\\
1 0.475 NaN\\
1 0.5 NaN\\
1 0.525 NaN\\
1 0.55 NaN\\
1 0.575 NaN\\
1 0.6 NaN\\
1 0.625 NaN\\
1 0.65 NaN\\
1 0.675 NaN\\
1 0.7 NaN\\
1 0.725 NaN\\
1 0.75 NaN\\
};
\addplot3 [
color=black, 
mark size=1pt,
only marks,
mark=*,
mark options={solid}]
table[row sep=crcr] {
0.8 0.5 1\\
0.896750966300958 0.286776127457757 0.839272723480317\\
0.902771977462955 0.443307096739818 0.93383368959585\\
0.610354789491249 0.590563611606405 0.64451033305754\\
0.717545631591306 0.636844667860001 0.727377537743687\\
0.852848313794446 0.435990016586211 0.913467103935741\\
0.660761305835014 0.452417241303641 0.705432859955636\\
0.820130475649893 0.320075943110966 0.823086003763105\\
0.610456067310923 0.319729894706417 0.618212623349431\\
0.996628822967953 0.281100166085246 0.875529882150331\\
0.829330911141074 0.36177041731311 0.849295082021718\\
0.738927787432516 0.702469856909833 0.732041922203926\\
0.626578235140411 0.702924333907962 0.626732349095003\\
0.981412845867243 0.705104676070429 0.879907806085337\\
0.679608034336165 0.657771390884537 0.6851181012674\\
0.7767791920296 0.333471872572624 0.800633283765651\\
0.905025376762873 0.553778152360609 0.927745059633021\\
0.94166257336887 0.316020689047971 0.872139819388935\\
0.818411855866936 0.603014587398564 0.847846998458529\\
0.666906186164037 0.438275770234889 0.704725459078505\\
0.813503588228468 0.595511914424481 0.847935785019619\\
0.625717367146037 0.668026662527658 0.646076574642231\\
0.906993019500661 0.662923413406259 0.8790223460283\\
0.572131847654909 0.48344100136562 0.55136374489833\\
0.903338950738432 0.571324421969113 0.904457648516763\\
0.986631397979414 0.628632455243067 0.921676710739234\\
0.890279572452008 0.746189348982611 0.836174584541974\\
0.804770328186065 0.384360635226165 0.855386401789428\\
0.556768447472082 0.688768012719625 0.580396205813352\\
0.597875400608463 0.402327162870215 0.601979570676646\\
0.734634322479852 0.517831329437117 0.82476977622095\\
0.588734931310846 0.731851647500458 0.586847519563577\\
0.805532816497391 0.641278569924126 0.816308025158311\\
0.761772587446987 0.495623970862494 0.900584911733728\\
0.73431033983156 0.374160408174907 0.774957441827437\\
0.957096547913993 0.505259946445197 0.975147675867003\\
0.780291105767305 0.454155929647519 0.910701245586381\\
0.892702383484364 0.563241363149134 0.915493394912678\\
0.612013984234454 0.339793008902145 0.629578674890848\\
0.50868945810227 0.624583002688413 0.5310125143195\\
0.598022997067236 0.648736113587461 0.623604375923224\\
0.712619762620216 0.373497888737766 0.749437338024905\\
0.64966640763613 0.423645275337246 0.661810362934186\\
0.957072824374163 0.684798012434607 0.887553041875415\\
0.670982508455675 0.451325421246072 0.711941857884977\\
0.87348890039871 0.636667687206959 0.865148725502853\\
0.91187476865064 0.72621176631384 0.854107011397613\\
0.689641681209292 0.635260414894502 0.695808234059464\\
0.557804714387723 0.625122397308514 0.58817110751464\\
0.800424794384846 0.697439096217222 0.795561685326973\\
0.504686378036765 0.419115029941283 0.480119592314448\\
0.529101751211831 0.641041544298777 0.541026369967113\\
0.953672341494932 0.329437527145471 0.889726725524599\\
0.556778717749001 0.44797090722905 0.540482258928318\\
0.600894508295261 0.465190351306326 0.58558850793151\\
0.833294096726429 0.603108885353926 0.848285439324318\\
0.989547621752577 0.736941860120139 0.873407493790533\\
0.697173860239204 0.637819207435455 0.703842887289515\\
0.508927167607604 0.532144542167342 0.506526866514799\\
0.681553544189616 0.572471948843176 0.739316337933137\\
0.557628145279771 0.527774442485995 0.546208519009345\\
0.950399058359332 0.735603954354949 0.869792490979582\\
0.510836847606421 0.345729971491733 0.50140404848119\\
0.841152553007951 0.478634600177314 0.938003877411206\\
0.560435455011861 0.666191663819268 0.59147435205699\\
0.944134517275711 0.734594407286704 0.863004201260534\\
0.900442140858193 0.332479457434833 0.860530395345068\\
0.632193997405328 0.369266158090538 0.653176964220209\\
0.516461734167434 0.42273759940176 0.492955209516943\\
0.69121462913145 0.377368393397148 0.723100892650939\\
0.649212958682711 0.599332068352685 0.689661128458198\\
0.57079593357119 0.63423768817409 0.601092134621222\\
0.723136464599151 0.576815374708839 0.775250900142052\\
0.928378223516919 0.349078429673212 0.889413792594436\\
0.505653869344781 0.748549542309387 0.494309533776813\\
0.545763188443633 0.48466628529025 0.526998830225958\\
0.825722294028466 0.515472651567212 0.939232871717726\\
0.71129147049187 0.33969018138519 0.735311938445748\\
0.961465217690023 0.315560651972318 0.888997412765631\\
0.952407770222575 0.643638364186628 0.90393178169303\\
0.747516926326245 0.313847884230513 0.762785171508784\\
0.796674491635933 0.388888360603717 0.854783725369774\\
0.672641662882902 0.634830830121167 0.687601837828491\\
0.603313292450109 0.65347035723999 0.62621514835668\\
0.796814421268911 0.263945218665011 0.773530738577388\\
0.735269940189609 0.383113193894543 0.774951066183014\\
0.656979767272955 0.620716536738546 0.681317362562259\\
0.748165576186312 0.259010531534803 0.730638550990171\\
0.939433799432894 0.665870143073621 0.891903091215741\\
0.728568244242342 0.551385087432613 0.796293082565311\\
0.650647237047691 0.275570801569803 0.658307114982358\\
0.806878408035618 0.62185710361723 0.831297546234613\\
0.544194662207204 0.302608389824434 0.529106287222048\\
0.622955861825668 0.577622950816906 0.657628532267427\\
0.818883711816406 0.734932788204561 0.789314350622776\\
0.846119357594698 0.6915399704329 0.810255496843789\\
0.553766009195282 0.575454247763055 0.577204088730682\\
0.555666158880712 0.568677618897226 0.575184959579343\\
0.604435974345879 0.640463967654522 0.630292576187905\\
0.526349642470887 0.56563144105305 0.556911608542055\\
0.721059175232061 0.297450940320884 0.729449796694457\\
0.84260103603702 0.345530155179932 0.848580683316567\\
0.759685718167462 0.557020628247086 0.829592102546951\\
0.923017767921153 0.69036343882904 0.877705824759301\\
0.557512611368792 0.261906665471005 0.534082166410148\\
0.544360948798569 0.468744207168826 0.528199703724798\\
0.629938261821685 0.606467334820059 0.657249845276013\\
0.917675500102484 0.733426039517901 0.854339443170051\\
0.68074167652991 0.688535317260708 0.672264082452536\\
0.617326694267432 0.382581145538756 0.633228002130265\\
0.681803363676878 0.518240854847237 0.768255828894616\\
0.885290171061452 0.558902348891923 0.916566347267922\\
0.999879067839387 0.448812646700452 0.964800265614896\\
0.556653679428374 0.591775511200523 0.591964849459059\\
0.511923460323181 0.731515256518571 0.495732776732701\\
0.595068302854429 0.720931728228616 0.598032971286192\\
0.653528993732435 0.524999851559058 0.705253320654822\\
0.791550251926774 0.361782993747501 0.834739392994923\\
0.713931367162337 0.405706449224023 0.752902980948934\\
0.805166281940449 0.43454148810554 0.890249507790754\\
0.517319137001344 0.503384924916418 0.50602924802476\\
0.79562218975414 0.68645392683332 0.785872728517091\\
0.591719690181946 0.695977721199807 0.60733396236567\\
0.659452620544525 0.731721188595494 0.635189836627249\\
0.993170985365447 0.741377113835917 0.872150056478664\\
0.791545234841724 0.475131522881481 0.929446203270343\\
0.579610587190408 0.591918525702021 0.611114392699507\\
0.919737896229818 0.434663358048463 0.944373721566923\\
0.555719911475165 0.551092672120086 0.564742908320994\\
0.757941722373145 0.373864137037787 0.807600299719858\\
0.654419236693616 0.275193740552363 0.658307114981473\\
0.889949054888152 0.294867099112871 0.8392480731738\\
0.668218844330459 0.326768892520437 0.690240601257308\\
0.918302268321324 0.647871097589267 0.889370117753943\\
0.652359263807669 0.637841326006166 0.670551404653193\\
0.553472481516696 0.476149594730903 0.537333779773439\\
0.916110193482973 0.45870182637889 0.944550277067333\\
0.78972060934824 0.60099276064908 0.833094980079698\\
0.966813818539717 0.434865117015749 0.95461216048948\\
0.853959402955557 0.485460244964119 0.951516517361965\\
0.540527794749159 0.431787412818572 0.526106152985447\\
0.505152550821154 0.515459444839436 0.500021453553122\\
0.578734333094401 0.690784956846228 0.607301337369057\\
0.897592273023424 0.57268934326616 0.897959344150643\\
0.80170852361882 0.74007759751037 0.782594617646723\\
0.956214471123851 0.637671814284754 0.909033808447692\\
0.514678816739023 0.505035171806744 0.504238702634991\\
0.665406548725784 0.429975246137746 0.699036558754099\\
0.972869269087839 0.53813192279225 0.961282303018511\\
0.71483991465287 0.570295440291489 0.77717729912404\\
0.523531322672027 0.503042847571715 0.507681925788384\\
0.848968807250136 0.526796175461105 0.952054183132048\\
0.763914494242464 0.416629221251442 0.833592660957771\\
0.701662840439768 0.538870654105546 0.768722605052172\\
0.979893450066785 0.454617438907175 0.963697089064412\\
0.803772445643009 0.682811751090593 0.806172361710719\\
0.615098939296913 0.268404367035586 0.61418348392245\\
0.798968557213536 0.38128354172878 0.849555776076709\\
0.890580710042851 0.738437300921625 0.838173698545216\\
0.967825084618818 0.536888363538073 0.961282302905535\\
0.69581419182488 0.496871966433736 0.756854657392401\\
0.544600752813757 0.40197995989309 0.525140578496941\\
0.936213193985016 0.444417592226183 0.953370760102227\\
0.598289326940254 0.635668147698973 0.626209779892736\\
0.851877230121635 0.653673406242905 0.839631596248418\\
0.752586851411443 0.472454038691497 0.883653903908508\\
0.837450216826856 0.522313656458288 0.958760870786082\\
0.889743706793846 0.28439031811912 0.835500672722579\\
0.981723549742713 0.345842180573096 0.907013524942862\\
0.862315767578514 0.26179302017712 0.814200548150073\\
0.529272376947821 0.405989155670227 0.508666055747985\\
0.967912669818733 0.592303195473019 0.934339025369056\\
0.830431113487101 0.539488684926696 0.930790341816855\\
0.621571966395797 0.414129110039419 0.644970269669494\\
0.57442603944283 0.735784725196401 0.568573876287934\\
0.687742464161246 0.455322621044195 0.713680580394482\\
0.590071029823612 0.618542337447304 0.618003920902795\\
0.617331885338194 0.405718710051996 0.640134004808607\\
0.72690845906154 0.70082605744322 0.720521631614683\\
0.634045924174415 0.63700555975707 0.656978319277619\\
0.752559219350107 0.257357444298136 0.730638544041826\\
0.75771271505066 0.577808150277288 0.815792146995026\\
0.838475119768797 0.465907037408674 0.934768643530128\\
0.844003912299951 0.510567992926764 0.956924708026947\\
0.976611594793114 0.685936179491611 0.889651517131377\\
0.981935894624428 0.493982272449778 0.97474611392749\\
0.923052998305207 0.457167450939811 0.950032560263399\\
0.652965308525666 0.607145655266464 0.687582169273398\\
0.860857500756139 0.384450601326728 0.879133046777193\\
0.622211396677159 0.355232097780553 0.642787639490026\\
0.791457906927876 0.436562126700811 0.887456421543083\\
0.674450201334774 0.635700701390973 0.690384639810512\\
0.864458138705059 0.494853097534895 0.943621016420854\\
0.665835532037635 0.457394875143017 0.711940222764198\\
0.621725699247968 0.711399117134646 0.612815659165949\\
0.770754244073012 0.650584031585053 0.784977322028015\\
0.730986697506542 0.295438385894555 0.729449762364987\\
0.800798167895389 0.356574648567488 0.834737117695244\\
0.592460076523037 0.456185420555379 0.583020412365842\\
0.842229522803865 0.277063905713421 0.810269823408535\\
0.526548414257021 0.267818589565266 0.503346775271177\\
0.835137254585206 0.628426907790893 0.844734419211725\\
0.908628810594412 0.741702177161977 0.853065354160847\\
0.8532643740468 0.42243260461999 0.908042492171366\\
0.522834301020647 0.670081960160418 0.538256237298426\\
0.946148794562511 0.749434273566672 0.860648935631012\\
0.604877099274976 0.320179925465855 0.618209154453782\\
0.821159119762907 0.450084003570264 0.92185847125331\\
0.655707947618119 0.741305278209221 0.619262274055211\\
0.818253846190852 0.578555867900184 0.869436813171756\\
0.952615054258229 0.427969357910103 0.949371907249102\\
0.564923381251329 0.489770209480822 0.547092014392716\\
0.676839393873037 0.303291533227662 0.696972991455807\\
0.848021602973121 0.2690429776502 0.810148774572527\\
0.734391454067498 0.417013474230646 0.798712643702618\\
0.976260164764427 0.517455232249628 0.975748894283769\\
0.646539376807889 0.574127157107886 0.704429694667209\\
0.708068392333254 0.689023454361018 0.693820036414735\\
0.565127848267236 0.457550272907665 0.545322054871627\\
0.923272787027268 0.693320898060286 0.871425322305068\\
0.64958740767562 0.368280328610595 0.670303358904558\\
0.755269923936985 0.592352884283737 0.81579213015038\\
0.829613366395353 0.709263676905948 0.800473039855859\\
0.687090495956269 0.560597930130718 0.751469007941859\\
0.708352078901677 0.459898933306801 0.770090766111296\\
0.60356654986875 0.524604656439891 0.606624070197419\\
0.633827593627922 0.728606225258041 0.614308284217199\\
0.806227835509854 0.722240128754432 0.791945779282046\\
0.509598833595111 0.472527691077761 0.493310987963627\\
0.795448592145339 0.40713967529347 0.865971683138114\\
0.857848427994472 0.698109151242465 0.814685794426743\\
0.531920624360976 0.741644972972042 0.51530810347383\\
0.8557431661952 0.30361060242885 0.829593835363759\\
0.96809664114014 0.518059368780743 0.975845598082311\\
0.563520581945619 0.625875894105425 0.590678142785055\\
0.653411476759634 0.314888736357681 0.65784335109975\\
0.680924942343551 0.525120814126042 0.768271203920019\\
0.991879636942119 0.462982362699806 0.96964826164967\\
0.595393067151581 0.427187978982608 0.587807160341862\\
0.963869982901046 0.519546686704649 0.975843749838208\\
0.966957924360731 0.745079342721426 0.868520096329731\\
0.553432214382848 0.320014753868527 0.562431814969311\\
0.638679974267844 0.470639301222573 0.658828357290775\\
0.873411670291924 0.505699577758516 0.940054913650896\\
0.732705912657507 0.530311976941279 0.816643395912019\\
0.964607985837518 0.31492514986018 0.888979504206689\\
0.94266717920143 0.589952678564109 0.917209466237829\\
0.576749240556923 0.571185597056978 0.596379645894414\\
0.957490711898216 0.528601842529843 0.966552422655988\\
0.830914484939197 0.694004136740118 0.808004296754529\\
0.535997183025795 0.474256633966791 0.512406746249571\\
0.543195953570968 0.589344789409276 0.57568796845776\\
0.869065550528584 0.382356373953292 0.880319951334876\\
0.791276864595015 0.468278281572604 0.929422494763623\\
0.578746863621295 0.582739120649017 0.607745070447274\\
0.956245399332731 0.706092484113529 0.87504814916121\\
0.622942837668894 0.644073384656418 0.651150136145143\\
0.565985873894862 0.364046567175199 0.575201380902627\\
0.789786835456242 0.405256210517513 0.860391665576196\\
0.672958216453504 0.34073261981703 0.695216850344088\\
0.802669041217025 0.596071307223995 0.839698373150363\\
0.935047292474637 0.647493043624819 0.896697776824583\\
0.908829310988788 0.624202335819216 0.890405325860163\\
0.961205847705437 0.412931052453473 0.948590151773737\\
0.670425202091271 0.589316936077182 0.716100301195496\\
0.809658643377415 0.432870142943479 0.891096198228246\\
0.6712806408244 0.38353675211867 0.692151176392364\\
0.902630827451677 0.321359962889079 0.860489360994115\\
0.703072531314131 0.574549578288336 0.766735952219584\\
0.547896315556114 0.342398023585522 0.518163850824825\\
0.972330940570122 0.582904678276249 0.941648431347019\\
0.678205364875292 0.564728636460807 0.735058836375957\\
0.642098093928631 0.379125673258032 0.663602657216547\\
0.751025442208974 0.58896013965036 0.805073694294529\\
0.721625935559945 0.252011485651552 0.706026848452305\\
0.942298156389932 0.406466802338889 0.945052553534018\\
0.979906883588691 0.530680066943437 0.966510430673195\\
0.783828554204179 0.575665967638753 0.835769426932146\\
0.866626713769668 0.514649457161133 0.938272345516284\\
0.760762604834385 0.659430632841898 0.769158392861616\\
0.543295865291924 0.416267580225373 0.532512762567709\\
0.982508217600462 0.38017318360995 0.93588191958865\\
0.942092747848833 0.382032002277144 0.923491831296274\\
0.817249939839354 0.264316819522339 0.79239897065449\\
0.861157368692872 0.457323928262554 0.927193011057767\\
0.859715064967857 0.613175275368393 0.863111647797975\\
0.800533289280708 0.319337808075139 0.810119557705202\\
0.734069621143138 0.380003045220266 0.774945572353662\\
0.81696018134481 0.662465188355179 0.810661655188596\\
0.665854300939101 0.715876169252962 0.64585308377486\\
0.879373911625107 0.373192954520311 0.875387293823503\\
0.97426709327304 0.395825166712584 0.935876048207923\\
0.573687924020635 0.50315311808187 0.549189000524449\\
0.696025311674954 0.601471771434695 0.732804783621415\\
0.549451515527481 0.522176640084032 0.54546461926932\\
0.733492408538993 0.568997944965594 0.794405370160141\\
0.634401544897786 0.30694619947175 0.65498412936384\\
0.875440888430685 0.607321435292432 0.876136027616009\\
0.971286121553407 0.405143118697132 0.948568789009619\\
0.668750702074159 0.411499750601688 0.690334832997657\\
0.827322055000691 0.634115377637947 0.841617657483547\\
0.721767632463402 0.329166134213662 0.746562617057691\\
0.57157316852064 0.59112884077839 0.604755296796067\\
0.785241574553222 0.554601571978825 0.877040589605434\\
0.939391442171645 0.516290193501678 0.95943646368972\\
0.868834348683977 0.441901831283072 0.93229100720996\\
0.525470794349296 0.503514562234064 0.510531427266999\\
0.744236395914103 0.664753411807545 0.742240282355334\\
0.979659004758977 0.275471362862015 0.865569845835385\\
0.59475169085539 0.607946713978303 0.622909486314154\\
0.673591549846092 0.730304578210596 0.645881527328834\\
0.601673188688195 0.632126487108433 0.633806494036328\\
0.763778061552928 0.502029186103611 0.896456307649101\\
};
\end{axis}
\end{tikzpicture}%

%% file: Realisation_Zoom_epsilon_0p0906_n1280_Model1_2D.tikz
\definecolor{mycolor1}{rgb}{0,0.447,0.741}

\begin{tikzpicture}

\begin{axis}[%
width=\figurewidth,
height=\figureheight,
unbounded coords=jump,
clip=false,
view={-20}{10},
scale only axis,
every outer x axis line/.append style={darkgray!60!black},
every x tick label/.append style={font=\color{darkgray!60!black}},
xmin=0.5,
xmax=1,
xmajorgrids,
every outer y axis line/.append style={darkgray!60!black},
every y tick label/.append style={font=\color{darkgray!60!black}},
ymin=0.25,
ymax=0.75,
ymajorgrids,
every outer z axis line/.append style={darkgray!60!black},
every z tick label/.append style={font=\color{darkgray!60!black}},
zmin=0.4,
zmax=1,
zmajorgrids,
hide axis,
grid style={solid},
]

\addplot3[%
mesh,
colormap={mymap}{[1pt] rgb(0pt)=(0.2081,0.1663,0.5292); rgb(1pt)=(0.211624,0.189781,0.577676); rgb(2pt)=(0.212252,0.213771,0.626971); rgb(3pt)=(0.2081,0.2386,0.677086); rgb(4pt)=(0.195905,0.264457,0.7279); rgb(5pt)=(0.170729,0.291938,0.779248); rgb(6pt)=(0.125271,0.324243,0.830271); rgb(7pt)=(0.0591333,0.359833,0.868333); rgb(8pt)=(0.0116952,0.38751,0.881957); rgb(9pt)=(0.00595714,0.408614,0.882843); rgb(10pt)=(0.0165143,0.4266,0.878633); rgb(11pt)=(0.0328524,0.443043,0.871957); rgb(12pt)=(0.0498143,0.458571,0.864057); rgb(13pt)=(0.0629333,0.47369,0.855438); rgb(14pt)=(0.0722667,0.488667,0.8467); rgb(15pt)=(0.0779429,0.503986,0.838371); rgb(16pt)=(0.0793476,0.520024,0.831181); rgb(17pt)=(0.0749429,0.537543,0.826271); rgb(18pt)=(0.0640571,0.556986,0.823957); rgb(19pt)=(0.0487714,0.577224,0.822829); rgb(20pt)=(0.0343429,0.596581,0.819852); rgb(21pt)=(0.0265,0.6137,0.8135); rgb(22pt)=(0.0238905,0.628662,0.803762); rgb(23pt)=(0.0230905,0.641786,0.791267); rgb(24pt)=(0.0227714,0.653486,0.776757); rgb(25pt)=(0.0266619,0.664195,0.760719); rgb(26pt)=(0.0383714,0.674271,0.743552); rgb(27pt)=(0.0589714,0.683757,0.725386); rgb(28pt)=(0.0843,0.692833,0.706167); rgb(29pt)=(0.113295,0.7015,0.685857); rgb(30pt)=(0.145271,0.709757,0.664629); rgb(31pt)=(0.180133,0.717657,0.642433); rgb(32pt)=(0.217829,0.725043,0.619262); rgb(33pt)=(0.258643,0.731714,0.595429); rgb(34pt)=(0.302171,0.737605,0.571186); rgb(35pt)=(0.348167,0.742433,0.547267); rgb(36pt)=(0.395257,0.7459,0.524443); rgb(37pt)=(0.44201,0.748081,0.503314); rgb(38pt)=(0.487124,0.749062,0.483976); rgb(39pt)=(0.530029,0.749114,0.466114); rgb(40pt)=(0.570857,0.748519,0.44939); rgb(41pt)=(0.609852,0.747314,0.433686); rgb(42pt)=(0.6473,0.7456,0.4188); rgb(43pt)=(0.683419,0.743476,0.404433); rgb(44pt)=(0.71841,0.741133,0.390476); rgb(45pt)=(0.752486,0.7384,0.376814); rgb(46pt)=(0.785843,0.735567,0.363271); rgb(47pt)=(0.818505,0.732733,0.34979); rgb(48pt)=(0.850657,0.7299,0.336029); rgb(49pt)=(0.882433,0.727433,0.3217); rgb(50pt)=(0.913933,0.725786,0.306276); rgb(51pt)=(0.944957,0.726114,0.288643); rgb(52pt)=(0.973895,0.731395,0.266648); rgb(53pt)=(0.993771,0.745457,0.240348); rgb(54pt)=(0.999043,0.765314,0.216414); rgb(55pt)=(0.995533,0.786057,0.196652); rgb(56pt)=(0.988,0.8066,0.179367); rgb(57pt)=(0.978857,0.827143,0.163314); rgb(58pt)=(0.9697,0.848138,0.147452); rgb(59pt)=(0.962586,0.870514,0.1309); rgb(60pt)=(0.958871,0.8949,0.113243); rgb(61pt)=(0.959824,0.921833,0.0948381); rgb(62pt)=(0.9661,0.951443,0.0755333); rgb(63pt)=(0.9763,0.9831,0.0538)},
shader=flat,
mesh/rows=21]
table[row sep=crcr,header=false] {
0.5 0.25 NaN\\
0.5 0.275 NaN\\
0.5 0.3 NaN\\
0.5 0.325 NaN\\
0.5 0.35 NaN\\
0.5 0.375 NaN\\
0.5 0.4 NaN\\
0.5 0.425 NaN\\
0.5 0.45 NaN\\
0.5 0.475 NaN\\
0.5 0.5 NaN\\
0.5 0.525 NaN\\
0.5 0.55 NaN\\
0.5 0.575 NaN\\
0.5 0.6 NaN\\
0.5 0.625 NaN\\
0.5 0.65 NaN\\
0.5 0.675 NaN\\
0.5 0.7 NaN\\
0.5 0.725 NaN\\
0.5 0.75 NaN\\
0.525 0.25 NaN\\
0.525 0.275 NaN\\
0.525 0.3 0.505730003763746\\
0.525 0.325 0.5065976445417\\
0.525 0.35 0.505884999277199\\
0.525 0.375 0.502839290111673\\
0.525 0.4 0.501085127420753\\
0.525 0.425 0.50418970265582\\
0.525 0.45 0.501533427436008\\
0.525 0.475 0.503276614137728\\
0.525 0.5 0.509179008317451\\
0.525 0.525 0.507957691059534\\
0.525 0.55 0.519859497357345\\
0.525 0.575 0.528548823011499\\
0.525 0.6 0.53577425329171\\
0.525 0.625 0.535789213761267\\
0.525 0.65 0.542049280734487\\
0.525 0.675 0.538983593010904\\
0.525 0.7 0.529573892309283\\
0.525 0.725 0.520690900491271\\
0.525 0.75 NaN\\
0.55 0.25 NaN\\
0.55 0.275 0.535037062797594\\
0.55 0.3 0.53940497387829\\
0.55 0.325 0.536080590335658\\
0.55 0.35 0.540021316256696\\
0.55 0.375 0.536340476035456\\
0.55 0.4 0.533418536963156\\
0.55 0.425 0.529273985307567\\
0.55 0.45 0.523434115230406\\
0.55 0.475 0.522240003766481\\
0.55 0.5 0.519973394700055\\
0.55 0.525 0.525419218808458\\
0.55 0.55 0.542480166146791\\
0.55 0.575 0.557261427268842\\
0.55 0.6 0.559796184231452\\
0.55 0.625 0.555988346103031\\
0.55 0.65 0.556689654500962\\
0.55 0.675 0.553500353853511\\
0.55 0.7 0.546641574327515\\
0.55 0.725 0.540690489076202\\
0.55 0.75 NaN\\
0.575 0.25 NaN\\
0.575 0.275 0.558868494591965\\
0.575 0.3 0.565200500365119\\
0.575 0.325 0.564613460187985\\
0.575 0.35 0.565326408665295\\
0.575 0.375 0.566049738217487\\
0.575 0.4 0.564866156687532\\
0.575 0.425 0.563447763446605\\
0.575 0.45 0.553587189294437\\
0.575 0.475 0.546484343959165\\
0.575 0.5 0.558278753226086\\
0.575 0.525 0.561929921864569\\
0.575 0.55 0.572169773200533\\
0.575 0.575 0.581792457995283\\
0.575 0.6 0.589326119587516\\
0.575 0.625 0.585826679157007\\
0.575 0.65 0.581554595254476\\
0.575 0.675 0.575697101173421\\
0.575 0.7 0.570322575603951\\
0.575 0.725 0.563163880084322\\
0.575 0.75 NaN\\
0.6 0.25 NaN\\
0.6 0.275 0.576205360393123\\
0.6 0.3 0.587096537593027\\
0.6 0.325 0.593727737704642\\
0.6 0.35 0.594533528416055\\
0.6 0.375 0.592007707911911\\
0.6 0.4 0.595392009632435\\
0.6 0.425 0.594123980317134\\
0.6 0.45 0.593702718099127\\
0.6 0.475 0.592994012036962\\
0.6 0.5 0.593228758096977\\
0.6 0.525 0.593857582429895\\
0.6 0.55 0.602554839474169\\
0.6 0.575 0.610670636728123\\
0.6 0.6 0.612197012318404\\
0.6 0.625 0.611634836425679\\
0.6 0.65 0.609441124317899\\
0.6 0.675 0.600364419033768\\
0.6 0.7 0.591425951204787\\
0.6 0.725 0.581551947156697\\
0.6 0.75 NaN\\
0.625 0.25 NaN\\
0.625 0.275 0.598290404613452\\
0.625 0.3 0.606746872073494\\
0.625 0.325 0.617785307702058\\
0.625 0.35 0.618725137102013\\
0.625 0.375 0.624306058522888\\
0.625 0.4 0.618662312371039\\
0.625 0.425 0.614643305146777\\
0.625 0.45 0.616189245237014\\
0.625 0.475 0.617618649064986\\
0.625 0.5 0.622685806077467\\
0.625 0.525 0.636421324704123\\
0.625 0.55 0.638856696284089\\
0.625 0.575 0.641072293461933\\
0.625 0.6 0.635222743479026\\
0.625 0.625 0.632167345027444\\
0.625 0.65 0.627320164850614\\
0.625 0.675 0.620904094650053\\
0.625 0.7 0.613319902591134\\
0.625 0.725 0.603182072712606\\
0.625 0.75 NaN\\
0.65 0.25 NaN\\
0.65 0.275 0.628835845067497\\
0.65 0.3 0.638990803991658\\
0.65 0.325 0.645372895340749\\
0.65 0.35 0.650346942417823\\
0.65 0.375 0.655554965293675\\
0.65 0.4 0.656667329723027\\
0.65 0.425 0.663304345205888\\
0.65 0.45 0.653076458356238\\
0.65 0.475 0.65111960618453\\
0.65 0.5 0.665927095853403\\
0.65 0.525 0.680832339589404\\
0.65 0.55 0.678620238335889\\
0.65 0.575 0.674263812428651\\
0.65 0.6 0.661875561065739\\
0.65 0.625 0.652807067810176\\
0.65 0.65 0.646347449546257\\
0.65 0.675 0.639758798196502\\
0.65 0.7 0.631330065718046\\
0.65 0.725 0.619223738894044\\
0.65 0.75 NaN\\
0.675 0.25 NaN\\
0.675 0.275 0.642562964587433\\
0.675 0.3 0.650475757107099\\
0.675 0.325 0.666601317340951\\
0.675 0.35 0.670444684456261\\
0.675 0.375 0.674220365828262\\
0.675 0.4 0.673778163107887\\
0.675 0.425 0.68177446707959\\
0.675 0.45 0.685914374821512\\
0.675 0.475 0.69544713007532\\
0.675 0.5 0.703240747854489\\
0.675 0.525 0.715567653676958\\
0.675 0.55 0.710642778284496\\
0.675 0.575 0.704604490765785\\
0.675 0.6 0.691428192720632\\
0.675 0.625 0.677647833139797\\
0.675 0.65 0.66458493773553\\
0.675 0.675 0.657706617901029\\
0.675 0.7 0.647891185860847\\
0.675 0.725 0.629558980863508\\
0.675 0.75 NaN\\
0.7 0.25 NaN\\
0.7 0.275 0.655866350030611\\
0.7 0.3 0.670628779985895\\
0.7 0.325 0.678001622406892\\
0.7 0.35 0.685894972905076\\
0.7 0.375 0.698399944525613\\
0.7 0.4 0.708518776222157\\
0.7 0.425 0.718628212634569\\
0.7 0.45 0.72903465256873\\
0.7 0.475 0.735951012490279\\
0.7 0.5 0.741812809711233\\
0.7 0.525 0.741317392328989\\
0.7 0.55 0.737278515819737\\
0.7 0.575 0.727703789877266\\
0.7 0.6 0.715728248023817\\
0.7 0.625 0.70810149650648\\
0.7 0.65 0.693402335830518\\
0.7 0.675 0.673092233047896\\
0.7 0.7 0.659713785342493\\
0.7 0.725 0.64942964237006\\
0.7 0.75 NaN\\
0.725 0.25 NaN\\
0.725 0.275 0.677201972877611\\
0.725 0.3 0.691690106337018\\
0.725 0.325 0.701297322229318\\
0.725 0.35 0.713072487308701\\
0.725 0.375 0.729560287602563\\
0.725 0.4 0.740185323160931\\
0.725 0.425 0.748376410912763\\
0.725 0.45 0.760994965237459\\
0.725 0.475 0.770962544031047\\
0.725 0.5 0.780841057310341\\
0.725 0.525 0.787247724318232\\
0.725 0.55 0.788622103010713\\
0.725 0.575 0.745260610542857\\
0.725 0.6 0.7398880420806\\
0.725 0.625 0.734033839216109\\
0.725 0.65 0.722396759542103\\
0.725 0.675 0.700912108491237\\
0.725 0.7 0.685963973411697\\
0.725 0.725 0.668277533521027\\
0.725 0.75 NaN\\
0.75 0.25 NaN\\
0.75 0.275 0.695415450304939\\
0.75 0.3 0.711493439317961\\
0.75 0.325 0.722535272291906\\
0.75 0.35 0.736394360756065\\
0.75 0.375 0.75249335723162\\
0.75 0.4 0.766077735146419\\
0.75 0.425 0.77622055942886\\
0.75 0.45 0.791833879542179\\
0.75 0.475 0.806048469293984\\
0.75 0.5 0.813434412958908\\
0.75 0.525 0.810506858429086\\
0.75 0.55 0.80620916315757\\
0.75 0.575 0.784700270275829\\
0.75 0.6 0.758562421497063\\
0.75 0.625 0.750649845509299\\
0.75 0.65 0.740423035527527\\
0.75 0.675 0.724447470878132\\
0.75 0.7 0.70590403324111\\
0.75 0.725 0.687125424671993\\
0.75 0.75 NaN\\
0.775 0.25 NaN\\
0.775 0.275 0.705449926359041\\
0.775 0.3 0.722344206695456\\
0.775 0.325 0.739160923504059\\
0.775 0.35 0.75422028895861\\
0.775 0.375 0.770860677671521\\
0.775 0.4 0.78714801972731\\
0.775 0.425 0.804678115646198\\
0.775 0.45 0.825510082350529\\
0.775 0.475 0.834281187671416\\
0.775 0.5 0.879199443862494\\
0.775 0.525 0.835140559716412\\
0.775 0.55 0.823987831032967\\
0.775 0.575 0.817633155830226\\
0.775 0.6 0.775003248353594\\
0.775 0.625 0.767045273290692\\
0.775 0.65 0.758063276594447\\
0.775 0.675 0.73725544916691\\
0.775 0.7 0.721612456393909\\
0.775 0.725 0.70597331582296\\
0.775 0.75 NaN\\
0.8 0.25 NaN\\
0.8 0.275 0.716817280015107\\
0.8 0.3 0.733877741068357\\
0.8 0.325 0.751612336297429\\
0.8 0.35 0.771416172882499\\
0.8 0.375 0.788155488495518\\
0.8 0.4 0.803109352349973\\
0.8 0.425 0.82618907481516\\
0.8 0.45 0.847049995047419\\
0.8 0.475 0.856470164258406\\
0.8 0.5 1\\
0.8 0.525 0.904747213752416\\
0.8 0.55 0.841569243361193\\
0.8 0.575 0.833416896817173\\
0.8 0.6 0.796411804762365\\
0.8 0.625 0.780583529132993\\
0.8 0.65 0.768054465011129\\
0.8 0.675 0.756680900212284\\
0.8 0.7 0.736934053003466\\
0.8 0.725 0.721651675797992\\
0.8 0.75 NaN\\
0.825 0.25 NaN\\
0.825 0.275 0.731610247873562\\
0.825 0.3 0.749562530349538\\
0.825 0.325 0.764778402294299\\
0.825 0.35 0.789584606364055\\
0.825 0.375 0.804419533875224\\
0.825 0.4 0.815038737277038\\
0.825 0.425 0.834122437390819\\
0.825 0.45 0.863593918657629\\
0.825 0.475 0.867149492162574\\
0.825 0.5 0.901856387951083\\
0.825 0.525 0.866268184391602\\
0.825 0.55 0.855053500249196\\
0.825 0.575 0.842243227355607\\
0.825 0.6 0.811226291740979\\
0.825 0.625 0.791868171496089\\
0.825 0.65 0.775219963916748\\
0.825 0.675 0.763158047909202\\
0.825 0.7 0.747338082448555\\
0.825 0.725 0.731874049626397\\
0.825 0.75 NaN\\
0.85 0.25 NaN\\
0.85 0.275 0.736792716631486\\
0.85 0.3 0.755430146678379\\
0.85 0.325 0.76756662605219\\
0.85 0.35 0.781461152005648\\
0.85 0.375 0.808092637605436\\
0.85 0.4 0.820314146798248\\
0.85 0.425 0.838348118203634\\
0.85 0.45 0.868308248032013\\
0.85 0.475 0.879870266551947\\
0.85 0.5 0.884195204521203\\
0.85 0.525 0.874523292740428\\
0.85 0.55 0.853101068193112\\
0.85 0.575 0.83620615252734\\
0.85 0.6 0.815128859223125\\
0.85 0.625 0.799039271996431\\
0.85 0.65 0.780149848669174\\
0.85 0.675 0.768087774008347\\
0.85 0.7 0.755304077555901\\
0.85 0.725 0.740086465174788\\
0.85 0.75 NaN\\
0.875 0.25 NaN\\
0.875 0.275 0.745554900852889\\
0.875 0.3 0.761272982788955\\
0.875 0.325 0.771614593846421\\
0.875 0.35 0.787816839160494\\
0.875 0.375 0.80998143263285\\
0.875 0.4 0.828295276920832\\
0.875 0.425 0.85007591854109\\
0.875 0.45 0.866869412095482\\
0.875 0.475 0.876818487179677\\
0.875 0.5 0.886545589001648\\
0.875 0.525 0.880469635652059\\
0.875 0.55 0.851384044322503\\
0.875 0.575 0.832505039697394\\
0.875 0.6 0.817212529142988\\
0.875 0.625 0.80385305288\\
0.875 0.65 0.785123674019633\\
0.875 0.675 0.76983686603604\\
0.875 0.7 0.75973996441779\\
0.875 0.725 0.747126951824328\\
0.875 0.75 NaN\\
0.9 0.25 NaN\\
0.9 0.275 0.749471113777848\\
0.9 0.3 0.764397428897073\\
0.9 0.325 0.776783966952591\\
0.9 0.35 0.793316407161644\\
0.9 0.375 0.810823997184585\\
0.9 0.4 0.824733091196836\\
0.9 0.425 0.840540537687175\\
0.9 0.45 0.855734188081738\\
0.9 0.475 0.86655579428163\\
0.9 0.5 0.87305998845023\\
0.9 0.525 0.866672967341736\\
0.9 0.55 0.851891351241823\\
0.9 0.575 0.836906760573379\\
0.9 0.6 0.820346174830889\\
0.9 0.625 0.804398321408984\\
0.9 0.65 0.78730855838956\\
0.9 0.675 0.772214497821146\\
0.9 0.7 0.763182025681712\\
0.9 0.725 0.753694549052621\\
0.9 0.75 NaN\\
0.925 0.25 NaN\\
0.925 0.275 0.754830851066586\\
0.925 0.3 0.769960451239142\\
0.925 0.325 0.783476512340978\\
0.925 0.35 0.79719039018442\\
0.925 0.375 0.811797641738035\\
0.925 0.4 0.821498061426922\\
0.925 0.425 0.838680715921462\\
0.925 0.45 0.851758451668761\\
0.925 0.475 0.854277515239211\\
0.925 0.5 0.859024693057625\\
0.925 0.525 0.856025605564711\\
0.925 0.55 0.845698069740811\\
0.925 0.575 0.832576212541958\\
0.925 0.6 0.817708659015868\\
0.925 0.625 0.805199943455212\\
0.925 0.65 0.793497812864876\\
0.925 0.675 0.776157596678216\\
0.925 0.7 0.76663981599101\\
0.925 0.725 0.758440775292443\\
0.925 0.75 NaN\\
0.95 0.25 NaN\\
0.95 0.275 0.760190588355324\\
0.95 0.3 0.775023887767379\\
0.95 0.325 0.789540879998353\\
0.95 0.35 0.802237382660757\\
0.95 0.375 0.813780103515684\\
0.95 0.4 0.821672078769581\\
0.95 0.425 0.832827458950064\\
0.95 0.45 0.846585274745776\\
0.95 0.475 0.849130104770524\\
0.95 0.5 0.851516260874012\\
0.95 0.525 0.852189413839313\\
0.95 0.55 0.840957590336063\\
0.95 0.575 0.827737547600812\\
0.95 0.6 0.817611559652169\\
0.95 0.625 0.806064495672393\\
0.95 0.65 0.796016172460998\\
0.95 0.675 0.779566615486988\\
0.95 0.7 0.770769141030506\\
0.95 0.725 0.761819935267235\\
0.95 0.75 NaN\\
0.975 0.25 NaN\\
0.975 0.275 0.765550325644062\\
0.975 0.3 0.776706922144215\\
0.975 0.325 0.789773532848773\\
0.975 0.35 0.803878361900613\\
0.975 0.375 0.813993455807846\\
0.975 0.4 0.821286218787528\\
0.975 0.425 0.830365144967029\\
0.975 0.45 0.838816848692432\\
0.975 0.475 0.843100850344941\\
0.975 0.5 0.846607211611326\\
0.975 0.525 0.844420462008928\\
0.975 0.55 0.837904690846226\\
0.975 0.575 0.826462798303389\\
0.975 0.6 0.815490890501091\\
0.975 0.625 0.801717259075703\\
0.975 0.65 0.790599163634862\\
0.975 0.675 0.780419127808687\\
0.975 0.7 0.773181613828203\\
0.975 0.725 0.763946612766412\\
0.975 0.75 NaN\\
1 0.25 NaN\\
1 0.275 NaN\\
1 0.3 NaN\\
1 0.325 NaN\\
1 0.35 NaN\\
1 0.375 NaN\\
1 0.4 NaN\\
1 0.425 NaN\\
1 0.45 NaN\\
1 0.475 NaN\\
1 0.5 NaN\\
1 0.525 NaN\\
1 0.55 NaN\\
1 0.575 NaN\\
1 0.6 NaN\\
1 0.625 NaN\\
1 0.65 NaN\\
1 0.675 NaN\\
1 0.7 NaN\\
1 0.725 NaN\\
1 0.75 NaN\\
};
\addplot3 [
color=black, 
mark size=1.0pt,
only marks,
mark=*,
mark options={solid}]
table[row sep=crcr] {
0.8 0.5 1\\
0.896750966300958 0.286776127457757 0.756600427310135\\
0.902771977462955 0.443307096739818 0.851085433350448\\
0.610354789491249 0.590563611606405 0.625491942624764\\
0.717545631591306 0.636844667860001 0.725714795502098\\
0.852848313794446 0.435990016586211 0.866343148009527\\
0.660761305835014 0.452417241303641 0.665345985582683\\
0.820130475649893 0.320075943110966 0.762084045995138\\
0.610456067310923 0.319729894706417 0.601850490032599\\
0.996628822967953 0.281100166085246 0.770084132730866\\
0.829330911141074 0.36177041731311 0.802175202229373\\
0.738927787432516 0.702469856909833 0.698004482132121\\
0.626578235140411 0.702924333907962 0.613843534927718\\
0.981412845867243 0.705104676070429 0.773058658654732\\
0.679608034336165 0.657771390884537 0.664312711987929\\
0.7767791920296 0.333471872572624 0.745281339251403\\
0.905025376762873 0.553778152360609 0.84749346767195\\
0.94166257336887 0.316020689047971 0.782589314556037\\
0.818411855866936 0.603014587398564 0.807011333230559\\
0.666906186164037 0.438275770234889 0.669803835943133\\
0.813503588228468 0.595511914424481 0.810504598727775\\
0.625717367146037 0.668026662527658 0.623660832707685\\
0.906993019500661 0.662923413406259 0.777901069693883\\
0.572131847654909 0.48344100136562 0.539154348091194\\
0.903338950738432 0.571324421969113 0.839132185922889\\
0.986631397979414 0.628632455243067 0.79543554779755\\
0.890279572452008 0.746189348982611 0.736689338316952\\
0.804770328186065 0.384360635226165 0.80237715562653\\
0.556768447472082 0.688768012719625 0.554952277463858\\
0.597875400608463 0.402327162870215 0.593531057183025\\
0.734634322479852 0.517831329437117 0.805211437277669\\
0.588734931310846 0.731851647500458 0.571205140133708\\
0.805532816497391 0.641278569924126 0.774233391499405\\
0.761772587446987 0.495623970862494 0.822057526612621\\
0.73431033983156 0.374160408174907 0.735631958767502\\
0.957096547913993 0.505259946445197 0.850313893478733\\
0.780291105767305 0.454155929647519 0.836078245364494\\
0.892702383484364 0.563241363149134 0.842714276207149\\
0.612013984234454 0.339793008902145 0.610166882498886\\
0.50868945810227 0.624583002688413 0.522420228159328\\
0.598022997067236 0.648736113587461 0.607694302796689\\
0.712619762620216 0.373497888737766 0.718594653863858\\
0.64966640763613 0.423645275337246 0.663447581829482\\
0.957072824374163 0.684798012434607 0.774570356922379\\
0.670982508455675 0.451325421246072 0.67690733140817\\
0.87348890039871 0.636667687206959 0.797940844124299\\
0.91187476865064 0.72621176631384 0.756185170252988\\
0.689641681209292 0.635260414894502 0.692288626653827\\
0.557804714387723 0.625122397308514 0.562324100094365\\
0.800424794384846 0.697439096217222 0.738750330435334\\
0.504686378036765 0.419115029941283 0.47870638607057\\
0.529101751211831 0.641041544298777 0.547018383220701\\
0.953672341494932 0.329437527145471 0.792974851231587\\
0.556778717749001 0.44797090722905 0.531267205238796\\
0.600894508295261 0.465190351306326 0.594117745481038\\
0.833294096726429 0.603108885353926 0.81183476921422\\
0.989547621752577 0.736941860120139 0.763295382352825\\
0.697173860239204 0.637819207435455 0.700108546278864\\
0.508927167607604 0.532144542167342 0.49794482745362\\
0.681553544189616 0.572471948843176 0.716693312357688\\
0.557628145279771 0.527774442485995 0.540128436305281\\
0.950399058359332 0.735603954354949 0.756206858043755\\
0.510836847606421 0.345729971491733 0.48664687447046\\
0.841152553007951 0.478634600177314 0.87665660211336\\
0.560435455011861 0.666191663819268 0.560679185795816\\
0.944134517275711 0.734594407286704 0.758003129004671\\
0.900442140858193 0.332479457434833 0.779139655482648\\
0.632193997405328 0.369266158090538 0.637222428799586\\
0.516461734167434 0.42273759940176 0.498092193696605\\
0.69121462913145 0.377368393397148 0.685333082565018\\
0.649212958682711 0.599332068352685 0.661042672248653\\
0.57079593357119 0.63423768817409 0.581290038758588\\
0.723136464599151 0.576815374708839 0.743930834209755\\
0.928378223516919 0.349078429673212 0.796964584742904\\
0.505653869344781 0.748549542309387 0.498677241736679\\
0.545763188443633 0.48466628529025 0.513715809218719\\
0.825722294028466 0.515472651567212 0.871348738507981\\
0.71129147049187 0.33969018138519 0.690998066922664\\
0.961465217690023 0.315560651972318 0.782298870836103\\
0.952407770222575 0.643638364186628 0.800695857834603\\
0.747516926326245 0.313847884230513 0.719479224275208\\
0.796674491635933 0.388888360603717 0.795953132034131\\
0.672641662882902 0.634830830121167 0.670735847288633\\
0.603313292450109 0.65347035723999 0.61225671162552\\
0.796814421268911 0.263945218665011 0.707168746664155\\
0.735269940189609 0.383113193894543 0.745171282649382\\
0.656979767272955 0.620716536738546 0.659335759187515\\
0.748165576186312 0.259010531534803 0.68483410857601\\
0.939433799432894 0.665870143073621 0.783406736401045\\
0.728568244242342 0.551385087432613 0.795756431735347\\
0.650647237047691 0.275570801569803 0.6316745831971\\
0.806878408035618 0.62185710361723 0.785259994784867\\
0.544194662207204 0.302608389824434 0.533833913780756\\
0.622955861825668 0.577622950816906 0.637916292223549\\
0.818883711816406 0.734932788204561 0.722414995647845\\
0.846119357594698 0.6915399704329 0.760675562912812\\
0.553766009195282 0.575454247763055 0.561580487701363\\
0.555666158880712 0.568677618897226 0.561002475759902\\
0.604435974345879 0.640463967654522 0.614384823887607\\
0.526349642470887 0.56563144105305 0.527534056359664\\
0.721059175232061 0.297450940320884 0.687387338009301\\
0.84260103603702 0.345530155179932 0.775716265931569\\
0.759685718167462 0.557020628247086 0.809950050985227\\
0.923017767921153 0.69036343882904 0.768676111902648\\
0.557512611368792 0.261906665471005 0.54103736113577\\
0.544360948798569 0.468744207168826 0.515553221317979\\
0.629938261821685 0.606467334820059 0.637871809217634\\
0.917675500102484 0.733426039517901 0.755998075263485\\
0.68074167652991 0.688535317260708 0.659237689303328\\
0.617326694267432 0.382581145538756 0.610789517255301\\
0.681803363676878 0.518240854847237 0.723632112882918\\
0.885290171061452 0.558902348891923 0.84320963247269\\
0.999879067839387 0.448812646700452 0.838558868964719\\
0.556653679428374 0.591775511200523 0.567417226869654\\
0.511923460323181 0.731515256518571 0.509005827376278\\
0.595068302854429 0.720931728228616 0.578059938851036\\
0.653528993732435 0.524999851559058 0.687101372856299\\
0.791550251926774 0.361782993747501 0.769927467036853\\
0.713931367162337 0.405706449224023 0.731168529630497\\
0.805166281940449 0.43454148810554 0.841394605189888\\
0.517319137001344 0.503384924916418 0.50286622023819\\
0.79562218975414 0.68645392683332 0.73966730219815\\
0.591719690181946 0.695977721199807 0.584790940122358\\
0.659452620544525 0.731721188595494 0.621385914604292\\
0.993170985365447 0.741377113835917 0.762630676560262\\
0.791545234841724 0.475131522881481 0.852952859950401\\
0.579610587190408 0.591918525702021 0.591244205062658\\
0.919737896229818 0.434663358048463 0.846003324119569\\
0.555719911475165 0.551092672120086 0.548106799574755\\
0.757941722373145 0.373864137037787 0.759536318662375\\
0.654419236693616 0.275193740552363 0.631672532659987\\
0.889949054888152 0.294867099112871 0.762215805431662\\
0.668218844330459 0.326768892520437 0.665126857436381\\
0.918302268321324 0.647871097589267 0.79531349562855\\
0.652359263807669 0.637841326006166 0.651943619781813\\
0.553472481516696 0.476149594730903 0.527016136861549\\
0.916110193482973 0.45870182637889 0.853600658471039\\
0.78972060934824 0.60099276064908 0.784397361720823\\
0.966813818539717 0.434865117015749 0.835325466938516\\
0.853959402955557 0.485460244964119 0.88892255366072\\
0.540527794749159 0.431787412818572 0.514172713612593\\
0.505152550821154 0.515459444839436 0.48767706013044\\
0.578734333094401 0.690784956846228 0.576768330318541\\
0.897592273023424 0.57268934326616 0.83848459359077\\
0.80170852361882 0.74007759751037 0.713242504519496\\
0.956214471123851 0.637671814284754 0.8006858223911\\
0.514678816739023 0.505035171806744 0.493328526123042\\
0.665406548725784 0.429975246137746 0.669666253132053\\
0.972869269087839 0.53813192279225 0.844198748921147\\
0.71483991465287 0.570295440291489 0.743238976737193\\
0.523531322672027 0.503042847571715 0.508701906610123\\
0.848968807250136 0.526796175461105 0.872672690034757\\
0.763914494242464 0.416629221251442 0.786169057035513\\
0.701662840439768 0.538870654105546 0.744103340440856\\
0.979893450066785 0.454617438907175 0.839066629323946\\
0.803772445643009 0.682811751090593 0.754634650531871\\
0.615098939296913 0.268404367035586 0.583802723471901\\
0.798968557213536 0.38128354172878 0.790674799217971\\
0.890580710042851 0.738437300921625 0.74333873304575\\
0.967825084618818 0.536888363538073 0.84419760125058\\
0.69581419182488 0.496871966433736 0.734169908803127\\
0.544600752813757 0.40197995989309 0.52653306031234\\
0.936213193985016 0.444417592226183 0.849701383669099\\
0.598289326940254 0.635668147698973 0.608521251924712\\
0.851877230121635 0.653673406242905 0.777581562149035\\
0.752586851411443 0.472454038691497 0.808678602100557\\
0.837450216826856 0.522313656458288 0.87308775474645\\
0.889743706793846 0.28439031811912 0.756571664447935\\
0.981723549742713 0.345842180573096 0.802398025667583\\
0.862315767578514 0.26179302017712 0.732615259519439\\
0.529272376947821 0.405989155670227 0.505619372343215\\
0.967912669818733 0.592303195473019 0.823923334832047\\
0.830431113487101 0.539488684926696 0.861767285115928\\
0.621571966395797 0.414129110039419 0.610488811891074\\
0.57442603944283 0.735784725196401 0.559490560938433\\
0.687742464161246 0.455322621044195 0.713135339300817\\
0.590071029823612 0.618542337447304 0.604801400125034\\
0.617331885338194 0.405718710051996 0.607006454594955\\
0.72690845906154 0.70082605744322 0.687777289535772\\
0.634045924174415 0.63700555975707 0.637609057594735\\
0.752559219350107 0.257357444298136 0.684826251841836\\
0.75771271505066 0.577808150277288 0.800674966939356\\
0.838475119768797 0.465907037408674 0.866383131762249\\
0.844003912299951 0.510567992926764 0.880050589175656\\
0.976611594793114 0.685936179491611 0.775430450428259\\
0.981935894624428 0.493982272449778 0.845535777397073\\
0.923052998305207 0.457167450939811 0.85443539270505\\
0.652965308525666 0.607145655266464 0.663254432326221\\
0.860857500756139 0.384450601326728 0.814260747667335\\
0.622211396677159 0.355232097780553 0.61539133186061\\
0.791457906927876 0.436562126700811 0.831381115488789\\
0.674450201334774 0.635700701390973 0.674750321084723\\
0.864458138705059 0.494853097534895 0.8888356982466\\
0.665835532037635 0.457394875143017 0.67773497556506\\
0.621725699247968 0.711399117134646 0.608739882868158\\
0.770754244073012 0.650584031585053 0.756100769853849\\
0.730986697506542 0.295438385894555 0.695099216719694\\
0.800798167895389 0.356574648567488 0.777072596351549\\
0.592460076523037 0.456185420555379 0.580570501178591\\
0.842229522803865 0.277063905713421 0.739796725974936\\
0.526548414257021 0.267818589565266 0.506872939369578\\
0.835137254585206 0.628426907790893 0.794160023665696\\
0.908628810594412 0.741702177161977 0.751001965150809\\
0.8532643740468 0.42243260461999 0.833041332044654\\
0.522834301020647 0.670081960160418 0.539392805080743\\
0.946148794562511 0.749434273566672 0.753328434336947\\
0.604877099274976 0.320179925465855 0.599270007983375\\
0.821159119762907 0.450084003570264 0.863095879787885\\
0.655707947618119 0.741305278209221 0.617201508686083\\
0.818253846190852 0.578555867900184 0.841133815749429\\
0.952615054258229 0.427969357910103 0.83394994458901\\
0.564923381251329 0.489770209480822 0.53176070209085\\
0.676839393873037 0.303291533227662 0.652496372514176\\
0.848021602973121 0.2690429776502 0.731907649744528\\
0.734391454067498 0.417013474230646 0.754208712595242\\
0.976260164764427 0.517455232249628 0.844713556955216\\
0.646539376807889 0.574127157107886 0.670365783889019\\
0.708068392333254 0.689023454361018 0.670718063938597\\
0.565127848267236 0.457550272907665 0.536297207730647\\
0.923272787027268 0.693320898060286 0.768585776597017\\
0.64958740767562 0.368280328610595 0.657308897859183\\
0.755269923936985 0.592352884283737 0.764422957516216\\
0.829613366395353 0.709263676905948 0.745192215286672\\
0.687090495956269 0.560597930130718 0.719115147357241\\
0.708352078901677 0.459898933306801 0.74550441788502\\
0.60356654986875 0.524604656439891 0.598307561452423\\
0.633827593627922 0.728606225258041 0.609071925291612\\
0.806227835509854 0.722240128754432 0.727238363904433\\
0.509598833595111 0.472527691077761 0.488577070902395\\
0.795448592145339 0.40713967529347 0.802303053784176\\
0.857848427994472 0.698109151242465 0.758001169783952\\
0.531920624360976 0.741644972972042 0.521091114628631\\
0.8557431661952 0.30361060242885 0.758757656293987\\
0.96809664114014 0.518059368780743 0.844262564212831\\
0.563520581945619 0.625875894105425 0.570962664437799\\
0.653411476759634 0.314888736357681 0.649024730444214\\
0.680924942343551 0.525120814126042 0.7234400198652\\
0.991879636942119 0.462982362699806 0.840463563585964\\
0.595393067151581 0.427187978982608 0.590631340469333\\
0.963869982901046 0.519546686704649 0.844200616220385\\
0.966957924360731 0.745079342721426 0.754339125323728\\
0.553432214382848 0.320014753868527 0.539354071407345\\
0.638679974267844 0.470639301222573 0.630685772413626\\
0.873411670291924 0.505699577758516 0.889424352324636\\
0.732705912657507 0.530311976941279 0.800449244452923\\
0.964607985837518 0.31492514986018 0.782293824728856\\
0.94266717920143 0.589952678564109 0.821281473521209\\
0.576749240556923 0.571185597056978 0.582871166962924\\
0.957490711898216 0.528601842529843 0.850820195237762\\
0.830914484939197 0.694004136740118 0.752234078830232\\
0.535997183025795 0.474256633966791 0.513167146523743\\
0.543195953570968 0.589344789409276 0.555518593967535\\
0.869065550528584 0.382356373953292 0.816116168893232\\
0.791276864595015 0.468278281572604 0.851760724915234\\
0.578746863621295 0.582739120649017 0.586572843413746\\
0.956245399332731 0.706092484113529 0.770418501916214\\
0.622942837668894 0.644073384656418 0.627066066163852\\
0.565985873894862 0.364046567175199 0.55522939705265\\
0.789786835456242 0.405256210517513 0.802301570084598\\
0.672958216453504 0.34073261981703 0.667690759908771\\
0.802669041217025 0.596071307223995 0.802659464360873\\
0.935047292474637 0.647493043624819 0.793895732486065\\
0.908829310988788 0.624202335819216 0.804992831572784\\
0.961205847705437 0.412931052453473 0.827927832104061\\
0.670425202091271 0.589316936077182 0.694137789594199\\
0.809658643377415 0.432870142943479 0.845096504863799\\
0.6712806408244 0.38353675211867 0.673042628945019\\
0.902630827451677 0.321359962889079 0.776384869578408\\
0.703072531314131 0.574549578288336 0.729986851181675\\
0.547896315556114 0.342398023585522 0.538867094888061\\
0.972330940570122 0.582904678276249 0.824182025234005\\
0.678205364875292 0.564728636460807 0.711976850132584\\
0.642098093928631 0.379125673258032 0.646731730081748\\
0.751025442208974 0.58896013965036 0.761833885291899\\
0.721625935559945 0.252011485651552 0.660098142508556\\
0.942298156389932 0.406466802338889 0.823704332602425\\
0.979906883588691 0.530680066943437 0.844510301319739\\
0.783828554204179 0.575665967638753 0.825540536298621\\
0.866626713769668 0.514649457161133 0.892176490122787\\
0.760762604834385 0.659430632841898 0.741240358432284\\
0.543295865291924 0.416267580225373 0.522874584323395\\
0.982508217600462 0.38017318360995 0.816124974405832\\
0.942092747848833 0.382032002277144 0.816522934757745\\
0.817249939839354 0.264316819522339 0.720923021472895\\
0.861157368692872 0.457323928262554 0.875675287098871\\
0.859715064967857 0.613175275368393 0.807171193245844\\
0.800533289280708 0.319337808075139 0.747426534505728\\
0.734069621143138 0.380003045220266 0.741903684922555\\
0.81696018134481 0.662465188355179 0.768542938226797\\
0.665854300939101 0.715876169252962 0.631563115493511\\
0.879373911625107 0.373192954520311 0.808977483574146\\
0.97426709327304 0.395825166712584 0.820032575318134\\
0.573687924020635 0.50315311808187 0.557494143834377\\
0.696025311674954 0.601471771434695 0.711409263078741\\
0.549451515527481 0.522176640084032 0.523019956840475\\
0.733492408538993 0.568997944965594 0.750507496051526\\
0.634401544897786 0.30694619947175 0.618492059461549\\
0.875440888430685 0.607321435292432 0.812770750722113\\
0.971286121553407 0.405143118697132 0.822090433086908\\
0.668750702074159 0.411499750601688 0.669079167184386\\
0.827322055000691 0.634115377637947 0.782555769834168\\
0.721767632463402 0.329166134213662 0.699751309707208\\
0.57157316852064 0.59112884077839 0.586572683039082\\
0.785241574553222 0.554601571978825 0.831952276869309\\
0.939391442171645 0.516290193501678 0.855053348123087\\
0.868834348683977 0.441901831283072 0.868118075926444\\
0.525470794349296 0.503514562234064 0.508701201195358\\
0.744236395914103 0.664753411807545 0.73032998379168\\
0.979659004758977 0.275471362862015 0.766862413365199\\
0.59475169085539 0.607946713978303 0.604838358140819\\
0.673591549846092 0.730304578210596 0.624993125160968\\
0.601673188688195 0.632126487108433 0.613984454415572\\
0.763778061552928 0.502029186103611 0.821834005719567\\
};
\end{axis}
\end{tikzpicture}%

%% file: Realisation_Zoom_epsilon_0p2_n1280_Model1_2D.tikz
\definecolor{mycolor1}{rgb}{0,0.447,0.741}

\begin{tikzpicture}

\begin{axis}[%
width=\figurewidth,
height=\figureheight,
unbounded coords=jump,
clip=false,
view={-20}{10},
scale only axis,
every outer x axis line/.append style={darkgray!60!black},
every x tick label/.append style={font=\color{darkgray!60!black}},
xmin=0.5,
xmax=1,
xmajorgrids,
every outer y axis line/.append style={darkgray!60!black},
every y tick label/.append style={font=\color{darkgray!60!black}},
ymin=0.25,
ymax=0.75,
ymajorgrids,
every outer z axis line/.append style={darkgray!60!black},
every z tick label/.append style={font=\color{darkgray!60!black}},
zmin=0.4,
zmax=1,
zmajorgrids,
hide axis,
grid style={solid},
]

\addplot3[%
mesh,
colormap={mymap}{[1pt] rgb(0pt)=(0.2081,0.1663,0.5292); rgb(1pt)=(0.211624,0.189781,0.577676); rgb(2pt)=(0.212252,0.213771,0.626971); rgb(3pt)=(0.2081,0.2386,0.677086); rgb(4pt)=(0.195905,0.264457,0.7279); rgb(5pt)=(0.170729,0.291938,0.779248); rgb(6pt)=(0.125271,0.324243,0.830271); rgb(7pt)=(0.0591333,0.359833,0.868333); rgb(8pt)=(0.0116952,0.38751,0.881957); rgb(9pt)=(0.00595714,0.408614,0.882843); rgb(10pt)=(0.0165143,0.4266,0.878633); rgb(11pt)=(0.0328524,0.443043,0.871957); rgb(12pt)=(0.0498143,0.458571,0.864057); rgb(13pt)=(0.0629333,0.47369,0.855438); rgb(14pt)=(0.0722667,0.488667,0.8467); rgb(15pt)=(0.0779429,0.503986,0.838371); rgb(16pt)=(0.0793476,0.520024,0.831181); rgb(17pt)=(0.0749429,0.537543,0.826271); rgb(18pt)=(0.0640571,0.556986,0.823957); rgb(19pt)=(0.0487714,0.577224,0.822829); rgb(20pt)=(0.0343429,0.596581,0.819852); rgb(21pt)=(0.0265,0.6137,0.8135); rgb(22pt)=(0.0238905,0.628662,0.803762); rgb(23pt)=(0.0230905,0.641786,0.791267); rgb(24pt)=(0.0227714,0.653486,0.776757); rgb(25pt)=(0.0266619,0.664195,0.760719); rgb(26pt)=(0.0383714,0.674271,0.743552); rgb(27pt)=(0.0589714,0.683757,0.725386); rgb(28pt)=(0.0843,0.692833,0.706167); rgb(29pt)=(0.113295,0.7015,0.685857); rgb(30pt)=(0.145271,0.709757,0.664629); rgb(31pt)=(0.180133,0.717657,0.642433); rgb(32pt)=(0.217829,0.725043,0.619262); rgb(33pt)=(0.258643,0.731714,0.595429); rgb(34pt)=(0.302171,0.737605,0.571186); rgb(35pt)=(0.348167,0.742433,0.547267); rgb(36pt)=(0.395257,0.7459,0.524443); rgb(37pt)=(0.44201,0.748081,0.503314); rgb(38pt)=(0.487124,0.749062,0.483976); rgb(39pt)=(0.530029,0.749114,0.466114); rgb(40pt)=(0.570857,0.748519,0.44939); rgb(41pt)=(0.609852,0.747314,0.433686); rgb(42pt)=(0.6473,0.7456,0.4188); rgb(43pt)=(0.683419,0.743476,0.404433); rgb(44pt)=(0.71841,0.741133,0.390476); rgb(45pt)=(0.752486,0.7384,0.376814); rgb(46pt)=(0.785843,0.735567,0.363271); rgb(47pt)=(0.818505,0.732733,0.34979); rgb(48pt)=(0.850657,0.7299,0.336029); rgb(49pt)=(0.882433,0.727433,0.3217); rgb(50pt)=(0.913933,0.725786,0.306276); rgb(51pt)=(0.944957,0.726114,0.288643); rgb(52pt)=(0.973895,0.731395,0.266648); rgb(53pt)=(0.993771,0.745457,0.240348); rgb(54pt)=(0.999043,0.765314,0.216414); rgb(55pt)=(0.995533,0.786057,0.196652); rgb(56pt)=(0.988,0.8066,0.179367); rgb(57pt)=(0.978857,0.827143,0.163314); rgb(58pt)=(0.9697,0.848138,0.147452); rgb(59pt)=(0.962586,0.870514,0.1309); rgb(60pt)=(0.958871,0.8949,0.113243); rgb(61pt)=(0.959824,0.921833,0.0948381); rgb(62pt)=(0.9661,0.951443,0.0755333); rgb(63pt)=(0.9763,0.9831,0.0538)},
shader=flat,
mesh/rows=21]
table[row sep=crcr,header=false] {
0.5 0.25 NaN\\
0.5 0.275 NaN\\
0.5 0.3 NaN\\
0.5 0.325 NaN\\
0.5 0.35 NaN\\
0.5 0.375 NaN\\
0.5 0.4 NaN\\
0.5 0.425 NaN\\
0.5 0.45 NaN\\
0.5 0.475 NaN\\
0.5 0.5 NaN\\
0.5 0.525 NaN\\
0.5 0.55 NaN\\
0.5 0.575 NaN\\
0.5 0.6 NaN\\
0.5 0.625 NaN\\
0.5 0.65 NaN\\
0.5 0.675 NaN\\
0.5 0.7 NaN\\
0.5 0.725 NaN\\
0.5 0.75 NaN\\
0.525 0.25 NaN\\
0.525 0.275 NaN\\
0.525 0.3 0.505472645511438\\
0.525 0.325 0.504179979058338\\
0.525 0.35 0.505654912073014\\
0.525 0.375 0.505120002874666\\
0.525 0.4 0.503408775039938\\
0.525 0.425 0.503300112481884\\
0.525 0.45 0.507537663954041\\
0.525 0.475 0.505701331781485\\
0.525 0.5 0.506729213316042\\
0.525 0.525 0.50893430841863\\
0.525 0.55 0.50970008110871\\
0.525 0.575 0.511157660704608\\
0.525 0.6 0.511101017063075\\
0.525 0.625 0.51167454616288\\
0.525 0.65 0.512245006813844\\
0.525 0.675 0.51231701524406\\
0.525 0.7 0.512768855964746\\
0.525 0.725 0.512578642103757\\
0.525 0.75 NaN\\
0.55 0.25 NaN\\
0.55 0.275 0.523072244269166\\
0.55 0.3 0.523702141261038\\
0.55 0.325 0.526710649119111\\
0.55 0.35 0.528299121418073\\
0.55 0.375 0.527554824061571\\
0.55 0.4 0.527799154552396\\
0.55 0.425 0.528541041297922\\
0.55 0.45 0.52862856146321\\
0.55 0.475 0.528238299347615\\
0.55 0.5 0.527995047010641\\
0.55 0.525 0.52712196665187\\
0.55 0.55 0.525757756202988\\
0.55 0.575 0.526310452692093\\
0.55 0.6 0.52909419011238\\
0.55 0.625 0.529937890469236\\
0.55 0.65 0.529852525376725\\
0.55 0.675 0.528565444425879\\
0.55 0.7 0.526666370572504\\
0.55 0.725 0.523808857608707\\
0.55 0.75 NaN\\
0.575 0.25 NaN\\
0.575 0.275 0.538127056479564\\
0.575 0.3 0.539071373244971\\
0.575 0.325 0.541161468310804\\
0.575 0.35 0.540679798907068\\
0.575 0.375 0.542784053085679\\
0.575 0.4 0.54529295871567\\
0.575 0.425 0.543877891462911\\
0.575 0.45 0.544571995808801\\
0.575 0.475 0.544043205635587\\
0.575 0.5 0.542673325269318\\
0.575 0.525 0.553627186284122\\
0.575 0.55 0.551839082671927\\
0.575 0.575 0.544053041852217\\
0.575 0.6 0.544999627076963\\
0.575 0.625 0.54579045711896\\
0.575 0.65 0.544423656973944\\
0.575 0.675 0.542673896575822\\
0.575 0.7 0.53934742492717\\
0.575 0.725 0.534899718797195\\
0.575 0.75 NaN\\
0.6 0.25 NaN\\
0.6 0.275 0.552282224164848\\
0.6 0.3 0.553880914886881\\
0.6 0.325 0.555117378012954\\
0.6 0.35 0.559410415319002\\
0.6 0.375 0.563029309413577\\
0.6 0.4 0.563284314339959\\
0.6 0.425 0.564571617779367\\
0.6 0.45 0.562877666012874\\
0.6 0.475 0.566477228254021\\
0.6 0.5 0.577564820389281\\
0.6 0.525 0.58867258829784\\
0.6 0.55 0.58248335249556\\
0.6 0.575 0.575403434335986\\
0.6 0.6 0.559012668085527\\
0.6 0.625 0.557571303287212\\
0.6 0.65 0.555157553018178\\
0.6 0.675 0.555065722599765\\
0.6 0.7 0.554089218023106\\
0.6 0.725 0.551381463500164\\
0.6 0.75 NaN\\
0.625 0.25 NaN\\
0.625 0.275 0.56589883889356\\
0.625 0.3 0.568188479070764\\
0.625 0.325 0.571996700139608\\
0.625 0.35 0.575874851379851\\
0.625 0.375 0.578905000479817\\
0.625 0.4 0.590559220053071\\
0.625 0.425 0.603671484806638\\
0.625 0.45 0.596764317492581\\
0.625 0.475 0.594409331947414\\
0.625 0.5 0.602616801723453\\
0.625 0.525 0.604337123564155\\
0.625 0.55 0.605793754730696\\
0.625 0.575 0.607227213489543\\
0.625 0.6 0.59970827145398\\
0.625 0.625 0.587419338507424\\
0.625 0.65 0.576400312885311\\
0.625 0.675 0.572256097213322\\
0.625 0.7 0.569426996624796\\
0.625 0.725 0.564804821972503\\
0.625 0.75 NaN\\
0.65 0.25 NaN\\
0.65 0.275 0.578608273486436\\
0.65 0.3 0.581773670224006\\
0.65 0.325 0.58566217825148\\
0.65 0.35 0.596693669174415\\
0.65 0.375 0.610864152234301\\
0.65 0.4 0.610214932819908\\
0.65 0.425 0.60986631567453\\
0.65 0.45 0.611692068465835\\
0.65 0.475 0.613291621099162\\
0.65 0.5 0.615029037492501\\
0.65 0.525 0.616887064189681\\
0.65 0.55 0.617853549733131\\
0.65 0.575 0.618482447736663\\
0.65 0.6 0.617855412446234\\
0.65 0.625 0.611529832931125\\
0.65 0.65 0.594506876878675\\
0.65 0.675 0.587470646072899\\
0.65 0.7 0.582452482908088\\
0.65 0.725 0.577733115878273\\
0.65 0.75 NaN\\
0.675 0.25 NaN\\
0.675 0.275 0.591115180797687\\
0.675 0.3 0.594107830886689\\
0.675 0.325 0.602397422530139\\
0.675 0.35 0.606487155403792\\
0.675 0.375 0.621689920661471\\
0.675 0.4 0.62409595375034\\
0.675 0.425 0.624910286997161\\
0.675 0.45 0.623756261428594\\
0.675 0.475 0.625334909297422\\
0.675 0.5 0.626948830239757\\
0.675 0.525 0.6287381827371\\
0.675 0.55 0.629099988279099\\
0.675 0.575 0.629070877159857\\
0.675 0.6 0.62947883674268\\
0.675 0.625 0.629385460620547\\
0.675 0.65 0.626806482543632\\
0.675 0.675 0.609571438848911\\
0.675 0.7 0.595178989814535\\
0.675 0.725 0.591553830636701\\
0.675 0.75 NaN\\
0.7 0.25 NaN\\
0.7 0.275 0.603103325217141\\
0.7 0.3 0.607902948850472\\
0.7 0.325 0.628851057528228\\
0.7 0.35 0.636240160664926\\
0.7 0.375 0.641120370862984\\
0.7 0.4 0.6417190917636\\
0.7 0.425 0.641854522583489\\
0.7 0.45 0.639704520640341\\
0.7 0.475 0.639076725796499\\
0.7 0.5 0.639511680430647\\
0.7 0.525 0.639727779435329\\
0.7 0.55 0.640306617880187\\
0.7 0.575 0.639613199723108\\
0.7 0.6 0.641429614307917\\
0.7 0.625 0.642313126041578\\
0.7 0.65 0.634871330889276\\
0.7 0.675 0.621106002489763\\
0.7 0.7 0.607276582368372\\
0.7 0.725 0.602311882526903\\
0.7 0.75 NaN\\
0.725 0.25 NaN\\
0.725 0.275 0.616884039781431\\
0.725 0.3 0.624738080361544\\
0.725 0.325 0.652512336622673\\
0.725 0.35 0.65687881022295\\
0.725 0.375 0.658280747021356\\
0.725 0.4 0.65827576010897\\
0.725 0.425 0.656735299300311\\
0.725 0.45 0.655895623459659\\
0.725 0.475 0.656774403118877\\
0.725 0.5 0.655013033569932\\
0.725 0.525 0.652816869958\\
0.725 0.55 0.652493447746977\\
0.725 0.575 0.652001894320537\\
0.725 0.6 0.65232598317083\\
0.725 0.625 0.652885194860816\\
0.725 0.65 0.649253423509748\\
0.725 0.675 0.635292205919027\\
0.725 0.7 0.61892793061767\\
0.725 0.725 0.612941308762305\\
0.725 0.75 NaN\\
0.75 0.25 NaN\\
0.75 0.275 0.636461637140775\\
0.75 0.3 0.654585598645127\\
0.75 0.325 0.665588092752055\\
0.75 0.35 0.669418037063921\\
0.75 0.375 0.672731396982973\\
0.75 0.4 0.675017977943458\\
0.75 0.425 0.674835517384416\\
0.75 0.45 0.674910463004376\\
0.75 0.475 0.674846431905522\\
0.75 0.5 0.671668895205466\\
0.75 0.525 0.669411596162158\\
0.75 0.55 0.668705970992671\\
0.75 0.575 0.66687315965461\\
0.75 0.6 0.66544212310498\\
0.75 0.625 0.663314902792863\\
0.75 0.65 0.660934456173345\\
0.75 0.675 0.64892139859806\\
0.75 0.7 0.633892083224402\\
0.75 0.725 0.623570734997707\\
0.75 0.75 NaN\\
0.775 0.25 NaN\\
0.775 0.275 0.642778673731734\\
0.775 0.3 0.660353893883464\\
0.775 0.325 0.674438647204868\\
0.775 0.35 0.678467768608246\\
0.775 0.375 0.682049927114702\\
0.775 0.4 0.68520895980144\\
0.775 0.425 0.688884926029725\\
0.775 0.45 0.690285387251844\\
0.775 0.475 0.690014424651055\\
0.775 0.5 0.783928794305322\\
0.775 0.525 0.703447199578057\\
0.775 0.55 0.687762253726163\\
0.775 0.575 0.682375249516265\\
0.775 0.6 0.67801164256617\\
0.775 0.625 0.674678813746585\\
0.775 0.65 0.671399786380155\\
0.775 0.675 0.667325208727879\\
0.775 0.7 0.655885715751959\\
0.775 0.725 0.634200161233109\\
0.775 0.75 NaN\\
0.8 0.25 NaN\\
0.8 0.275 0.649195327254884\\
0.8 0.3 0.666164276662542\\
0.8 0.325 0.680267012274078\\
0.8 0.35 0.683835724772039\\
0.8 0.375 0.689120545660808\\
0.8 0.4 0.693912204087028\\
0.8 0.425 0.698315007763262\\
0.8 0.45 0.701161814123545\\
0.8 0.475 0.703292817593592\\
0.8 0.5 1\\
0.8 0.525 0.815558683020543\\
0.8 0.55 0.701732366542172\\
0.8 0.575 0.698507961557964\\
0.8 0.6 0.696720503709305\\
0.8 0.625 0.691329387846529\\
0.8 0.65 0.68553414288035\\
0.8 0.675 0.681645451488278\\
0.8 0.7 0.675930452299051\\
0.8 0.725 0.644877935262181\\
0.8 0.75 NaN\\
0.825 0.25 NaN\\
0.825 0.275 0.657696863290759\\
0.825 0.3 0.67469817522547\\
0.825 0.325 0.690779276608957\\
0.825 0.35 0.693942006753002\\
0.825 0.375 0.699151163585749\\
0.825 0.4 0.703482594450304\\
0.825 0.425 0.708602379109689\\
0.825 0.45 0.714166586020513\\
0.825 0.475 0.713610883112134\\
0.825 0.5 0.780461278461772\\
0.825 0.525 0.713987424155232\\
0.825 0.55 0.71274665902261\\
0.825 0.575 0.712193296483207\\
0.825 0.6 0.710636664009398\\
0.825 0.625 0.703042428422149\\
0.825 0.65 0.696513017573759\\
0.825 0.675 0.690208452965853\\
0.825 0.7 0.677085767929139\\
0.825 0.725 0.653511818095715\\
0.825 0.75 NaN\\
0.85 0.25 NaN\\
0.85 0.275 0.660482100201738\\
0.85 0.3 0.670003393076301\\
0.85 0.325 0.686775528485864\\
0.85 0.35 0.705293755931173\\
0.85 0.375 0.711006031186179\\
0.85 0.4 0.715149736851303\\
0.85 0.425 0.718163906221248\\
0.85 0.45 0.721808656638602\\
0.85 0.475 0.722681209195213\\
0.85 0.5 0.724059901521644\\
0.85 0.525 0.725618242251153\\
0.85 0.55 0.722112885080832\\
0.85 0.575 0.721231279789079\\
0.85 0.6 0.718873526530712\\
0.85 0.625 0.712635919834461\\
0.85 0.65 0.706221338025217\\
0.85 0.675 0.697926899363672\\
0.85 0.7 0.670957134070823\\
0.85 0.725 0.65952974068517\\
0.85 0.75 NaN\\
0.875 0.25 NaN\\
0.875 0.275 0.667511479663821\\
0.875 0.3 0.673653018949114\\
0.875 0.325 0.697929255785691\\
0.875 0.35 0.713370473617977\\
0.875 0.375 0.721209318515701\\
0.875 0.4 0.726402988544962\\
0.875 0.425 0.73134022900032\\
0.875 0.45 0.732419031595428\\
0.875 0.475 0.733014002770592\\
0.875 0.5 0.735135521704359\\
0.875 0.525 0.733229313752117\\
0.875 0.55 0.732135699739647\\
0.875 0.575 0.729627116776133\\
0.875 0.6 0.726144022565251\\
0.875 0.625 0.721855144222411\\
0.875 0.65 0.715484342783991\\
0.875 0.675 0.698295097262818\\
0.875 0.7 0.674690630559089\\
0.875 0.725 0.66618058959346\\
0.875 0.75 NaN\\
0.9 0.25 NaN\\
0.9 0.275 0.671888659081964\\
0.9 0.3 0.67912974071055\\
0.9 0.325 0.696344475817362\\
0.9 0.35 0.721129509550106\\
0.9 0.375 0.727763128838666\\
0.9 0.4 0.734948866511291\\
0.9 0.425 0.740236644950266\\
0.9 0.45 0.742833703356078\\
0.9 0.475 0.74345815840299\\
0.9 0.5 0.743671264792403\\
0.9 0.525 0.742217928214385\\
0.9 0.55 0.739728112474102\\
0.9 0.575 0.736308617414371\\
0.9 0.6 0.73236854648558\\
0.9 0.625 0.728111342459697\\
0.9 0.65 0.723448200150154\\
0.9 0.675 0.706297235833327\\
0.9 0.7 0.680723681486428\\
0.9 0.725 0.67282324013069\\
0.9 0.75 NaN\\
0.925 0.25 NaN\\
0.925 0.275 0.675742008523967\\
0.925 0.3 0.683033007487227\\
0.925 0.325 0.696471042849547\\
0.925 0.35 0.728285336587685\\
0.925 0.375 0.733767005106664\\
0.925 0.4 0.742090220582246\\
0.925 0.425 0.748281204840598\\
0.925 0.45 0.752264222644997\\
0.925 0.475 0.75234474369684\\
0.925 0.5 0.752008470610487\\
0.925 0.525 0.75098074023689\\
0.925 0.55 0.747870835706009\\
0.925 0.575 0.741861119927463\\
0.925 0.6 0.737783802710446\\
0.925 0.625 0.725023454764451\\
0.925 0.65 0.72565206698038\\
0.925 0.675 0.700471182942635\\
0.925 0.7 0.686631422183904\\
0.925 0.725 0.678638842312809\\
0.925 0.75 NaN\\
0.95 0.25 NaN\\
0.95 0.275 0.67959535796597\\
0.95 0.3 0.686680067116815\\
0.95 0.325 0.696446196199938\\
0.95 0.35 0.714993067148598\\
0.95 0.375 0.732817832506431\\
0.95 0.4 0.739155408441924\\
0.95 0.425 0.753239749866806\\
0.95 0.45 0.757556009413587\\
0.95 0.475 0.758955894530504\\
0.95 0.5 0.760187242379049\\
0.95 0.525 0.760137004748179\\
0.95 0.55 0.755202529148636\\
0.95 0.575 0.749089351816448\\
0.95 0.6 0.737192765307298\\
0.95 0.625 0.717435672789084\\
0.95 0.65 0.704693929414651\\
0.95 0.675 0.697703293817272\\
0.95 0.7 0.691821942916391\\
0.95 0.725 0.683732613374027\\
0.95 0.75 NaN\\
0.975 0.25 NaN\\
0.975 0.275 0.683448707407973\\
0.975 0.3 0.692282617870923\\
0.975 0.325 0.699603001816264\\
0.975 0.35 0.710461043024809\\
0.975 0.375 0.716079590486136\\
0.975 0.4 0.722489285944224\\
0.975 0.425 0.75150225222867\\
0.975 0.45 0.763442617105365\\
0.975 0.475 0.765460667337421\\
0.975 0.5 0.768990262070842\\
0.975 0.525 0.769148904048168\\
0.975 0.55 0.762106293318419\\
0.975 0.575 0.753309566508431\\
0.975 0.6 0.738887069852503\\
0.975 0.625 0.715924998341198\\
0.975 0.65 0.705902857948013\\
0.975 0.675 0.701181245369461\\
0.975 0.7 0.695094415680569\\
0.975 0.725 0.685560914318484\\
0.975 0.75 NaN\\
1 0.25 NaN\\
1 0.275 NaN\\
1 0.3 NaN\\
1 0.325 NaN\\
1 0.35 NaN\\
1 0.375 NaN\\
1 0.4 NaN\\
1 0.425 NaN\\
1 0.45 NaN\\
1 0.475 NaN\\
1 0.5 NaN\\
1 0.525 NaN\\
1 0.55 NaN\\
1 0.575 NaN\\
1 0.6 NaN\\
1 0.625 NaN\\
1 0.65 NaN\\
1 0.675 NaN\\
1 0.7 NaN\\
1 0.725 NaN\\
1 0.75 NaN\\
};
\addplot3 [
color=black, 
mark size=1.0pt,
only marks,
mark=*,
mark options={solid}]
table[row sep=crcr] {
0.8 0.5 1\\
0.896750966300958 0.286776127457757 0.675309021228622\\
0.902771977462955 0.443307096739818 0.744836468149218\\
0.610354789491249 0.590563611606405 0.563738447957281\\
0.717545631591306 0.636844667860001 0.650155945238496\\
0.852848313794446 0.435990016586211 0.72214249788317\\
0.660761305835014 0.452417241303641 0.617324922289526\\
0.820130475649893 0.320075943110966 0.688343795990109\\
0.610456067310923 0.319729894706417 0.563224615829313\\
0.996628822967953 0.281100166085246 0.690045813560512\\
0.829330911141074 0.36177041731311 0.697124515786589\\
0.738927787432516 0.702469856909833 0.623317488530984\\
0.626578235140411 0.702924333907962 0.570028291796796\\
0.981412845867243 0.705104676070429 0.694308716428758\\
0.679608034336165 0.657771390884537 0.632088246304667\\
0.7767791920296 0.333471872572624 0.677414779170404\\
0.905025376762873 0.553778152360609 0.741089496835124\\
0.94166257336887 0.316020689047971 0.689828748303971\\
0.818411855866936 0.603014587398564 0.706250552185257\\
0.666906186164037 0.438275770234889 0.618816128338631\\
0.813503588228468 0.595511914424481 0.706620932746701\\
0.625717367146037 0.668026662527658 0.573468959215293\\
0.906993019500661 0.662923413406259 0.722785175277229\\
0.572131847654909 0.48344100136562 0.542121740197994\\
0.903338950738432 0.571324421969113 0.737258088267436\\
0.986631397979414 0.628632455243067 0.709988327204149\\
0.890279572452008 0.746189348982611 0.662629687735164\\
0.804770328186065 0.384360635226165 0.692580163386609\\
0.556768447472082 0.688768012719625 0.530881466281883\\
0.597875400608463 0.402327162870215 0.561533615132784\\
0.734634322479852 0.517831329437117 0.658064668491838\\
0.588734931310846 0.731851647500458 0.544922691831549\\
0.805532816497391 0.641278569924126 0.690018852915826\\
0.761772587446987 0.495623970862494 0.680982433685929\\
0.73431033983156 0.374160408174907 0.664255227959639\\
0.957096547913993 0.505259946445197 0.762223343870089\\
0.780291105767305 0.454155929647519 0.69313996891296\\
0.892702383484364 0.563241363149134 0.736461140265627\\
0.612013984234454 0.339793008902145 0.567895093097158\\
0.50868945810227 0.624583002688413 0.499747064178206\\
0.598022997067236 0.648736113587461 0.554415809479207\\
0.712619762620216 0.373497888737766 0.650562973601767\\
0.64966640763613 0.423645275337246 0.609598238101671\\
0.957072824374163 0.684798012434607 0.696533967318057\\
0.670982508455675 0.451325421246072 0.620729799632379\\
0.87348890039871 0.636667687206959 0.719307792836015\\
0.91187476865064 0.72621176631384 0.675490958981554\\
0.689641681209292 0.635260414894502 0.638154645929706\\
0.557804714387723 0.625122397308514 0.535643294915641\\
0.800424794384846 0.697439096217222 0.680051355152973\\
0.504686378036765 0.419115029941283 0.483917193903741\\
0.529101751211831 0.641041544298777 0.515130143874614\\
0.953672341494932 0.329437527145471 0.69958350926354\\
0.556778717749001 0.44797090722905 0.533155033808704\\
0.600894508295261 0.465190351306326 0.56340006780058\\
0.833294096726429 0.603108885353926 0.713581976767414\\
0.989547621752577 0.736941860120139 0.682998345505369\\
0.697173860239204 0.637819207435455 0.641633978388963\\
0.508927167607604 0.532144542167342 0.497609495956456\\
0.681553544189616 0.572471948843176 0.633478379382559\\
0.557628145279771 0.527774442485995 0.530539030772093\\
0.950399058359332 0.735603954354949 0.679524403159992\\
0.510836847606421 0.345729971491733 0.492143904644896\\
0.841152553007951 0.478634600177314 0.720271683775339\\
0.560435455011861 0.666191663819268 0.537082635902509\\
0.944134517275711 0.734594407286704 0.679595820872517\\
0.900442140858193 0.332479457434833 0.715822428890343\\
0.632193997405328 0.369266158090538 0.580368079025297\\
0.516461734167434 0.42273759940176 0.493080656996309\\
0.69121462913145 0.377368393397148 0.634655571882189\\
0.649212958682711 0.599332068352685 0.617378173075356\\
0.57079593357119 0.63423768817409 0.543658782476681\\
0.723136464599151 0.576815374708839 0.651058852233531\\
0.928378223516919 0.349078429673212 0.728966489828421\\
0.505653869344781 0.748549542309387 0.498126816813527\\
0.545763188443633 0.48466628529025 0.526366197826546\\
0.825722294028466 0.515472651567212 0.714984637378247\\
0.71129147049187 0.33969018138519 0.646776168311117\\
0.961465217690023 0.315560651972318 0.695495573327347\\
0.952407770222575 0.643638364186628 0.707007424789972\\
0.747516926326245 0.313847884230513 0.663838800420421\\
0.796674491635933 0.388888360603717 0.690134891259959\\
0.672641662882902 0.634830830121167 0.627749149954137\\
0.603313292450109 0.65347035723999 0.556242883319727\\
0.796814421268911 0.263945218665011 0.640126444715355\\
0.735269940189609 0.383113193894543 0.664412924367409\\
0.656979767272955 0.620716536738546 0.622162794916472\\
0.748165576186312 0.259010531534803 0.623499756625529\\
0.939433799432894 0.665870143073621 0.697910119302324\\
0.728568244242342 0.551385087432613 0.654287651814854\\
0.650647237047691 0.275570801569803 0.578696328150886\\
0.806878408035618 0.62185710361723 0.696375158056645\\
0.544194662207204 0.302608389824434 0.520197040840609\\
0.622955861825668 0.577622950816906 0.606377060264836\\
0.818883711816406 0.734932788204561 0.649414399883414\\
0.846119357594698 0.6915399704329 0.69615128847075\\
0.553766009195282 0.575454247763055 0.52847955256031\\
0.555666158880712 0.568677618897226 0.532072697044367\\
0.604435974345879 0.640463967654522 0.558265394585855\\
0.526349642470887 0.56563144105305 0.512113779990246\\
0.721059175232061 0.297450940320884 0.619373777548769\\
0.84260103603702 0.345530155179932 0.701431949683509\\
0.759685718167462 0.557020628247086 0.675970229981688\\
0.923017767921153 0.69036343882904 0.689327389982002\\
0.557512611368792 0.261906665471005 0.527388292044186\\
0.544360948798569 0.468744207168826 0.525051508403691\\
0.629938261821685 0.606467334820059 0.610013447328318\\
0.917675500102484 0.733426039517901 0.671865658122397\\
0.68074167652991 0.688535317260708 0.600157758549229\\
0.617326694267432 0.382581145538756 0.577559172545856\\
0.681803363676878 0.518240854847237 0.631146415292364\\
0.885290171061452 0.558902348891923 0.734862488107348\\
0.999879067839387 0.448812646700452 0.727401342130736\\
0.556653679428374 0.591775511200523 0.533768996349688\\
0.511923460323181 0.731515256518571 0.505799371823714\\
0.595068302854429 0.720931728228616 0.548981064911093\\
0.653528993732435 0.524999851559058 0.61865860201313\\
0.791550251926774 0.361782993747501 0.68278854035169\\
0.713931367162337 0.405706449224023 0.651902070708594\\
0.805166281940449 0.43454148810554 0.702480128778132\\
0.517319137001344 0.503384924916418 0.502280002033299\\
0.79562218975414 0.68645392683332 0.678239268896186\\
0.591719690181946 0.695977721199807 0.550224315876001\\
0.659452620544525 0.731721188595494 0.580120785553871\\
0.993170985365447 0.741377113835917 0.680560651364204\\
0.791545234841724 0.475131522881481 0.699803137652766\\
0.579610587190408 0.591918525702021 0.547969590226936\\
0.919737896229818 0.434663358048463 0.749171090895793\\
0.555719911475165 0.551092672120086 0.529539821111044\\
0.757941722373145 0.373864137037787 0.677085077356758\\
0.654419236693616 0.275193740552363 0.581269366193757\\
0.889949054888152 0.294867099112871 0.675396038980631\\
0.668218844330459 0.326768892520437 0.595967269336275\\
0.918302268321324 0.647871097589267 0.728182273286895\\
0.652359263807669 0.637841326006166 0.592377831349789\\
0.553472481516696 0.476149594730903 0.530102452886497\\
0.916110193482973 0.45870182637889 0.748810824127394\\
0.78972060934824 0.60099276064908 0.68606043423479\\
0.966813818539717 0.434865117015749 0.762093748356929\\
0.853959402955557 0.485460244964119 0.724372954106271\\
0.540527794749159 0.431787412818572 0.52207101613061\\
0.505152550821154 0.515459444839436 0.493772051085467\\
0.578734333094401 0.690784956846228 0.543130223080311\\
0.897592273023424 0.57268934326616 0.736327023861278\\
0.80170852361882 0.74007759751037 0.642575041031485\\
0.956214471123851 0.637671814284754 0.707007273029842\\
0.514678816739023 0.505035171806744 0.500945083073484\\
0.665406548725784 0.429975246137746 0.618270699688186\\
0.972869269087839 0.53813192279225 0.763021204707754\\
0.71483991465287 0.570295440291489 0.64771130416623\\
0.523531322672027 0.503042847571715 0.505733330354179\\
0.848968807250136 0.526796175461105 0.725362570365034\\
0.763914494242464 0.416629221251442 0.685106029564832\\
0.701662840439768 0.538870654105546 0.640731030598895\\
0.979893450066785 0.454617438907175 0.764440791492274\\
0.803772445643009 0.682811751090593 0.682563195583193\\
0.615098939296913 0.268404367035586 0.560409569780512\\
0.798968557213536 0.38128354172878 0.690361143289597\\
0.890580710042851 0.738437300921625 0.66574687597849\\
0.967825084618818 0.536888363538073 0.763662347698319\\
0.69581419182488 0.496871966433736 0.637261007266206\\
0.544600752813757 0.40197995989309 0.524219720937104\\
0.936213193985016 0.444417592226183 0.754187416251243\\
0.598289326940254 0.635668147698973 0.554921602890308\\
0.851877230121635 0.653673406242905 0.705943811577635\\
0.752586851411443 0.472454038691497 0.676891838208833\\
0.837450216826856 0.522313656458288 0.720207703574648\\
0.889743706793846 0.28439031811912 0.672820801564725\\
0.981723549742713 0.345842180573096 0.706260925874137\\
0.862315767578514 0.26179302017712 0.661682652642288\\
0.529272376947821 0.405989155670227 0.507163967227763\\
0.967912669818733 0.592303195473019 0.749972190622589\\
0.830431113487101 0.539488684926696 0.715201113476676\\
0.621571966395797 0.414129110039419 0.602728306086779\\
0.57442603944283 0.735784725196401 0.532651595604926\\
0.687742464161246 0.455322621044195 0.632589635526563\\
0.590071029823612 0.618542337447304 0.552619748375517\\
0.617331885338194 0.405718710051996 0.577770753496953\\
0.72690845906154 0.70082605744322 0.619682882941244\\
0.634045924174415 0.63700555975707 0.585139996790123\\
0.752559219350107 0.257357444298136 0.624877252632397\\
0.75771271505066 0.577808150277288 0.671874743860647\\
0.838475119768797 0.465907037408674 0.719186688398525\\
0.844003912299951 0.510567992926764 0.722261228753131\\
0.976611594793114 0.685936179491611 0.699122702982295\\
0.981935894624428 0.493982272449778 0.771543556991402\\
0.923052998305207 0.457167450939811 0.752686444744642\\
0.652965308525666 0.607145655266464 0.620260326615805\\
0.860857500756139 0.384450601326728 0.71876365208381\\
0.622211396677159 0.355232097780553 0.575115262781263\\
0.791457906927876 0.436562126700811 0.69441088410237\\
0.674450201334774 0.635700701390973 0.628614584313959\\
0.864458138705059 0.494853097534895 0.730749269667513\\
0.665835532037635 0.457394875143017 0.61939136049261\\
0.621725699247968 0.711399117134646 0.5639306875787\\
0.770754244073012 0.650584031585053 0.66890850052086\\
0.730986697506542 0.295438385894555 0.62381936272907\\
0.800798167895389 0.356574648567488 0.684913854525035\\
0.592460076523037 0.456185420555379 0.554762345567348\\
0.842229522803865 0.277063905713421 0.659125457257369\\
0.526548414257021 0.267818589565266 0.508335320903809\\
0.835137254585206 0.628426907790893 0.706321188458983\\
0.908628810594412 0.741702177161977 0.668025428431373\\
0.8532643740468 0.42243260461999 0.71876954913036\\
0.522834301020647 0.670081960160418 0.510979527754725\\
0.946148794562511 0.749434273566672 0.67334643409243\\
0.604877099274976 0.320179925465855 0.557932819593626\\
0.821159119762907 0.450084003570264 0.712971968212817\\
0.655707947618119 0.741305278209221 0.576779452822177\\
0.818253846190852 0.578555867900184 0.709606985359673\\
0.952615054258229 0.427969357910103 0.754238552086707\\
0.564923381251329 0.489770209480822 0.536056090384328\\
0.676839393873037 0.303291533227662 0.595383883950727\\
0.848021602973121 0.2690429776502 0.658061759689753\\
0.734391454067498 0.417013474230646 0.663262506907909\\
0.976260164764427 0.517455232249628 0.769772220628357\\
0.646539376807889 0.574127157107886 0.617052273164587\\
0.708068392333254 0.689023454361018 0.612824007755618\\
0.565127848267236 0.457550272907665 0.538825784733203\\
0.923272787027268 0.693320898060286 0.688403487729057\\
0.64958740767562 0.368280328610595 0.610046482483158\\
0.755269923936985 0.592352884283737 0.66806589801012\\
0.829613366395353 0.709263676905948 0.658781266657541\\
0.687090495956269 0.560597930130718 0.634719476988533\\
0.708352078901677 0.459898933306801 0.643098526713472\\
0.60356654986875 0.524604656439891 0.593554547712214\\
0.633827593627922 0.728606225258041 0.569314904316656\\
0.806227835509854 0.722240128754432 0.648093846629936\\
0.509598833595111 0.472527691077761 0.493897030675512\\
0.795448592145339 0.40713967529347 0.693614124668796\\
0.857848427994472 0.698109151242465 0.671135773347803\\
0.531920624360976 0.741644972972042 0.514361294879432\\
0.8557431661952 0.30361060242885 0.670628806367391\\
0.96809664114014 0.518059368780743 0.767298335282901\\
0.563520581945619 0.625875894105425 0.540711886416362\\
0.653411476759634 0.314888736357681 0.585951995496444\\
0.680924942343551 0.525120814126042 0.631533688868062\\
0.991879636942119 0.462982362699806 0.770774079733044\\
0.595393067151581 0.427187978982608 0.55643188734346\\
0.963869982901046 0.519546686704649 0.766723007126275\\
0.966957924360731 0.745079342721426 0.676502803142463\\
0.553432214382848 0.320014753868527 0.529217595190137\\
0.638679974267844 0.470639301222573 0.607459747857872\\
0.873411670291924 0.505699577758516 0.734836526217227\\
0.732705912657507 0.530311976941279 0.655035654469484\\
0.964607985837518 0.31492514986018 0.695493388544131\\
0.94266717920143 0.589952678564109 0.743287129726907\\
0.576749240556923 0.571185597056978 0.546095871445357\\
0.957490711898216 0.528601842529843 0.762629282734032\\
0.830914484939197 0.694004136740118 0.689733389092391\\
0.535997183025795 0.474256633966791 0.514027475652709\\
0.543195953570968 0.589344789409276 0.523729431038071\\
0.869065550528584 0.382356373953292 0.720807183579411\\
0.791276864595015 0.468278281572604 0.698144782372788\\
0.578746863621295 0.582739120649017 0.545227904157538\\
0.956245399332731 0.706092484113529 0.691603959345842\\
0.622942837668894 0.644073384656418 0.576101987040201\\
0.565985873894862 0.364046567175199 0.535377153642951\\
0.789786835456242 0.405256210517513 0.689571607438769\\
0.672958216453504 0.34073261981703 0.599400257435966\\
0.802669041217025 0.596071307223995 0.700707358840591\\
0.935047292474637 0.647493043624819 0.73311845007489\\
0.908829310988788 0.624202335819216 0.730469157247578\\
0.961205847705437 0.412931052453473 0.752663083835372\\
0.670425202091271 0.589316936077182 0.626567870645837\\
0.809658643377415 0.432870142943479 0.704457714049562\\
0.6712806408244 0.38353675211867 0.624237478763425\\
0.902630827451677 0.321359962889079 0.684880077285516\\
0.703072531314131 0.574549578288336 0.640589291436286\\
0.547896315556114 0.342398023585522 0.527538969950843\\
0.972330940570122 0.582904678276249 0.752560903530411\\
0.678205364875292 0.564728636460807 0.630813652481144\\
0.642098093928631 0.379125673258032 0.60651992940854\\
0.751025442208974 0.58896013965036 0.662532679644722\\
0.721625935559945 0.252011485651552 0.610721655584641\\
0.942298156389932 0.406466802338889 0.748632945823733\\
0.979906883588691 0.530680066943437 0.770440088595673\\
0.783828554204179 0.575665967638753 0.687440558525653\\
0.866626713769668 0.514649457161133 0.731139898773891\\
0.760762604834385 0.659430632841898 0.663808860277865\\
0.543295865291924 0.416267580225373 0.524255098379243\\
0.982508217600462 0.38017318360995 0.711995967296736\\
0.942092747848833 0.382032002277144 0.739394287101452\\
0.817249939839354 0.264316819522339 0.650420481227829\\
0.861157368692872 0.457323928262554 0.725055342046589\\
0.859715064967857 0.613175275368393 0.718619427991183\\
0.800533289280708 0.319337808075139 0.679552021157372\\
0.734069621143138 0.380003045220266 0.663596872326873\\
0.81696018134481 0.662465188355179 0.688874705956616\\
0.665854300939101 0.715876169252962 0.587588879580629\\
0.879373911625107 0.373192954520311 0.722418458640827\\
0.97426709327304 0.395825166712584 0.715542370927374\\
0.573687924020635 0.50315311808187 0.541504401179271\\
0.696025311674954 0.601471771434695 0.639709676443671\\
0.549451515527481 0.522176640084032 0.526910190358761\\
0.733492408538993 0.568997944965594 0.656167144508744\\
0.634401544897786 0.30694619947175 0.573322160334179\\
0.875440888430685 0.607321435292432 0.725252191713433\\
0.971286121553407 0.405143118697132 0.752648839481407\\
0.668750702074159 0.411499750601688 0.619697352833719\\
0.827322055000691 0.634115377637947 0.703081931588635\\
0.721767632463402 0.329166134213662 0.653944257627703\\
0.57157316852064 0.59112884077839 0.543090497398338\\
0.785241574553222 0.554601571978825 0.695203229064323\\
0.939391442171645 0.516290193501678 0.756899755047945\\
0.868834348683977 0.441901831283072 0.732484535693276\\
0.525470794349296 0.503514562234064 0.507289290270958\\
0.744236395914103 0.664753411807545 0.657643631744071\\
0.979659004758977 0.275471362862015 0.684323770096047\\
0.59475169085539 0.607946713978303 0.556823610020793\\
0.673591549846092 0.730304578210596 0.590034854607251\\
0.601673188688195 0.632126487108433 0.557106714989733\\
0.763778061552928 0.502029186103611 0.681664103165553\\
};
\end{axis}
\end{tikzpicture}%

%% file: SmoothEpsConnVError_2DModel1_n1280_p4.tikz
\begin{tikzpicture}
\def\xl{-0.035}
\def\xu{0.26}
\def\yl{0}
\def\yu{0.4}

\begin{axis}[%
width=\figurewidth,
height=\figureheight,
scale only axis,
xmin={\xl-(\xu-\xl)/5},
xmax={\xu+(\xu-\xl)/5},
ymin={\yl-(\yu-\yl)/5},
ymax={\yu+(\yu-\yl)*0.05},
hide axis,
axis background/.style={fill=white!100}
]
\draw [->] (axis cs: \xl,\yl) -- (axis cs: {\xu+(\xu-\xl)*0.05},\yl);
\draw [->] (axis cs: \xl,\yl) -- (axis cs: \xl,{\yu+(\yu-\yl)*0.05});
\node at (axis cs: {(\xu+(\xu-\xl)*0.05+\xl)/2},{\yl-0.16*(\yu-\yl)}) {$\eps-\eps_{\mathrm{conn}}$};
\node[rotate=90] at (axis cs: {\xl-0.16*(\xu-\xl)},{(\yu+(\yu-\yl)*0.05+\yl)/2}) {$\err_n^{(4)}(f_n)$};
\foreach \yValue in {0,0.1,0.2,0.3,0.4} {
    \edef\temp{\noexpand\draw [-] ({\xl+(\xu-\xl)/100},\yValue) -- (\xl,\yValue) node[left] {\yValue};} 
    \temp
}
\foreach \xValue in {0,0.05,0.1,0.15,0.2,0.25} {
	\edef\temp{\noexpand\draw [-] (\xValue,{\yl+(\yu-\yl)/100}) -- (\xValue,\yl) node[below] {\xValue};}    
    \temp
}
\addplot [
color=black,
solid,
mark options={solid},
line width=1.5pt,
forget plot
]
table[row sep=crcr]{
-0.0330612244897959 0.319333328843084\\
-0.0257142857142857 0.282404288810456\\
-0.0183673469387755 0.261158993315612\\
-0.0110204081632653 0.151920700812303\\
-0.00367346938775511 0.0805710651581672\\
0.00367346938775509 0.0396704786073163\\
0.0110204081632653 0.0526225079593209\\
0.0183673469387755 0.0661489834977423\\
0.0257142857142857 0.0788847941798471\\
0.0330612244897959 0.0936514058140233\\
0.0404081632653061 0.101834247231017\\
0.0477551020408163 0.111420581870816\\
0.0551020408163265 0.119711420433269\\
0.0624489795918367 0.125248602240685\\
0.0697959183673469 0.130281483421952\\
0.0771428571428571 0.135656684579896\\
0.0844897959183673 0.13861298635686\\
0.0918367346938776 0.143334777933284\\
0.0991836734693878 0.14837196715967\\
0.106530612244898 0.152407673331202\\
0.113877551020408 0.155126728177351\\
0.121224489795918 0.159048434836214\\
0.128571428571429 0.161583104654381\\
0.135918367346939 0.165491935705983\\
0.143265306122449 0.169262748524062\\
0.150612244897959 0.17240486317562\\
0.157959183673469 0.175435246105262\\
0.16530612244898 0.179557427993269\\
0.17265306122449 0.182311889426967\\
0.18 0.18626887793178\\
0.18734693877551 0.190165437973264\\
0.19469387755102 0.193468750496589\\
0.202040816326531 0.197007021250273\\
0.209387755102041 0.200211192537866\\
0.216734693877551 0.20266880786828\\
0.224081632653061 0.206263566149743\\
0.231428571428571 0.209827923507983\\
0.238775510204082 0.212998674731556\\
0.246122448979592 0.215555744094105\\
0.253469387755102 0.216173495696284\\
};
\addplot [
color=black,
dashed,
line width=1pt,
forget plot
]
table[row sep=crcr]{
-0.0330612244897959 0.328096392892804\\
-0.0257142857142857 0.333511368202853\\
-0.0183673469387755 0.370138516350193\\
-0.0110204081632653 0.233515087839971\\
-0.00367346938775511 0.124719405624537\\
0.00367346938775509 0.0516782997137671\\
0.0110204081632653 0.0684037509437666\\
0.0183673469387755 0.0869113631408184\\
0.0257142857142857 0.0926167573084404\\
0.0330612244897959 0.105085628683409\\
0.0404081632653061 0.112335212233349\\
0.0477551020408163 0.120473064654056\\
0.0551020408163265 0.126787984810976\\
0.0624489795918367 0.132644817188272\\
0.0697959183673469 0.136984678668434\\
0.0771428571428571 0.141188312540389\\
0.0844897959183673 0.144009627849754\\
0.0918367346938776 0.148919529358348\\
0.0991836734693878 0.152748878807354\\
0.106530612244898 0.15773648123065\\
0.113877551020408 0.160077517320811\\
0.121224489795918 0.163006073273156\\
0.128571428571429 0.165977323402532\\
0.135918367346939 0.169437412271668\\
0.143265306122449 0.173164705537809\\
0.150612244897959 0.177749573334737\\
0.157959183673469 0.179262917522618\\
0.16530612244898 0.183146895359947\\
0.17265306122449 0.185841825234973\\
0.18 0.190126284830811\\
0.18734693877551 0.193803330711826\\
0.19469387755102 0.196876530295128\\
0.202040816326531 0.200867126086956\\
0.209387755102041 0.203776720567651\\
0.216734693877551 0.206421961840791\\
0.224081632653061 0.210093361665516\\
0.231428571428571 0.213384066093644\\
0.238775510204082 0.216394627548196\\
0.246122448979592 0.217642911413715\\
0.253469387755102 0.217612910478752\\
};
\addplot [
color=black,
dashed,
line width=1pt,
forget plot
]
table[row sep=crcr]{
-0.0330612244897959 0.312631787438358\\
-0.0257142857142857 0.188837970706642\\
-0.0183673469387755 0.114352271225189\\
-0.0110204081632653 0.0642380914400652\\
-0.00367346938775511 0.0382959956023033\\
0.00367346938775509 0.029677964523765\\
0.0110204081632653 0.0384756208139897\\
0.0183673469387755 0.0482597121433615\\
0.0257142857142857 0.065190186309635\\
0.0330612244897959 0.0813733704398543\\
0.0404081632653061 0.0910149199823701\\
0.0477551020408163 0.10176023586781\\
0.0551020408163265 0.112070554604923\\
0.0624489795918367 0.117761159679668\\
0.0697959183673469 0.123450978566745\\
0.0771428571428571 0.129970253717508\\
0.0844897959183673 0.133183705589932\\
0.0918367346938776 0.137408400220266\\
0.0991836734693878 0.143453380350383\\
0.106530612244898 0.146752331705343\\
0.113877551020408 0.150288944311049\\
0.121224489795918 0.154583115741296\\
0.128571428571429 0.157306404921797\\
0.135918367346939 0.161224014737599\\
0.143265306122449 0.165626902744966\\
0.150612244897959 0.168363386534809\\
0.157959183673469 0.171211463793102\\
0.16530612244898 0.175629808519532\\
0.17265306122449 0.178258348044894\\
0.18 0.182324986715554\\
0.18734693877551 0.1864483433428\\
0.19469387755102 0.189908508847346\\
0.202040816326531 0.192530888717584\\
0.209387755102041 0.197353917409116\\
0.216734693877551 0.199668917685657\\
0.224081632653061 0.202878398270384\\
0.231428571428571 0.206888205589\\
0.238775510204082 0.209793784531592\\
0.246122448979592 0.213215027989101\\
0.253469387755102 0.214574725927619\\
};
\end{axis}
\end{tikzpicture}%

%% file: Error_n1280_Model2_2D_p4.tikz
\definecolor{purple}{rgb}{0.5020,0,0.5020}

\begin{tikzpicture}
\def\xl{0.03}
\def\xu{0.12}
\def\yl{0}
\def\yu{0.4}

\begin{axis}[%
width=\figurewidth,
height=\figureheight,
scale only axis,
xmin={\xl-(\xu-\xl)/5},
xmax={\xu+(\xu-\xl)/5},
ymin={\yl-(\yu-\yl)/5},
ymax={\yu+(\yu-\yl)*0.05},
hide axis,
axis background/.style={fill=white!100}
]
\draw [->] (axis cs: \xl,\yl) -- (axis cs: {\xu+(\xu-\xl)*0.05},\yl);
\draw [->] (axis cs: \xl,\yl) -- (axis cs: \xl,{\yu+(\yu-\yl)*0.05});
\draw [blue,->] (axis cs: \xu,\yl) -- (axis cs: \xu,{\yu+(\yu-\yl)*0.05});
\node at (axis cs: {(\xu+(\xu-\xl)*0.05+\xl)/2},{\yl-0.16*(\yu-\yl)}) {$\eps$};
\node[rotate=90] at (axis cs: {\xl-0.16*(\xu-\xl)},{(\yu+(\yu-\yl)*0.05+\yl)/2}) {$\err_n^{(4)}(f_n)$};
\node[blue,rotate=90] at (axis cs: {\xu+0.16*(\xu-\xl)},{(\yu+(\yu-\yl)*0.05+\yl)/2}) {$\%$ Graphs Connected};
\foreach \yValue in {0,0.1,0.2,0.3,0.4,0.5,0.6} {
    \edef\temp{\noexpand\draw [-] ({\xl+(\xu-\xl)/100},\yValue) -- (\xl,\yValue) node[left] {\yValue};} 
    \temp
}
\foreach \yValue in {20,40,60,80,100} {
    \edef\temp{\noexpand\draw [blue,-] ({\xu-(\xu-\xl)/100},\yValue*\yu/100) -- (\xu,\yValue*\yu/100) node[right] {\yValue};} 
    \temp
}
\foreach \xValue in {0.03,0.06,0.09,0.12} {
	\edef\temp{\noexpand\draw [-] (\xValue,{\yl+(\yu-\yl)/100}) -- (\xValue,\yl) node[below] {\xValue};}    
    \temp
}

\addplot [
color=yellow,
mark size=2.0pt,
only marks,
mark=diamond*,
mark options={solid},
forget plot
]
table[row sep=crcr]{
0.0998 0.0639\\
};
\addplot [
color=yellow,
solid,
mark options={solid},
line width=1.5pt,
forget plot
]
table[row sep=crcr]{
0.03 0.318914572948825\\
0.0318367346938776 0.347196312499084\\
0.0336734693877551 0.30878894013019\\
0.0355102040816326 0.23693090491882\\
0.0373469387755102 0.160358208046852\\
0.0391836734693878 0.135434748002313\\
0.0410204081632653 0.115960584777222\\
0.0428571428571429 0.116375735858256\\
0.0446938775510204 0.128650035627268\\
0.046530612244898 0.119380749445566\\
0.0483673469387755 0.102427553212769\\
0.0502040816326531 0.100926165801522\\
0.0520408163265306 0.0877659128434039\\
0.0538775510204082 0.0906236524512034\\
0.0557142857142857 0.087641403153616\\
0.0575510204081633 0.0896956421992476\\
0.0593877551020408 0.0888940901658814\\
0.0612244897959184 0.090660230089644\\
0.0630612244897959 0.0921948182208781\\
0.0648979591836735 0.0942076229727873\\
0.066734693877551 0.0963822189977471\\
0.0685714285714286 0.0986597675206774\\
0.0704081632653061 0.101484487905546\\
0.0722448979591837 0.104220573511792\\
0.0740816326530612 0.10651813611101\\
0.0759183673469388 0.108777452460548\\
0.0777551020408163 0.111351817273623\\
0.0795918367346939 0.114053180695064\\
0.0814285714285714 0.11676950407523\\
0.083265306122449 0.119339875853053\\
0.0851020408163265 0.122034747346621\\
0.0869387755102041 0.124729834836394\\
0.0887755102040816 0.127453585671757\\
0.0906122448979592 0.130242396443778\\
0.0924489795918367 0.133029901800148\\
0.0942857142857143 0.135840376597612\\
0.0961224489795918 0.138712994018755\\
0.0979591836734694 0.141740713100497\\
0.0997959183673469 0.144896007892697\\
0.101632653061224 0.147779128180994\\
0.103469387755102 0.151003261687961\\
0.10530612244898 0.153861633477483\\
0.107142857142857 0.156923860087184\\
0.108979591836735 0.159888858748047\\
0.110816326530612 0.162983979653193\\
0.11265306122449 0.165872522732078\\
0.114489795918367 0.168983077066924\\
0.116326530612245 0.172246368956991\\
0.118163265306122 0.175419858101836\\
0.12 0.178567805871698\\
};
\addplot [
color=blue,
solid,
mark options={solid},
line width=1.5pt,
forget plot
]
table[row sep=crcr]{
0.030000000000000 0\\
0.031836734693878 0\\
0.033673469387755 0\\
0.035510204081633 0\\
0.037346938775510 0\\
0.039183673469388 0\\
0.041020408163265 0\\
0.042857142857143 0.012000000000000\\
0.044693877551020 0.056000000000000\\
0.046530612244898 0.100000000000000\\
0.048367346938776 0.168000000000000\\
0.050204081632653 0.212000000000000\\
0.052040816326531 0.272000000000000\\
0.053877551020408 0.316000000000000\\
0.055714285714286 0.336000000000000\\
0.057551020408163 0.356000000000000\\
0.059387755102041 0.360000000000000\\
0.061224489795918 0.368000000000000\\
0.063061224489796 0.372000000000000\\
0.064897959183673 0.376000000000000\\
0.066734693877551 0.380000000000000\\
0.068571428571429 0.392000000000000\\
0.070408163265306 0.396000000000000\\
0.072244897959184 0.396000000000000\\
0.074081632653061 0.400000000000000\\
0.075918367346939 0.400000000000000\\
0.077755102040816 0.400000000000000\\
0.079591836734694 0.400000000000000\\
0.081428571428571 0.400000000000000\\
0.083265306122449 0.400000000000000\\
0.085102040816327 0.400000000000000\\
0.086938775510204 0.400000000000000\\
0.088775510204082 0.400000000000000\\
0.090612244897959 0.400000000000000\\
0.092448979591837 0.400000000000000\\
0.094285714285714 0.400000000000000\\
0.096122448979592 0.400000000000000\\
0.097959183673469 0.400000000000000\\
0.099795918367347 0.400000000000000\\
0.101632653061224 0.400000000000000\\
0.103469387755102 0.400000000000000\\
0.105306122448980 0.400000000000000\\
0.107142857142857 0.400000000000000\\
0.108979591836735 0.400000000000000\\
0.110816326530612 0.400000000000000\\
0.112653061224490 0.400000000000000\\
0.114489795918367 0.400000000000000\\
0.116326530612245 0.400000000000000\\
0.118163265306122 0.400000000000000\\
0.120000000000000 0.400000000000000\\
};
\addplot [
color=orange,
solid,
mark options={solid},
line width=1.5pt,
forget plot
]
table[row sep=crcr]{
0.0300 0.3214\\
0.0318 0.3462\\
0.0337 0.3425\\
0.0355 0.2813\\
0.0373 0.1867\\
0.0392 0.1442\\
0.0410 0.1158\\
0.0429 0.1112\\
0.0447 0.1275\\
0.0465 0.1195\\
0.0484 0.0838\\
0.0502 0.0636\\
0.0520 0.0526\\
0.0539 0.0485\\
0.0557 0.0479\\
0.0576 0.0479\\
0.0594 0.0465\\
0.0612 0.0460\\
0.0631 0.0465\\
0.0649 0.0472\\
0.0667 0.0475\\
0.0686 0.0473\\
0.0704 0.0481\\
0.0722 0.0490\\
0.0741 0.0495\\
0.0759 0.0504\\
0.0778 0.0513\\
0.0796 0.0522\\
0.0814 0.0532\\
0.0833 0.0541\\
0.0851 0.0552\\
0.0869 0.0564\\
0.0888 0.0575\\
0.0906 0.0587\\
0.0924 0.0599\\
0.0943 0.0607\\
0.0961 0.0617\\
0.0980 0.0629\\
0.0998 0.0639\\
0.1016 0.0653\\
0.1035 0.0665\\
0.1053 0.0678\\
0.1071 0.0691\\
0.1090 0.0704\\
0.1108 0.0715\\
0.1127 0.0728\\
0.1145 0.0741\\
0.1163 0.0755\\
0.1182 0.0768\\
0.1200 0.0780\\
};
\addplot [
color=black,
solid,
mark options={solid},
line width=1.5pt,
forget plot
]
table[row sep=crcr]{
0.03 0.32147522070433\\
0.0355102040816326 0.280189256869023\\
0.0410204081632653 0.123528841683459\\
0.046530612244898 0.0695725337187247\\
0.0520408163265306 0.0433756116094781\\
0.0575510204081633 0.0429961638397036\\
0.0630612244897959 0.0499855877451332\\
0.0685714285714286 0.0616524757228678\\
0.0740816326530612 0.0727504834527598\\
0.0795918367346939 0.0832828853717174\\
0.0851020408163265 0.0925616995006072\\
0.0906122448979592 0.100936594101105\\
0.0961224489795918 0.108045071178075\\
0.101632653061224 0.113602824116318\\
0.107142857142857 0.119193575984676\\
0.11265306122449 0.123807340101691\\
0.118163265306122 0.128170088478517\\
};
\addplot [
color=darkred,
solid,
mark options={solid},
line width=1.5pt,
forget plot
]
table[row sep=crcr]{
0.0300 0.3214\\
0.0318 0.3462\\
0.0337 0.3425\\
0.0355 0.2813\\
0.0373 0.1868\\
0.0392 0.1444\\
0.0410 0.1118\\
0.0429 0.0953\\
0.0447 0.0870\\
0.0465 0.0733\\
0.0484 0.0621\\
0.0502 0.0553\\
0.0520 0.0467\\
0.0539 0.0442\\
0.0557 0.0434\\
0.0576 0.0418\\
0.0594 0.0414\\
0.0612 0.0409\\
0.0631 0.0410\\
0.0649 0.0410\\
0.0667 0.0408\\
0.0686 0.0400\\
0.0704 0.0399\\
0.0722 0.0401\\
0.0741 0.0401\\
0.0759 0.0403\\
0.0778 0.0404\\
0.0796 0.0402\\
0.0814 0.0405\\
0.0833 0.0408\\
0.0851 0.0411\\
0.0869 0.0412\\
0.0888 0.0414\\
0.0906 0.0416\\
0.0924 0.0417\\
0.0943 0.0419\\
0.0961 0.0421\\
0.0980 0.0421\\
0.0998 0.0422\\
0.1016 0.0424\\
0.1035 0.0426\\
0.1053 0.0427\\
0.1071 0.0429\\
0.1090 0.0430\\
0.1108 0.0433\\
0.1127 0.0435\\
0.1145 0.0437\\
0.1163 0.0437\\
0.1182 0.0439\\
0.1200 0.0440\\
};
\end{axis}
\end{tikzpicture}%